\documentstyle[amsfonts, amssymb, epsf, makeidx]{amsbook}
\theoremstyle{plain}
 \newtheorem{thm}{Theorem}[chapter]
 \newtheorem{lem}[thm]{Lemma}
 \newtheorem{cor}[thm]{Corollary}
 \newtheorem{prop}[thm]{Proposition}

\theoremstyle{definition}
  \newtheorem{defn}{Definition}[chapter]
  \newtheorem{ex}{Example}[chapter]
  \newtheorem{assumption}{Assumption}[chapter]

\theoremstyle{remark}
  \newtheorem{rem}{Remark}[chapter]
  \newtheorem{claim}{Claim}

\newcommand{\ci}[2]{\cite[#1]{#2}}
\newcommand{\ta}[2]{\renewcommand{\arraystretch}{0.7}
                    \begin{tabular}{c} 
                      $\scriptstyle #1$ \\
                      $\scriptstyle #2$
                    \end{tabular}}

\numberwithin{equation}{chapter}
\numberwithin{figure}{chapter}

\begin{document}
\frontmatter
\title{The mapping class group from the viewpoint of measure equivalence theory}
\author{Yoshikata Kida}
\address{Department of Mathematics, Graduate School of Science, Kyoto University, Kyoto 606-8502, Japan}
\email{kida@@math.kyoto-u.ac.jp}

\maketitle

\tableofcontents

\chapter*{Abstract}
We obtain some classification result for the mapping class groups of compact orientable surfaces in terms of measure equivalence. In particular, the mapping class groups of different closed surfaces can not be measure equivalent. Moreover, we give various examples of discrete groups which are not measure equivalent to the mapping class groups. In the course of the proof, we investigate amenability in a measurable sense for the actions of the mapping class group on the boundary at infinity of the curve complex and on the Thurston boundary. Using this investigation, we prove that the mapping class group of a compact orientable surface is exact.

\footnote[0]{{\it Date}: October 7, 2005, revised on March 26, 2006.}
\footnote[0]{2000 {\it Mathematics Subject Classification}. Primary 20F38, 37A20.}
\footnote[0]{{\it Key words and phrases}. The mapping class group, the curve complex, measure equivalence, orbit equivalence.}

\mainmatter
\chapter{Introduction}\label{introduction}

Given a connected compact orientable surface $M$ which may have non-empty boundary, we define the {\it mapping class group}\index{mapping class group} as the group consisting of all isotopy classes of orientation-preserving self-diffeomorphisms of $M$. In this paper, we study the mapping class groups from the viewpoint of measure equivalence, which is a notion introduced by Gromov as a measure-theoretical analogue of a coarse geometric notion of quasi-isometry among finitely generated groups (see Remark \ref{rem-qi-me} for this analogy). Let us recall the definition of measure equivalence. Throughout the paper, discrete groups are assumed to be countable.

\begin{defn}[\ci{0.5.E}{gromov2}]
Two discrete groups $\Gamma$ and $\Lambda$ are said to be {\it measure equivalent}\index{measure equivalent} if there exist commuting, measure-preserving, essentially free actions of $\Gamma$ and $\Lambda$ on some standard Borel space $(\Omega, m)$ with a $\sigma$-finite non-zero positive measure such that the action of each of the groups $\Gamma$ and $\Lambda$ admits a fundamental domain with finite measure.
\end{defn}

In fact, this defines an equivalence relation among discrete groups \cite[Section 2]{furman1}. Note that two isomorphic groups modulo finite kernels and cokernels are measure equivalent. In this case, we say that the two groups are {\it almost isomorphic}\index{almost isomorphic}. It is easy to see that all finite groups form one complete class of measure equivalent groups.

Two lattices (i.e., discrete subgroups with cofinite volume) in the same locally compact second countable group are measure equivalent (see the comment right after Definition \ref{defn-me}). This example is one of the geometric motivations for introducing the notion of measure equivalence. On the other hand, we have another equivalent formulation for measure equivalence of two discrete groups, which is expressed in terms of orbit equivalence.

\begin{defn}\label{defn-oe}
Suppose that we have two essentially free, measure-preserving Borel actions of discrete groups $\Gamma_{1}$ and $\Gamma_{2}$ on standard probability spaces $(X_{1}, \mu_{1})$ and $(X_{2}, \mu_{2})$, respectively. The two actions are said to be {\it weakly orbit equivalent}\index{weakly!orbit equivalent} if there exist Borel subsets $A_{1}\subseteq X_{1}$ and $A_{2}\subseteq X_{2}$ with positive measure satisfying $\Gamma_{1}A_{1}=X_{1}$ and $\Gamma_{2}A_{2}=X_{2}$ and we have a Borel isomorphism $f\colon A_{1}\rightarrow A_{2}$ such that
\begin{enumerate}
\item[(i)] two measures $f_{*}(\mu_{1}|_{A_{1}})$ and $\mu_{2}|_{A_{2}}$ are equivalent;
\item[(ii)] we have $f(\Gamma_{1}x\cap A_{1})=\Gamma_{2}f(x)\cap A_{2}$ for a.e.\ $x\in A_{1}$.
\end{enumerate}
If we can take both $A_{1}$ and $A_{2}$ to have full measure, then the two actions are said to be {\it orbit equivalent}\index{orbit equivalent}.
\end{defn}

The study of this notion has much longer history than that of measure equivalence and it is related to ergodic theory and operator algebras, especially the construction of von Neumann algebras by Murray-von Neumann \cite{mvn}. It contributes to the recent solution of a longstanding problem on von Neumann algebras by Popa \cite{popa} through an excellent work of Gaboriau \cite{gab-l2} on $\ell^{2}$-Betti numbers of a measure-preserving group actions. We know that two discrete groups are measure equivalent if and only if they have essentially free, measure-preserving ergodic actions on standard probability spaces which are weakly orbit equivalent (see Theorem \ref{me-eq-oe}). We shall recall a few deep consequences related to measure equivalence and orbit equivalence. We recommend the reader to consult \cite{gab2} and \cite{gab-survey} for further topics about measure equivalence and orbit equivalence.

It is Dye \cite{dye}, \cite{dye2} who initiated study of orbit equivalence in the measure-preserving case. He treated some classes of infinite amenable groups (e.g., abelian ones and finitely generated ones with polynomial growth) from the viewpoint of orbit equivalence. As an ultimate result in this direction, Ornstein-Weiss' theorem \cite{ow} concludes that all essentially free, measure-preserving ergodic actions of infinite amenable groups are orbit equivalent. In particular, the class of groups measure equivalent to ${\Bbb Z}$ consists precisely of all infinite amenable groups since amenability is preserved under measure equivalence. (The result of Ornstein-Weiss was generalized to a statement in terms of discrete measured equivalence relations by Connes-Feldman-Weiss \cite{cfw}.) This conclusion gives a notable difference between measure equivalence and quasi-isometry since there exist infinitely many different quasi-isometry classes contained in the class of all infinite amenable groups. For example, ${\Bbb Z}^{n}$ and ${\Bbb Z}^{m}$ are not quasi-isometric for different $n$ and $m$ (e.g., see \cite[Chapitre 1, Th\'eor\`eme 17]{ghys-harpe}). In general, little relationships between the two equivalence relations among groups seem to be known.

In the case of non-amenable groups, Zimmer \cite{zim2} proved that if lattices in two connected semisimple Lie groups with finite center and no compact factors and of higher ${\Bbb R}$-rank (e.g., $SL(n, {\Bbb R})$ with $n\geq 3$) have essentially free, measure-preserving actions which are orbit equivalent, then the ambient Lie groups are locally isomorphic. In fact, he showed a stronger rigidity result, which is called Zimmer's cocycle superrigidity and can be seen as a generalization of the Mostow-Margulis rigidity theorem. More generally, Furman \cite{furman1} showed that if a lattice $\Gamma$ in a connected simple Lie group $G$ with finite center and higher ${\Bbb R}$-rank is measure equivalent to a discrete group $\Lambda$, then $\Lambda$ is almost isomorphic to some lattice in the Lie group $G$. This rigidity means that the class of all discrete groups measure equivalent to $\Gamma$ consists precisely of all groups almost isomorphic to lattices in $G$. On the other hand, notably, the quasi-isometry classes of irreducible lattices in semisimple Lie groups can be completely understood through combination of several papers by several authors. The reference \cite{farb} is a nicely written survey about this classification.

We shall recall another important result on measure equivalence. Gaboriau \cite{gab-l2} showed that $\ell^{2}$-{\it Betti numbers}\index{l2-Betti number@$\ell^{2}$-!Betti number} for discrete groups, which are introduced by Cheeger-Gromov \cite{cheeger-gromov}, are invariant under measure equivalence in the following sense: if two discrete groups $\Gamma_{1}$ and $\Gamma_{2}$ are measure equivalent, then there exists a positive constant $c$ such that $\beta_{n}(\Gamma_{1})=c\beta_{n}(\Gamma_{2})$ for all $n$, where $\beta_{n}(\Lambda)$\index{$\ b n \ C$@$\beta_{n}(\Gamma)$} denotes the $n$-th $\ell^{2}$-Betti number of a discrete group $\Lambda$. One consequence of this theorem is that if lattices in $SO(2n, 1)$ and $SO(2m, 1)$ are measure equivalent, then $n=m$.

Let us mention a certain analogy between the mapping class groups and discrete groups in the above examples. The group $SL(2, {\Bbb Z})$ can be understood from various points of view. For example, it is a lattice in $SL(2, {\Bbb R})$ and acts properly discontinuously and isometrically on the hyperbolic plane ${\Bbb H}^{2}$. The quotient group $PSL(2, {\Bbb Z})$ is a lattice in the orientation-preserving isometry group $PSO(2, 1)$ of ${\Bbb H}^{2}$. Besides these examples, $SL(2, {\Bbb Z})$ can be identified with the mapping class group of the torus. Hence, $SL(2, {\Bbb Z})$ is in the intersection of three natural classes of discrete groups, lattices in $SL(n, {\Bbb R})$ for $n\geq 2$, lattices in the orientation-preserving isometry group $PSO(n, 1)$ of the real hyperbolic space ${\Bbb H}^{n}$ for $n\geq 2$ and the mapping class groups for various surfaces. When we study one of the above three classes of groups, it is natural to find properties which all groups in the class possess in common, or conversely, to find sufficient conditions to distinguish groups for different $n$ or surfaces in the class. From another point of view, one should note that these classes have many analogous common properties. In fact, there are several works on the mapping class groups conducted by ideas coming from study in the other two classes of groups. It is widely believed that this approach of seeking analogy leads to fundamental discoveries of the structure of the mapping class groups. However, it goes without saying that there exist many differences among the three classes, which often give rise to difficulty and mystery and they attract us strongly (e.g., see \cite{FLM}, \cite{ivanov1} and \cite{ivanov2}).

As mentioned above, lattices in the groups in the two classes $\{ SL(n, {\Bbb R})\}_{n}$ and $\{ SO(n, 1)\}_{n}$ have been already studied and (still not completely) classified from the viewpoint of measure equivalence substantially. In this paper, we classify groups in the remaining class, that is, the mapping class groups up to measure equivalence in the following form.

Let $M=M_{g, p}$\index{$M = M g p $@$M=M_{g, p}$} be a connected compact orientable surface of {\it type} $(g, p)$\index{type!$(g, p)$}, that is, of genus $g$ and with $p$ boundary components. (Throughout the paper, we assume a surface to be connected, compact and orientable unless otherwise stated.) Let $\Gamma(M)$\index{$\ C M $@$\Gamma(M)$} be the mapping class group of $M$. We denote $\kappa(M)=3g+p-4$\index{$\ j M $@$\kappa(M)$}, which is called the {\it complexity}\index{complexity} of $M$, and denote
\begin{equation*}
g_{0}(M)= \begin{cases}
            2 & {\rm if} \ g\leq 2, \\
            g & {\rm if} \ g>2.
           \end{cases}
\end{equation*}
\index{$g 0 M$@$g_{0}(M)$}
One of our main theorems in this paper is the following:

\begin{thm}\label{thm-int-classification}
Suppose that $M^{1}$ and $M^{2}$ are compact orientable surfaces with $\kappa(M^{1}), \kappa(M^{2})\geq 0$ and that the mapping class groups $\Gamma(M^{1})$ and $\Gamma(M^{2})$ are measure equivalent. Then $\kappa(M^{1})=\kappa(M^{2})$ and $g_{0}(M^{1})=g_{0}(M^{2})$.
\end{thm}

If a surface $M$ satisfies $\kappa(M)<0$ and is not the torus, then $\Gamma(M)$ is finite. Note that all $\Gamma(M_{0, 4})$, $\Gamma(M_{1, 0})$, $\Gamma(M_{1, 1})$ and $SL(2, {\Bbb Z})$ are almost isomorphic (see the comment right after Theorem \ref{commuting-dehn}). Moreover, $\Gamma(M_{0, 6})$ and $\Gamma(M_{2, 0})$ are almost isomorphic (see the end of Section 6.8 in \cite{ivanov2}).

It is worth mentioning the notion of virtual cohomological dimension here, which is also an invariant under almost isomorphism, and comparing our result with it. Given a discrete group $\Gamma$ which has a torsion-free subgroup $\Gamma'$ of finite index, we define the {\it virtual cohomological dimension} ${\rm vcd}(\Gamma)$\index{virtual!cohomological dimension} of $\Gamma$ by the cohomological dimension of $\Gamma'$, which is independent of the choice of $\Gamma'$. Remark that two discrete groups with different virtual cohomological dimensions cannot be almost isomorphic. The virtual cohomological dimensions for the mapping class groups can be computed as follows (e.g., see \cite[Section 6.4]{ivanov2}): let $M$ be a compact orientable surface of type $(g, p)$.
\begin{enumerate}
\item[(i)] If $g\geq 2$ and $p=0$, then ${\rm vcd}(\Gamma(M))=4g-5$.
\item[(ii)] If $g\geq 2$ and $p\geq 1$, then ${\rm vcd}(\Gamma(M))=4g-4+p$.
\item[(iii)] If $g=0$, then ${\rm vcd}(\Gamma(M))=\max \{0, p-3\}$.
\item[(iv)] If $g=1$, then ${\rm vcd}(\Gamma(M))=\max \{ 1, p\}$.
\end{enumerate}
We shall give a few remarks on some relationships between the above fact and our classification result.
\begin{enumerate}
\item The mapping class groups of two surfaces of types $(g+1, 0)$ and $(g, 3)$ with $g\geq 2$ have the same complexity and the same virtual cohomological dimension, but are not measure equivalent.
\item The mapping class groups of two surfaces of types $(0, p)$ and $(1, p-3)$ with $p\geq 4$ (resp. of types $(1, 3)$ and $(2, 0)$) have the same complexity and the same virtual cohomological dimension.
\end{enumerate}
It is not known whether the two mapping class groups in the second remark are measure equivalent except for those of two surfaces of types $(0, 4)$ and $(1, 1)$.

The above mentioned fact that the complexity is a measure equivalence invariant can be proved by using the $\ell^{2}$-Betti number as well, though our proof is completely different. The $\ell^{2}$-Betti numbers of the mapping class groups can be calculated by the results due to Gromov \cite{gromov-kahler} and McMullen \cite{mcmullen} as follows (see Appendix \ref{app-cost}): if $M$ is a compact orientable surface of type $(g, p)$ with $\kappa(M)\geq 0$, then we see that
\[\beta_{3g-3+p}(\Gamma(M))>0, \ \ \beta_{n}(\Gamma(M))=0 {\rm \ for \ each \ } n\neq 3g-3+p. \]
(Moreover, the explicit value of $\beta_{3g-3+p}(\Gamma(M))$ can be obtained.) It follows from this result that if the mapping class groups $\Gamma(M^{1})$ and $\Gamma(M^{2})$ of two surfaces $M^{1}$ and $M^{2}$ with $\kappa(M^{1}), \kappa(M^{2})\geq 0$ are measure equivalent, then $\kappa(M^{1})=\kappa(M^{2})$.

Mosher-Whyte announced in \cite{mosher} that the mapping class group $\Gamma(M)$ of once-punctured surface $M$ has rigidity for quasi-isometry, that is, any finitely generated group quasi-isometric to $\Gamma(M)$ is almost isomorphic to $\Gamma(M)$. It is conjectured that this rigidity is satisfied for the mapping class groups of all surfaces $M$ with $\kappa(M)>0$. Remark that $SL(2, {\Bbb Z})$ has rigidity for quasi-isometry (see \cite[Chapitre 1, Th\'eor\`eme 21]{ghys-harpe}). From these announcement and conjecture, one might expect, or dream of, rigidity for measure equivalence of the mapping class groups, that is, any discrete group measure equivalent to $\Gamma(M)$ is almost isomorphic to $\Gamma(M)$.

In addition to the above classification, we give various types of discrete groups not measure equivalent to the mapping class groups. Let us denote
\[n(M)=g+\left[ \frac{g+p-2}{2}\right]\]\index{$n M$@$n(M)$}
for a compact orientable surface $M$ of type $(g, p)$, where $[a]$\index{$a for a R$@$[a]$ for $a\in {\Bbb R}$} denotes the maximal integer less than or equal to $a$ for $a\in {\Bbb R}$.

\begin{thm}\label{thm-main-main-free-groups}
Let $M$ be a compact orientable surface with $\kappa(M)\geq 0$. Suppose that an infinite subgroup $\Gamma$ of the mapping class group $\Gamma(M)$ is measure equivalent to a discrete group containing a subgroup isomorphic to the direct product of $n$ free groups of rank $2$. Then $n\leq n(M)$.
\end{thm}

In fact, $\Gamma(M)$ itself contains the direct product of $n(M)$ free groups of rank $2$ as a subgroup. Therefore, the inequality in Theorem \ref{thm-main-main-free-groups} is sharp. Geometrically, the number $n(M)$ is interpreted as follows (see Lemma \ref{main-geometric-lem}): given a family of non-peripheral, disjoint simple closed curves on $M$ in non-trivial, mutually distinct isotopy classes, we cut $M$ along all the curves in the family. The number $n(M)$ is equal to the maximal number of components not of type $(0, 3)$ of the cut surface among all such families of curves.

We recall a few results related to Theorem \ref{thm-main-main-free-groups}. Using the notion of geometric dimension of discrete measure-preserving equivalence relations, Gaboriau \cite{gab-l2} showed that if the direct product of $k$ free groups of rank more than $1$ is measure equivalent to a discrete group containing the direct product of $l$ free groups of rank more than $1$, then $k\geq l$. (We can deduce this fact as a corollary of Theorem \ref{thm-main-main-free-groups} in Corollary \ref{cor-gab-free-groups}.) Developing the theory of bounded cohomology, Monod-Shalom \cite{ms} gave the following result: let $\Gamma_{1},\ldots, \Gamma_{k}$ and $\Lambda_{1},\ldots, \Lambda_{l}$ be non-trivial finitely generated torsion-free groups with $k\leq l$. Suppose that all the $\Gamma_{i}$'s are non-elementary hyperbolic groups and that the direct products $\Gamma_{1}\times \cdots \times \Gamma_{k}$ and $\Lambda_{1}\times \cdots \times \Lambda_{l}$ are measure equivalent. Then $k=l$. More strongly, after permutation of the indices, $\Gamma_{i}$ is measure equivalent to $\Lambda_{i}$ for all $i$. Hjorth-Kechris \cite{hjorth-kechris} showed a similar result. Theorem \ref{thm-main-main-free-groups} can be regarded as an analogue of these results. However, the mapping class groups and the groups of the form $\Gamma_{1}\times \cdots \times \Gamma_{k}$ with $k\geq 2$ are different from the viewpoint of measure equivalence. This will be clarified in Theorem \ref{intro-indecomposability} (i) below.

In the course of the proof of the above theorems, the curve complex plays the most important role.

\begin{defn}
Let $M$ be a compact orientable surface with $\kappa(M)>0$. The {\it curve complex}\index{curve complex!for a non-exceptional surface} $C = C(M)$\index{$C = C M$@$C=C(M)$} of $M$ is the simplicial complex whose vertex set $V(C)$\index{$V C $@$V(C)$} is the set of non-trivial isotopy classes of non-peripheral simple closed curves on $M$, and a (finite) subset of $V(C)$ forms a simplex of $C$ if the set of curves representing it can be realized disjointly on $M$. (The curve complex can be defined also for surfaces $M$ with $\kappa(M)=0$ in a slightly different way. We recall its definition in Definition \ref{defn-curve-complex-exceptional}.)
\end{defn}

This complex was introduced by Harvey \cite{harvey}, who intended to investigate its combinatorial structure for the study of the asymptotic geometry of the 
Teichm\"uller space in analogy with Tits buildings for symmetric spaces. Note that this complex is locally infinite, which makes its geometry much more complicated. A remarkable progress in its study was made by Masur-Minsky \cite{masur-minsky1}, who studied the structure of the curve complex as a metric space and showed that it is connected, (Gromov-)hyperbolic and has infinite diameter. Moreover, they established a certain relationship between the Teichm\"uller space (with the Teichm\"uller metric) and the curve complex as metric spaces. These results give an analogy between the Teichm\"uller space and the hyperbolic spaces. In fact, the Teichm\"uller space has a lot of properties in common with the hyperbolic spaces other than the above result. Note that the mapping class group is (almost) equal to the isometry group of the Teichm\"uller space by a famous theorem of Royden \cite{royden} for closed surfaces and of Earle-Kra \cite{ek} for surfaces with non-empty boundary (see \cite{ivanov3}, \cite{ivanov4} for Ivanov's different proof of the theorem, in which the curve complex is used). Hence, it seems natural to seek for an analogy between the mapping class group and (lattices in) the isometry groups of the hyperbolic spaces, or more generally, (Gromov-)hyperbolic groups. We follow this strategy in the first half of the paper. For this purpose, it is helpful to know some results about an essentially free, measure-preserving action of a hyperbolic group on a standard probability space.

Following Zimmer's idea in \cite{zim3}, Adams \cite{adams2} proceeded to study essentially free Borel actions of hyperbolic groups. His method shows that any non-elementary hyperbolic group cannot be measure equivalent to the direct product $\Gamma_{1}\times \Gamma_{2}$, where both $\Gamma_{1}$ and $\Gamma_{2}$ are infinite and either $\Gamma_{1}$ or $\Gamma_{2}$ has an infinite amenable subgroup (see Theorem \ref{thm-hyp-main}). Inspired by Adams' method and using the hyperbolicity of the curve complex, we show the following result:

\begin{thm}\label{intro-indecomposability}
Let $M$ be a compact orientable surface with $\kappa(M)\geq 0$ and $\Gamma(M)$ be the mapping class group of $M$.
\begin{enumerate}
\item[(i)] If $\Gamma_{1}$ and $\Gamma_{2}$ are two infinite discrete groups and either $\Gamma_{1}$ or $\Gamma_{2}$ has an infinite amenable subgroup, then the direct product $\Gamma_{1}\times \Gamma_{2}$ and a sufficiently large subgroup of $\Gamma(M)$ are not measure equivalent.
\item[(ii)] If a discrete group $\Gamma$ has an infinite amenable normal subgroup, then $\Gamma$ and a sufficiently large subgroup of $\Gamma(M)$ are not measure equivalent.
\end{enumerate}
\end{thm}

We refer the reader to Theorem \ref{subgroup-classification} for the definition of sufficiently large subgroups of $\Gamma(M)$, whose class contains $\Gamma(M)$ itself.

Let us comment on Hamenst\"adt's recent work in \cite{ham2}. She constructs non-trivial bounded cohomology classes for various subgroups of the mapping class groups. Combining this with Monod-Shalom's result \cite{ms}, we can obtain a consequence containing Theorem \ref{intro-indecomposability} (see Remark \ref{rem-ms-ham}).

We can not apply Adams' method directly to our setting because the curve complex $C$ is locally infinite and its boundary $\partial C$ at infinity is not compact. To overcome these difficulties, we need to use stronger geometric properties of the curve complex than its hyperbolicity. Masur-Minsky \cite{masur-minsky2} introduced a certain type of geodesics, which are called {\it tight} ones\index{tight geodesic}, on the curve complex and showed that they have some finiteness properties. (This discovery plays an important role in the study of ends of hyperbolic 3-manifolds by Minsky and his collaborators, culminating in the proof of the ending lamination conjecture \cite{bcm}.) Moreover, Bowditch \cite{bowditch} refined their statement about the finiteness properties (see Theorem \ref{thm:bow}). Using Bowditch's formulation, we overcome the local infiniteness of the curve complex. On the other hand, to avoid difficulties arising from the non-compactness of the boundary $\partial C$, we use the Thurston boundary ${\cal PMF}$, which is a compact ideal boundary of the Teichm\"uller space on which the mapping class group acts continuously. It has a Borel subset ${\cal MIN}$, the space of minimal measured foliations, which is invariant under the action of the mapping class group. Klarreich \cite{kla} constructed a continuous map from ${\cal MIN}$ onto $\partial C$ equivariant for the actions of the mapping class group. This map plays a crucial role in our modification of Adams' method to our setting.

One of Adams' important observations is amenability in a topological sense of the action of a hyperbolic group on its boundary at infinity \cite{adams1}, which goes back to Spatzier and Zimmer's results \cite{spa}, \cite{spa-zim} in a special case. Since we follow Adams' method, the following similar result need to be shown:

\begin{thm}\label{intro-amenable-action}
Let $M$ be a compact orientable surface with $\kappa(M)\geq 0$ and $\Gamma(M)$ be the mapping class group of $M$. Let $\partial C$ and ${\cal PMF}$ denote the boundary at infinity of the curve complex and the Thurston boundary of $M$, respectively.
\begin{enumerate}
\item[(i)] The Borel space $\partial C$ is standard.
\item[(ii)] The action of $\Gamma(M)$ on $(\partial C, \mu)$ is amenable in a measurable sense for any quasi-invariant probability measure $\mu$ on $\partial C$.
\item[(iii)] The action of $\Gamma(M)$ on $({\cal PMF}, \mu)$ is amenable in a measurable sense for any quasi-invariant probability measure $\mu$ on ${\cal PMF}$ with $\mu({\cal MIN})=1$.
\end{enumerate}
\end{thm}

We refer the reader to Appendix \ref{general-amenable} for the definition of amenability of an action of a discrete group. Note that there exist several natural measures on ${\cal PMF}$. Masur \cite{masur} and Rees \cite{rees} constructed a natural measure on ${\cal PMF}$, for which the action of $\Gamma(M)$ on ${\cal PMF}$ is ergodic. Moreover, the measure satisfies the hypothesis in the assertion (iii) of Theorem \ref{intro-amenable-action}. Kaimanovich-Masur \cite{kai-mas} verified that ${\cal PMF}$ can be identified with the Poisson boundary of random walks on the mapping class group and that the corresponding harmonic measures satisfy the same hypothesis. Remark that for such a harmonic measure $\mu$, we can deduce amenability of the action of $\Gamma(M)$ on $({\cal PMF}, \mu)$ in a measurable sense by Zimmer's work \cite{zim1}.   

It is significant to investigate amenability of the action of a given discrete group on a (compact) space not only in the study of measure equivalence as above. Let us explain another significance of it. In general, given a continuous action of a discrete group $G$ on a compact Hausdorff space $X$, we know the following fact \cite{ana} (see also Theorem \ref{amenable-and-fixed}): the action is amenable in a topological sense if and only if the action of $G$ on $(X, \mu)$ is amenable in a measurable sense for any quasi-invariant measure $\mu$ on $X$. If a discrete group $G$ admits such an amenable action for some compact Hausdorff space $X$, then we say that $G$ is {\it exact}\index{exact}. If $G$ is finitely generated, then this property is described as some geometric condition on the Cayley graph of $G$, which is called {\it property A}\index{property A} \cite{higson-roe}. Property A for metric spaces was introduced by Yu \cite{yu} in the work on the Baum-Connes conjecture. Thanks to his result, all finitely generated groups with property A satisfy the Novikov conjecture. Notably, the class of exact groups is believed to be huge since exactness is closed under various operations of groups except for taking quotients and is satisfied for all amenable groups and free groups (see \cite{dykema}, \cite{kirch-wass} and \cite{tu}). On the other hand, exactness has the following characterization within the framework of operator algebras (in fact, the notion of exactness originally comes from $C^{*}$-algebras): a discrete group $G$ is exact if and only if the reduced group $C^{*}$-algebra of $G$ is exact \cite{ozawa1}. In the theory of operator algebras, the exactness of $C^{*}$-algebras is one of the important notions and has various applications \cite{kirch}, \cite{rordam}. Hence, it is an interesting challenge to prove that a given discrete group is exact.

It follows from Adams' result that all hyperbolic groups are exact. On the other hand, although we can obtain Theorem \ref{intro-amenable-action} (ii) and (iii) related to amenability of the actions of the mapping class group, we see immediately that $\partial C$ is not compact (see Proposition \ref{prop-boundary-non-compact}) and that its action on ${\cal PMF}$ is not amenable in a topological sense since there exist non-amenable stabilizers of the action. However, in the course of the proof of Theorem \ref{intro-amenable-action}, we obtain the following result (see also Remark \ref{rem-mcg-exact}). In the proof of the assertion (ii) in the following theorem, we use a result of Ozawa \cite{ozawa5} (see Appendix \ref{app-exact}).

\begin{thm}\label{intro-property-A}
Let $M$ be a compact orientable surface with $\kappa(M)\geq 0$.
\begin{enumerate}
\item[(i)] The curve complex of $M$ has property A as a metric space.
\item[(ii)] The mapping class group of $M$ is exact.
\end{enumerate}
\end{thm}

Although we have stated many kinds of analogous properties of the mapping class groups to those of hyperbolic groups so far, we can obtain the following difference between them:

\begin{thm}\label{intro-not-me-hyp}
The mapping class group of a compact orientable surface $M$ with $\kappa(M)>0$ is not measure equivalent to any hyperbolic group.
\end{thm}

This theorem can be shown by using Adams' results in \cite{adams2} and investigating a certain group-theoretic difference between the mapping class groups and hyperbolic groups.

The paper is organized as follows. We begin with the proof of Theorem \ref{intro-property-A} (i) in Chapter \ref{chapter-property-A}, where we discuss the definitions of the curve complex, tight geodesics and property A for metric spaces. We give the proof of Theorem \ref{intro-property-A} (ii) in Appendix \ref{app-exact}. In Chapter \ref{chapter:amenable-action}, we review some important facts about the Thurston boundary, the boundary of the curve complex and the mapping class group, which will be needed in the subsequent chapters. Theorem \ref{intro-amenable-action} is verified in this chapter. In Chapter \ref{chapter-indec}, we start to study discrete measured equivalence relations generated by the mapping class group and prove Theorems \ref{intro-indecomposability} and \ref{intro-not-me-hyp}. Here, we recall fundamentals of discrete measured equivalence relations. In particular, the notion of normal subrelations introduced by Feldman-Sutherland-Zimmer \cite{fsz} plays an important role in the proof of the theorems and Theorem \ref{thm-int-classification}. In this chapter, two types of subrelations of relations generated by the mapping class group will be introduced, which are called irreducible and amenable ones and reducible ones. Moreover, we introduce the canonical reduction system for reducible subrelations, which is an analogue of that for subgroups of the mapping class group introduced by Birman-Lubotzky-McCarthy \cite{blm} and Ivanov \cite{ivanov1}. In Chapter \ref{chapter-best}, we analyze reducible subrelations further and give a partial classification result of Theorem \ref{thm-int-classification} and the proof of Theorem \ref{thm-main-main-free-groups}. In Chapter \ref{chapter-best2}, we prove a complete classification result of Theorem \ref{thm-int-classification}. In the beginnings of Chapters \ref{chapter-indec}, \ref{chapter-best} and \ref{chapter-best2}, we give outlines of the proofs of the theorems in the chapters. 

Most of the notation we use will be standard and will be explained as needed. In particular, throughout the paper, ${\Bbb N}$ denotes the set $\{ 0, 1, 2,\ldots \}$ of natural numbers\index{$N = 0 1 2 $@${\Bbb N}=\{ 0, 1, 2,\ldots \}$}. Given a set $A$, we denote by $|A|$ \index{$A for a set A$@$"|A"|$ for a set $A$} its cardinality.

After we finished to write the first draft of this paper, we found the paper \cite{ham3} related to Theorem \ref{intro-property-A} (ii).

\vspace{1em}

\noindent{\it Acknowledgements.} The author wishes to thank everyone who gave his precious opinions for this paper and Professor Narutaka Ozawa for showing him Proposition \ref{prop-ozawa}. Especially, the author would like to express his deep gratitude to his advisor, Professor Masaki Izumi. He read carefully the whole manuscript of the paper and gave many valuable suggestions throughout the paper. The author could not complete the paper without his continuous and warm encouragement.


\chapter{Property A for the curve complex}\label{chapter-property-A}

In this chapter, we show that the curve complex of a compact orientable surface $M$ satisfies property A. Property A was introduced by Yu in his work on the Baum-Connes conjecture \cite{yu}, and has many equivalent, geometric and operator algebraic, conditions \cite{higson-roe}, \cite{ozawa1}, \cite{ozawa3}. It is a challenging problem to decide whether a given discrete metric space has property A or not. 

For the proof, we use an excellent method treating the curve complex, developed by Masur-Minsky \cite{masur-minsky1}, \cite{masur-minsky2} and refined by Bowditch \cite{bowditch}. In Section \ref{sec:geom}, we prepare the geometric properties of the curve complex which are necessary for the proof. In Section \ref{sec:pro}, we introduce property A for a discrete metric space and recall some background of the property, and in the following sections, we prove that the curve complex has property A.

\section{Geometry of the curve complex}\label{sec:geom}

\subsection{Hyperbolicity of the curve complex}

Let $M$ be a compact orientable surface of {\it type} $(g, p)$\index{type!$(g, p)$} (i.e., of genus $g$ and with $p$ boundary components). The {\it complexity}\index{complexity} of $M$ is defined by $\kappa(M) = 3g+p-4$\index{$\ j M $@$\kappa(M)$}. Let $\Gamma(M)$\index{$\ C M $@$\Gamma(M)$} be the {\it mapping class group}\index{mapping class group} of $M$, that is, the group of orientation-preserving self-diffeomorphisms of $M$ up to isotopy. 

\begin{defn}\label{defn-curve-complex-non-exceptional}
Let $M$ be a surface with $\kappa(M)>0$. The {\it curve complex}\index{curve complex!for a non-exceptional surface} $C = C(M)$\index{$C = C M$@$C=C(M)$} of $M$ is the simplicial complex whose vertex set $V(C)$\index{$V C $@$V(C)$} is the set of non-trivial isotopy classes of non-peripheral simple closed curves on $M$, and a (finite) subset of $V(C)$ forms a simplex of $C$ if the set of curves representing it can be realized disjointly on $M$. 
\end{defn}

This complex was introduced by Harvey \cite{harvey}. The mapping class group $\Gamma(M)$ acts on $C(M)$ simplicially and thus, isometrically when $C(M)$ is equipped with the natural combinatorial metric.

If a surface $M$ satisfies $\kappa(M)>0$, then it is said to be {\it non-exceptional}\index{non-exceptional surface}, otherwise it is said to be {\it exceptional}\index{exceptional surface}. The definition of the curve complex of an exceptional surface is slightly different, therefore it should be investigated independently (see \cite{minsky} and Section \ref{section-exceptional}).

Let $M$ be a non-exceptional surface. Then the dimension of $C(M)$ is equal to $\kappa(M)$ and $C(M)$ is locally infinite because there exist infinitely many curves on $M$ not intersecting $\alpha$ for any curve $\alpha$ on $M$. Sometimes this property makes the geometry of $C(M)$ difficult. However, Masur-Minsky \cite{masur-minsky1} have shown the following significant theorem: 

\begin{thm}[\ci{Theorem 1.1}{masur-minsky1}]\label{thm:mm1}

Let $M$ be a surface with $\kappa(M)>0$. Then the curve complex $C(M)$ is a connected Gromov-hyperbolic metric space and has infinite diameter.
\end{thm}

We shall recall the definition of Gromov-hyperbolic metric spaces.

For a metric space $X$, we say that $X$ is {\it geodesic}\index{geodesic!metric space} if any two points $x, y\in X$ can be joined by a geodesic $l$, that is, there exists a map
\[l\colon [0, d(x, y)]\rightarrow X\]
such that $l(0)=x$, $l(d(x, y))=y$ and for any $s, t\in [0, d(x, y)]$, we have $d(l(s), l(t))=|s-t|$. A map $l$ from the interval $[0, \infty)$ (resp. $(-\infty, \infty)$) to $X$ is called a {\it geodesic ray}\index{geodesic!ray} (resp. {\it bi-infinite geodesic}\index{bi-infinite geodesic}) if $d(l(s), l(t))=|s-t|$ for any $s$, $t$ in the domain of $l$. If $l\colon [0, \infty)\rightarrow X$ is a geodesic ray, then we call $l(0)$ the {\it origin}\index{origin} of $l$. Sometimes we identify a geodesic $l\colon I\rightarrow X$, a map from some interval $I$ into $X$, with its image in $X$, and also write $l$ as the image in $X$. When $X$ is a graph with vertex set $V(X)$ and equipped with the path metric, we often identify a geodesic whose endpoints are in $V(X)$ with the sequence of vertices through which $l$ passes. A {\it geodesic triangle}\index{geodesic!triangle} on $X$ consists of three points $x_{1}, x_{2}, x_{3}\in X$ and three geodesics between $x_{i}$ and $x_{j}$ for $i\neq j$. The geodesics are called {\it edges}\index{edge} of the triangle.  

For a geodesic metric space $X$, we denote
\[(x|y)_{z}=(1/2)(d(x, z)+d(y, z)-d(x, y))\]
\index{$x y z$@$(x"|y)_{z}$}for $x, y, z\in X$, which is called the {\it Gromov-product}\index{Gromov-!product} of two points $x$, $y$ with a base point $z$. We say that $X$ is {\it (Gromov-)hyperbolic}\index{Gromov-!hyperbolic space}\index{hyperbolic!space} if there exists a constant $\delta \geq 0$ satisfying 
\[(x|z)_{w}\geq \min \{ (x|y)_{w}, (y|z)_{w}\} -\delta \]
for any $w, x, y, z\in X$. In this case, there exists $\delta'\geq 0$ depending only on $\delta$ such that for any geodesic triangle $T$ on $X$, any edge of $T$ is contained in the $\delta'$-neighborhood of the other two edges of $T$ (such a triangle is said to be $\delta'$-{\it slim}\index{delta-slim@$\delta$-!slim}). Conversely, if there exists a constant $\delta \geq 0$ such that any geodesic triangle on $X$ is $\delta$-slim, then $X$ is hyperbolic \cite[Chapitre 2, Proposition 21]{ghys-harpe}. For a geodesic metric space $X$ and a constant $\delta \geq 0$, we say that $X$ is $\delta$-{\it hyperbolic}\index{delta-hyperbolic@$\delta$-!hyperbolic} if $\delta$ satisfies the above inequality for any $w, x, y, z\in X$ and every geodesic triangle on $X$ is $\delta$-slim. 

The classical hyperbolic spaces, a simplicial tree and the universal cover (or the Cayley graph of the fundamental group) of a closed negatively curved manifold are important examples of hyperbolic spaces. Note that every geodesic space with finite diameter is hyperbolic. For the general theory of hyperbolic spaces, we recommend the reader to consult \cite{bridson-haefliger}, \cite{ghys-harpe} and \cite{gromov1}.

\subsection{Finiteness properties of tight geodesics}\label{sec:finite}

Let $M$ be a non-exceptional surface. In the rest of this section, we identify $C=C(M)$ with its 1-skeleton, which is quasi-isometric to the curve complex. Let $V(C)$ be the vertex set.

Masur-Minsky showed in \cite{masur-minsky2} that the set of geodesics of some type on $C(M)$, called {\it tight} geodesics, has finiteness properties in some sense. For describing them, we need to give the definition of tight geodesics.

\begin{defn}
\begin{enumerate}
\item[(i)] A {\it multicurve}\index{multicurve} is a subset of $V(C)$ consisting of pairwise disjoint elements (i.e., forming a simplex of the curve complex).

\item[(ii)] A {\it multigeodesic}\index{multigeodesic} is a sequence $\{ A_{i}\}_{i=0}^{n}$ of multicurves such that $d(\alpha, \beta)=|i-j|$ for all $i\neq j$ and all $\alpha \in A_{i}$ and $\beta \in A_{j}$. 

\item[(iii)] We say that a multigeodesic $\{ A_{i}\}_{i=0}^{n}$ is {\it tight} at index $i\neq 0, n$ if for all $\alpha \in A_{i}$, each curve intersecting $\alpha$ also intersects some element of $A_{i-1}\cup A_{i+1}$. (This definition is due to Bowditch \cite{bowditch} and weaker than that given in \cite{masur-minsky2}.)

\item[(iv)] We say that a multigeodesic $\{ A_{i}\}_{i=0}^{n}$ is {\it tight} if it is tight at all indices in $\{ 1,\ldots, n-1\}$.

\item[(v)] A geodesic $\{ \alpha_{i}\}_{i=0}^{n}$ consisting of vertices is said to be {\it tight}\index{tight geodesic} if there exists some tight multigeodesic $\{ A_{i}\}_{i=0}^{n}$ such that $\alpha_{i} \in A_{i}$ for all $i$. 
\end{enumerate}
\end{defn}

By definition, the set of tight geodesics is $\Gamma(M)$-invariant. For describing the finiteness properties of tight geodesics, we need some preparations.

Let $X$ be a (not necessarily locally finite) $\delta$-hyperbolic graph for $\delta \geq 0$ with vertex set $V(X)$.

\begin{defn}
Suppose that $f$ is a geodesic segment on $X$ identified with the sequence $x_{0}x_{1}\ldots x_{n}$ of vertices. For $r\in \Bbb{N}$, we call the geodesic segment $x_{r}x_{r+1}\ldots x_{n-r}$ the $r$-{\it central segment}\index{r-central@$r$-!central segment} of $f$. 

If $g=y_{0}y_{1}\ldots y_{n}$ is another geodesic segment of the same length, then we say that $f$ and $g$ are $r$-{\it close}\index{r-close@$r$-!close} when $d(x_{i}, y_{i})\leq r$ for all $i$. For $x, y\in V(X)$, let ${\mathcal L}(x, y)$\index{$L x y $@${\cal L}(x, y)$} denote the set of all geodesics on $X$ from $x$ to $y$.
\end{defn}

\begin{lem}[\ci{Lemma 2.1}{bowditch}]\label{lem:bow}
There exists a constant $\delta_{0}\geq 0$ depending only on $\delta$ such that if $f\in {\mathcal L}(x, y)$ and $g\in {\mathcal L}(x^{\prime}, y^{\prime})$, then the $r$-central segment of $f$ is $\delta_{0}$-close to some geodesic segment of $g$, where $r=\max \{ d(x, x^{\prime}), d(y, y^{\prime})\}$. 
\end{lem}

Let $C$ be the 1-skeleton of the curve complex of a non-exceptional surface $M$. Thanks to Theorem \ref{thm:mm1}, the graph $C$ is $\delta$-hyperbolic for some $\delta \geq 0$. For $\alpha$, $\beta \in V(C)$, let ${\mathcal L}_{T}(\alpha, \beta)\subseteq {\mathcal L}(\alpha, \beta)$\index{$L T x y $@${\cal L}_{T}(x, y)$} be the subset of tight geodesics. Denote

\begin{align*}
G(\alpha, \beta)&=\{ \gamma \in V(C): \gamma \in f, \ f\in {\mathcal L}_{T}(\alpha, \beta)\},\\ 
{\mathcal L}_{T}(\alpha, \beta ; r)&=\bigcup\{ {\mathcal L}_{T}(\alpha', \beta'): \alpha'\in B(\alpha ;r), \ \beta'\in B(\beta ;r)\},\\
G(\alpha, \beta; r)&=\bigcup\{G(\alpha^{\prime}, \beta^{\prime}): \alpha^{\prime}\in B(\alpha; r), \ \beta^{\prime}\in B(\beta; r)\}
\end{align*} 
\index{$G x y $@$G(x, y)$}\index{$L T x y r$@${\cal L}_{T}(x, y;r)$}\index{$G x y r$@$G(x, y; r)$}for $r\geq 0$, where $B(\gamma ;r)$\index{$B x r$@$B(x ;r)$} is (the intersection of the vertex set $V(C)$ and) the closed ball with center $\gamma$ and radius $r$.

In the next theorem, the assertion (i) and the finiteness in (ii) were first proved by Masur-Minsky and the existence of constants appearing below was shown by Bowditch:
  
\begin{thm}[\ci{Corollary 6.14}{masur-minsky2}, \ci{Theorems 1.1, 1.2}{bowditch}]\label{thm:bow}
In the above notation,
\begin{enumerate}
\renewcommand{\labelenumi}{\rm(\roman{enumi})}
\item the set ${\mathcal L}_{T}(\alpha, \beta)$ is non-empty for any $\alpha, \beta \in V(C)$.

\item there exists a constant $K_{0}$ depending only on $\kappa(M)$ such that if $\alpha, \beta \in V(C)$ and $\gamma \in G(\alpha, \beta)$, then the set $G(\alpha, \beta)\cap B(\gamma; \delta_{0})$ has at most $K_{0}$ elements.

\item there exist constants $\delta_{1}>\delta_{0}$ and $K_{1}$ depending only on $\kappa(M)$ such that if $\alpha, \beta \in V(C)$, $r\in {\Bbb N}$ and $\gamma \in G(\alpha, \beta)$ with $d(\gamma, \{ \alpha, \beta \} )\geq r+\delta_{1}$, then the set $G(\alpha, \beta; r)\cap B(\gamma; \delta_{0})$ has at most $K_{1}$ elements.

\end{enumerate} 
\end{thm}
Remark that in the assertion (ii), any geodesic from $\alpha$ to $\beta$ passes through $B(\gamma; \delta_{0})$ by Lemma \ref{lem:bow}, and similarly in the assertion (iii), any geodesic segment with endpoints in $B(\alpha; r)$ and $B(\beta; r)$ respectively passes through $B(\gamma; \delta_{0})$. Thus, the quantity of $\delta_{0}$ is not important as far as it satisfies Lemma \ref{lem:bow} because if $\delta_{0}$ increases, then the universal constants $\delta_{1}$, $K_{0}$, $K_{1}$ in the above theorem increase simply. 

These properties are the key ingredients for the proof of property A for the curve complex. In Section \ref{sec:Acc}, we will consider a sufficient condition for property A in a slightly more general setting, that is, we show that the vertex set of a hyperbolic graph satisfying the similar conditions in the above theorem has property A. It will be formulated explicitly in Section \ref{sec:Acc}.


\section{Generalities for property A}\label{sec:pro}

In \cite{yu}, Yu introduced property A for any metric space $X$. However, in this chapter, we adopt, as our definition of property A, a condition considered by Higson-Roe, which is equivalent to Yu's definition for discrete metric spaces with bounded geometry \cite[Lemma 3.5]{higson-roe}. Here, a metric space $X$ is said to have {\it bounded geometry}\index{bounded geometry} if for every $R>0$, there exists a natural number $N=N(R)$ such that the cardinality of $B(x; R)$ is less than $N$ for all $x\in X$. For example, the Cayley graph of a finitely generated group has bounded geometry. 

Given a countable set $Z$, we denote by ${\rm Prob}(Z)$\index{$P r ob Z$@${\rm Prob}(Z)$} the space of probability measures on $Z$, i.e., $\{\phi \colon Z\rightarrow [0, 1]: \Vert \phi \Vert_{1}=1\}$, where $\Vert \cdot \Vert_{1}$ is the $\ell^{1}$-norm.

\begin{defn}\label{defn-higson-roe} 
Let $X$ be a discrete metric space. Then we say that $X$ has {\it property A}\index{property A} if there exists a sequence of maps
\[a^{n}\colon X\rightarrow {\rm Prob}(X)\]
satisfying the following two conditions:
\begin{enumerate}
\renewcommand{\labelenumi}{\rm(\roman{enumi})}
\item for every $n\in {\Bbb N}$, there exists a constant $R=R(n)>0$ such that 
\[{\rm supp}(a^{n}_{x})\subseteq B(x; R)\]
for all $x\in X$, where ${\rm supp}(a)$ denotes the support of $a\in {\rm Prob}(X)$; 
\item for every $K>0$, we have
\[\lim_{n\rightarrow \infty}\sup_{d(x, y)<K}\Vert a^{n}_{x}-a^{n}_{y}\Vert_{1}=0.\]
\end{enumerate}

For a general (non-discrete) metric space $Z$, it is said to have {\it property A} if there exists a discrete metric subspace $X$ of $Z$ such that
\begin{enumerate}
\renewcommand{\labelenumi}{\rm(\roman{enumi})}
\item there exists a constant $c>0$ for which $d(z, X)\leq c$ for all $z\in Z$;
\item $X$ has property A.
\end{enumerate}  
\end{defn}

Whether a metric space has property A depends only on the coarse equivalence class of the metric (see \cite{roe} for a discussion of coarse geometry). Thus, for a finitely generated group, if we equip it with the word-length metric associated to some finite system of generators, then whether the metric space has property A does not depend on the choice of the generating system.

One of the important examples of spaces with property A is the class of discrete metric spaces with bounded geometry which are Gromov-hyperbolic. Germain \cite[Appendix B]{ana} showed this fact for the Cayley graph of a hyperbolic group (see also Example \ref{ex-adams-exact}) and Tu \cite[Proposition 8.1]{tu} generalized the result to any hyperbolic space with bounded geometry, using Germain's proof.

For finitely generated groups, property A has many equivalent conditions. We collect them in the following theorem. For the amenability of a group action, the reader is referred to \cite{ana} or Appendix \ref{general-amenable}. For the exactness\index{exact} of a discrete group and a $C^{*}$-algebra, see \cite{guentner-kaminker} and \cite{was}. 

\begin{thm}[\cite{higson-roe}, \cite{ozawa1}]\label{exact-property-A}
Let $\Gamma$ be a finitely generated group. The following conditions are equivalent:
\begin{enumerate}
\renewcommand{\labelenumi}{\rm(\roman{enumi})}
\item The group $\Gamma$ has property A;
\item The action of $\Gamma$ on its Stone-\v Cech compactification is  amenable in a topological sense, or equivalently, there exists a continuous action of $\Gamma$ on some compact Hausdorff space which is topologically amenable;
\item The reduced group $C^{*}$-algebra $C^{*}_{r}(\Gamma)$ is exact.
\end{enumerate}
\end{thm}

\begin{ex}\label{ex-adams-exact} 
A finitely generated group $\Gamma$ is said to be {\it (Gromov-)hyperbolic}\index{Gromov-!hyperbolic group}\index{hyperbolic!group} if the Cayley graph of $\Gamma$ is a hyperbolic metric space. This definition is independent of the choice of a generating system. Adams \cite{adams1} showed that the action of a hyperbolic group on its boundary (which is a compact Hausdorff space) is topologically amenable. Germain \cite[Appendix B]{ana} gave a simpler proof of the fact. We recommend the reader to consult \cite{bridson-haefliger}, \cite{ghys-harpe} and \cite{gromov1} for fundamentals of hyperbolic groups and their boundaries. 
\end{ex}

\begin{rem}\label{rem-mcg-exact}
In the case of the mapping class group, there exists a compact Hausdorff space on which it acts naturally, called the Thurston boundary. It is an ideal boundary of the Teichm\"uller space and is also identified with the Poisson boundary for random walks on the mapping class group \cite{kai-mas}. However, the action of the mapping class group on the Thurston boundary is not topologically amenable since there exist non-amenable stabilizers of the action (see Chapter \ref{chapter:amenable-action}, Section \ref{mcg}). We show exactness of the mapping class group in Appendix \ref{app-exact}, constructing another compact space on which the mapping class group acts amenably in a topological sense.

On the other hand, thanks to \cite{bb} and \cite{kor}, the mapping class groups $\Gamma(M_{2, 0})$ and $\Gamma(M_{0, p})$ for each $p$ are linear, where $M_{g, p}$ denotes a compact orientable surface of type $(g, p)$. Thus, they are exact by \cite{ghw}. (In \cite{bb} and \cite{kor}, it is shown that the mapping class groups ${\cal M}_{2, 0}$ and ${\cal M}_{0, p}$ for each $p$ are linear, where ${\cal M}_{g, p}$ denotes the group of isotopy classes of orientation-preserving self-diffeomorphisms of the closed surface $M_{g, 0}$ which preserve the set of $p$ marked points on $M_{g, 0}$. The two groups $\Gamma(M_{g, p})$ and ${\cal M}_{g, p}$ are naturally isomorphic \cite[Section 5.1]{ivanov2}.)    

Recently, Bell-Fujiwara \cite{bell-fuji} show that the asymptotic dimension of the mapping class group for a surface of genus less than or equal to $2$ is finite. This implies that the mapping class group of such a surface is exact by \cite[Lemma 4.3]{higson-roe}.  
\end{rem}


\section{Property A for the curve complex}\label{sec:Acc}

In this section, we prove that the curve complex $C=C(M)$ of a non-exceptional surface $M$ has property A. In the first subsection, we show property A in Definition \ref{defn-higson-roe} and in the second subsection, we recall Yu's definition of property A and show that $C$ satisfies also this property A.

\subsection{Higson-Roe's property A}

We shall recall the {\it Gromov-boundary} (or {\it boundary at infinity})\index{Gromov-!boundary}\index{boundary at infinity} $\partial X$ for a (not necessarily proper) hyperbolic space $X$. We say that a sequence $\{ x_{n}\}$ in $X$ {\it converges to infinity}\index{converge at infinity} if 
\[\liminf_{n, m\rightarrow \infty}(x_{n}|x_{m})_{e}=\infty,\]
where $e\in X$ is some base point. This definition is independent of the choice of $e$. Two sequences $\{ x_{n}\}$ and $\{ y_{m}\}$ converging to infinity are {\it equivalent}\index{equivalent} if 
\[\liminf_{n, m\rightarrow \infty}(x_{n}|y_{m})_{e}=\infty.\]
In fact, this defines an equivalence relation on the set of all sequences in $X$ converging to infinity. We define the boundary $\partial X$\index{$\ z X$@$\partial X$} of $X$ as the set of equivalence classes of sequences which converge to infinity. The reader should be referred to \cite{bridson-haefliger} and \cite{ghys-harpe} for more details. Remark that $\partial X$ may be empty even if $X$ has infinite diameter (see \cite[Chapitre 7, Exercise 7]{ghys-harpe}). We study the boundary of the curve complex in Chapter \ref{chapter:amenable-action}, Section \ref{boundary-of-curve-complex}, using the description of it by Klarreich \cite{kla}.

\begin{lem}\label{Gromov-product}
Let $X$ be a geodesic $\delta$-hyperbolic space. Then for each geodesic triangle $[x, y]\cup [y, z]\cup [z, x]$ in $X$, we have
\[(y|z)_{x}\leq d(x, [y, z])\leq (y|z)_{x}+4\delta.\]
\end{lem}

\begin{pf}
This follows from \cite[Chapitre 2, Lemme 17, Proposition 21]{ghys-harpe}. 
\end{pf}

\begin{lem}\label{hyp-fund}
Let $X$ be a $\delta$-hyperbolic graph. For any two geodesic rays $f, g\colon {\Bbb N}\rightarrow V(X)$ with origins in $X$, the following four conditions are equivalent:
\begin{enumerate}
\renewcommand{\labelenumi}{\rm(\roman{enumi})}
\item The rays $f$ and $g$ define the same point in $\partial X$, that is, 
\[\liminf_{i, j\rightarrow \infty}(f(i)|g(j))=\infty,\]
where $(x|y)=(x|y)_{e}$ for some fixed point $e\in X$;  
\item There exists $K>0$ such that $\inf \{d(f(i), g(j)): j\in {\Bbb N}\}\leq K$ for each $i\in {\Bbb N}$;
\item $\sup_{i\in {\Bbb N}}d(f(i), g(i))<\infty$;
\item There exist $l\in {\Bbb Z}$ and $i_{0}\geq \max \{ 0, l\}$ such that 
\[d(f(i), g(i-l))\leq 16\delta\]
for any $i\geq i_{0}$. 
\end{enumerate}
\end{lem} 

\begin{pf}
We show that the condition (iii) implies the condition (i). For $i\geq j$, we have
\begin{align*}
2(f(i)|g(j))_{e}&=d(f(i), e)+d(g(j), e)-d(f(i), g(j))\\
                &\geq d(f(i), f(0))-d(f(0), e)+d(g(j), g(0))-d(g(0), e)\\
                &\ \ \ \ \ \ \ \ \ \ \ \ \ -d(f(i), f(j))-d(f(j), g(j))\\
                &\geq i+j-d(f(0), e)-d(g(0), e)-(i-j)-d(f(j), g(j))\\
                &\geq 2j-d(f(0), e)-d(g(0), e)-d(f(j), g(j)). 
\end{align*}
This shows the assertion.

The equivalence of the three conditions (ii), (iii) and (iv) is shown in \cite[Chapitre 7, Proposition 2]{ghys-harpe}.

We show that the condition (i) implies the condition (ii).  First, we prove that there exists a geodesic ray $f'$ such that $f(0)=f'(0)$ and $M=\sup_{i\in {\Bbb N}}d(f'(i), g(i))<\infty$. 

For each $n\in {\Bbb N}$, we have $n-d(f(0), g(n))\leq d(f(0), g(0))$ by the triangle inequality for the triangle $f(0)g(0)g(n)$. Let $n_{0}\in {\Bbb N}$ be an integer such that 
\[n_{0}-d(f(0), g(n_{0}))=\max_{n\in {\Bbb N}}(n-d(f(0), g(n))).\]
Choose a geodesic $f''$ from $f(0)$ to $g(n_{0})$ and define $f'\colon {\Bbb N}\rightarrow V(X)$ by connecting $f''$ with the geodesic ray $g|_{[n_{0}, \infty)}$. We prove that $f'$ is a geodesic ray. 

It suffices to show that 
\[d(f''(l), g(n_{0}+k))=L-l+k\]
for any $l\in [0, L]\cap {\Bbb N}$ and $k\in {\Bbb N}$, where $L=d(f(0), g(n_{0}))$. It is clear that the left hand side is not bigger than the right hand side. By the maximal property of $n_{0}$, we see that
\[n_{0}+k-d(f(0), g(n_{0}+k))\leq n_{0}-d(f(0), g(n_{0}))=n_{0}-L.\]
Hence,
\begin{align*}
d(f''(l), g(n_{0}+k)) & \geq d(f(0), g(n_{0}+k))-d(f(0), f''(l))\\
                      & \geq L+k-l,
\end{align*}
which shows the claim.

Suppose that the condition (ii) is not satisfied for $f$ and $g$. Then there would exist $n\in {\Bbb N}$ such that if $x=f(n)$, then $B(x; M+\delta)\cap g=\emptyset$ (recall that $M=\sup_{i\in {\Bbb N}}d(f'(i), g(i))<\infty$). For each $k\in {\Bbb N}$, choose a geodesic segment $[f(n+k), f'(n+k)]$ from $f(n+k)$ to $f'(n+k)$. By $\delta$-slimness of the triangle $f(0)f(n+k)f'(n+k)$, the set $B(x; \delta)\cap ([f(0), f'(n+k)]\cup [f'(n+k), f(n+k)])$ is non-empty. If $B(x;\delta)\cap [f(0), f'(n+k)]$ were non-empty, then $B(x; M+\delta)\cap g$ would be non-empty. It is a contradiction. Thus, $B(x;\delta)\cap [f'(n+k), f(n+k)]$ must be non-empty. By Lemma \ref{Gromov-product}, we have
\begin{align*}
(f(n+k)|f'(n+k))_{f(0)}& \leq d(f(0), [f(n+k), f'(n+k)])\\
                       & \leq d(f(0), x)+d(x,[f(n+k), f'(n+k)])\\
                       & \leq n+\delta.
\end{align*}
This inequality is satisfied for each $k\in {\Bbb N}$, which contradicts the condition (i). 
\end{pf}

\begin{lem}\label{lem:bow2}
Let $X$ be a $\delta$-hyperbolic graph with $\partial X\neq \emptyset$ . Then there exists a constant $\delta_{0}'\geq 0$ depending only on $\delta$ which satisfies the following condition:
let $a$ be a point on the boundary $\partial X$ and $f$, $g$ be two geodesic rays to $a$ and $r\geq d(f(0), g(0))$ be a natural number. Then the subray $[f(r), a]\subseteq f$ is $\delta_{0}'$-close to some subray of $g$. 
\end{lem}

\begin{pf}
By Lemma \ref{hyp-fund}, we have 
\[M=\sup_{i\in {\Bbb N}}d(f(i), g(i))<\infty.\]
First, we show that there exists $r'\in {\Bbb N}$ such that $d(f(r), g(r'))\leq 4\delta$. 

Choose geodesics $[g(0), f(M+r)]$, $[f(0), g(0)]$ and $[f(M+r), g(M+r)]$, where $[p, q]$ denotes some geodesic from $p$ to $q$ for $p, q\in X$. By $\delta$-slimness of the triangle $f(0)g(0)f(M+r)$, there exists $x\in [f(0), g(0)]\cup [g(0), f(M+r)]$ such that $d(f(r), x)\leq \delta$. If $x\in [f(0), g(0)]$, then
\begin{align*}
d(f(r), g(0))& \leq d(f(r), x)+d(x, g(0))\\
             & \leq \delta +d(f(0), g(0))-d(f(0), x)\\
             & \leq \delta +r-(d(f(0), f(r))-d(f(r), x))\\
             & \leq 2\delta.
\end{align*}
If $x\in [g(0), f(M+r)]$, then there exists $y\in [g(0), g(M+r)]\cup [f(M+r), g(M+r)]$ such that $d(x, y)\leq \delta$ by $\delta$-slimness of the triangle $g(0)g(M+r)f(M+r)$. If $y\in [g(0), g(M+r)]$, then 
\[d(f(r), y)\leq d(f(r), x)+d(x, y)\leq 2\delta.\]
If $y\in [f(M+r), g(M+r)]$, then 
\begin{align*}
d(f(r), g(M+r))& \leq d(f(r), y)+d(y, g(M+r))\\
               & \leq d(f(r), x)+d(x, y)+d(f(M+r), g(M+r))-d(f(M+r), y)\\
               & \leq 2\delta +d(f(M+r), g(M+r))-(d(f(M+r), x)-d(x, y))\\
               & \leq 2\delta +M-(d(f(M+r), f(r))-d(f(r), x))+\delta \\
               & \leq 4\delta.
\end{align*}

Finally, we show that two subrays $[f(r), a]\subseteq f$ and $[g(r'), a]\subseteq g$ are $12\delta$-close. We may assume that $\delta$ is a natural number. It suffices to show that two sub-rays $[f(r+4\delta), a]\subseteq f$ and $[g(r'+4\delta), a]\subseteq g$ are $12\delta$-close since $d(f(r), g(r'))\leq 4\delta$. 

For any $k\in {\Bbb N}$, applying the above argument to two geodesic rays $[f(r), a]$ and $[g(r'), a]$, we see that there exists a natural number $r''\geq r'$ such that $d(f(r+4\delta +k), g(r''))\leq 4\delta$. Considering the quadrilateral $f(r)f(r+4\delta +k)g(r'')g(r')$, we can see that 
\[|4\delta +k-d(g(r'), g(r''))|\leq 8\delta \]
and that 
\begin{align*}
&d(f(r+4\delta +k), g(r'+4\delta +k))\\
\leq &d(f(r+4\delta +k), g(r''))+d(g(r''), g(r'+4\delta +k))\\
\leq &4\delta +|d(g(r'), g(r'+4\delta +k))-d(g(r'), g(r''))|\\
\leq &12\delta.
\end{align*}
Thus, $[f(r), a]$ and $[g(r'), a]$ are $12\delta$-close. 
\end{pf}

In what follows, we use the symbol $\delta_{0}$ instead of $\delta_{0}'$ for simplicity.

Let $X$ be a $\delta$-hyperbolic graph for some $\delta \geq 0$ with $\partial X\neq \emptyset$ and $\delta_{0}$ be the constant in Lemma \ref{lem:bow}. We suppose that for any $x, y\in V(X)$, we are given a non-empty subset ${\mathcal L}_{T}(x, y)\subseteq {\mathcal L}(x, y)$\index{$L T x y $@${\cal L}_{T}(x, y)$} of geodesics, which we call {\it tight} ones\index{tight geodesic}, regardless of whether ${\mathcal L}_{T}(x, y)$ consists of genuine tight geodesics defined in Subsection \ref{sec:finite}. We define $G(x, y)$, ${\mathcal L}_{T}(x, y; r)$, $G(x, y; r)$\index{$G x y $@$G(x, y)$}\index{$L T x y r$@${\cal L}_{T}(x, y;r)$}\index{$G x y r$@$G(x, y; r)$} for $r\geq 0$ as before. We suppose that every subpath of a tight geodesic is also tight and that ${\mathcal L}_{T}(x, y)$ satisfies the following finiteness properties: 

\begin{enumerate}
\item[(F1)]\index{$F 1$@(F1)} There exists a constant $P_{0}\in {\Bbb N}$ such that for any $x, y\in V(X)$ and $z\in G(x, y)$, the intersection $G(x, y)\cap B(z; \delta_{0})$ has at most $P_{0}$ elements;

\item[(F2)]\index{$F 2$@(F2)} There exist constants $P_{1}\in {\Bbb N}$, $\delta_{1}>0$ such that for any $r\in {\Bbb N}$, $x, y\in V(X)$, and $z\in G(x, y)$ with $d(z, \{ x, y\} )\geq r+\delta_{1}$, the intersection $G(x, y; r)\cap B(z; \delta_{0})$ has at most $P_{1}$ elements.
\end{enumerate}

We can define a tight geodesic ray or bi-infinite geodesic as one for which every finite subpath is tight. We extend the notations ${\mathcal L}(x, y)$ etc. to allow $x, y\in X\cup \partial X$, provided that for any $x\in \partial X$, $y\in X\cup \partial X$ with $x\neq y$, we have $d(x, y)=d(y, x)=\infty$. 

\begin{lem}[\ci{Lemmas 3.1, 3.2}{bowditch}]\label{non-empty-lem}
With the notation as above,
\begin{enumerate}
\renewcommand{\labelenumi}{\rm(\roman{enumi})}
\item the set ${\mathcal L}_{T}(x, y)$ is non-empty for all $x, y\in V(X)\cup \partial X$.
\item the finiteness properties $({\rm F1})$ and $({\rm F2})$ hold for $x, y\in V(X)\cup \partial X$.
\end{enumerate}
\end{lem}
If a $\delta$-hyperbolic graph $X$ has bounded geometry and ${\mathcal L}_{T}(x, y)={\mathcal L}(x, y)$, then $X$ satisfies the properties (F1) and (F2). The following theorem is the main result in this section:

\begin{thm}\label{main1}
Let $X$ be a $\delta$-hyperbolic graph for some $\delta \geq 0$ with $\partial X\neq \emptyset$ and suppose that we can find a collection of tight geodesics $\{ {\mathcal L}_{T}(x, y)\}_{x, y\in V(X)}$ satisfying the properties $({\rm F1})$ and $({\rm F2})$. Then $X$ has property A.
\end{thm}

\begin{pf}
This is a modification of Germain and Tu's proof \cite[Appendix B]{ana}, \cite[Proposition 8.1]{tu} for a hyperbolic graph with bounded geometry, though when $X$ does not have bounded geometry, then we cannot directly apply the same argument. To overcome the difficulty, we use the finiteness properties (F1) and (F2).  

We identify $X$ with the vertex set of the graph $X$.

Let $a$ be a fixed point on the boundary $\partial X$, and $n, k\in {\Bbb N}$ be natural numbers such that
\begin{equation}
n-2k-\delta_{0}-\delta_{1}>0.
\label{2.1}
\end{equation}
For $x\in X$, define a function $F_{a}(x, k, n)$ on $X$ by the characteristic function on the set
\[\bigcup_{f\in {\mathcal L}_{T}(x, a; k)}f([n, 2n]) \]
and define 
\[H_{a}(x, n)=n^{-3/2}\sum_{k<\sqrt{n}}F_{a}(x, k, n)\]
for any $n\in {\Bbb N}$ such that $n\geq N_{0}$, where $N_{0}\in {\Bbb N}$ is a constant satisfying $N_{0}-2\sqrt{N_{0}}-\delta_{0}-\delta_{1}>0$. We will show that $H$ satisfies the following properties:

\begin{enumerate}
\renewcommand{\labelenumi}{\rm(\roman{enumi})}
\item for any $R>0$, we have $\Vert H_{a}(x, n)-H_{a}(y, n)\Vert_{1}\rightarrow 0$ as $n\rightarrow \infty$ uniformly on the set $\{ (a, x, y)\in \partial X\times X\times X: d(x, y)<R\}$;
\item the $\ell^{1}$-norm of $H_{a}(x, n)$ is bounded above and away from 0 uniformly on the set $\{ (a, x, n)\in \partial X\times X\times {\Bbb N}\}$.
\end{enumerate}

\begin{lem}
With the above definition, we have
\begin{equation}
n\leq \Vert F_{a}(x, k, n)\Vert_{1}\leq (n+2k+2\delta_{0}+1)P_{1}
\label{2.2}
\end{equation}
for any $a\in \partial X$, $x\in X$, and $k, n\in {\Bbb N}$ satisfying the inequality $(\ref{2.1})$.
\end{lem}

\begin{pf}
It suffices to show the inequality in the right hand side. Fix $g_{0}\in {\mathcal L}_{T}(x, a)$. We will show that the support $S$ of $F_{a}(x, k, n)$ is contained in the set 
\[\bigcup_{m=n-k-\delta_{0}}^{2n+k+\delta_{0}}B(g_{0}(m); \delta_{0}).\]
If it is shown, since $n-k-\delta_{0}>k+\delta_{1}$, the proof will be completed by (F2). For any $y\in S$, we can find a tight geodesic $g\in {\mathcal L}_{T}(x, a;k)$ containing $y$ in the segment $g|_{[n, 2n]}$ and a point $t\in g_{0}$ such that $d(y, t)\leq \delta_{0}$ by Lemma \ref{lem:bow2}. Considering the quadrilateral $g(0)xty$ and using the triangle inequality, we see that $n-k-\delta_{0}\leq d(x, t)\leq 2n+k+\delta_{0}$. 
\end{pf}

\begin{lem} 
We have 
\begin{equation}
F(x, k, n)\leq F(y, k+R, n)
\label{2.3}
\end{equation}
for every $R>0$, $n, k\in {\Bbb N}$ with $n-2k-2R-\delta_{0}-\delta_{1}>0$, and $x, y\in X$ with $d(x, y)<R$.
\end{lem}
This is clear by definition.

Return to the proof of Theorem \ref{main1}. Fix any $R>0$. Then using the inequality (\ref{2.3}), for $x, y\in X$ with $d(x, y)<R$ and $n\in {\Bbb N}$ with $n-2\sqrt{n}-\delta_{0}-\delta_{1}>0$, we have
\begin{align*}
\sum_{0\leq k <\sqrt{n}}F_{a}(x, k, n) &\leq \sum_{\sqrt{n}-R\leq k<\sqrt{n}}F_{a}(x, k, n)+\sum_{0\leq k<\sqrt{n}-R}F_{a}(y, k+R, n)\\
&\leq \sum_{\sqrt{n}-R\leq k<\sqrt{n}}F_{a}(x, k, n)+\sum_{0\leq k<\sqrt{n}}F_{a}(y, k, n).
\end{align*}
By symmetry,
\begin{equation*}
\left\vert \sum_{0\leq k <\sqrt{n}}(F_{a}(x, k, n)-F_{a}(y, k, n))\right\vert  \leq \sum_{\sqrt{n}-R\leq k<\sqrt{n}}(F_{a}(x, k, n)+F_{a}(y, k, n)),
\end{equation*}
which implies
\begin{equation*}
\vert H_{a}(x, n)-H_{a}(y, n)\vert \leq n^{-3/2}\sum_{\sqrt{n}-R\leq k<\sqrt{n}}(F_{a}(x, k, n)+F_{a}(y, k, n)),
\end{equation*}
and 
\begin{equation*}
\Vert H_{a}(x, n)-H_{a}(y, n)\Vert_{1}\leq 2n^{-3/2}R(n+2\sqrt{n}+2\delta_{0}+1)P_{1}
\end{equation*}
by the inequality (\ref{2.2}). The right hand side converges to 0 as $n$ increases and thus, the left hand side converges to 0 uniformly on the set $\{(a, x, y)\in \partial X\times X\times X: d(x, y)<R\}$. On the other hand, by the inequality (\ref{2.2}), the norm of $H_{a}(x, n)$ is bounded above and away from 0 uniformly on the set $\{ (a, x, n)\in \partial X\times X\times {\Bbb N}: n\geq N_{0}\}$. Using these facts, it is easy to prove the existence of a sequence of maps from $X$ to Prob$(X)$ implementing property A of $X$.  
\end{pf}

\begin{rem}\label{constructed-function-action}
Note that the curve complex of a non-exceptional surface has the non-empty boundary at infinity by Theorem \ref{thm-kla-boundary}. When $X$ is (the 1-skeleton of) the curve complex $C$ of a non-exceptional surface $M$, the following condition is also satisfied for the action of $\Gamma(M)$:
\begin{enumerate}
\item[(iii)] $H_{ga}(gx, n)=g\cdot H_{a}(x, n)$ for $g\in \Gamma(M)$, $a\in \partial C$, $x\in C$ and $n\geq N_{0}$, where the action of $\Gamma(M)$ on the space of functions $\varphi$ on $C$ is defined by the formula $(g\cdot \varphi)(x)=\varphi(g^{-1}x)$ for $x\in C$ and $g\in \Gamma(M)$.
\end{enumerate}
This follows by definition. 
\end{rem}

Combining the above theorem with Theorem \ref{thm:bow}, we get the following corollary:

\begin{cor}
The curve complex of a non-exceptional surface has property A.
\end{cor}

\begin{rem}
Recently, Bell-Fujiwara \cite{bell-fuji} obtained the same conclusion, using Bowditch's formulation of the finiteness property of tight geodesics. 
\end{rem}


\subsection{Yu's property A}  

We add some remarks on Yu's definition of property A \cite{yu}. 

\begin{defn}
For a set $Z$, we denote by ${\mathcal F}(Z)$ the set of all non-empty finite subsets of $Z$. Let $X$ be a discrete metric space, we say that $X$ has {\it property A in the sense of Yu}\index{property A!in the sense of Yu} if there exists a sequence of maps $A_{n}\colon X \rightarrow {\mathcal F}(X\times {\Bbb N})$ such that
\begin{enumerate}
\renewcommand{\labelenumi}{\rm(\roman{enumi})}
\item for every $n\in {\Bbb N}$, there exists a constant $R>0$ for which
\[A_{n}(x)\subseteq \{ (x^{\prime}, j)\in X\times {\Bbb N}: d(x, x^{\prime})<R\}\]
for all $x\in X$;
\item for every $K>0$, we have
\[\lim_{n\rightarrow \infty}\sup_{d(x, y)<K}\frac{\vert A_{n}(x)\triangle A_{n}(y)\vert}{\vert A_{n}(x)\cap A_{n}(y)\vert }=0.\]
\end{enumerate}
\end{defn}

Higson-Roe \cite[Lemma 3.5]{higson-roe} shows that Yu's definition of property A is equivalent to that given in Section \ref{sec:pro} for any discrete metric space with bounded geometry. For a general discrete metric space, their argument can be applied to the proof that Yu's property A implies our version of property A. However, it cannot be applied for the converse. The author does not know whether there exists a discrete metric space with our version of property A, but without Yu's property A. 

Yet, in the case of a metric space satisfying the hypothesis in Theorem \ref{main1}, their argument can be applied and we can show that it also has Yu's property A as follows: 

Let $H_{a}(x, n)$ be the function considered in the proof of Theorem \ref{main1}. It follows from its definition and (\ref{2.2}) that for each $n\in {\Bbb N}$ for which $H_{a}(x, n)$ can be defined, the cardinality of the support of $H_{a}(x, n)$ is not bigger than $(n+2\sqrt{n}+2\delta_{0}+1)P_{1}$. Since the norm of $H_{a}(x, n)$ is bounded above and away from 0, by normalization and approximation argument, we can find a sequence of maps $a^{n}$ from $X$ to Prob$(X)$ satisfying the property implementing our version of property A with the additional property that for each $n\in {\Bbb N}$, there exists a natural number $P$ such that for any $x\in X$, the function $a^{n}_{x}\in \ell^{1}(X)$ takes only values in the range
\[0/P, 1/P,\ldots, P/P.\]
Define the subset
\[A_{n}(x)=\{ (x^{\prime}, j)\in X\times {\Bbb N}: j\geq 1,\ j/P\leq a^{n}_{x}(x^{\prime})\}.\]
By definition, it follows the condition (i) of Yu's definition and that $\vert A_{n}(x)\vert =P$ for every $x\in X$. In addition, we have
\[\vert A_{n}(x)\triangle A_{n}(y)\vert = P\Vert a^{n}_{x}-a^{n}_{y}\Vert_{1}=\vert A_{n}(x)\vert \Vert a_{x}^{n}-a_{y}^{n}\Vert_{1},\]
and the condition (ii) of Yu's definition follows from
\[\lim_{n\rightarrow \infty}\sup_{d(x, y)<K}\frac{\vert A_{n}(x)\triangle A_{n}(y)\vert}{\vert A_{n}(x)\vert}=\lim_{n\rightarrow \infty}\sup_{d(x, y)<K}\Vert a_{x}^{n}-a_{y}^{n}\Vert_{1}=0.\]
  
By \cite[Theorem 2.7]{yu}, we get

\begin{cor}
Any discrete metric space satisfying the hypothesis of Theorem 2.3.2 has property A in the sense of Yu and admits a uniform embedding into a Hilbert space. In particular, the curve complex of a non-exceptional surface has this property. 
\end{cor}

\begin{rem}
Kaimanovich-Masur \cite[Theorem 1.3.2]{kai-mas} showed that the Teichm\"uller space with the Teichm\"uller metric is quasi-isometric to some graph with uniformly finite degrees. It is interesting to know whether the graph has property A or not. (This graph is not hyperbolic \cite{masur-wolf}, \cite{mc-pa2}, \cite{mc-pa3}.)  
\end{rem}


\section{Exceptional surfaces}\label{section-exceptional}

In this section, we study the curve complex of an exceptional surface. Denote a compact orientable surface of type $(g, p)$ by $M_{g, p}$. Recall that the complexity $\kappa(M_{g, p})$ of $M_{g, p}$ is defined by $\kappa(M_{g, p})=3g+p-4$. 

Let $M$ be a surface with $\kappa(M)=0$. If we defined the curve complex for $M$ similarly to that for a non-exceptional surface, then the simplicial complex would be $0$-dimensional since any two non-peripheral simple closed curves on $M$ in distinct non-trivial isotopy classes always intersect. Hence, we modify the definition slightly as follows:

\begin{defn}\label{defn-curve-complex-exceptional}
If $\kappa(M)=0$, that is, $M$ is either of type $(1, 1)$ or $(0, 4)$, then we define the {\it curve complex}\index{curve complex!for an exceptional surface} $C=C(M)$ of $M$ by the simplicial complex whose vertex set $V(C)$ is the set of non-trivial isotopy classes of non-peripheral simple closed curves on $M$, and whose edges are the set of pairs $\{ \alpha, \beta \}$ of two distinct elements in $V(C)$ where $\alpha$ and $\beta$ have the lowest possible intersection number, that is, 1 for $M_{1, 1}$ and 2 for $M_{0, 4}$.

If $M=M_{0, 3}$, define $C(M)$ to be the empty set.
\end{defn}
Since other exceptional surfaces do not arise as a subsurface of a hyperbolic surface, we do not need to treat them.

\begin{figure}
\centering
\epsfxsize=\textwidth
\epsfbox{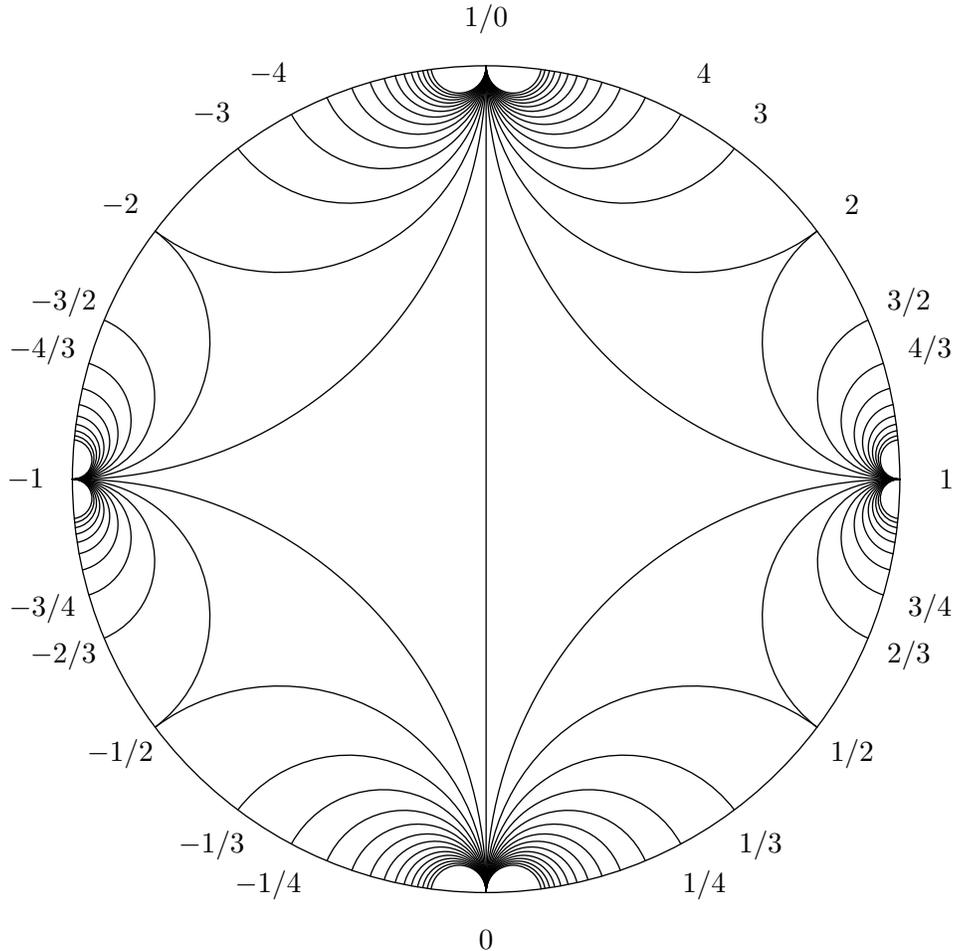}
\caption{Edges incident to $0$, $1$, $1/0$, or $-1$}
\label{picture-ccc}
\end{figure}

If $\kappa(M)=0$, then $C=C(M)$ is a graph and in fact, is isomorphic to the familiar Farey graph in the hyperbolic plane (see \cite{minsky} and Figure \ref{picture-ccc}). We explain this and show that $C$ has property A. 

Assume that $M$ is of type $(1, 1)$. The case of a surface of type $(0, 4)$ can be treated in a similar way. If we choose an identification of the first homology group $H_{1}(M)$ with ${\Bbb Z}^{2}$, then non-trivial homotopy classes of non-peripheral simple closed curves on $M$ are classically well-known to be in one-to-one correspondence with the slopes of elements in ${\Bbb Z}^{2}$, namely the rational numbers $p/q$ including $1/0=\infty$. Thus, $V(C)$ is identified with $\Hat{{\Bbb Q}}={\Bbb Q}\cup \{ \infty \}$\index{$Q $@$\Hat{{\Bbb Q}}$} on the circle $S^{1}=\Hat{{\Bbb R}}={\Bbb R}\cup \{ \infty \}$\index{$R $@$\Hat{{\Bbb R}}$}, and the identification is unique up to the action of $SL(2, {\Bbb Z})$, the group of automorphisms of ${\Bbb Z}^{2}$, on the circle. For more details on this identification, we refer the reader to \cite[Expos\'e 1, \S VI]{FLP} in the case of $M_{1, 1}$ and to \cite[Section 7]{ivanov0} in the case of $M_{0, 4}$. Although $M_{1, 0}$ is treated in \cite[Expos\'e 1, \S VI]{FLP}, we can proceed the above identification similarly in the case of $M_{1, 1}$ with a bit of modification.

The graph $C$ can be embedded as an ideal triangulation of the unit disk bounded by $S^{1}$ (where the word ``ideal'' means that the vertices are on the boundary) by realizing each edge as a hyperbolic geodesic in the Poincar\'{e} disk. Moreover, the graph $C$ is connected \cite[Theorem 1.1]{hat-thu}. Note that the geometric intersection number of curves represented by $p/q$ and $r/s$ is simply the absolute value of the determinant $|ps-qr|$ if the pairs $p$, $q$ and $r$, $s$ have no non-trivial common divisors, respectively. (For the case of $M_{0, 4}$, twice the determinant is the intersection number.) We can write each vertex other than $1/0$ uniquely by $p/q$ so that $(p, q)\in {\mathbb Z}\times ({\mathbb N}\setminus \{ 0\})$ and $p$ and $q$ have no non-trivial common divisors. In what follows, we write each vertex in this way.

Minsky showed the following theorem:

\begin{thm}[\ci{Section 3}{minsky}]\label{excep-hyp}
If $\kappa(M)=0$, then the curve complex $C(M)$ is hyperbolic.
\end{thm}

For $x, y\in V(C)$, let ${\mathcal L}_{T}(x, y)={\mathcal L}(x, y)$, which is non-empty since $C$ is connected. We want to show that $C$ has the finiteness properties (F1) and (F2). 

By the embedding of $C$ in the disk, the complex $C$ has the property that any edge $e$ separates $C$ into two components. Let $e$ be an edge connecting two vertices $v_{1}, v_{2}\in V(C)$ and $x$, $y$ be two points in $V(C)$. We say that $e$ {\it separates}\index{separate} $x$ and $y$ if $x$ and $y$ belong to different connected components of the complement of $\overline{e}=e\cup v_{1}\cup v_{2}$. Given two vertices $x$ and $y$, let $E(x, y)$\index{$E x y $@$E(x, y)$} denote the set of edges that separate $x$ from $y$. Then not only geodesics but every continuous path in $C$ from $x$ to $y$ must pass through some vertex of every edge $e\in E(x, y)$. Conversely, we can show the following:

\begin{lem}\label{edge-pass}
Let $x, y\in V(C)$ with $d(x, y)\geq 2$ and $z\in G(x, y)\setminus \{ x, y\}$. Then $z$ is a vertex of some edge in $E(x, y)$. In particular, $E(x, y)$ is non-empty.
\end{lem}

For the proof, we need the following lemma:

\begin{lem}\label{edge-pass-lem}
Let $f\in {\mathcal L}(x, y)$ be a geodesic such that $d(x, y)\geq 2$ and let $n_{0}\in {\mathbb N}$ with $1\leq n_{0}\leq d(x, y)-1$. Denote $x'=f(n_{0}-1)$, $y'=f(n_{0}+1)$. Suppose that there exists an edge $e\in E(x', y')$ such that $f(n_{0})$ is a vertex of $e$. Then $e$ is also in $E(x, y)$.
\end{lem}

\begin{pf}
First, note that neither $x$ nor $y$ are the vertices of $e$ by the following argument: it is clear that neither $x$ nor $y$ are equal to $f(n_{0})$. If $x$ were the vertex of $e$ different from $f(n_{0})$, then $d(x, x')\geq 1$, $d(x', f(n_{0}))=1$ and $d(x, f(n_{0}))=1$, which is impossible since $f$ is a geodesic. Similarly, $y$ is not a vertex of $e$.     

Denote $n_{1}=d(x, y)$. Suppose that both $x$ and $y$ are in the same side of $e$. Then either $f|_{[0, n_{0}-1]}$ or $f|_{[n_{0}+1, n_{1}]}$ would have to pass through at least one vertex of $e$ since either $x'$ or $y'$ are in the opposite side of $e$ to $x$, $y$. It is impossible. 
\end{pf}

\begin{pf*}{{\sc Proof of Lemma \ref{edge-pass}.}} Let $f\in {\mathcal L}(x, y)$ be a geodesic passing through $z$ and $n_{0}\in {\mathbb N}$ be an integer with $f(n_{0})=z$. Identifying the set of vertices of $C$ with the set $\Hat{{\Bbb Q}}={\Bbb Q}\cup \{ \infty \}$ of rational points on the circle, denote $p/q=f(n_{0}-1)$, $r/s=f(n_{0}+1)$ with $q$, $s\geq 0$. Applying the action of $SL(2, {\mathbb Z})$, we may assume that $p/q=1/0=\infty$ and $z=f(n_{0})=0$. Then $r=1$ or $-1$.

Since $f$ is a geodesic, $r/s$ is neither equal to $1$ nor $-1$ because $d(1/0, 1)=d(1/0, -1)=1$. If $r=1$, then the edge between 0 and 1 is in $E(1/0, r/s)$ and if $r=-1$, then the edge between 0 and $-1$ is in $E(1/0, r/s)$. Applying the previous lemma, the edge is also in $E(x, y)$. 
\end{pf*}

Let $V(x, y)$ be the set of vertices of all edges in $E(x, y)$. By Lemma \ref{edge-pass}, vertices of any geodesic in ${\cal L}(x, y)$ except for $x$ and $y$ are contained in $V(x, y)$. We can see Theorem \ref{excep-hyp} immediately from this fact.

Fix two distinct points $x, y\in V(C)$. If $E(x, y)$ is empty, then $d(x, y)=1$ by Lemma \ref{edge-pass} and there exists an edge between $x$ and $y$.

From now on, we assume $E(x, y)\neq \emptyset$. 

It is easy to see that $E(x, y)$ admits a totally order where $e<f$ if $e$ separates $x$ from the interior of $f$. A vertex $v$ is called a {\it pivot}\index{pivot} for $x$ and $y$ if there exist at least two edges in $E(x, y)$ incident to $v$. If $v$ is a pivot, let $w(v)$\index{$w(v)$} be the number of edges in $E(x, y)$ incident to $v$. 

\begin{lem}
Suppose $E(x, y)\neq \emptyset$ and let $v$ be a pivot for two distinct vertices $x$, $y$.
\begin{enumerate}
\item[(i)] The number $w(v)$ is finite.
\item[(ii)] Let $e_{1}< \cdots < e_{w(v)}$ be the sequence of edges in $E(x, y)$ incident to $v$ and $v_{j}$ be the vertex of $e_{j}$ different from $v$. Then there exists an edge between $v_{j}$ and $v_{j+1}$ for each $j$.
\end{enumerate}
\end{lem}

\begin{pf}
We may assume $v=1/0$. Since $v$ is a pivot, there exists an edge in $E(x, y)$ incident to $v$. We may assume that the edge is $(1/0, 0)$ and $x=p/q<0$, $y=s/t>0$. The set of all edges incident to $1/0$ is $\{ (1/0, n): n\in {\mathbb Z}\}$. Thus, the set of all edges in $E(x, y)$ incident to $1/0$ is $E=\{ (1/0, n): n\in {\mathbb Z},\ p/q<n<s/t\}$, which is finite.

The induced order on $E$ is given by $(1/0, n)<(1/0, m)$ if and only if $n<m$. Since there exists an edge between $n$ and $n+1$ for each $n\in {\mathbb Z}$, the assertion (ii) holds. 
\end{pf}

For $e\in E(x, y)$, we consider the possible ways in which a geodesic in ${\cal L}(x, y)$ may cross $e$. 

\begin{lem}\label{distance-two}
Let $e\in E(x, y)$. If $e$ is incident to no pivots, then $E(x, y)=\{ e\}$ and $x$ and $y$ are opposite vertices of a quadrilateral consisting of two triangles in $C$. Thus $d(x, y)=2$ and there exist two geodesics, going around the quadrilateral in the two possible ways (see Figure \ref{picture-dis} (i)).
\end{lem}

\begin{pf}
We may assume that $e=(1/0, 0)$ and $x=p/q<0$, $y=s/t>0$. Since $1/0$ is not a pivot, we have $\{ (1/0, n): n\in {\mathbb Z},\ p/q<n<s/t\} =\{ (1/0, 0)\}$ and thus, $-1\leq p/q<0$ and $0<s/t\leq 1$. Similarly, since $0$ is not a pivot, both the sets $\{ (0, -1/n): n\in {\mathbb N}\setminus \{ 0\},\ -1/n<p/q \}$ and $\{ (0, 1/n): n\in {\mathbb N}\setminus \{ 0\},\ s/t<1/n\}$ are empty. Remark that the set of all edges incident to $0$ is
\[\{ (0, -1/n): n\in {\mathbb N}\setminus \{ 0\} \} \cup \{ (0, 1/n): n\in {\mathbb N}\setminus \{ 0\} \} \cup \{ (0, 1/0)\}.\]    
Thus, $p/q\leq -1$ and $1\leq s/t$. Combining these inequalities, we get $x=p/q=-1$ and $y=s/t=1$. 
\end{pf}

If $e$ is incident to a pivot $v$, let $e_{1}< \cdots < e_{w(v)}$ be the sequence of edges (including $e$) incident to $v$ and let $v_{j}$ be the vertex of $e_{j}$ different from $v$. 

Any of the shortest paths from $x$ to $y$ must enter the union $\bigcup \overline{e_{i}}$ at one of the vertices of $e_{1}$ and exit through one of the vertices of $e_{w(v)}$. If $w(v)\geq 4$, then the path must pass through $v$ and not around through any of $v_{2},\ldots , v_{w(v)-1}$. If $w(v)=2$ or $3$, then there exist uniformly finitely many choices of the shortest ways to pass the union $\bigcup \overline{e_{i}}$. In any case, any of the shortest paths cannot pass through any vertices on the other side of an edge $(v_{j}, v_{j+1})$ from $v$ since any such digression would be longer than simply traversing the edge.

\begin{lem}\label{lem-distinguish}
Let $x, y\in V(C)$ with $d(x, y)\geq 2$ and $z\in G(x, y)\setminus \{ x, y\}$. Let $e$ be an edge in $E(x, y)$ incident to $z$ (see Lemma \ref{edge-pass}). Let $v$ be the vertex of $e$ different from $z$. If $z$ is not a pivot for $x$ and $y$, then one of the following three cases occurs (see Figure \ref{picture-dis}):
\begin{enumerate}
\item[(i)] The vertex $v$ is not a pivot and then $d(x, y)=2$;
\item[(ii)] We have $w(v)=2$. Let $e_{1}<e_{2}$ be edges in $E(x, y)$ incident to $v$ and $v_{i}$ be the vertex of $e_{i}$ different from $v$. In this case, either $v_{1}$ or $v_{2}$ is $z$ and if $v_{1}=z$ (resp. $v_{2}=z$), then $d(v_{1}, x)=1$ (resp. $d(v_{2}, y)=1$);
\item[(iii)] We have $w(v)=3$. Let $e_{1}<e_{2}<e_{3}$ be edges in $E(x, y)$ incident to $v$ and $v_{i}$ be the vertex of $e_{i}$ different from $v$. In this case, $v_{2}=z$.
\end{enumerate}
\end{lem}

\begin{pf}
If $v$ is not a pivot, then $d(x, y)=2$ by Lemma \ref{distance-two}. 
 
We assume that $v$ is a pivot, that is, $w(v)\geq 2$. If $w(v)=2$ or $w(v)\geq 4$, then $z$ is either $v_{1}$ or $v_{w(v)}$ by the observation before this lemma because $z$ is in $G(x, y)$. Suppose that $v_{1}=z$. We may assume that $v=0$, $v_{1}=1/0$ and $0<v_{w(v)}< +\infty$ by applying the action of $SL(2, {\mathbb Z})$. Since $e_{1}$ is the smallest edge in $E(x, y)$ incident to $v$, we have $-\infty <x\leq -1$. On the other hand, since $v_{1}=1/0$ is not a pivot, the vertex $x$ must be $-1$, which shows $d(v_{1}, x)=1$. We can show similarly that if $v_{w(v)}=z$, then $d(v_{w(v)}, y)=1$. 

Moreover, if $w(v)\geq 4$ and $d(v_{1}, x)=1$, then any geodesic $f$ in ${\cal L}(x, y)$ must pass through $v$ by the observation before this lemma. It follows from $d(x, z)=d(x, v)=1$ that $f$ can not pass through $z$, which contradicts $z\in G(x, y)$. Similarly, we can deduce a contradiction in the case where $w(v)\geq 4$ and $d(v_{w(v)}, y)=1$. Thus, we have $w(v)=2$.  

Assume $w(v)=3$ and $v_{1}=z$. Then we can also deduce $d(v_{1}, x)=1$ similarly. Any geodesic $f\in {\cal L}(x, y)$ which passes through $z$ can not pass through $v$ because $d(x, z)=d(x, v)=1$. Since $f$ passes through one of the vertices of any edge in $E(x, y)$, it must pass through $v_{2}$ and $v_{3}$. It is a contradiction because $d(x, v_{2})=d(x, v_{3})=2$. Similarly, we can deduce a contradiction in the case where $w(v)=3$ and $v_{3}=z$.    
\end{pf}

\begin{figure}[h]
\centering
\epsfxsize=\textwidth
\epsfbox{diss.1}
\caption{}
\label{picture-dis}
\end{figure}

\begin{lem}\label{excep-f1}
There exists $P_{0}\in {\Bbb N}$ such that for any $x, y\in V(C)$ and $z\in G(x, y)\setminus \{ x, y\}$, the intersection $G(x, y)\cap B(z; 1)$ has at most $P_{0}$ elements. 
\end{lem}

\begin{pf}
We may assume that $x\neq y$.

If $z$ is not a pivot for $x$ and $y$, then we have the possibilities in Lemma \ref{lem-distinguish}. In the case (i), the cardinality $|G(x, y)|$ is $4$. In the other cases, we have $|G(x, y)\cap B(z;1)|\leq 4$.

We assume that $z$ is a pivot for $x$ and $y$. Let $e_{1}<\cdots <e_{w(z)}$ be the edges in $E(x, y)$ incident to $z$. We may assume that $z=0$ and $e_{1}$ is equal to the edge $(0, 1/0)$ and that $-\infty <x<0$ and $0<y<+\infty$. In this case, the vertex $v_{i}$ of $e_{i}$ different from $z$ is $1/(i-1)$. It follows that $-\infty <x\leq -1$ and $1/w(z)\leq y<1/(w(z)-1)$. 

Let $f$ be a geodesic in ${\mathcal L}(x, y)$. First, assume that $f$ does not pass through $z$. Then $f$ must pass through all $v_{1}, \ldots, v_{w(z)}$. It follows that $w(z)=2$ or $3$. Moreover, we claim that $f$ can not pass through any vertex $-1/k$ with $k\geq 2$. If $f$ passed through $-1/k$ with $k\geq 2$, then the restriction of $f$ from $x$ to $-1/k$ (resp. from $-1/k$ to $y$) would pass through the vertex $-1$ since the edge $(0, -1)$ is in $E(x, -1/k)\cap E(-1/k, y)$, which is a contradiction. Similarly, we can verify that $f$ can not pass through any vertex $1/k$ with $k\geq w(z)+1$. It follows that 
\[f\cap B(z;1)\subseteq \{ -1, v_{1},\ldots, v_{w(z)}, 1/w(z)\}.\]

Next, assume that $f$ passes through $z$. Let $t_{0}\in {\mathbb N}$ be a natural number such that $f(t_{0})=z=0$. Suppose that there exists $t\in {\mathbb N}$ such that $f(t)=-1/k$ with $k\geq 2$. If $t\geq t_{0}$, then $t=t_{0}+1$ since $d(0, -1/k)=1$. The geodesic segment $f|_{[t, d(x, y)]}$ must pass through $0$ or $1/0$ because $1/w(z)\leq y<1/(w(z)-1)$ and thus, $(0, 1/0)\in E(-1/k, y)$. This is a contradiction. Thus, we have $t< t_{0}$. On the other hand, since the geodesic segment $f|_{[0, t]}$ must pass through $-1$ or $1/0$, there exists $t'< t$ such that $f(t')=-1$ or $=1/0$. Thus, $1=d(-1, 0)=d(1/0, 0)=t_{0}-t'>1$, which is a contradiction. Hence, we see that $f$ can not pass through any vertex $-1/k$ with $k\geq 2$. Similarly, we can show that $f$ can not pass either through any vertex $1/k$ with $k\geq w(z)+1$. 

Hence, if $w(z)\leq 3$, then 
\[f\cap B(z;1)\subseteq \{ z, v_{1},\ldots, v_{w(z)}, -1, 1/w(z)\}\]
and if $w(z)\geq 4$, then
\[f\cap B(z;1)\subseteq \{ z, v_{1}, v_{w(z)}, -1, 1/w(z)\}.\]

It follows that $G(x, y)\cap B(z;1)$ has at most six elements in any case. 
\end{pf}

\begin{prop}
The curve complex $C=C(M)$ with $\kappa(M)=0$ satisfies the property ${\rm (F1)}$.
\end{prop}

\begin{pf}
Let $x, y\in V(C)$ and we may assume that $d(x, y)\geq 2$. Then 
\[G(x, y)\subseteq \bigcup_{n=1}^{d(x, y)-1}B(f(n); 1)\]
for a fixed geodesic $f\in {\mathcal L}(x, y)$ by Lemma \ref{edge-pass} and the fact that any geodesic in ${\mathcal L}(x, y)$ passes through at least one vertex of each edge in $E(x, y)$. 

Let $z\in G(x, y)$ and $\delta_{0}\geq 0$ be the constant in Lemma \ref{lem:bow}. Let $1\leq n\leq d(x, y)-1$ be a natural number such that $d(f(n), z)\leq 1$. Then we have
\[G(x, y)\cap B(z; \delta_{0})\subseteq G(x, y)\cap \left( \bigcup_{i=n-\delta_{0}-2}^{n+\delta_{0}+2}B(f(i);1)\right).\]
Hence, the intersection $G(x, y)\cap B(z;\delta_{0})$ has at most $(2\delta_{0}+5)P_{0}$ elements.  
\end{pf}

\begin{prop}
The curve complex $C=C(M)$ with $\kappa(M)=0$ satisfies the property ${\rm (F2)}$. 
\end{prop}

\begin{pf}
The proof is almost the same as that for the property (F1). 

Let $x, y\in V(C)$, $r\in {\mathbb N}$ and $\delta_{0}$ be the constant appearing in Lemma \ref{lem:bow}. We consider only the case $d(x, y)\geq 2r+2\delta_{0}+10$. Let $z\in G(x, y; r)$ be a point with $d(z, \{ x, y\} )\geq r+2$ and let $f\in {\mathcal L}(x, y; r)$ be a geodesic with $z=f(n_{0})$ for some $n_{0}\in {\mathbb N}$. We denote by $x_{0}\in B(x; r)$ and $y_{0}\in B(y; r)$ the initial and terminal points of $f$, respectively. It follows from the proof of Lemma \ref{edge-pass} that there exists an edge $e\in E(f(n_{0}-1), f(n_{0}+1))$ which has $z$ as a vertex. By Lemma \ref{edge-pass-lem}, we have $e\in E(x_{0}, y_{0})$. 

Since any two points in $B(x;r)$ can be joined by a path in $B(x;r)$ (resp. $B(y;r)$) through no vertices of $e$, the ball $B(x;r)$ (resp. $B(y;r)$) must be in one side of $e$. Since $e\in E(x_{0}, y_{0})$, the points $x_{0}\in B(x;r)$ and $y_{0}\in B(y;r)$ are in different sides of $e$. Thus, $B(x;r)$ and $B(y;r)$ must be in different sides of $e$. Hence, we have $e\in E(x, y)$.

It follows from the above observation that the set 
\[S=\{ z\in G(x, y;r): d(z, \{ x, y\})\geq r+2\}\]
is contained in the set of vertices of all edges in $E(x, y)$. We define $E(x, y;r)$ to be the set of edges separating $B(x;r)$ and $B(y;r)$ and we can equip a similar total order on $E(x, y;r)$ induced from $E(x, y)$. It follows that any point in $S$ is a vertex of some edge in $E(x, y;r)$.

Let $g\in {\mathcal L}(x, y)$ be a fixed geodesic. By the fact that $g$ passes through at least one vertex of any edge in $E(x, y; r)$, we obtain
\[S\subseteq \bigcup_{i=r+1}^{d(x, y)-r-1}B(g(i);1).\]

Let $z\in S$. We assume that $z=1/0$. Since $d(z, \{ x, y\})\geq r+2$, there exist $n_{1}, n_{2}\in {\Bbb Z}$ such that $n_{1}<x'<n_{1}+1$ and $n_{2}<y'<n_{2}+1$ for any $x'\in B(x;r)$, $y'\in B(y;r)$.

By this fact and a similar argument as in the proof of Lemmas \ref{lem-distinguish} and \ref{excep-f1}, we can show that there exists $P_{1}'\in {\mathbb N}$ such that $S\cap B(z;1)$ has at most $P_{1}'$ elements for any $z\in S$.

Let $z\in G(x, y;r)$ with $d(z, \{ x, y\})\geq r+\delta_{0}+5$. Since $z\in S$, we can find a natural number $r+\delta_{0}+4\leq n\leq d(x, y)-r-\delta_{0}-4$ such that $d(g(n), z)\leq 1$ and 
\[G(x, y;r)\cap B(z;\delta_{0})\subseteq S\cap \bigcup_{i=n-\delta_{0}-2}^{n+\delta_{0}+2}B(g(i);1).\]
The right hand side has at most $(2\delta_{0}+5)P_{1}'$ elements. Thus, the property (F2) holds for $P_{1}=(2\delta_{0}+5)P_{1}'$ and $\delta_{1}=\delta_{0}+5$. 
\end{pf}

We have shown the following theorem: 

\begin{thm}\label{thm-main-excep}
If $\kappa(M)=0$, then the curve complex $C(M)$ satisfies the properties ${\rm (F1)}$ and ${\rm (F2)}$. In particular, $C(M)$ has property A and can be uniformly embedded into a Hilbert space.     
\end{thm}

\chapter[Amenability for the action]{Amenability for the action of the mapping class group on the boundary of the curve complex}\label{chapter:amenable-action}

In this chapter, we study the actions of the mapping class group $\Gamma(M)$ for a surface $M$ with $\kappa(M)\geq 0$ on the boundary $\partial C$ of the curve complex $C$ and on the Thurston boundary ${\cal PMF}$ focusing on amenability in a measurable sense introduced by Zimmer \cite{zim1}. 

Adams \cite{adams1} showed that the action of a hyperbolic group $G$ on the boundary $\partial G$ is amenable for any quasi-invariant measure on $\partial G$. He used this significant property in the process of solving a certain problem concerning equivalence relations generated by hyperbolic groups \cite{adams2}. 

In this chapter, we show that the action of $\Gamma(M)$ on $\partial C$ is amenable for any quasi-invariant measure on $\partial C$ and see as a corollary that the action of $\Gamma(M)$ on ${\cal PMF}$ is amenable for some class of quasi-invariant measures on ${\cal PMF}$. Following Adams' idea, in the next chapter, we use these results for attacking some problems about equivalence relations generated by the mapping class groups.

In Sections \ref{mcg} and \ref{boundary-of-curve-complex}, we prepare fundamental notions about the objects $\Gamma(M)$, ${\cal PMF}$ and $\partial C$. The materials of Sections \ref{mcg} and \ref{boundary-of-curve-complex} can be found in the references \cite{FLP}, \cite{gardiner}, \cite{ivanov1}, \cite{ivanov2}, \cite{kai-mas}, \cite{kla} and \cite{mc-pa}. In particular, the reference \cite[Section 1]{kai-mas} is a detailed survey treating not only the above objects but also the Teichm\"uller space and a connection between them. 

In Section \ref{section-amenable-action}, we verify the amenability for the action of $\Gamma(M)$ on $\partial C$ with any quasi-invariant measure and on ${\cal PMF}$ with certain measures. In the proof, we use the functions constructed in the proof of Theorem \ref{main1}. We recommend the reader to see Appendix \ref{general-amenable} for generalities of amenability for actions of discrete groups.

\section{The mapping class group and the Thurston boundary}\label{mcg}

Let $M$ be a compact orientable surface of type $(g, p)$ with $\kappa(M)=3g+p-4\geq 0$. Let $\Gamma(M)$ be the mapping class group of $M$. It is finitely generated and moreover, finitely presented \cite[Theorem 4.2.D]{ivanov2}.  

Let $C=C(M)$ be the curve complex of $M$ as defined in Definitions \ref{defn-curve-complex-non-exceptional} and \ref{defn-curve-complex-exceptional}. If $\kappa(M)>0$, then let $S(M)$\index{$S M$@$S(M)$} denote the set of all simplices of $C(M)$, otherwise, denote $S(M)=V(C(M))$. The {\it geometric intersection number}\index{geometric intersection number} of two elements $\alpha$, $\beta \in V(C)$ is defined by the minimal number of intersections of any two their representatives, and is denoted by $i(\alpha, \beta)$\index{$i(\alpha, \beta)$}. Let ${\Bbb R}^{V(C)}_{\geq 0}$\index{$R V C 0$@${\Bbb R}^{V(C)}_{\geq 0}$} be the space of non-negative functions on $V(C)$ given the product topology. The quotient of ${\Bbb R}^{V(C)}_{\geq 0}\setminus \{ 0\}$ by the multiplicative action of ${\Bbb R}_{>0}$, the multiplicative group of all positive real numbers, is denoted by $P{\Bbb R}^{V(C)}_{\geq 0}$\index{$P R V C 0$@$P{\Bbb R}^{V(C)}_{\geq 0}$}.

The map $\alpha \mapsto i(\alpha, \cdot)$ defines an embedding of $V(C)$ into ${\Bbb R}^{V(C)}_{\geq 0}$ and moreover, its projection into $P{\Bbb R}^{V(C)}_{\geq 0}$ is also an embedding. The closure of the set $\{r\alpha : r\geq 0,\ \alpha \in V(C)\}$ in ${\Bbb R}^{V(C)}_{\geq 0}$ is denoted by ${\cal MF}$\index{$M F$@${\cal MF}$}, and the closure of the embedding of $V(C)$ into $P{\Bbb R}^{V(C)}_{\geq 0}$ is denoted by ${\cal PMF}$\index{$P M F$@${\cal PMF}$}. They are called the {\it space of measured foliations}\index{space!of measured foliations} and the {\it space of projective measured foliations}\index{space!of projective measured foliations}, respectively (see \cite[Expos\'e 3, \S III]{FLP}). The space ${\cal PMF}$ is also called the {\it Thurston boundary}\index{Thurston boundary}.

We denote by $V(C)'$\index{$V C p$@$V(C)'$} the set of all isotopy classes of (non-oriented and not necessarily connected) one-dimensional closed submanifolds in $M$ any of whose connected components is not isotopic to a point (see \cite[Expos\'e 4, \S II]{FLP}). In particular, $S(M)$ can be identified with a subset of $V(C)'$. We can define $i(\alpha, \beta)$ for $\alpha \in V(C)'$, $\beta \in V(C)$ similarly and regard $V(C)'$ as a subset of ${\Bbb R}^{V(C)}_{\geq 0}$. Then the closure of the projection of $V(C)'$ in $P{\Bbb R}^{V(C)}_{\geq 0}$ is equal to ${\cal PMF}$. Moreover, two elements $\alpha_{1}, \alpha_{2}\in V(C)'$ are equal in $P{\Bbb R}_{\geq 0}^{V(C)}$ if and only if there exist $n_{1}, n_{2}\in {\Bbb N}$ and $\alpha_{0}\in V(C)'$ such that $\alpha_{1}=n_{1}\alpha_{0}$ and $\alpha_{2}=n_{2}\alpha_{0}$. In particular, $S(M)$ can be regarded as a subset of ${\cal PMF}$.    

Roughly speaking, a measured foliation on $M$ is a one-dimensional foliation of $M$ with singularities of specific types together with a transverse measure, which assigns to each arc transverse to the foliation a non-negative real number subject to certain conditions (see \cite[Expos\'e 5]{FLP} for a precise definition). If $M$ is closed, then a measured foliation $F$ is defined by a finite number of singularities $p_{k}\in M$ of some type on $M$ and an atlas of coordinate charts $(x_{i}, y_{i})\colon U_{i}\rightarrow {\Bbb R}^{2}$ on the complement $M\setminus \{ p_{k}\}$ such that we have
\[x_{j}=f_{ij}(x_{i}, y_{i}), \ \ y_{j}=\pm y_{i}+C\]
for any two overlapping charts $(x_{i}, y_{i})$, $(x_{j}, y_{j})$. The leaves of $F$ is determined by the line where the $y$-coordinate is constant, and the transverse measure of $F$ is $|dy|$. Each measured foliation $F$ defines the function $i(F, \cdot)\colon V(C)\rightarrow {\Bbb R}_{\geq 0}$ by
\[i(F, \alpha)=\inf_{c\in \alpha}\int_{c}|dy|,\]
where the infimum is taken over all representatives $c$ of $\alpha \in V(C)$.   

Let $M$ be a surface with $\kappa(M)\geq 0$. Two measured foliations $F$ and $G$ on $M$ satisfies 
\[i(F, \alpha)=i(G, \alpha)\]
for any $\alpha \in V(C)$ if and only if $F$ and $G$ can be deformed to each other by the two topological operations, isotopies and Whitehead operations. Hence, each equivalence class of a measured foliation for the two topological operations can be regarded as a point in ${\Bbb R}_{\geq 0}^{V(C)}$. 

In fact, ${\cal MF}$ can be identified with the quotient of the set of all measured foliations on $M$ by the two topological operations (see \cite{FLP}, \cite{gardiner} etc.\ for more details). We shall work not with foliations, but with their equivalence classes and also use the term foliations for these equivalence classes.

Topologically, ${\cal MF}$ is homeomorphic to ${\Bbb R}^{6g-6+2p}$ and ${\cal PMF}$ is homeomorphic to the sphere of dimension $6g-7+2p$. Moreover, ${\cal PMF}$ can be identified with the boundary of the Teichm\"uller space of $M$ as described by Thurston, and the union of the Teichm\"uller space and ${\cal PMF}$ is homeomorphic to the closed unit ball of dimension $6g-6+2p$. Under this homeomorphism, ${\cal PMF}$ corresponds to the boundary of the closed unit ball \cite[Expos\'e 1, Th\'eor\`eme 4]{FLP}. 

The intersection number $i(\cdot, \cdot)$ can be extended continuously on ${\cal MF}\times {\cal MF}$ with 
\[i(r_{1}F_{1}, r_{2}F_{2})=r_{1}r_{2}i(F_{1}, F_{2})\]
for any $r_{1}, r_{2}> 0$ and $F_{1}, F_{2}\in {\cal MF}$. Thus, for any two elements $F_{1}, F_{2}\in {\cal PMF}$, whether $i(F_{1}, F_{2})=0$ or $\neq 0$ makes sense. We say that a measured foliation $F\in {\cal MF}$ is {\it minimal}\index{minimal measured foliation} if $i(F, \alpha)>0$ for any $\alpha \in V(C)$, and denote by ${\cal MIN}$\index{$MIN$@${\cal MIN}$} the subset of ${\cal PMF}$ that consists of the projections of all minimal foliations. A measured foliation $F$ is said to be {\it uniquely ergodic}\index{uniquely ergodic measured foliation} if for any $G\in {\cal MF}$ with $i(F, G)=0$, we have $G=rF$ for some $r>0$. The subset of ${\cal PMF}$ consisting of the projections of all uniquely ergodic foliations is denoted by ${\cal UE}$\index{$UE$@${\cal UE}$}, which is a subset of ${\cal MIN}$. Almost all measured foliations are uniquely ergodic in some sense (see Remark \ref{masur-measure}).

The mapping class group $\Gamma(M)$ acts naturally on ${\cal MF}$ and on ${\cal PMF}$ respectively. The intersection number is $\Gamma(M)$-invariant:
\[i(gF_{1}, gF_{2})=i(F_{1}, F_{2})\]
for any $g\in \Gamma(M)$ and $F_{1}, F_{2}\in {\cal MF}$. The subsets ${\cal MIN}$ and ${\cal UE}$ are $\Gamma(M)$-invariant. It follows from 
\[{\cal MIN}=\{ F\in {\cal PMF}: i(F, \alpha)>0 \ {\rm for \ any \ }\alpha \in V(C)\}\]
that ${\cal MIN}$ is a Borel subset of ${\cal PMF}$. We can verify that ${\cal UE}$ is also measurable in ${\cal PMF}$ as follows: let $p\colon ({\cal PMF})^{2}\rightarrow {\cal PMF}$ be the projection onto the first coordinate and $D\subseteq ({\cal PMF})^{2}$ be the diagonal subset. Then we have    
\begin{align*}
{\cal PMF}\setminus {\cal UE}&=\{ F\in {\cal PMF}: i(F, G)=0 \ {\rm for \ some \ }G\in {\cal PMF}\setminus \{ F\} \} \\
         &=p(\{ (F_{1}, F_{2})\in ({\cal PMF})^{2}: i(F_{1}, F_{2})=0\}\cap (({\cal PMF})^{2}\setminus D)).
\end{align*} 
Since the set $\{ (F_{1}, F_{2})\in ({\cal PMF})^{2}: i(F_{1}, F_{2})=0\}\cap (({\cal PMF})^{2}\setminus D)$ can be expressed as the union of an increasing sequence of compact subsets in $({\cal PMF})^{2}\setminus D$, the set ${\cal PMF}\setminus {\cal UE}$ is measurable.
 
Since the union of the Teichm\"uller space and ${\cal PMF}$ is homeomorphic to the closed unit ball of dimension $6g-6+2p$, any element $g$ in $\Gamma(M)$ has a fixed point on the union by Brouwer's fixed point theorem. Each element $g\in \Gamma(M)$ can be classified into one of finite order, a reducible one of infinite order and a pseudo-Anosov one by considering the fixed point set of $g$ on the union of the Teichm\"uller space and ${\cal PMF}$.  

If $g\in \Gamma(M)$ has a fixed point in the Teichm\"uller space, then the isotopy class $g$ contains a diffeomorphism of finite order and thus, $g$ has finite order \cite[Expos\'e 3, Th\'eor\`eme 18]{FLP}. If $g$ has a fixed point in ${\cal PMF}$, then we have the following three possibilities:
\begin{enumerate}
\item[(i)] there exists a non-minimal measured foliation $F$ and $\lambda >0$ such that $gF=\lambda F$;
\item[(ii)] there exists a minimal measured foliation $F$ such that $gF=F$;
\item[(iii)] there exists a minimal measured foliation $F$ and $\lambda >0$, $\lambda \neq 1$ such that $gF=\lambda F$;
\end{enumerate}
If $g$ satisfies (i), then it is said to be {\it reducible}\index{reducible!element}. In this case, there exists an element in $S(M)$ fixed by $g$. (This can also be seen by using the map $H\colon {\cal PMF}\setminus {\cal MIN}\rightarrow S(M)$ constructed in Definition \ref{construction-of-H}.) On the other hand, any element in $S(M)$ can be identified with some non-minimal foliation in ${\cal PMF}$. Thus, the set of all reducible elements can be characterized by the set of all elements fixing some element in $S(M)$. 

If $g$ satisfies (ii), then it has finite order \cite[Expos\'e 9, Lemme 7]{FLP}.

If $g$ satisfies (iii), then it is said to be {\it pseudo-Anosov}\index{pseudo-Anosov!element}. In this case, $g$ has exactly two fixed points $F_{+}(g), F_{-}(g)\in {\cal MIN}$ such that $g^{n}F\rightarrow F_{+}(g)$ for any $F\neq F_{-}(g)$ and $g^{-n}F\rightarrow F_{-}(g)$ for any $F\neq F_{+}(g)$ as $n\rightarrow +\infty$. We remark that $g$ has a stronger convergence property than that stated above (see \cite[Theorem 3.5]{ivanov1}). The name ``pseudo-Anosov'' comes from the following definition. A diffeomorphism $\varphi$ on a surface $M$ is said to be {\it pseudo-Anosov}\index{pseudo-Anosov!diffeomorphism} if there exist two transversing measured foliations $(F^{u}, \mu^{u})$, $(F^{s}, \mu^{s})$ on $M$ and $\lambda >1$ such that
\begin{align*}
\varphi(F^{u})&=F^{u}, \ \varphi_{*}(\mu^{u})=\lambda \mu^{u},\\
\varphi(F^{s})&=F^{s}, \ \varphi_{*}(\mu^{s})=(1/\lambda)\mu^{s},
\end{align*}
where $\mu^{u}$ and $\mu^{s}$ are the transverse measures on the foliations $F^{u}$ and $F^{s}$, respectively. Two measured foliations are said to be transverse if the sets of singularities of the foliations are equal and they are transverse on the complement of the set of singularities. In the above case, the measured foliations $(F^{u}, \mu^{u})$, $(F^{s}, \mu^{s})$ are called the unstable, stable foliations for $\varphi$, respectively. The isotopy class of each pseudo-Anosov element $g\in \Gamma(M)$ contains a pseudo-Anosov diffeomorphism $\varphi$ and then two measured foliations $(F^{u}, \mu^{u})$, $(F^{s}, \mu^{s})$ are in the classes $F_{+}(g)$, $F_{-}(g)$, respectively \cite[Expos\'e 1, Th\'eor\`eme 5]{FLP}.

Moreover, the two foliations $F_{+}(g)$, $F_{-}(g)$ for a pseudo-Anosov element $g$ are uniquely ergodic \cite[Expos\'e 12, Th\'eor\`eme I]{FLP}. The fixed points of pseudo-Anosov elements are called {\it pseudo-Anosov foliations}\index{pseudo-Anosov!foliation}.

For more details of the classification of diffeomorphisms, we refer the reader to \cite[Expos\'es 9, 11]{FLP}.

McCarthy and Papadopoulos classified subgroups of $\Gamma(M)$, using the action on ${\cal PMF}$ as follows:

\begin{thm}[\ci{Theorem 4.6}{mc-pa}]\label{subgroup-classification}
Every subgroup of the mapping class group $\Gamma(M)$ of a surface $M$ with $\kappa(M)\geq 0$ is classified into one of the following four types:
\begin{enumerate}
\renewcommand{\labelenumi}{\rm(\roman{enumi})}
\item a subgroup containing a pair of independent pseudo-Anosov elements (such a subgroup is said to be {\rm sufficiently large})\index{sufficiently large subgroup};
\item an infinite subgroup fixing the set $\{ F_{+}(g), F_{-}(g)\}$ for a pseudo-Anosov element $g\in \Gamma(M)$ (such a subgroup is virtually infinite cyclic, that is, contains an infinite cyclic subgroup of finite index);
\item an infinite subgroup fixing an element in $S(M)$ (such a subgroup is said to be {\rm reducible})\index{reducible!subgroup};
\item a finite subgroup.
\end{enumerate}
\end{thm}
Here, we say that a pair of two pseudo-Anosov elements is {\it independent}\index{independent pair} when their fixed point sets are disjoint. Every sufficiently large subgroup contains a free subgroup of rank $2$ and thus, is non-amenable (see \cite[Corollary 8.4]{ivanov1}).

If $\kappa(M)=0$, then all infinite reducible subgroups are virtually infinite cyclic. On the other hand, if $\kappa(M)>0$, then there always exist non-amenable reducible subgroups. For example, the subgroup generated by Dehn twists around two intersecting curves each of which are disjoint from a third curve is a desired one (see Chapter \ref{chapter-indec}, Section \ref{section-reducible-elements} for the definition of Dehn twists). This fact shows the existence of non-amenable stabilizers of the action of $\Gamma(M)$ on ${\cal PMF}$, which implies the action is not topologically amenable (see Remark \ref{topologically-amenable-action}). Non-amenability of the subgroup can be seen from the following two general theorems. The first one is often called the Tits alternative for the mapping class group.

\begin{thm}[\ci{Theorem A}{mc2}, \ci{Corollary 8.10}{ivanov1}]\label{tits-alternative-mcg}
If $M$ is a surface with $\kappa(M)\geq 0$, then every subgroup of $\Gamma(M)$ either contains a free subgroup of rank $2$ or is virtually abelian, that is, contains an abelian subgroup of finite index.
\end{thm}

\begin{thm}[\ci{Theorem 7.5.C}{ivanov2}]\label{commuting-dehn}
Let $M$ be a surface with $\kappa(M)\geq 0$. Let $t_{\alpha}$, $t_{\beta}$ be two Dehn twists about $\alpha, \beta \in V(C)$, respectively. Let $m$, $n$ be two non-zero integers. Then $t_{\alpha}^{n}t_{\beta}^{m}=t_{\beta}^{m}t_{\alpha}^{n}$ if and only if $\alpha$ and $\beta$ can be represented by disjoint curves.
\end{thm}

Hence, infinite reducible subgroups are often complicated. On the other hand, a hyperbolic group does not have such subgroups, that is, any subgroup of a hyperbolic group is classified into one of the following three types (see \cite[Chapitre 8, Th\'eor\`eme 37]{ghys-harpe}):
\begin{enumerate}
\item[(i)] a subgroup containing a pair of independent hyperbolic elements (such a subgroup contains a free subgroup of rank 2);
\item[(ii)] an infinite subgroup fixing the fixed point set of some hyperbolic element (such a subgroup is virtually infinite cyclic);
\item[(iii)] a finite subgroup.
\end{enumerate}
Here, an element of a hyperbolic group is said to be hyperbolic if it is torsion-free and has exactly two fixed points on the boundary. We can say that a pseudo-Anosov element of the mapping class group is an analogue of a hyperbolic element of a hyperbolic group. Note that any element of a hyperbolic group is hyperbolic or of finite order. This fact is one of the most significant difference between the mapping class group of a non-exceptional surface and a hyperbolic group. Remark that $\Gamma(M)$ is isomorphic to $SL(2, {\Bbb Z})$ if $M=M_{1, 1}$ (see \cite[Expos\'e 1, \S VI]{FLP} in the case of the torus). If $M=M_{0, 4}$, then $\Gamma(M)$ is almost isomorphic to $SL(2, {\Bbb Z})$ (see \cite[Theorem 7.9]{ivanov0}), where two groups are said to be {\it almost isomorphic}\index{almost isomorphic} if they are isomorphic up to taking a subgroup of finite index or a quotient by a finite normal subgroup. Thus, $\Gamma(M)$ is hyperbolic if $\kappa(M)=0$. 

On the other hand, a subgroup of the mapping class group of a surface $M$ with $\kappa(M)\geq 0$ which fixes two curves filling $M$ is almost trivial by the next lemma. We say that two elements $\alpha, \beta\in V(C)$ {\it fill}\index{fill} $M$ if either $i(\alpha, \gamma)\neq 0$ or $i(\beta, \gamma)\neq 0$ for any $\gamma \in V(C)\setminus \{ \alpha, \beta \}$. Remark that if $M$ is non-exceptional, then $\alpha$, $\beta \in V(C)$ fills $M$ if and only if $d(\alpha, \beta)\geq 3$ on the curve complex $C$. 

\begin{lem}\label{fill-curves-fix}
Let $M$ be a surface with $\kappa(M)\geq 0$. 
\begin{enumerate}
\item[(i)] If $\kappa(M)>0$ and $\alpha, \beta \in V(C)$ fill $M$, then the subgroup
\[\Gamma_{\alpha, \beta}=\{ g\in \Gamma(M): g\alpha =\alpha, \ g\beta =\beta \}\]
is finite \cite[Lemma 10]{bestvina-fujiwara}.
\item[(ii)] If $\kappa(M)=0$ and $\alpha, \beta \in V(C)$ satisfies $d(\alpha, \beta)\geq 2$ on the graph $C$, then the subgroup $\Gamma_{\alpha, \beta}$ is finite.
\end{enumerate}  
\end{lem}    

The latter assertion can be easily seen from the description of $C$ in Chapter \ref{chapter-property-A}, Section \ref{section-exceptional} (see also Remark \ref{rem-kla-excep}).

Finally, we give some information on a pair of pants $M_{0, 3}$:

\begin{lem}[\ci{Expos\'e 2, \S III}{FLP}]\label{lem-pants-fund}
The mapping class group $\Gamma(M_{0, 3})$ consists of six elements and is isomorphic to the symmetric group of three letters. 
\end{lem}

In other words, each element $g$ in $\Gamma(M_{0, 3})$ can be distinguished by the behavior of the permutation of the boundary components of $M_{0, 3}$ by $g$.


\section{The boundary at infinity of the curve complex}\label{boundary-of-curve-complex}

Let $M$ be a compact orientable surface of type $(g, p)$ with $\kappa(M)\geq 0$. Let $C=C(M)$ be the associated curve complex. Since it is hyperbolic as a metric space, we can consider the Gromov-boundary $\partial C$ (see Chapter \ref{chapter-property-A}, Section \ref{sec:Acc}). 

For describing it explicitly, we recall the topology of the boundary $\partial X$ of a hyperbolic metric space $X$. Let $e\in X$ be a base point. Define the Gromov product\index{Gromov-!product} on $\partial X$ by
\[(a|b)_{e}=\sup \liminf_{n, m\rightarrow \infty}(x_{n}|y_{m})_{e},\]
where the supremum is taken over all sequences $\{ x_{n}\}$ and $\{ y_{m}\}$ in $X$ such that $\{ x_{n}\}$, $\{ y_{m}\}$ are in the classes $a, b\in \partial X$, respectively. The open neighborhood basis of $a\in \partial X$ is defined by
\[\{ b\in \partial X: (b|a)_{e}>r\}\]
for $r\geq 0$. This topology is independent of the choice of $e\in X$ and metrizable (see \cite[Chapitre 7, Section 3]{ghys-harpe} or \cite[Proposition 3.21]{bridson-haefliger}). If $X$ is not proper, then $\partial X$ is not compact in general.

In fact, $\partial C$ is not compact (see Proposition \ref{prop-boundary-non-compact} and Theorem \ref{thm-excep-identified}). Klarreich \cite{kla} studied $\partial C$ for a non-exceptional surface, and gave the following characterization of it (see also \cite{ham}):

\begin{thm}[\ci{Theorem 1.1}{kla}]\label{thm-kla-boundary}
Let $M$ be a non-exceptional surface. The boundary at infinity of the curve complex is homeomorphic to the space of minimal topological foliations on $M$.
\end{thm}

In particular, $\partial C$ is non-empty. The space of minimal topological foliations is the quotient space of the space ${\cal MIN}$ of minimal measured foliations by identifying topologically equivalent foliations, where two minimal measured foliations $F$, $G$ are said to be {\it topologically equivalent}\index{topologically equivalent!minimal measured foliations} if $i(F, G)=0$. Note that this is actually an equivalence relation of ${\cal MIN}$ and that two minimal measured foliations are topologically equivalent if and only if they can be deformed to each other by some topological operations (see \cite[Theorem 1.12]{rees}). The topology of the space is induced by the measure forgetting quotient map from ${\cal MIN}$ onto it. We denote this map by $\pi \colon {\cal MIN}\rightarrow \partial C$\index{$\ p MIN \ z C$@$\pi \colon {\cal MIN}\rightarrow \partial C$}. Moreover, Klarreich showed the following fact:

\begin{prop}[\ci{Proposition 7.1}{kla}]\label{prop-kla-pi}
The continuous map $\pi$ is closed, and the preimage of any point of $\partial C$ is compact. 
\end{prop}

By \cite[Section 5, Example 1]{mc-pa}, the set $\Lambda_{0}$ of all pseudo-Anosov foliations is dense in ${\cal PMF}$. Moreover, $\Lambda_{0}$ is contained in ${\cal UE}$ as mentioned in Section \ref{mcg}. Thus, the projection of $\Lambda_{0}$ by $\pi$ is dense in $\partial C$ and in particular, $\partial C$ is separable.

\begin{prop}\label{prop-boundary-non-compact}
The boundary $\partial C$ of the curve complex of a non-exceptional surface is non-compact.
\end{prop}

\begin{pf}
We assume that $\partial C$ is compact and deduce a contradiction. Take $F\in {\cal PMF}\setminus {\cal MIN}$. Since ${\cal MIN}$ is dense in ${\cal PMF}$, we can choose a sequence $\{ F_{n}\}$ in ${\cal MIN}$ converging to $F$. Taking a subsequence of $\{ F_{n}\}$, we may assume that $\{ \pi(F_{n})\}$ converges to $a$ in $\partial C$ because $\partial C$ is compact and metrizable. Since $\pi^{-1}(a)$ is compact in ${\cal MIN}$ by Proposition \ref{prop-kla-pi} and thus, in ${\cal PMF}$ and $F\in {\cal PMF}\setminus {\cal MIN}$, there exist open sets $U_{1}$, $U_{2}$ of ${\cal PMF}$ such that $\pi^{-1}(a)\subseteq U_{1}$, $F\in U_{2}$ and $U_{1}\cap U_{2}=\emptyset$. Taking a subsequence, we may assume that $F_{n}\in U_{2}$ for all $n$. It follows from Proposition \ref{prop-kla-pi} that $\partial C\setminus \pi({\cal MIN}\setminus U_{1})$ is an open neighborhood of $a$ in $\partial C$. However, $\pi(F_{n})\in \pi({\cal MIN}\setminus U_{1})$ for all $n$ , which is a contradiction since $\{ \pi(F_{n})\}$ converges to $a$.    
\end{pf}

Next, we study properties of $\partial C$ as a Borel space. A topological space is called a {\it Polish space}\index{Polish space} if it is homeomorphic to a separable complete metric space. A Borel space is said to be {\it standard}\index{standard!Borel space} if it is Borel isomorphic to the Borel space of a Polish space generated by the topology. It is often useful to assume that Borel spaces are standard in order to avoid pathology. We collect some results about standard Borel spaces as follows:

\begin{thm}\label{thm-standard-borel-space}
\begin{enumerate}
\item[(i)] Any Borel subset of a standard Borel space is standard \cite[Corollary 13.4]{kechris}.
\item[(ii)] If $f\colon X\rightarrow Y$ is a Borel map between standard Borel spaces and the restriction of $f$ to a Borel subset $A\subseteq X$ is injective, then the set $f(A)$ is Borel in $Y$ and the restriction $f|_{A}$ from $A$ to $f(A)$ is a Borel isomorphism \cite[Corollary 15.2]{kechris}.
\item[(iii)] Let $f\colon X\rightarrow Y$ be a Borel map between standard Borel spaces such that the preimage $f^{-1}(y)$ is countable for any $y\in Y$. Then the image $f(X)$ is a Borel subset of $Y$ and there is a countable Borel partition $X=\bigsqcup_{n\in {\Bbb N}}X_{n}$ such that the restriction $f|_{X_{n}}$ is injective for each $n\in {\Bbb N}$ and $f(X_{0})=f(X)$. If there exists $N\in {\Bbb N}$ such that the preimage $f^{-1}(y)$ has at most $N$ elements for any $y\in Y$, then the partition can be chosen to have at most $N$ sets \cite[Theorem 1.3]{sauer}.
\item[(iv)] Let $f\colon X\rightarrow Y$ be a Borel map between standard Borel spaces such that the preimage $f^{-1}(y)$ is countable for any $y\in Y$. Then the image $f(A)$ of a Borel subset $A$ of $X$ is a Borel subset of $Y$.
\end{enumerate}
\end{thm} 

The assertion (iv) follows from (ii) and (iii). We shall return to our case.

\begin{prop}\label{boundary-standard}
Let $M$ be a non-exceptional surface. The map $\pi \colon {\cal MIN}\rightarrow \partial C$ has a Borel section, that is, there exists a Borel subset $S$ of ${\cal MIN}$ such that the restriction map $\pi \colon S\rightarrow \partial C$ is bijective. Moreover, it is a Borel isomorphism. In particular, $\partial C$ is a standard Borel space.
\end{prop} 
This proposition follows from the next two propositions:

\begin{prop}\label{prop-tak-similar}
Suppose that $X$ is a metric space and second countable. Let ${\cal R}$ be an equivalence relation on $X$ such that each equivalence class is compact and the ${\cal R}$-saturation of every closed subset of $X$ is measurable. Then there exists a Borel subset $S$ of $X$ such that each equivalence class of ${\cal R}$ intersects $S$ at exactly one point.
\end{prop}

\begin{pf}
This proof is almost the same as that for \cite[Theorem A15]{tak1}. 

Let $d$ be the metric on $X$. Let 
\[\{ X(n_{1}, n_{2}, \ldots, n_{k}): n_{i}\in {\Bbb N}, \ 1\leq i\leq k, \ k\in {\Bbb N}\}\]
be a system of open covering of $X$ such that 
\begin{equation*}
{\rm diam}(X(n_{1},\ldots, n_{k}))\leq 1/2^{k},
\end{equation*}
\begin{equation*}
X(n_{1},\ldots, n_{k})=\bigcup_{j=1}^{\infty}X(n_{1},\ldots, n_{k}, j),
\end{equation*}
\begin{equation*}
\overline{X(n_{1},\ldots, n_{k}, j)}\subseteq X(n_{1},\ldots, n_{k}), 
\end{equation*}
where ${\rm diam}(Z)$ denotes the diameter for a metric space $Z$. We consider the lexicographic ordering in each ${\Bbb N}^{k}$, $k=2, 3,\ldots$. For each subset $A$ of $X$, we denote by ${\cal R}A$ the saturation of $A$ under ${\cal R}$. Put $S(1)=\overline{X(1)}$ and 
\[S(n)=\overline{X(n)}\setminus \left( \bigcup_{j=1}^{n-1}{\cal R}\overline{X(j)}\right).\]
By assumption, $S(n)$ is a Borel subset (but may be empty). For each $k$, we set $S(\underbrace{1,\ldots, 1}_{k})=\overline{X(\underbrace{1,\ldots, 1}_{k})}$ and 
\begin{multline*}
S(n_{1},\ldots, n_{k})\\
=\overline{X(n_{1},\ldots, n_{k})}\setminus \bigcup\{ {\cal R}\overline{X(m_{1},\ldots, m_{k})}:(m_{1},\ldots, m_{k})<(n_{1},\ldots, n_{k})\}.
\end{multline*}
It follows that $S(n_{1},\ldots, n_{k})$ is a Borel set. Set
\begin{equation*}
S_{k}=\bigcup_{{\Bbb N}^{k}}S(n_{1},\ldots, n_{k}), \ \ \ S=\bigcap_{k=1}^{\infty}S_{k}.
\end{equation*}
Let $H$ be an equivalence class under ${\cal R}$. For each $k\in {\Bbb N}$, we have
\begin{equation*}
\emptyset \neq H\cap S_{k}=H\cap S(n_{1},\ldots, n_{k})=H\cap \overline{X(n_{1},\ldots, n_{k})}
\end{equation*}
for a unique $(n_{1},\ldots, n_{k})$. Hence, we see that $H\cap S_{k}$ is compact and that diam$(H\cap S_{k})\leq 1/2^{k}$ and $H\cap S_{k+1}\subseteq H\cap S_{k}$. It follows that $H\cap S=\bigcap_{k=1}^{\infty}(H\cap S_{k})$ consists of exactly a single point.  
\end{pf}

Considering the equivalence relation on ${\cal MIN}$ induced by the map $\pi \colon {\cal MIN}\rightarrow \partial C$ and applying Proposition \ref{prop-tak-similar}, we see that there exists a Borel subset $S$ of ${\cal MIN}$ which is a fundamental domain for the equivalence relation.

We say that a topological space is a {\it Souslin space}\index{Souslin!space} if it is a metrizable space which is a continuous image of a Polish space. A subset of a topological space is said to be a {\it Souslin set}\index{Souslin!set} if it is a Souslin space as a topological space.

\begin{prop}[\ci{Corollaries A.8, A.10}{tak1}]\label{Souslin}
\begin{enumerate}
\renewcommand{\labelenumi}{\rm(\roman{enumi})}
\item Every Borel subset of a Souslin space is a Souslin set.
\item Let $X$ be a Souslin space and $Y$ be a separable metrizable topological space. If $f$ is a Borel map of $X$ into $Y$, then $f(X)$ is a Souslin set in $Y$. If $f$ is injective in addition, then $f$ is a Borel isomorphism from $X$ onto $f(X)$.
\end{enumerate}
\end{prop}
It follows from this proposition that $S$ is a Souslin set of ${\cal PMF}$ and that the restriction $\pi \colon S\rightarrow \partial C$ is a Borel isomorphism, which shows Proposition \ref{boundary-standard}.

\begin{rem}\label{rem-kla-excep}
In Section \ref{section-excep-boundary}, we verify that if $M$ is a surface with $\kappa(M)=0$, then the boundary $\partial C$ of the curve complex $C$ is topologically identified with the set $\Hat{{\Bbb R}}\setminus \Hat{{\Bbb Q}}$ of irrational numbers on the circle $\Hat{{\Bbb R}}=S^{1}$ when $C$ is embedded in the Poincar\'e disk (see Chapter \ref{chapter-property-A}, Section \ref{section-exceptional}).

Under the identification between $\Hat{{\Bbb Q}}$ and $V(C)$ in Chapter \ref{chapter-property-A}, Section \ref{section-exceptional}, if $M$ is of type $(1, 1)$, then we have an isomorphism $p_{M}\colon SL(2, {\Bbb Z})\rightarrow \Gamma(M)$. If $M$ is of type $(0, 4)$, then we can construct an injective homomorphism $p_{M}\colon PSL(2, {\Bbb Z})\rightarrow \Gamma(M)$ and the image of this homomorphism is the subgroup of $\Gamma(M)$ consisting of elements preserving one fixed boundary component of $M$ \cite[Theorem 7.9]{ivanov0}. The above topological identification between $\Hat{{\Bbb R}}\setminus \Hat{{\Bbb Q}}$ and $\partial C$ can be taken as to be equivariant under $p_{M}$.  

Moreover, it is well-known that the action of (the image of $p_{M}$ in) $\Gamma(M)$ on the Thurston boundary can be topologically identified with the standard action of $(P)SL(2, {\Bbb Z})$ on the circle $\Hat{{\Bbb R}}$ under $p_{M}$ and that $\Hat{{\Bbb R}}\setminus \Hat{{\Bbb Q}}$ is identified with the set of the projections of all minimal measured foliations. This can be verified as well as in the case of the torus in \cite[Expos\'e 1, \S VI]{FLP}.

Hence, it follows that Theorem \ref{thm-kla-boundary} and Propositions \ref{prop-kla-pi}, \ref{boundary-standard} hold also for a surface $M$ with $\kappa(M)=0$ if we define $\pi$ as the identity map on $\Hat{{\Bbb R}}\setminus \Hat{{\Bbb Q}}$. 
\end{rem}



\section{Amenability for the actions of the mapping class group}\label{section-amenable-action}

Let $X$ be a $\delta$-hyperbolic graph with $\partial X\neq \emptyset$ satisfying the properties (F1) and (F2) in Chapter \ref{chapter-property-A}, Section \ref{sec:Acc}. We identify $X$ with the vertex set of the graph $X$. We recall the functions constructed in the proof of Theorem \ref{main1}. For $x\in X$, $a\in \partial X$ and $n, k\in {\Bbb N}$ satisfying $n-2k-\delta_{0}-\delta_{1}>0$, we have defined a function $F_{a}(x, k, n)$ on $X$ by the characteristic function on the set  
\[\bigcup_{f\in {\cal L}_{T}(x, a; k)}f([n, 2n]),\]
and defined
\[H_{a}(x, n)=n^{-3/2}\sum_{k<\sqrt{n}}F_{a}(x, k, n)\]
for $n\in {\Bbb N}$ such that $n\geq N_{0}$, where $N_{0}$ is a certain fixed natural number. We must verify the measurability of the two functions for $a\in \partial X$. 

\begin{prop}\label{measurable}
The function $\partial X \ni a\mapsto F_{a}(x, k, n)(y)$ is the pointwise limit of an increasing sequence of upper semi-continuous functions on $\partial X$ for fixed $x, y\in X$ and $k, n\in {\Bbb N}$ for which we define $F$. In particular, it is measurable.
\end{prop}

For the proof, we need some lemmas.

\begin{lem}\label{hyp-basic1}
Let $X$ be a $\delta$-hyperbolic graph. Then there exists a constant $\delta'\geq 0$ depending only on $\delta$ such that any geodesic triangle whose vertices in $X\cup \partial X$ is $\delta'$-slim.
\end{lem}

\begin{pf}
We show that any geodesic triangle $abc$ with $a, b\in X$, $c\in \partial X$ is $18\delta$-slim. Denote $f=[a, c]$ and $g=[b, c]$. By Lemma \ref{hyp-fund}, there exist $n$, $m\in {\Bbb N}$ such that 
\[d(f(n+k), g(m+k))\leq 16\delta \]
for each $k\in {\Bbb N}$. Choose geodesic segments $[f(n), b]$ and $[f(n), g(m)]$. It follows from $\delta$-slimness of the triangles $f(n)g(m)b$, $abf(n)$ that $abc$ is $18\delta$-slim. 

In the other cases, we can show similarly that there exists some constant $\delta'$ depending only on $\delta$ such that any triangle whose vertices in $X\cup \partial X$ is $\delta'$-slim. 
\end{pf}

\begin{lem}\label{hyp-basic2}
Let $X$ be a $\delta$-hyperbolic graph with $\partial X\neq \emptyset$ satisfying the properties ${\rm (F1)}$ and ${\rm (F2)}$ and $\delta'\geq 0$ be the constant given in Lemma \ref{hyp-basic1}. Let $a$, $b$ be any two distinct points in $\partial X$ and $[a, b]$ be a bi-infinite geodesic from $a$ to $b$. Let $x\in X$. Then
\[(a|b)_{x}\leq d(x, [a, b])+2\delta'.\]
\end{lem} 

\begin{pf}
It suffices to show that for any $y\in [x, a]$ and $z\in [x, b]$, we have
\[(y|z)_{x}\leq d(x, [a, b])+2\delta'.\]
There exists $t\in [a, b]$ such that $d(x, t)=d(x, [a, b])$. Choose a geodesic ray $[y, b]$ and consider the triangle $yba$. By $\delta'$-slimness, there exists a point $t'\in [y, b]\cup [y, a]$ such that $d(t, t')\leq \delta'$. 

If $t'\in [y, a]$, then
\begin{align*}
(y|z)_{x}&\leq d(x, [y, z])\leq d(x, y)\\
         &\leq d(x, t')\leq d(x, t)+d(t, t')\\
         &\leq d(x, [a, b])+\delta'.
\end{align*}
If $t'\in [y, b]$, then choose a geodesic $[y, z]$. By $\delta'$-slimness of the triangle $yzb$, we have a point $t''\in [y, z]\cup [z, b]$ such that $d(t', t'')\leq \delta'$. 

If $t''\in [z, b]$, then 
\begin{align*}
(y|z)_{x}&\leq d(x, [y, z])\leq d(x, z)\leq d(x, t'')\\
         &\leq d(x, t)+d(t, t')+d(t', t'')\\
         &\leq d(x, [a, b])+2\delta'.
\end{align*}

If $t''\in [y, z]$, then
\begin{align*}
(y|z)_{x}&\leq d(x, [y, z])\leq d(x, t'')\\
         &\leq d(x, t)+d(t, t')+d(t', t'')\\
         &\leq d(x, [a, b])+2\delta'.
\end{align*}
This completes the proof. 
\end{pf}

\begin{lem}\label{hyp-basic3}
Let $X$ be as in Lemma \ref{hyp-basic2}. Let $x\in X$ and $a, b\in \partial X$ be two distinct points. Let $g_{x, a}\in {\cal L}(x, a)$, $g_{x, b}\in {\cal L}(x, b)$ and $t_{0}\in {\Bbb N}$. If $g_{x, b}\cap B(g_{x, a}(t_{0}); \delta'+1)=\emptyset$, then $(a|b)_{x}\leq t_{0}+3\delta'$.
\end{lem}

\begin{pf}
By $\delta'$-slimness of the triangle $xab$, there exists $y\in [a, b]$ such that $d(g_{x, a}(t_{0}), y)\leq \delta'$. Thus, we have
\[(a|b)_{x}\leq d(x, [a, b])+2\delta' \leq t_{0}+\delta' +2\delta'=t_{0}+3\delta'.\]
It completes the proof. 
\end{pf}

For $x\in X$, $a\in \partial X$ and $k\in {\Bbb N}$, we denote
\[G(x, a)_{k}=G(x, a)\cap \{ y\in X: d(x, y)=k\}=\bigcup_{f\in {\cal L}_{T}(x, a)}f(k).\]
\index{$G x a k$@$G(x, a)_{k}$}

\begin{lem}\label{upper-semi-cont}
For fixed $r\in {\Bbb N}$, $a\in \partial X$ and $x\in X$, there exists a neighborhood $V$ in $\partial X$ of $a$ satisfying the following condition: for any $b\in V$, $g\in {\cal L}_{T}(x, b)$ and $k\in {\Bbb N}$ with $k\leq r$, we have $g(k)\in  G(x, a)_{k}$.  
\end{lem}

\begin{pf}
This can be verified similarly to \cite[Appendix B, Lemma 3.3]{ana}, using the properties (F1) and (F2). 

Suppose that the lemma is not true. There would exist $a_{n}\in \partial C$, $g_{x, a_{n}}\in {\cal L}_{T}(x, a_{n})$ and $0\leq k\leq r$ such that $a_{n}\rightarrow a$ with $(a|a_{n})_{x}\nearrow \infty$ and 
\[g_{x, a_{n}}(k)\notin G(x, a)_{k}\]
for any $n$ (remark that $\partial X$ is metrizable). Choosing a subsequence, we may assume that
\[(a|a_{n})_{x}> n+3\delta'\]
for any $n\in {\Bbb N}$. Since $a_{n}\neq a$, it follows from Lemma \ref{hyp-basic3} that
\[g_{x, a_{n}}\cap B(f_{x, a}(m); \delta'+1)\neq \emptyset \]
for any $0\leq m\leq n$, where $f_{x, a}\in {\cal L}_{T}(x, a)$ is a fixed geodesic ray. 

We will show that there exists a subsequence of $\{ g_{x, a_{n}}\}$ converging to some geodesic ray in ${\cal L}_{T}(x, a)$ uniformly on any bounded set.

Let $n_{0}\in {\Bbb N}$ be any integer such that $n_{0}>2(\delta'+1)+\delta_{1}$, where $\delta_{1}$ is the constant appearing in the property (F2). Then for any $n\geq 2n_{0}$, there exists a point $y\in g_{x, a_{n}}\cap B(f_{x, a}(2n_{0}); \delta'+1)$ and we see that $d(g_{x, a_{n}}(2n_{0}), y)\leq \delta'+1$ because 
\[d(g_{x, a_{n}}(2n_{0}), y)=|d(g_{x, a_{n}}(2n_{0}), x)-d(y, x)|=|2n_{0}-d(x, y)|\] 
and
\[2n_{0}-(\delta'+1)\leq d(x, y)\leq 2n_{0}+\delta'+1\] 
by considering the triangle $xf_{x, a}(2n_{0})y$. Moreover, we have
\begin{align*}
d(g_{x, a_{n}}(2n_{0}), f_{x, a}(2n_{0}))& \leq d(g_{x, a_{n}}(2n_{0}), y)+d(y, f_{x, a}(2n_{0}))\\
     & \leq 2(\delta'+1).
\end{align*} 
Thus, $g_{x, a_{n}}|_{[0, 2n_{0}]}\in {\cal L}_{T}(x, f_{x, a}(2n_{0}); 2(\delta'+1))$. It follows from Lemma \ref{lem:bow} that the geodesic $g_{x, a_{n}}|_{[0, 2n_{0}]}$ passes through the ball $B(f_{x, a}(n_{0}); \delta_{0})$ since $n_{0}>2(\delta'+1)$. Taking a subsequence, we may assume that all $g_{x, a_{n}}$ coincide on $B(f_{x, a}(n_{0}); \delta_{0})$ by (F2). Applying (F1) for $x$ and $g_{x, a_{n}}(n_{0})$ (which is the same point for all $n\in {\Bbb N}$), we can choose a subsequence $\{ g_{x, a_{n}}\}$ such that all $g_{x, a_{n}}(t)$ coincide for any $0\leq t\leq n_{0}$.

This operation can be applied for any large $n_{0}\in {\Bbb N}$. Thus, we can show the existence of a desired subsequence. The limit geodesic $g_{\infty}$ is in ${\cal L}_{T}(x, a)$. On the other hand, since $g_{x, a_{n}}(k)\notin G(x, a)_{k}$ for any $n$, we conclude $g_{\infty}(k)\notin G(x, a)_{k}$. It is a contradiction. 
\end{pf}

\begin{cor}\label{measurable-cor}
\begin{enumerate}
\renewcommand{\labelenumi}{\rm(\roman{enumi})}
\item For fixed $x, y\in X$, the map $\partial X\ni a\mapsto \chi_{G(x, a)}(y)$ is upper semi-continuous, where $\chi_{A}$ is the characteristic function on a subset $A$ of $X$.
\item For fixed $x, y\in X$ and $k\in {\Bbb N}$, the map $\partial X\ni a\mapsto \chi_{G(x, a)_{k}}(y)$ is upper semi-continuous.
\end{enumerate}
\end{cor}

For the proof of Proposition \ref{measurable}, it is enough to remark that the support of the function $F_{a}(x, k, n)$ is equal to the set
\[\bigcup_{y\in B(x; k)}\bigcup_{m=n}^{2n}G(y, a)_{m}.\] 

In the rest of this section, we will show the following theorem:

\begin{thm}\label{amenable-action-cc-non-excep}
Let $M$ be a compact orientable surface with $\kappa(M)>0$ and $C$ and $\Gamma$ be the curve complex and the mapping class group of $M$, respectively. 
Let $\mu$ be a quasi-invariant probability measure on $\partial C$. Then the action of $\Gamma$ on $(\partial C, \mu)$ is amenable.
\end{thm}

We refer the reader to Appendix \ref{general-amenable} for the definition of amenability for an action of a discrete group. 
\vspace{1em}
 
{\sc Proof}. Recall the map $H$ stated in the beginning of this section, which is defined on the set $\partial C\times C\times \{ n\geq N_{0}\}$ and valued in the space of non-zero, non-negative valued functions in $\ell^{1}(C)$ with finite norm. It follows from Proposition \ref{measurable} that the map $H$ is measurable. By normalizing, we may assume that $H$ satisfies the following conditions (see the proof of Theorem \ref{main1} and Remark \ref{constructed-function-action}):

\begin{enumerate}
\renewcommand{\labelenumi}{\rm(\roman{enumi})}
\item for any $R>0$, we have $\Vert H_{a}(x, n)-H_{a}(y, n)\Vert_{1}\rightarrow 0$ as $n\rightarrow \infty$ uniformly on the set $\{ (a, x, y)\in \partial C\times C\times C: d(x, y)<R\}$;
\item the $\ell^{1}$-norm of $H_{a}(x, n)$ is 1 for any $a\in \partial C$, $x\in X$ and $n\geq N_{0}$;
\item $H_{ga}(gx, n)=g\cdot H_{a}(x, n)$ for $g\in \Gamma$, $a\in \partial C$, $x\in C$ and $n\geq N_{0}$, where the action of $\Gamma$ on the space of functions $\varphi$ on $C$ is defined by the formula $(g\cdot \varphi)(x)=\varphi(g^{-1}x)$ for $x\in C$ and $g\in \Gamma$.
\end{enumerate}
Define a non-negative valued Borel function $f_{n}$ on $\partial C\times C$ by
\[f_{n}(a, x)=H_{a}(x_{0}, n)(x),\]
where $x_{0}\in C$ is a fixed point. Then
\[\sum_{x\in C}f_{n}(a, x)=1,\]
and since $f_{n}(g^{-1}a, g^{-1}x)=H_{g^{-1}a}(x_{0}, n)(g^{-1}x)=H_{a}(gx_{0}, n)(x)$ for any $g\in \Gamma$, it follows from (i) that
\[\lim_{n\rightarrow \infty}\sup_{a\in \partial C}\sum_{x\in C}\vert f_{n}(a, x)-f_{n}(g^{-1}a, g^{-1}x)\vert =0.\]  
Thus, the sequence $\{ f_{n}\}$ satisfies 
\[\lim_{n\rightarrow \infty}\sum_{x\in C}\int_{\partial C}h(a)\varphi(a, x)(f_{n}(g^{-1}a, g^{-1}x)-f_{n}(a, x))d\mu(a)=0\]
for any $h\in L^{1}(\partial C)$, $\varphi \in L^{\infty}(\partial C\times C)$ and $g\in \Gamma$. This means that $\{ f_{n}\}$ satisfies the condition in Definition \ref{app-inv-mean}, that is, $\{ f_{n}\}$ is an approximate weakly invariant mean for the $(\partial C\rtimes \Gamma)$-space $\partial C\times C$.

It follows from Theorem \ref{ame-eq-cond} that there exists an invariant measurable system of means  
\[m\colon \partial C\rightarrow S(\ell^{\infty}(C)),\]
for the $(\partial C\rtimes \Gamma)$-space $\partial C\times C$, where $S$ indicates the set of states.

Since the action of $\Gamma$ on $C$ is cocompact \cite{harvey}, there exists a finite family $\{ \alpha\}$ of elements in $C$ such that $\{ \Gamma \alpha \}$ is a complete set of $\Gamma$-orbits on $C$. If $\chi_{\Gamma \alpha}$ denotes the characteristic function on the orbit $\Gamma \alpha$, then
\[1=m_{a}(1)=m_{a}(\sum_{\alpha}\chi_{\Gamma \alpha})=\sum_{\alpha}m_{a}(\chi_{\Gamma \alpha}),\]
for a.e. $a\in \partial C$. Let
\[X_{\alpha}=\{ a\in \partial C: m_{a}(\chi_{\Gamma \alpha})\geq 1/K \},\]
where $K$ denotes the cardinality of the set $\{ \alpha \}$. Then $X_{\alpha}$ is $\Gamma$-invariant, and $\partial C=\bigcup X_{\alpha}$. 

For the proof of the amenability for the action of $\Gamma$ on $(\partial C, \mu)$, it suffices to show that for each $X_{\alpha}$ with $\mu(X_{\alpha})>0$, the action of $\Gamma$ on $(X_{\alpha}, \mu)$ is amenable. By normalizing the mean $m_{a}$ restricted to $\Gamma \alpha$, we get an invariant measurable system of means 
\[m^{\alpha}\colon X_{\alpha}\rightarrow S(\ell^{\infty}(\Gamma \alpha))\]
for the $(X_{\alpha}\rtimes \Gamma)$-space $X_{\alpha}\times (\Gamma \alpha)$. The orbit $\Gamma \alpha$ can be identified with the $\Gamma$-space $\Gamma /\Gamma_{\alpha}$, where $\Gamma_{\alpha}$ is the stabilizer of $\alpha$. Thus, the existence of the invariant system $m^{\alpha}$ of means implies that the amenability of the $(X_{\alpha}\rtimes \Gamma)$-space $X_{\alpha}\times (\Gamma /\Gamma_{\alpha})$. By Proposition \ref{rel-ame}, it is enough to show that the action of $\Gamma_{\alpha}$ on $(X_{\alpha}, \mu)$ is amenable. This follows from the next lemma:

\begin{lem}\label{lem-amenable-action-stabilizer}
Let $\alpha$ be a vertex in the complex $C$. We denote the $\Gamma$-stabilizer of $\alpha$ by $\Gamma_{\alpha}$. Then the action of $\Gamma_{\alpha}$ on $(\partial C, \mu)$ is amenable for any quasi-invariant measure $\mu$ for the action of $\Gamma_{\alpha}$. 
\end{lem}

\begin{pf}
Denote
\[S= \{ x\in C: d(x, \alpha)=3\},\]
which is a $\Gamma_{\alpha}$-invariant set. For a set $T$, let ${\cal F}(T)$ be the set of all non-empty finite subsets in $T$. Define $F\colon \partial C\rightarrow {\cal F}(S)$ by
\[F(a)=S\cap G(\alpha, a)=G(\alpha, a)_{3}\]
for $a\in \partial C$. The right hand side is non-empty. Let $z$ be a point in $G(\alpha, a)_{3}$ for $a\in \partial C$ and $f\in {\cal L}_{T}(\alpha, a)$ be a geodesic passing through $z$. Since we have
\[G(\alpha, a)\subseteq \bigcup_{n\in {\Bbb N}}B(f(n);\delta_{0})\]
by Lemma \ref{lem:bow2}, it follows from (F1) that the set $G(\alpha, a)_{3}$ is finite and its cardinality is uniformly bounded for $a\in \partial C$. Moreover, we have
\[F(ga)=gF(a)\]
for any $g\in \Gamma_{\alpha}$ and $a\in \partial C$. We give the set ${\cal F}(S)$ the discrete Borel structure. Then for any $x\in C$, the map $\phi_{x}\colon \partial C\rightarrow \{ 0, 1\}$ defined by the formula $\phi_{x}(a)=\chi_{F(a)}(x)$ is measurable by Corollary \ref{measurable-cor} (i). Moreover, for any $A\in {\cal F}(S)$, we have
\[F^{-1}(A)=\left( \bigcap_{x\in A}\phi_{x}^{-1}(1)\right) \cap \left( \bigcap_{x\in S\setminus A}\phi_{x}^{-1}(0)\right).\]
This means that the map $F$ is measurable. 

There exists a subset $\{ \beta \}$ of $S$ such that $S$ can be identified with the disjoint union
\[\bigsqcup_{\beta}(\Gamma_{\alpha}/\Gamma_{\alpha, \beta})\]
as $\Gamma_{\alpha}$-spaces, where 
\[\Gamma_{\alpha, \beta}=\{ g\in \Gamma: g\alpha =\alpha, \ g\beta =\beta \}.\]
Define
\[X_{\alpha, \beta}=\{ a\in \partial C: F(a)\cap (\Gamma_{\alpha}/\Gamma_{\alpha, \beta})\neq \emptyset \}.\]
It is $\Gamma_{\alpha}$-invariant and a Borel subset of $\partial C$ since $F$ is measurable. Then we have 
\[\partial C=\bigcup_{\beta}X_{\alpha, \beta}.\]      
It suffices to show that the action of $\Gamma_{\alpha}$ on $(X_{\alpha, \beta}, \mu)$ is amenable for any $\beta$ with $\mu(X_{\alpha, \beta})>0$. 

Let $\beta$ be such an element. Define $F'\colon X_{\alpha, \beta}\rightarrow {\cal F}(\Gamma_{\alpha}/\Gamma_{\alpha, \beta})$ by 
\[F'(a)=F(a)\cap (\Gamma_{\alpha}/\Gamma_{\alpha, \beta}).\]
The right hand side is a non-empty finite set and its cardinality is uniformly bounded for $a\in X_{\alpha, \beta}$. Moreover, we have
\[F'(ga)=gF'(a)\]
for any $g\in \Gamma_{\alpha}$ and $a\in X_{\alpha, \beta}$. The map $F'$ is also measurable when we give the set ${\cal F}(\Gamma_{\alpha}/\Gamma_{\alpha, \beta})$ the discrete Borel structure.

We can construct an invariant mean for the $(X_{\alpha, \beta}\rtimes \Gamma_{\alpha})$-space $X_{\alpha, \beta}\times (\Gamma_{\alpha}/\Gamma_{\alpha, \beta})$ from $F'$ easily as follows: we denote $X'=X_{\alpha, \beta}$, $Z=\Gamma_{\alpha}/\Gamma_{\alpha, \beta}$ and define a map $f$ from $X'$ into the space of non-negative valued functions in $\ell^{1}(Z)$ with norm $1$ by the formula
\[(f(a))(z)=(\chi_{F'(a)}(z))/(\sum_{z'\in Z}\chi_{F'(a)}(z'))\]
for $a\in X'$ and $z\in Z$. The function on $X'$ defined by the formula $a\mapsto (f(a))(z)$ is measurable for any $z\in Z$ since $F'$ is measurable. Recall that $L^{\infty}(X', \ell^{1}(Z))$ is the space of measurable function $\varphi$ on $X'\times Z$ such that the essential supremum of $\sum_{z\in Z}|\varphi(a, z)|$ for $a\in X'$ is finite and whose norm is given by this quantity (see Appendix \ref{general-amenable}, Section \ref{section-amenability}). Any non-negative function $\varphi \in L^{\infty}(X', \ell^{1}(Z))$ is the pointwise limit of an increasing sequence of sums of functions 
\[X'\times Z\ni (a, z)\mapsto t\chi_{Y\times \{ z'\}}(a, z), \]
where $Y$ is a Borel subset of $X'$ and $t\geq 0$, $z'\in Z$ and $\chi_{Y\times \{ z'\}}$ denotes the characteristic function on $Y\times \{ z'\}$. It follows that the function $X'\ni a\mapsto \sum_{z\in Z}(f(a))(z)\varphi(a, z)$ is measurable. By the invariance of $F'$, we see that
\[(f(g^{-1}a))(z)=(f(a))(gz)\]
for any $a\in X'$, $z\in Z$ and $g\in \Gamma_{\alpha}$. For $a\in X'$, define a mean $m_{a}$ on $\ell^{\infty}(Z)$ by 
\[m_{a}(\psi)=\sum_{z\in Z}(f(a))(z)\psi(z)\]
for $\psi \in \ell^{\infty}(Z)$. Then for $a\in X'$, $g\in \Gamma_{\alpha}$ and $\varphi \in L^{\infty}(X', \ell^{1}(Z))$, we obtain the equality
\begin{align*}
m_{g^{-1}a}(\varphi(g^{-1}a, \cdot))&=\sum_{z\in Z}(f(g^{-1}a))(z)\varphi(g^{-1}a, z)\\
      &=\sum_{z\in Z}(f(a))(gz)\varphi(g^{-1}a, z)\\
      &=\sum_{z\in Z}(f(a))(z)\varphi(g^{-1}a, g^{-1}z)\\
      &=\sum_{z\in Z}(f(a))(z)(g\cdot \varphi)(a, z)\\
      &=m_{a}((g\cdot \varphi)(a, \cdot)).
\end{align*}    
It means that $\{ m_{a}\}_{a\in X'}$ is an invariant measurable system of means for the $(X'\rtimes \Gamma_{\alpha})$-space $X'\times Z$. 

On the other hand, it follows from Lemma \ref{fill-curves-fix} that $\Gamma_{\alpha, \beta}$ is finite since $d(\alpha, \beta)=3$. Thus, the action of  $\Gamma_{\alpha}$ on $(X', \mu)$ is amenable by Proposition \ref{rel-ame}. 
\end{pf}

\begin{rem}\label{rem-normal-typeI}
Let $\alpha \in V(C)$ and $\mu$ be a quasi-invariant measure on $\partial C$ for the $\Gamma_{\alpha}$-action. It follows from the above proof that there exists a non-negative Borel function $k \in L^{\infty}(\partial C, \ell^{1}(\Gamma_{\alpha}))$ such that 
\[\sum_{g\in \Gamma_{\alpha}}k(a, g)=1\]
for a.e.\ $a\in \partial C$ and the family $\{ m_{a}\}_{a\in \partial C}$ of means defined by the formula
\[m_{a}(\psi)=\sum_{g\in \Gamma_{\alpha}}k(a, g)\psi(g)\]
for $a\in \partial C$ and $\psi \in \ell^{\infty}(\Gamma_{\alpha})$ is an invariant measurable system of means for the $\Gamma_{\alpha}$-action on $(\partial C, \mu)$. By restriction and normalization, we may assume that the function $\Gamma_{\alpha}\ni g \mapsto k(a, g)$ has a finite support for a.e.\ $a\in \partial C$. Thus, we have a $\Gamma_{\alpha}$-equivariant Borel map $\phi: \partial C\rightarrow {\cal F}(\Gamma_{\alpha})$, where ${\cal F}(\Gamma_{\alpha})$ denotes the set of all non-empty finite subsets of $\Gamma_{\alpha}$ and is equipped with the discrete Borel structure. Let $F'$ be a set of all representatives for the $\Gamma_{\alpha}$-orbits on ${\cal F}(\Gamma_{\alpha})$. Then the Borel subset $F=\phi^{-1}(F')$ of $\partial C$ satisfies $\partial C=\Gamma_{\alpha}F$ up to null sets and for any distinct points $a, b\in F$, we have $\Gamma_{\alpha}a\neq \Gamma_{\alpha}b$. Such a Borel subset $F$ is called a {\it fundamental domain}\index{fundamental domain} for the $\Gamma_{\alpha}$-action on $(\partial C, \mu)$.  

In \cite[Section 7]{mc-pa}, the dynamical behavior for the action of an infinite reducible subgroup $G$ on ${\cal PMF}$ is studied. Note that $G$ acts properly discontinuously on a certain open subset of ${\cal PMF}$ containing ${\cal MIN}$ \cite[Theorem 7.17]{mc-pa}.   
\end{rem}

Recall that we have the continuous $\Gamma$-equivariant map $\pi$ from ${\cal MIN}$ to $\partial C$ (see Section \ref{boundary-of-curve-complex}). By Theorem \ref{aeg-thm}, we get the following corollary:

\begin{cor}\label{amenable-action-non-excep}
Let $M$ be a compact orientable surface with $\kappa(M)>0$ and let $\mu$ be a quasi-invariant probability measure on ${\cal PMF}$ whose support is contained in ${\cal MIN}$. Then the action of the mapping class group $\Gamma$ on $({\cal PMF}, \mu)$ is amenable.
\end{cor}

\begin{rem}\label{masur-measure}
Masur \cite{masur} and Rees \cite{rees} constructed a natural finite measure $\mu$ on ${\cal PMF}$ satisfying the following conditions:
\begin{enumerate}
\item[(i)] The support of $\mu$ is contained in ${\cal UE}$.
\item[(ii)] The action of the mapping class group $\Gamma$ on $({\cal PMF}, \mu)$ is non-singular and ergodic.
\end{enumerate}
(Masur's measure satisfies more notable properties.) It is natural to ask the type of the measure. For some measures on the boundary of free groups and fundamental groups of negatively curved manifolds, their types have been calculated in \cite{oka} and \cite{spa}, respectively. 
\end{rem}



\section{The boundary of the curve complex for an exceptional surface}\label{section-excep-boundary} 

In this section, as briefly mentioned in \cite[Section 5]{minsky}, we identify the boundary at infinity of the curve complex $C$ for an exceptional surface $M$ with the set $\Hat{{\Bbb R}}\setminus \Hat{{\Bbb Q}}$ of irrational numbers on the circle $\Hat{{\Bbb R}}=S^{1}$ when $C$ is embedded in the Poincar\'e disk as described in Chapter \ref{chapter-property-A}, Section \ref{section-exceptional}. Under this embedding, the vertex set $V(C)$ is identified with the set $\Hat{{\Bbb Q}}$ of rational numbers. If $M$ is of type $(1, 1)$, then we denote by $p_{M}$ the isomorphism $SL(2, {\Bbb Z})\rightarrow \Gamma(M)$ and if $M$ is type $(0, 4)$, then we denote by $p_{M}$ the injective homomorphism $PSL(2, {\Bbb Z})\rightarrow \Gamma(M)$ in Remark \ref{rem-kla-excep}.

\begin{thm}\label{thm-excep-identified}
There exists a homeomorphism $\varphi$ between the set $\Hat{{\Bbb R}}\setminus \Hat{{\Bbb Q}}$ of irrational numbers on the circle $\Hat{{\Bbb R}}=S^{1}$ and the boundary $\partial C$ at infinity of the curve complex $C$ for a surface $M$ with $\kappa(M)=0$. Moreover, $\varphi$ can be taken to be equivariant under the homomorphism $p_{M}$ when $V(C)$ is identified with $\Hat{{\Bbb Q}}$ as described in Chapter \ref{chapter-property-A}, Section \ref{section-exceptional}.
\end{thm}

\begin{lem}\label{lem-number-theory}
For any irrational number $a$ and positive number $\varepsilon$, there exist two rational numbers $p/q$, $r/s$ such that $p/q<a<r/s$, $r/s-p/q<\varepsilon$ and $|ps-qr|=1$. 
\end{lem}

For the proof, we recommend the reader to consult any standard textbook on the continued fraction of a real number. For example, see \cite[Chapter 1]{lang}.

Remark that any edge $e$ in $C$ separates the closure $\overline{D}$ of the Poincar\'e disk $D$ into two components, that is, $\overline{D}\setminus \overline{e}$ consists of two components. For any two subsets $A$, $B$ of the closed disk, we denote by $E(A, B)$\index{$E A B $@$E(A, B)$} the set of all edges $e$ in $C$ such that one component of $\overline{D}\setminus \overline{e}$ contains $A$ and the other contains $B$. This notation is compatible with that in Chapter \ref{chapter-property-A}, Section \ref{section-exceptional}. 

It follows from Lemma \ref{lem-number-theory} that for any irrational number $a$ and edge $e$ in $C$, the set $E(e, a)$ is non-empty. Define $\varphi \colon \Hat{{\Bbb R}}\setminus \Hat{{\Bbb Q}}\rightarrow \partial C$ as follows: let $a\in \Hat{{\Bbb R}}\setminus \Hat{{\Bbb Q}}$. Then there exists a sequence $\{ e_{n}\}_{n\in {\Bbb N}}$ of edges in $C$ such that $e_{n}\in E(e_{n-1}, a)$ for each $n$ and that both vertices of the edge $e_{n}$ converge to $a$ on the circle. Let $x_{n}\in V(C)$ be any vertex of the edge $e_{n}$. Define $\varphi(a)$ as the element in $\partial C$ determined by the sequence $\{ x_{n}\}_{n\in {\Bbb N}}$. We need the following lemma to show that the map $\varphi$ is well-defined:

\begin{lem}\label{lem-gromov-product}
Fix $0\in \Hat{{\Bbb Q}}=V(C)$. 
\begin{enumerate}
\item[(i)] Let $x, y\in V(C)$ and $e\in E(0, \{ x, y\})$. Then we have $(x|y)_{0}\geq d(0, e)-2$.
\item[(ii)] Let $x, y\in V(C)$ and $e\in E(\{ 0, x\}, y)$. Then we have $(x|y)_{0}\leq d(0, e)+1$.
\end{enumerate}
\end{lem}  

\begin{pf}
Let $x, y\in V(C)$ and $e\in E(0, \{ x, y\})$. Let $z$ be a vertex of $e$ such that $d(0, e)=d(0, z)$. Using the fact that any geodesic from $0$ to $x$ (resp. $y$) passes through one of vertices of $e$, we see that
\begin{align*}
2(x|y)_{0}&=d(x, 0)+d(y, 0)-d(x, y)\\
          &\geq d(0, z)+d(z, x)-2+d(0, z)+d(z, y)-2-d(x, y)\\
          &\geq 2d(0, e)-4,
\end{align*}
and that the assertion (i) follows. 

Let $x, y\in V(C)$ and $e\in E(\{ 0, x\}, y)$. Let $z$ be a vertex of $e$ such that $d(0, e)=d(0, z)$. Similarly, we have
\begin{align*}
2(x|y)_{0}&=d(x, 0)+d(y, 0)-d(x, y)\\
          &\leq d(x, 0)+d(0, z)+d(z, y)-(d(x, z)+d(z, y)-2)\\
          &\leq 2d(0, z)+2=2d(0, e)+2,
\end{align*}
which shows the assertion (ii).
\end{pf}

\begin{lem}
The map $\varphi \colon \Hat{{\Bbb R}}\setminus \Hat{{\Bbb Q}}\rightarrow \partial C$ is well-defined, that is, for any irrational number $a$, the sequence $\{ x_{n}\}_{n\in {\Bbb N}}$ taken as above converges to infinity and the element in $\partial C$ determined by $\{ x_{n}\}_{n\in {\Bbb N}}$ depends only on the irrational number $a$.  
\end{lem}

\begin{pf}
Let $0\in \Hat{{\Bbb Q}}=V(C)$ be a fixed point. Let $a\in \Hat{{\Bbb R}}\setminus \Hat{{\Bbb Q}}$ and $\{ e_{n}\}_{n\in {\Bbb N}}$ be a sequence of edges in $C$ such that $e_{n}\in E(e_{n-1}, a)$ for each $n$ and that both vertices of the edge $e_{n}$ converge to $a$ on the circle $S^{1}$. Let $x_{n}\in V(C)$ be any vertex of the edge $e_{n}$. We may assume that $e_{n}\in E(0, a)$ for any $n$. Since $e_{n}\in E(e_{n-1}, a)$ and $e_{n}\in E(0, a)$, we have $d(0, e_{n})\geq n$ for each $n$. It follows that $d(0, e_{n})\rightarrow \infty$ as $n\rightarrow \infty$ and that the sequence $\{ x_{n}\}$ converges to infinity by Lemma \ref{lem-gromov-product} (i).

Let $\{ e_{m}'\}$ and $\{ x_{m}'\}$ be another choice of sequences for $a$. For each $n$, there exists $K$ such that $e_{m}'\in E(e_{n}, a)$ for any $m\geq K$ because both vertices of $e_{m}'$ converge to $a$ on the circle. Then $(x_{n}|x_{m}')_{0}\geq d(0, e_{n})-2$ for any $m\geq K$ by Lemma \ref{lem-gromov-product} (i), which means that two sequences $\{ e_{n}\}$ and $\{ e_{m}'\}$ are equivalent.  
\end{pf}

\begin{lem}\label{lem-converge-varphi}
Let $a\in \Hat{{\Bbb R}}\setminus \Hat{{\Bbb Q}}$ and $\{ x_{n}\}_{n\in {\Bbb N}}$ be a sequence in $V(C)$ converging to $\varphi(a)\in \partial C$. Then the sequence $\{ x_{n}\}_{n\in {\Bbb N}}$ in $\Hat{{\Bbb Q}}$ converges to $a$ on the circle.
\end{lem}

\begin{pf}
Let $\{ e_{n}\}_{n\in {\Bbb N}}$ be a sequence of edges in $C$ such that $e_{n}\in E(e_{n-1}, a)$ for each $n$ and that both vertices of the edge $e_{n}$ converge to $a$ on the circle. Let $y_{n}\in V(C)$ be any vertex of $e_{n}$. If the lemma were not true, then there would exist $N$ and a subsequence $\{ n_{i}\}_{i}$ of ${\Bbb N}$ such that $e_{N}\in E(\{ 0, x_{n_{i}}\}, a)$ for each $i$. Then we have $(x_{n_{i}}|y_{m})_{0}\leq d(0, e_{N})+1$ for each $i$ and $m\geq N+1$ by Lemma \ref{lem-gromov-product} (ii), which is a contradiction since both sequences $\{ x_{n}\}$, $\{ y_{m}\}$ converge to $\varphi(a)\in \partial C$.
\end{pf}

\begin{lem}\label{lem-bijective-varphi}
The map $\varphi \colon \Hat{{\Bbb R}}\setminus \Hat{{\Bbb Q}}\rightarrow \partial C$ is bijective.
\end{lem}

\begin{pf}
Note that the set ${\cal L}(0, a)$ of all geodesic rays from $0\in V(C)$ to $a\in \partial C$ is non-empty for any $a\in \partial C$ by Lemma \ref{non-empty-lem} and Theorem \ref{thm-main-excep}.

Let $a\neq b\in \Hat{{\Bbb R}}\setminus \Hat{{\Bbb Q}}$. There exist two rays $f\in {\cal L}(0, \varphi(a))$, $g\in {\cal L}(0, \varphi(b))$. Let $e\in E(\{ 0, b\}, a)$. It follows from Lemma \ref{lem-converge-varphi} that there exist $N$ such that $f(n)$ (resp. $g(n)$) is in the same side of $e$ as $a$ (resp. $b$) for any $n\geq N$. Moreover, since the edge $e$ is in $E(\{ f(n)\}_{n\geq N}, \{ g(n)\}_{n\geq N})$, we have
\[d(f(n), g|_{[N, \infty )})\geq d(f(n), e)\]  
for any $n\geq N$ and $d(f(n), e)\rightarrow \infty$ as $n\rightarrow \infty$. This means that $f$ and $g$ can not define the same point in $\partial C$ by Lemma \ref{hyp-fund} (ii). Hence, $\varphi$ is injective.

Next, we show the surjectivity. Let $\alpha \in \partial C$ and $f\in {\cal L}(0, \alpha)$. Since the circle is compact, there exists a subsequence of $\{ n_{i}\}_{i}$ of ${\Bbb N}$ such that the sequence $\{ f(n_{i})\}_{i}$ converges to some $a\in \Hat{{\Bbb R}}$ on the circle. 

First, we prove that $a$ is irrational. For deducing a contradiction, we may assume that $a=0$. Since $\{ f(n_{i})\}$ converges to $0$ on the circle, there exists an edge incident to $0$ through which $f$ passes. Recall that the set of all edges incident to $0$ is
\[\{ (0, 1/k): k\in {\Bbb N}\setminus \{ 0\}\}\cup \{ (0, -1/k): k\in {\Bbb N}\setminus \{ 0\}\}\cup \{ (0, 1/0)\}.\]
Suppose that $f$ passes through $0$. Let $n_{0}\in {\Bbb N}$ be a natural number such that $f(n_{0})=0$. Assume that $f(n_{0}+1)=1/k$, $k\in {\Bbb N}\setminus \{ 0\}$. Then $1/(k+1)<f(n)<1/(k-1)$ for any $n\geq n_{0}+1$, otherwise $f|_{[n_{0}+1, \infty)}$ would pass through either edge $(0, 1/(k+1))$ or $(0, 1/(k-1))$, which contradicts the fact that $f$ is a geodesic. Thus, the sequence $\{ f(n_{i})\}_{i}$ can not converge to $0$ on the circle, which is a contradiction. In the cases of $f(n_{0}+1)=-1/k$, $k\in {\Bbb N}\setminus \{ 0\}$ and $f(n_{0}+1)=1/0$, we can deduce a contradiction similarly. 

Suppose that $f$ does not pass through $0$. Since $\{ f(n_{i})\}$ converges to $0$, there exists $k\in {\Bbb N}$ such that $f$ passes either through four edges $(0, 1/k)$, $(0, 1/(k+1))$, $(0, 1/(k+2))$ and $(0, 1/(k+3))$ or through four edges $(0, -1/k)$, $(0, -1/(k+1))$, $(0, -1/(k+2))$ and $(0, -1/(k+3))$. It follows from the assumption that $f$ must pass either through any of $1/k$, $1/(k+1)$, $1/(k+2)$ and $1/(k+3))$ or through any of $-1/k$, $-1/(k+1)$, $-1/(k+2)$ and $-1/(k+3)$, which contradicts the fact that $f$ is a geodesic. 

Therefore, we have shown that the sequence $\{ f(n_{i})\}_{i}$ converges to an irrational number $a$ on the circle. 

Finally, we show that $f$ defines the boundary point $\varphi(a)\in \partial C$, that is, $\alpha =\varphi(a)$. Let $\{ e_{n}\}_{n\in {\Bbb N}}$ be a sequence of edges in $E(0, a)$ such that $e_{n}\in E(e_{n-1}, a)$ for each $n$ and that both vertices of the edge $e_{n}$ converge to $a$ on the circle. Let $x_{n}$ be any vertex of $e_{n}$. The geodesic ray $f$ must pass through any edge $e_{n}$ because the sequence $\{ f(n_{i})\}_{i}$ converges to $a$ on the circle and for each $n$, we have $e_{n}\in E(0, f(n_{i}))$ for large $i$. For each $n$, let $m_{n}\in {\Bbb N}$ be a natural number such that $f(m_{n})\in \overline{e_{n}}$. It is clear that $m_{n}<m_{n'}$ if $n<n'$. Since $d(f(m_{n}), x_{n})\leq 1$ for each $n$, two sequences $\{ f(m_{n})\}_{n\in {\Bbb N}}$, $\{ x_{n}\}_{n\in {\Bbb N}}$ define the same point on $\partial C$.      
\end{pf}

\begin{lem}
The map $\varphi \colon \Hat{{\Bbb R}}\setminus \Hat{{\Bbb Q}}\rightarrow \partial C$ is a homeomorphism.
\end{lem}

\begin{pf}
Let $U$ be an open neighborhood on $\Hat{{\Bbb R}}\setminus \Hat{{\Bbb Q}}$ of an irrational number $a\in \Hat{{\Bbb R}}\setminus \Hat{{\Bbb Q}}$. Then there exist three edges $e_{1}$, $e_{2}$, $e_{3}$ in $C$ such that the intersection of the circle and the component of $\overline{D}\setminus \overline{e_{i}}$ to which $a$ belongs is contained in $U$ for each $i$ and $e_{i+1}\in E(e_{i}, a)$ for $i=1, 2$. It follows from Lemma \ref{upper-semi-cont} that there exists an open neighborhood $V$ on $\partial C$ of $\varphi(a)$ such that 
\[G(0, V)_{r}\subseteq G(0, \varphi(a))_{r},\]
where $r=d(0, e_{3})$. We show $\varphi^{-1}(V)\subseteq U$, which means that $\varphi^{-1}$ is continuous.

Let $v\in V$ and $f\in {\cal L}(0, v)$. Then $f(r)\in G(0, \varphi(a))_{r}$ and $d(f(r), e_{3})\leq 1$ since any geodesic ray in ${\cal L}(0, \varphi(a))$ passes through $e_{3}$. Thus, $f$ must pass through $e_{1}$ since we have either $e_{2}\in E(0, f(r))$ or $f(r)\in \overline{e_{2}}$. It follows that after passing through $e_{1}$, the geodesic ray $f$ keeps to be in the component of $\overline{D}\setminus \overline{e_{1}}$ to which $a$ belongs. It follows from the proof of the surjectivity in Lemma \ref{lem-bijective-varphi} that the image by $\varphi$ of any cluster point of the sequence $\{ f(n)\}_{n\in {\Bbb N}}$ on the circle is equal to the boundary point determined by $f$. Since the intersection of the circle and the component of $\overline{D}\setminus \overline{e_{1}}$ to which $a$ belongs is contained in $U$, we get $\varphi(v)\in U$.   

Let $\alpha \in \partial C$ and $V$ be an open neighborhood in $\partial C$ of $\alpha$. There exists a unique $a\in \Hat{{\Bbb R}}\setminus \Hat{{\Bbb Q}}$ such that $\varphi(a)=\alpha$. Let $N\in {\Bbb N}$ be a natural number such that $(\alpha |v)_{0}\geq N-4\delta$ for $v\in \partial C$ implies $v\in V$, where $\delta$ is a constant such that $C$ is $\delta$-hyperbolic. Let $e\in E(0, a)$ be an edge with $d(0, e)\geq N$. We can take such an edge as follows: if $e_{1}, \ldots, e_{N}$ be a sequence of edges in $E(0, a)$ such that $e_{i}\in E(e_{i-1}, a)$ for each $i$, then $e=e_{N}$ satisfies $d(0, e)\geq N$. Let $U$ be the intersection of the circle and the component of $\overline{D}\setminus \overline{e}$ to which $a$ belongs. We show $\varphi(U)\subseteq V$, which means that $\varphi$ is continuous.

Let $u\in U$ and $f\in {\cal L}(0, \varphi(u))$, $g\in {\cal L}(0, \alpha)$. Since both $f$ and $g$ pass through $e$, for all sufficiently large $n$ and $m$, we see that any geodesic segment $[f(n), g(m)]$ between $f(n)$ and $g(m)$ is contained in the component of $\overline{D}\setminus \overline{e}$ to which $a$ belongs and that 
\[(f(n)|g(m))_{0}+4\delta \geq d(0, [f(n), g(m)])\geq d(0, e)\geq N.\]
The first inequality comes from Lemma \ref{Gromov-product}. It follows that $(\varphi(u)|\alpha)_{0}\geq N-4\delta$ and that $\varphi(u)\in V$.    
\end{pf}

Since the equivariance of $\varphi$ under the homomorphism $p_{M}$ is clear by definition, Theorem \ref{thm-excep-identified} follows from the above lemmas.

It is well-known that the standard action of $SL(2, {\Bbb Z})$ on $\Hat{{\Bbb R}}$ is amenable for any quasi-invariant measure on $\Hat{{\Bbb R}}$ (e.g., see \cite{spa-zim}). Therefore, the following theorem follows from Theorem \ref{amenable-action-cc-non-excep} and Theorem \ref{thm-excep-identified}.

\begin{thm}\label{amenable-action-cc}
Let $M$ be a compact orientable surface with $\kappa(M)\geq 0$ and $C$ and $\Gamma$ be the curve complex and the mapping class group of $M$, respectively. 
Let $\mu$ be a quasi-invariant probability measure on $\partial C$. Then the action of $\Gamma$ on $(\partial C, \mu)$ is amenable.
\end{thm}

Moreover, by Remark \ref{rem-kla-excep}, we obtain

\begin{cor}\label{amenable-action}
Let $M$ be a compact orientable surface with $\kappa(M)\geq 0$ and let $\mu$ be a quasi-invariant probability measure on ${\cal PMF}$ whose support is contained in ${\cal MIN}$. Then the action of $\Gamma$ on $({\cal PMF}, \mu)$ is amenable. 

If $\kappa(M)=0$, the same statement holds for any quasi-invariant measure $\mu$ on ${\cal PMF}$.
\end{cor}

\begin{rem}\label{rem-excep-funct-const}
We can obtain Theorem \ref{amenable-action-cc} for exceptional surfaces by the same argument as in Theorem \ref{amenable-action-cc-non-excep} for non-exceptional surfaces since the curve complex $C$ of an exceptional surface $M$ with $\kappa(M)=0$ satisfies (F1) and (F2) by Theorem \ref{thm-main-excep} and the action of the mapping class group $\Gamma$ on $C$ is cocompact. In particular, we can construct non-negative valued Borel functions $f_{n}$ on $\partial C \times C$ satisfying that
\[\sum_{x\in C}f_{n}(a, x)=1\]
for any $a\in \partial C$ and
\[\lim_{n\rightarrow \infty}\sup_{a\in \partial C}\sum_{x\in C}|f_{n}(a, x)-f_{n}(g^{-1}a, g^{-1}x)|=0\]
for any $g\in \Gamma$ (see the proof of Theorem \ref{amenable-action-cc-non-excep}).
\end{rem}

\chapter[Indecomposability of equivalence relations]{Indecomposability of equivalence relations generated by the mapping class group}\label{chapter-indec}

In this chapter, we study equivalence relations generated by an essentially free, non-singular Borel action of the mapping class group on a standard Borel space with a finite positive measure, and we prove Theorems \ref{intro-indecomposability} and \ref{intro-not-me-hyp}. This work is inspired by the following two theorems shown by Zimmer \cite{zim3} and Adams \cite{adams2}:

\begin{thm}[\ci{Theorems 1.1, 1.3}{zim3}]\label{zim3-main}
Let $\Gamma$ be a lattice in a connected simple Lie group with finite center or the fundamental group of a finite volume manifold with negative sectional curvature bounded away from $0$. Let $(S, \mu)$ and $(S_{i}, \mu_{i})$ for $i=1, 2$ be standard Borel spaces with finite positive measures. Suppose that we have an essentially free, measure-preserving Borel action of $\Gamma$ on $(S, \mu)$. Let ${\cal R}$ be the discrete measured equivalence relation generated by the $\Gamma$-action. 

If there exists a discrete measured equivalence relation ${\cal S}_{i}$ on $(S_{i}, \mu_{i})$ for $i=1, 2$ such that the product relation ${\cal S}_{1}\times {\cal S}_{2}$ is weakly isomorphic to the relation ${\cal R}$, then either $S_{1}$ or $S_{2}$ is finite (up to null sets). 
\end{thm}

\begin{thm}[\ci{Theorem 6.1}{adams2}]\label{adams2-main}
Let $\Gamma$ be a hyperbolic group. Let $(X, \mu)$ and $(Y_{i}, \nu_{i})$ for $i=1, 2$ be standard Borel spaces with finite positive measures. Suppose that we have an essentially free, non-singular Borel action of $\Gamma$ on $(X, \mu)$. Let ${\cal R}$ be the discrete measured equivalence relation generated by the $\Gamma$-action.   

Assume that we have Borel subsets $X'\subseteq X$, $Y'\subseteq Y_{1}\times Y_{2}$ with positive measure and a recurrent discrete measured equivalence relation ${\cal S}_{i}$ on $(Y_{i}, \nu_{i})$ for $i=1, 2$ such that the restricted relation $({\cal S}_{1}\times {\cal S}_{2})_{Y'}$ of the product relation is isomorphic to the restricted relation $({\cal R})_{X'}$. Then the relation $({\cal R})_{X'}$ is amenable.
\end{thm}

Here, the symbol $({\cal R})_{A}$ means the restriction to a Borel subset $A\subseteq X$ of a discrete measured equivalence relation ${\cal R}$ on $(X, \mu)$ (see Definition \ref{defn-restriction}).   

In the proofs of the both theorems, a certain property of the ``action'' of the relation on a relevant boundary object was of basic importance. The technique using this property has repeatedly appeared in attacking many other problems about equivalence relations studied mainly by Zimmer (for example, see \cite{zim2} and the references therein). We will show that the mapping class group also has the same property as the groups in the above theorems, using the boundary at infinity of the curve complex and developing the notion of reducibility for subrelations of relations generated by the mapping class group.

In general, in order to study the structure of a given group, if it has a continuous action on some (compact) topological space, then it is often effective to investigate the fixed points for its subgroups. We can say a similar thing also for equivalence relations. For example, if $\Gamma$ is a hyperbolic group and ${\cal R}$ is an equivalence relation generated by an essentially free, non-singular action of $\Gamma$, then we have the natural Borel cocycle $\rho \colon {\cal R}\rightarrow \Gamma$ defined by
\[\rho(gx, x)=g\]
for $x\in X$ and $g\in \Gamma$. This map satisfies the cocycle identity as follows:
\[\rho(x, y)\rho(y, z)=\rho(x, z)\]
for $(x, y), (y, z)\in {\cal R}$. The cocycle induces the ``action'' of ${\cal R}$ on $M(\partial \Gamma)$, the space of probability measures on the boundary $\partial \Gamma$. Adams \cite{adams2} investigated invariant Borel maps $\varphi \colon X\rightarrow M(\partial \Gamma)$ for a recurrent subrelation ${\cal S}$ of ${\cal R}$, where the word ``invariant'' means that the map $\varphi$ satisfies
\[\rho(x, y)\varphi(y)=\varphi(x)\]
for $(x, y)\in {\cal S}$. This invariant Borel map corresponds to a fixed point for the ``action'' of ${\cal S}$. If ${\cal S}$ is amenable, then there exists an invariant Borel map by definition. Adams verified that if a recurrent subrelation ${\cal S}$ has an invariant Borel map $\varphi \colon X\rightarrow M(\partial \Gamma)$, then the support of $\varphi(x)$ consists of at most two points for a.e.\ $x\in X$ and that $\varphi$ induces the invariant Borel map $X\rightarrow \partial_{2}\Gamma$ (see \cite[Lemma 3.5]{adams2}), where $\partial_{2}\Gamma$ denotes the quotient space of $\partial \Gamma \times \partial \Gamma$ by the coordinate interchanging action of the symmetric group of two letters. The space $\partial_{2}\Gamma$ can be naturally identified with a Borel subset of $M(\partial \Gamma)$.

Using amenability of the $\Gamma$-action on $\partial \Gamma$, we can show that if a recurrent subrelation ${\cal S}$ has an invariant Borel map, then ${\cal S}$ is amenable (see \cite[Lemma 5.3]{adams2} and Proposition \ref{irreducible-amenable-hyp}). Thus, we see that a recurrent subrelation ${\cal S}$ is amenable if and only if it has an invariant Borel map $X\rightarrow M(\partial \Gamma)$. Moreover, in this case, ${\cal S}$ has a unique invariant Borel map $\varphi_{0}\colon X\rightarrow \partial_{2}\Gamma \subseteq M(\partial \Gamma)$ with the maximal property in the following sense: if $\varphi \colon X\rightarrow M(\partial \Gamma)$ is an invariant Borel map for ${\cal S}$, then we have
\[{\rm supp}(\varphi(x))\subseteq {\rm supp}(\varphi_{0}(x))\]
for $x\in X$ \cite[Lemma 3.3]{adams2}, where ${\rm supp}(\nu)$\index{$supp \ m$@${\rm supp}(\nu)$} denotes the support of a measure $\nu$.

With these claims, we can sketch an outline of the proof of Theorem \ref{adams2-main} as follows: for simplicity, assume that ${\cal R}$ and ${\cal S}_{1}\times {\cal S}_{2}$ are isomorphic. Then we have a Borel cocycle
\[\rho \colon {\cal S}_{1}\times {\cal S}_{2}\rightarrow \Gamma.\]
Let ${\cal T}_{1}$ be an amenable recurrent subrelation of ${\cal S}_{1}$ (see Lemma \ref{lem-ame-rec-exist}). Then there exists a unique invariant Borel map $\varphi_{0}\colon Y_{1}\times Y_{2}\rightarrow \partial_{2}\Gamma$ for the relation ${\cal T}_{1}\times {\cal D}_{Y_{2}}$ with the maximal property, where ${\cal D}_{Z}$ denotes the trivial relation on a standard Borel space $Z$. For any Borel isomorphism $g$ on $Y_{1}\times Y_{2}$ such that $(g(x), x)\in {\cal D}_{Y_{1}}\times {\cal S}_{2}$ for $x\in Y_{1}\times Y_{2}$, we can see easily that the Borel map
\[Y_{1}\times Y_{2}\ni x\mapsto \rho(x, g(x))\varphi_{0}(g(x))\in \partial_{2}\Gamma\]
is also invariant for ${\cal T}_{1}\times {\cal D}_{Y_{2}}$, which implies 
\[\rho(x, g(x))\varphi_{0}(g(x))=\varphi_{0}(x)\]
for $x\in Y_{1}\times Y_{2}$ by the maximality of $\varphi_{0}$ and that $\varphi_{0}$ is invariant for the relation ${\cal D}_{Y_{1}}\times {\cal S}_{2}$. Thus, ${\cal S}_{2}$ is amenable. In a similar way, we can show that ${\cal S}_{1}$ is also amenable, which proves Theorem \ref{adams2-main}.

It is Zimmer \cite{zim3} who first got the idea of using the fact that any amenable recurrent subrelation ${\cal S}$ has a unique invariant Borel map with a certain maximal property in order to prove indecomposability of equivalence relations.

If $\Gamma(M)$ is the mapping class group of a surface $M$ with $\kappa(M)\ge 0$ and we have an essentially free, non-singular Borel action on a standard Borel space $(X, \mu)$, then it induces an equivalence relation ${\cal R}$ and a Borel cocycle $\rho \colon {\cal R}\rightarrow \Gamma(M)$. Following Adams' argument, we study the action of ${\cal R}$ on $M({\cal PMF})$, the space of probability measures on ${\cal PMF}$, and investigate its invariant Borel maps for a recurrent subrelation ${\cal S}$ of ${\cal R}$. 

Let us denote $\Phi ={\cal PMF}\setminus {\cal MIN}$. We will show that if there exists an invariant Borel map for ${\cal S}$, then it is impossible for ${\cal S}$ to have an invariant Borel map $\varphi \colon X\rightarrow M({\cal PMF})$ such that $\varphi(x)({\cal MIN})>0$ and $\varphi(x)(\Phi)>0$ for $x$ in a Borel subset with positive measure (see Theorem \ref{alternative} and Corollary \ref{alternative-cor}). By this assertion, it is natural to classify recurrent subrelations ${\cal S}$ as follows: if ${\cal S}$ has an invariant Borel map $\varphi \colon X\rightarrow M({\cal PMF})$ such that $\varphi(x)({\cal MIN})=1$ (resp. $\varphi(x)(\Phi)=1$) for a.e.\ $x\in X$, then it is said to be {\it irreducible and amenable} (resp. {\it reducible})\index{irreducible and amenable!subrelation}\index{reducible!subrelation}.  

If ${\cal S}$ is irreducible and amenable, then we can apply Adams' argument for hyperbolic groups to our case and show that ${\cal S}$ is amenable, using the $\Gamma(M)$-equivariant map $\pi \colon {\cal MIN}\rightarrow \partial C$ in Chapter \ref{chapter:amenable-action}, Section \ref{boundary-of-curve-complex} and the hyperbolicity of the curve complex $C$ (see Proposition \ref{irreducible-amenable}).

If ${\cal S}$ is reducible, then we can define the canonical reduction system for ${\cal S}$ (see Definition \ref{crs}), using the $\Gamma(M)$-equivariant Borel map $H\colon \Phi \rightarrow S(M)$ in Definition \ref{construction-of-H}, where $S(M)$ is the set of all simplices of $C$ if $\kappa(M)>0$, and $S(M)=V(C(M))$ if $\kappa(M)=0$. The canonical reduction system $\sigma(G)\in S(M)$ for a reducible subgroup $G$ of $\Gamma(M)$ was introduced in \cite{blm}, \cite{ivanov1} and defined as the set of all invariant elements in $V(C(M))$ for $G$ with a certain additional property. It plays an important role in the study of reducible subgroups.

Using these notions, for any amenable recurrent subrelation ${\cal S}$, we can construct an invariant Borel map $X\rightarrow M({\cal PMF})$ with some maximal property. Then we can verify Theorem \ref{intro-indecomposability} by the fact that ${\cal R}$ can be neither amenable nor reducible if ${\cal R}$ is of type ${\rm II}_{1}$ (see Proposition \ref{typeII}).

For Theorem \ref{intro-not-me-hyp}, we notice a group-theoretic difference between the mapping class group $\Gamma_{1}$ of a non-exceptional surface and a hyperbolic group $\Gamma_{2}$, that is, there exists a non-amenable subgroup of $\Gamma_{1}$ which contains an infinite amenable normal subgroup, but there exist no such subgroups in $\Gamma_{2}$. We apply this difference to the proof of Theorem \ref{intro-not-me-hyp} as follows: let ${\cal R}_{1}$ and ${\cal R}_{2}$ be equivalence relations generated by essentially free, measure-preserving actions of $\Gamma_{1}$ and $\Gamma_{2}$, respectively. Since there exists a non-amenable subgroup of $\Gamma_{1}$ which contains an infinite amenable normal subgroup, we see that there exists a non-amenable subrelation ${\cal S}$ of ${\cal R}_{1}$ which has an amenable recurrent normal subrelation. On the other hand, if ${\cal T}$ is a subrelation of ${\cal R}_{2}$ and has an amenable recurrent normal subrelation, then ${\cal T}$ is amenable. This fact has been essentially shown by Adams \cite{adams2}.


\section{Construction of Busemann functions and the MIN set map}

Let $C$ be (the vertex set of) the curve complex for a surface $M$ with $\kappa(M)\geq 0$. We construct a $\Gamma(M)$-equivariant Borel map $MS'\colon \delta C\rightarrow {\cal F}'$, where 
\[\delta C=\{ (a, b, c)\in (\partial C)^{3}: a\neq b\neq c\neq a\}\]
and ${\cal F}'$ is the set of all non-empty finite subsets of $C$ whose diameters are more than or equal to $3$. The $\Gamma(M)$-action on $\delta C$ is given by the formula $g(a, b, c)=(ga, gb, gc)$ for $g\in \Gamma(M)$ and $(a, b, c)\in \delta C$. 

Adams \cite{adams1} constructed an analogous map to $MS'$ and call it the MIN set map, which plays a crucial role in the proof of Theorem \ref{adams2-main}. Needless to say, the local infiniteness of $C$ makes the same construction as Adams' one impossible. The finiteness properties (F1) and (F2) still play an indispensable role in the following construction of the map $MS'$. We need to introduce Busemann functions on $C$, which are treated in the first subsection. In the second subsection, we construct the map $MS'$.

\subsection{Busemann functions}   

In this subsection, we introduce Busemann functions on a $\delta$-hyperbolic graph $X$ with $\partial X\neq \emptyset$ satisfying the properties (F1) and (F2) appearing in Chapter \ref{chapter-property-A}, Section \ref{sec:Acc}. Our definition is slightly different from that in \cite{ghys-harpe} for a proper hyperbolic space. However, since their difference is covered by some universal constant associated with a hyperbolic constant $\delta$, their properties are almost the same. The advantage of our definition is that we can easily show measurability of the Busemann functions on the boundary $\partial X$. 

In this subsection, we denote by $X$ (the vertex set of) a $\delta$-hyperbolic graph satisfying the properties (F1), (F2) and $\partial X\neq \emptyset$.

\begin{defn}
We define a function $\alpha \colon \partial X\times X\times X\rightarrow {\Bbb Z}$ by the formula
\[\alpha_{a}(x, y)=\limsup_{t\rightarrow \infty}(d(x, G(y, a)_{t})-t),\]
\index{$\ a a x y$@$\alpha_{a}(x, y)$}where $G(y, a)_{t}=\{z\in G(y, a):d(y, z)=t\}$ for $y\in X$, $a\in \partial X$ and $t\in {\Bbb N}$. We call the function {\it (our version of) the Busemann function}\index{Busemann function} at the point $a\in \partial X$.

On the other hand, the Busemann function (for a proper hyperbolic space) appearing in \cite[Chapitre 8]{ghys-harpe} is defined as follows: first, for $a\in \partial X$, $x, y\in X$ and a geodesic $h\in {\cal L}(y, a)$, define
\[\beta_{a}(x, h)=\limsup_{t\rightarrow \infty}(d(x, h(t))-t)\leq d(x, y).\]
\index{$\ b a x h$@$\beta_{a}(x, h)$}Secondly, define a function $\beta \colon \partial X\times X\times X\rightarrow {\Bbb Z}$ by the formula
\[\beta_{a}(x, y)=\sup \{\beta_{a}(x, h): h\in {\cal L}(y, a)\}.\]
\end{defn}
\index{$\ b a x y$@$\beta_{a}(x, y)$}
Remark that $\beta$ is well-defined even if the graph $X$ satisfying (F1), (F2) and $\partial X\neq \emptyset$ is not proper, since ${\cal L}(y, a)$ is non-empty for any $a\in \partial X$ and $y\in X$ by Lemma \ref{non-empty-lem}. We collect fundamental properties of $\alpha$ and $\beta$.

\begin{lem}\label{hyp-basic4}
\begin{enumerate}
\renewcommand{\labelenumi}{\rm(\roman{enumi})}
\item If $g$ and $h$ are two equivalent geodesic rays with the same origin in $X$, then we have
\[\sup_{t\in {\Bbb N}}d(g(t), h(t))\leq 8\delta. \]
\item For $a\in \partial X$, $x, y\in C$ and $h, h'\in {\cal L}(y, a)$, we have
\begin{align*}
\vert \beta_{a}(x, h)-\beta_{a}(x, h')\vert &\leq 8\delta,\\
\vert \beta_{a}(x, y)-\beta_{a}(x, h)\vert &\leq 8\delta.
\end{align*}
\item For $a\in \partial X$, $x, y\in X$ and $h\in {\cal L}(y, a)$, there exists $t_{0}\in {\Bbb N}$ such that 
\[|\beta_{a}(x, y)-(d(x, h(t))-t)|\leq 40\delta\]
for any $t\geq t_{0}$. 
\item For $a\in \partial X$ and $x, y, z\in X$, we have
\[|\beta_{a}(x, z)-\beta_{a}(y, z)|\leq d(x, y)+80\delta.\]
\end{enumerate}
\end{lem}

\begin{pf}
The assertion (i) is shown in \cite[Chapitre 7, Corollaire 3]{ghys-harpe}. The assertion (ii) follows from (i) immediately. We can verify the assertion (iii) in the same way as \cite[Chapitre 8, Lemme 1]{ghys-harpe}, using Lemma \ref{hyp-basic1}. The assertion (iv) follows from (iii). 
\end{pf}

\begin{lem}\label{beta-fund}
Let $a\in \partial X$ and $x, x', y, y', z\in X$.
\begin{enumerate}
\renewcommand{\labelenumi}{\rm(\roman{enumi})}
\item For any $b\in \partial X\setminus \{ a\}$, we have
\[\lim_{x\rightarrow b}\beta_{a}(x, y)=+\infty,\]
and 
\[\lim_{x\rightarrow a}\beta_{a}(x, y)=-\infty.\]
\item $|\beta_{a}(x, y)+\beta_{a}(y, x)|\leq 120\delta$.
\item $|\beta_{a}(x, y)+\beta_{a}(y, z)-\beta_{a}(x, z)|\leq 200\delta$.
\item $|\beta_{a}(x, y)-\beta_{a}(x', y')|\leq d(x, x')+d(y, y')+400\delta$.
\end{enumerate}
\end{lem}

We can prove this lemma similarly to \cite[Chapitre 8, Proposition 2]{ghys-harpe}, using Lemma \ref{hyp-basic4}.

\begin{lem}\label{beta}
For any $a\in \partial X$ and $x, y\in X$, we have
\[\vert \alpha_{a}(x, y)-\beta_{a}(x, y)\vert \leq 8\delta.\]
\end{lem}

\begin{pf}
Let $\varepsilon$ be any positive number. There exists a geodesic ray $h\in {\cal L}(y, a)$ such that
\[0\leq \beta_{a}(x, y)-\beta_{a}(x, h) <\varepsilon.\]
By Lemma \ref{hyp-basic4} (i), we have
\[d(h(t), z)\leq 8\delta\]
for any $t\in {\Bbb N}$ and $z\in G(y, a)_{t}$. Thus,
\[\vert (d(x, h(t))-t)-(d(x, z)-t)\vert \leq d(h(t), z)\leq 8\delta \]
for any $z\in G(y, a)_{t}$, which implies
\[\vert \beta_{a}(x, h)-\alpha_{a}(x, y)\vert \leq 8\delta,\]
and
\[\vert \beta_{a}(x, y)-\alpha_{a}(x, y)\vert \leq 8\delta 
+\varepsilon. \]
This completes the proof. 
\end{pf}

\begin{prop}\label{alpha-measurable}
The function $\alpha \colon \partial X\times X\times X\rightarrow {\Bbb Z}$ is measurable.
\end{prop}

\begin{pf}
We will show that the function $\alpha^{t}:\partial X\times X\times X\rightarrow {\Bbb Z}$ defined by the formula
\[\alpha^{t}_{a}(x, y)=d(x, G(y, a)_{t})-t\]
is lower semi-continuous for each $t\in {\Bbb N}$, which proves the proposition.

Fix $t\in {\Bbb N}$. For any $a\in \partial X$, there exists a neighborhood $V$ of $a$ in $\partial X$ such that 
\[G(y, V)_{t}=\{ z\in G(y, V): d(y, z)=t\} \subseteq G(y, a)_{t}\]
by Lemma \ref{upper-semi-cont}, which implies the above claim. 
\end{pf}


\subsection{The MIN set map}\label{MIN-set-map} 

For any hyperbolic group $\Gamma$, Adams \cite{adams1} constructed a certain $\Gamma$-equivariant measurable map from the set
\[\delta \Gamma =\{(a, b, c)\in (\partial \Gamma)^{3} :a\neq b\neq c\neq a\} \]
to ${\cal F}(\Gamma)$, the set of all non-empty finite subsets of $\Gamma$. The $\Gamma$-action on $\delta \Gamma$ is given by the formula $g(a, b, c)=(ga, gb, gc)$ for $g\in \Gamma$ and $(a, b, c)\in \delta \Gamma$. In the construction, he studied the sum 
\[F_{abc}^{y}(x)=\beta_{a}(x, y)+\beta_{b}(x, y)+\beta_{c}(x, y)\]
of Busemann functions of three distinct points $a, b, c\in \partial \Gamma$ for $x, y\in \Gamma$. When $y$ is fixed, the function $F_{abc}^{y}$ is proper and bounded below. Thus, it gives the non-empty finite subset $MS_{abc}^{y}$ of $\Gamma$ where it attains its minimum. Moreover, it can be shown that the set 
\[MS_{abc}=\bigcup_{y\in \Gamma}MS_{abc}^{y}\]
is also finite. Therefore, we can define a $\Gamma$-equivariant Borel map $\delta \Gamma \rightarrow {\cal F}(\Gamma)$ by the formula 
\[(a, b, c)\rightarrow MS_{abc},\]
which is called the {\it MIN set map}\index{MIN set map}.

In this subsection, we construct a similar function for the mapping class group and the curve complex.

Let $X$ be a $\delta$-hyperbolic graph satisfying the properties (F1), (F2) and $|\partial X|\geq 3$. We identify $X$ with its vertex set. Denote
\[\delta X=\{(a, b, c)\in (\partial X)^{3}: a\neq b\neq c\neq a\}.\]
\index{$\ d X$@$\delta X$}

\begin{lem}\label{chi-is-measurable} 
The function $X^{4}\times \delta X\rightarrow \{ 0, 1\}$ defined by the formula
\[(x, y, z, w, (a, b, c))\mapsto \chi_{abc}^{xyz}(w)\]
is upper semi-continuous, where $\chi_{abc}^{xyz}$ is the characteristic function on the set $G(x, b)\cup G(y, c)\cup G(z, a)$.
\end{lem}

\begin{pf}
Fix $x, y, z\in X$ and $(a, b, c)\in \delta X$. Let $w\in X\setminus (G(x, b)\cup G(y, c)\cup G(z, a))$. By  Lemma \ref{upper-semi-cont}, for each $i\in \{ a, b, c\}$, there exists neighborhood $V_{i}$ of $i$ in $\partial X$ such that 
\[w\notin G(x, V_{b}),\ \  w\notin G(y, V_{c}),\ \  w\notin G(z, V_{a}).\]
It means that the function in the lemma is upper semi-continuous.  
\end{pf} 

Define a function $F\colon X\times \delta X\rightarrow \{ 0, 1\}$ by the formula
\[F(w, (a, b, c))=\sup_{x, y, z}\chi_{abc}^{xyz}(w),\]
where the supremum is taken over all $x\in G(a, b)$, $y\in G(b, c)$, $z\in G(c, a)$. This function is measurable, and when $X$ is the curve complex $C$ of a surface $M$ with $\kappa(M)\geq 0$, it is invariant under the action of the mapping class group $\Gamma(M)$, that is, 
\[F(gw, g(a, b, c))=F(w, (a, b, c))\]
for any $g\in \Gamma(M)$, where the action of $\Gamma(M)$ on $\delta C$ is given by the formula $g(a, b, c)=(ga, gb, gc)$. Let us denote by $X(a, b, c)$\index{$X a b c$@$X(a, b, c)$} the support of the function $F(\cdot , (a, b, c))$, that is,
\[X(a, b, c)=\bigcup_{x\in G(a, b),\ y\in G(b, c),\ z\in G(c, a)}G(x, b)\cup G(y, c)\cup G(z, a).\]
(Although it seems more natural to treat the function
\[(w, (a, b, c))\mapsto \chi_{G_{abc}}(w),\]
where $G_{abc}=G(a, b)\cup G(b, c)\cup G(c, a)$ and $\chi_{G_{abc}}$ is the characteristic function on the set $G_{abc}$, the author does not know whether it is measurable.)

\begin{lem}\label{lem:bow3}
Let $x, y\in X\cup \partial X$.  
\begin{enumerate}
\renewcommand{\labelenumi}{\rm(\roman{enumi})}
\item For any $g\in {\cal L}_{T}(x, y)$, we have
\[G(x, y)\subseteq \bigcup_{n\in {\rm dom}(g)}B(g(n); \delta_{0}).\]
\item The intersection $G(x, y)\cap B$ is finite for any bounded ball $B$ in $X$.
\end{enumerate} 
\end{lem}

\begin{pf}
The assertion (i) in the case where either $x$ or $y$ is in $X$ follows from Lemmas \ref{lem:bow} and \ref{lem:bow2} immediately. Assume that both $x$ and $y$ are in $\partial X$ and let $g$, $h\in {\cal L}_{T}(x, y)$. Applying Lemma \ref{lem:bow2} to two geodesic rays $g|_{[0, \infty)}$, $h|_{[0, \infty)}$, we see that there exist $n_{0}, m_{0}\in {\Bbb N}$ such that 
\[d(g(n_{0}+t), h(m_{0}+t))\leq \delta_{0}\]
for any $t\in {\Bbb N}$. Let $t_{0}\geq \delta_{0}$ be a natural number. Applying Lemma \ref{lem:bow2} to two geodesic rays $g|_{(-\infty, n_{0}+t_{0}]}$, $h|_{(-\infty, m_{0}+t_{0}]}$, we see that 
\[d(g(n_{0}-t), h(m_{0}-t))\leq \delta_{0}\]
for any $t\in {\Bbb N}$. Thus, it follows from Lemma \ref{lem:bow} that
\[h\subseteq \bigcup_{n\in \Bbb Z}B(g(n);\delta_{0}).\]  

Let $B$ be a bounded ball in $X$ and $g\in {\cal L}_{T}(x, y)$. Since there exist at most finitely many $n\in {\Bbb Z}$ such that $B\cap B(g(n); \delta_{0})$ is non-empty, the assertion (ii) follows from (i) and the property (F1). 
\end{pf}

\begin{prop}\label{bounded}   
For any bounded ball $B$ in $X$ and any $(a, b, c)\in \delta X$, the intersection of $B$ and $X(a, b, c)$ is finite. 
\end{prop}

\begin{pf}
It is enough to show that the set
\[\left( \bigcup_{x\in G(a, b)}G(x, b)\right)\cap B\]
is finite for any $a, b\in \partial X$ with $a\neq b$.

Let $B$ be the ball of radius $R$ with center $w_{0}\in X$. For our purpose, we may assume that $w_{0}\in G(a, b)$. We fix a geodesic ray $g\in {\cal L}_{T}(a, b)$ with $g(-\infty)=a$, $g(0)=w_{0}$, $g(+\infty)=b$ in the proof of this proposition. 

\begin{claim}
For any bounded ball $A$ in $X$, the intersection
\[\left( \bigcup_{x\in A\cap G(a, b)}G(x, b)\right) \cap B\]
is a finite set. 
\end{claim}

\begin{pf}
By Lemma \ref{lem:bow3} (ii), both the intersections $A\cap G(a, b)$ and $G(x, b)\cap B$ are finite. 
\end{pf}  

\begin{claim}
We have
\[G(x, b)\cap B=\emptyset \]
for any $x\in G(a, b)$ with  $\beta_{b}(w_{0}, x)>R+8\delta$. 
\end{claim} 

\begin{pf}
For such $x\in G(a, b)$, let $y\in G(x, b)$ and $h\in {\cal L}_{T}(x, b)$ be a geodesic ray with $h(t_{0})=y$, $t_{0}\in {\Bbb N}$. Then
\begin{align*}
\beta_{b}(w_{0}, y)&\geq \beta_{b}(w_{0}, h|_{[t_{0}, \infty )})\\
                   &=\beta_{b}(w_{0}, h)+d(x, y)\\
                   &\geq \beta_{b}(w_{0}, x)-8\delta >R.
\end{align*}
The second inequality comes from Lemma \ref{hyp-basic4} (ii). On the other hand, if $y\in B$, then for any $h'\in {\cal L}(y, b)$ and $t\in {\Bbb N}$, we have
\[d(w_{0}, h'(t))-t\leq d(w_{0}, y)+d(y, h'(t))-t=d(w_{0}, y)\leq R,\] 
which implies $\beta_{b}(w_{0}, y)\leq R$. 
\end{pf}

\begin{claim}
Let $R_{0}$ be a positive number satisfying $R_{0}>R+3\delta_{0}+\delta_{1}$. Then we have
\[G(x, b)\cap B\subseteq G(g(-R-2\delta_{0}-\delta_{1}), b;\delta_{0})\cap B\]for any $x\in G(a, b)\cap B(g(-R_{0});\delta_{0})$.
\end{claim} 

\begin{pf}
Let $x\in G(a, b)\cap B(g(-R_{0}); \delta_{0})$. It follows from Lemma \ref{lem:bow2} that any geodesic ray in ${\cal L }_{T}(x, b)$ passes through the ball $B(g(-R-2\delta_{0}-\delta_{1}); \delta_{0})$ since $d(x, g(-R_{0}))\leq \delta_{0}$. For any point $y$ in $G(x, b)\cap B$, there exists a geodesic ray $h$ in ${\cal L}_{T}(x, b)$ passing through $y$. Since $h$ passes through the ball $B(g(-R-2\delta_{0}-\delta_{1}); \delta_{0})$, there exists $n_{0}\in {\Bbb N}$ such that $h(n_{0})\in B(g(-R-2\delta_{0}-\delta_{1});\delta_{0})$. Let $h'\in {\cal L}_{T}(h(n_{0}), b)$ be the restriction of $h$. We need to show that $h'$ passes through $y$. 

Suppose that the geodesic ray $h'$ does not pass through the point $y$, that is, $y\in h\setminus h'$. By Lemma \ref{lem:bow2}, there exists $n\in {\Bbb N}$ such that $d(h'(n), w_{0})\leq \delta_{0}$. Since $y\in h\setminus h'$ and $h'(0)\in B(g(-R-2\delta_{0}-\delta_{1}); \delta_{0})$, we would have
\begin{align*}
d(y, w_{0}) &\geq d(y, h'(n))-d(h'(n), w_{0}) \\
            &\geq d(h'(0), h'(n))+1-\delta_{0}\\
            &\geq d(g(-R-2\delta_{0}-\delta_{1}), w_{0})-d(h'(0), g(-R-2\delta_{0}-\delta_{1}))\\
            &\ \ \ \ \ \ \ \ \ \ \ \ \ \ \ \ \ \ \ \ \ \ \ \ \ \ \ \ \ \ \ \ \ -d(h'(n), w_{0})+1-\delta_{0}\\
            &\geq R+2\delta_{0}+\delta_{1}-\delta_{0}-\delta_{0}+1-\delta_{0}\\
            &>R
\end{align*}
(remark that $\delta_{1}>\delta_{0}$). This contradicts $y\in B$. 
\end{pf}

\begin{claim}For any $r>0$, the set
\[\{ y\in X: -r\leq \beta_{b}(w_{0}, y)\leq r\}\cap \bigcup_{x\in G(a, b)}G(x, b)\]
is bounded. 
\end{claim}

\begin{pf}
For any point $y$ in the above set, let $x\in G(a, b)$ with $y\in G(x, b)$. Considering a geodesic in ${\cal L}_{T}(a, b)$ through $x$ and a geodesic in ${\cal L}_{T}(x, b)$ through $y$, we can find $y_{0}\in G(a, b)$ such that $d(y, y_{0})\leq \delta_{0}$. Moreover, there exists a point $y_{1}$ in the fixed bi-infinite geodesic $g\in {\cal L}_{T}(a, b)$ such that $d(y_{0}, y_{1})\leq \delta_{0}$ by Lemma \ref{lem:bow3} (i). Combining these facts with Lemma \ref{beta-fund} (iv), we have 
\[|\beta_{b}(w_{0}, y)-\beta_{b}(w_{0}, y_{1})|\leq 2\delta_{0}+400\delta.\]
On the other hand, by Lemma \ref{hyp-basic4} (iii), we see that 
\[||\beta_{b}(w_{0}, y_{1})|-d(w_{0}, y_{1})|\leq 8\delta \]
since $y_{1}\in g$ and the geodesic $g\in {\cal L}_{T}(a, b)$ passes through $w_{0}$. Hence, we have
\begin{align*}
d(y, w_{0})&\leq d(y, y_{1})+d(y_{1}, w_{0})\\
           &\leq 2\delta_{0}+|\beta_{b}(w_{0}, y_{1})|+8\delta \\
           &\leq |\beta_{b}(w_{0}, y)|+4\delta_{0}+408\delta \\
           &\leq r+4\delta_{0}+408\delta.
\end{align*}
The claim follows. 
\end{pf}
  
Return to the proof of Proposition \ref{bounded}. We introduce two sets
\[E=\{ x\in G(a, b): \beta_{b}(w_{0}, x)>R+8\delta\}, \]
\[F=G(a, b)\cap \left( \bigcup_{n> R+3\delta_{0}+\delta_{1}}B(g(-n); \delta_{0})\right).\]  
Let $x\in G(a, b)\setminus E$. By Lemma \ref{lem:bow3} (i), there exists $t_{1}\in {\Bbb Z}$ such that $d(x, g(t_{1}))\leq \delta_{0}$. It follows from Lemma \ref{beta-fund} (iv) that
\[R+8\delta \geq \beta_{b}(w_{0}, x)\geq \beta_{b}(w_{0}, g(t_{1}))-d(x, g(t_{1}))-400\delta.\]
Since $|t_{1}-\beta_{b}(w_{0}, g(t_{1}))|=|\beta_{b}(w_{0}, g|_{[t_{1}, \infty)})-\beta_{b}(w_{0}, g(t_{1}))|\leq 8\delta$ by Lemma \ref{hyp-basic4} (ii), we have
\[t_{1}\leq R+416\delta +\delta_{0}.\]
Thus, if $x\in G(a, b)\setminus (E\cup F)$, then
\[x\in \bigcup_{-R-3\delta_{0}-\delta_{1}\leq n\leq R+416\delta+\delta_{0}}B(g(n); \delta_{0}).\]
Let $A_{0}$ be a bounded ball in $X$ with center $w_{0}$ which contains the bounded subset
\[\bigcup_{-R-3\delta_{0}-\delta_{1}\leq n\leq R+416\delta+\delta_{0}}B(g(n); \delta_{0}).\]
We have
\[\left( \bigcup_{x\in G(a, b)\setminus (E\cup F)}G(x, b)\right) \cap B\subseteq \left( \bigcup_{x\in A_{0}\cap G(a, b)}G(x, b)\right) \cap B.\]
The right hand side is finite by Claim 1.

Next, the set 
\[\left( \bigcup_{x\in E}G(x, b)\right) \cap B\]
is empty by Claim 2.

Finally, we show that 
\begin{multline*} 
\left( \bigcup_{x\in F}G(x, b)\right) \cap B \\
\subseteq G(g(-R-2\delta_{0}-\delta_{1}), b; \delta_{0})\cap \left( \bigcup_{-R-\delta_{0}\leq n\leq R+\delta_{0}}B(g(n); \delta_{0})\right).
\end{multline*}
Since the right hand side is finite by the property (F2), so is the left hand side. Thus, we can show that the set 
\[\left( \bigcup_{x\in G(a, b)}G(x, b)\right) \cap B\]
is finite if the above inclusion is verified.  
 
By Claim 3, it suffices to show that
\begin{multline*}
G(g(-R-2\delta_{0}-\delta_{1}), b; \delta_{0})\cap B\\
\subseteq G(g(-R-2\delta_{0}-\delta_{1}), b; \delta_{0})\cap \left( \bigcup_{-R-\delta_{0}\leq n\leq R+\delta_{0}}B(g(n); \delta_{0})\right).
\end{multline*}
Let $f\in {\cal L}_{T}(g(-R-2\delta_{0}-\delta_{1}), b;\delta_{0})$ and let $s_{0}\in {\Bbb N}$ be a natural number such that $f(s_{0})\in B$. Since $d(f(0), f(s_{0}))\geq \delta_{0}+\delta_{1}>\delta_{0}$ by the triangle inequality on the quadrilateral $f(0)f(s_{0})w_{0}g(-R-2\delta_{0}-\delta_{1})$, there exists $n_{0}\in {\Bbb Z}$ such that $d(f(s_{0}), g(n_{0}))\leq \delta_{0}$ by Lemma \ref{lem:bow2}. It follows that
\[d(g(n_{0}), w_{0})\leq d(g(n_{0}), f(s_{0}))+d(f(s_{0}), w_{0})\leq \delta_{0}+R.\]
Since $g(0)=w_{0}$, we have $|n_{0}|\leq \delta_{0}+R$. This completes the proof since $f(s_{0})\in B(g(n_{0}); \delta_{0})$. 
\end{pf}

\begin{defn}
For $(a, b, c)\in \delta X$ and $x$, $y\in X$, we define
\[F^{y}_{abc}(x)=\alpha_{a}(x, y)+\alpha_{b}(x, y)+\alpha_{c}(x, y).\]
\end{defn}

By Proposition \ref{alpha-measurable}, this function is measurable for $(a, b, c)\in \delta X$.

\begin{prop}\label{F-is-proper}
Fix $(a, b, c)\in \delta X$ and $y\in X$. Then there exists a constant $M_{0}>0$ such that 
\[F^{y}_{abc}(w)\geq d(w, y)-M_{0}\]
for any point $w\in X(a, b, c)$.
\end{prop}

For the proof, we need the following lemma:

\begin{lem}\label{lem-constant-lem}
Fix $x, y\in X$ and $a, b\in \partial X$ with $a\neq b$. Then there exists a constant $M_{1}>0$ such that 
\[(f(n)|g(m))_{y}\leq M_{1}\]
for any geodesic rays $f\in {\cal L}(x, a)$, $g\in {\cal L}(y, b)$ and natural numbers $n, m\in {\Bbb N}$.
\end{lem}

\begin{pf}
This proof is essentially due to \cite[Lemma 6.3]{adams1}.

Let $f'\in {\cal L}(x, a)$ and $g'\in {\cal L}(y, b)$ be any geodesic rays. First, we prove that there exists a constant $M_{2}>0$ such that
\[(f'(n)|g'(m))_{y}\leq M_{2}\]
for any $n, m\in {\Bbb N}$.

Assume that it does not hold. Then there would exist sequences $\{ i_{k}\}_{k\in {\Bbb N}}$, $\{ j_{k}\}_{k\in {\Bbb N}}$ of positive integers such that
\[(f'(i_{k})|g'(j_{k}))_{y}\rightarrow \infty \]
as $k\rightarrow \infty$. Either $i_{k}\rightarrow \infty$ or $j_{k}\rightarrow \infty$, for otherwise we could replace both $\{ i_{k}\}$ and $\{ j_{k}\}$ by constant subsequences and $(f'(i_{k})|g'(j_{k}))_{y}\rightarrow \infty$ would be impossible. Hence, we may assume that $j_{k}\rightarrow \infty$. 

We first consider the special case where $i_{k}\rightarrow \infty$. Since $f'(n)\rightarrow a\in \partial X$ as $n\rightarrow \infty$, we have $(f'(i_{k})|f'(i_{l}))_{y}\rightarrow \infty$ as $k, l\rightarrow \infty$. Thus, by $\delta$-hyperbolicity of $X$, we see that
\[(f'(i_{k})|g'(j_{l}))_{y}\geq \min \{ (f'(i_{k})|f'(i_{l}))_{y},\ (f'(i_{l})|g'(j_{l}))_{y}\} -\delta \rightarrow \infty, \] 
as $k, l\rightarrow \infty$. Since $f'(i_{k})\rightarrow a$ and $g'(j_{l})\rightarrow b$, we have $a=b$, which is a contradiction.

Next, we consider the case where the sequence $\{i_{k}\}$ does not diverge to infinity. Taking a subsequence, we may assume that $\{ i_{k}\}$ is a constant sequence such that $i_{k}=n_{0}$ for any $k\in {\Bbb N}$ and $(f'(n_{0})|g'(j_{k}))_{y}\rightarrow \infty$. On the other hand,
\begin{align*}
2(f'(n_{0})|g'(j_{k}))_{y}&=d(f'(n_{0}), y)+d(g'(j_{k}), y)-d(f'(n_{0}), g'(j_{k}))\\
                          &\leq d(f'(n_{0}), y)+d(y, f'(n_{0}))=2d(f'(n_{0}), y)<\infty.
\end{align*}
It is a contradiction. Thus, there exists a constant $M_{2}>0$ such that 
\[(f'(n)|g'(m))_{y}\leq M_{2}\]
for any $n, m\in {\Bbb N}$.

For the proof of the lemma, note that
\[d(f'(n), f(n))\leq 8\delta, \ \ d(g'(n), g(n))\leq 8\delta \]
for any geodesic rays $f\in {\cal L}(x, a)$, $g\in {\cal L}(y, b)$ and $n\in {\Bbb N}$ by Lemma \ref{hyp-basic4} (i). It follows that
\[(f(n)|g(m))_{y}\leq M_{2}+32\delta\]
for any $n$, $m\in {\Bbb N}$. Hence, $M_{1}=M_{2}+32\delta$ is a desired constant.
\end{pf}

\begin{pf*}{{\sc Proof of Proposition \ref{F-is-proper}}.}
Fix $z\in X$. By Lemma \ref{lem-constant-lem}, there exists a constant $M>0$ such that 
\[(g(n)|g_{1}(m))_{y}\leq M, \ \ (g(n)|g_{2}(m))_{y}\leq M\]
for any geodesic rays $g\in {\cal L}(z, a)$, $g_{1}\in {\cal L}(y, b)$, $g_{2}\in {\cal L}(y, c)$ and $n$, $m\in {\Bbb N}$. This means that
\[d(g(n), g_{1}(m))\geq d(g(n), y)+d(g_{1}(m), y)-2M,\]
\[d(g(n), g_{2}(m))\geq d(g(n), y)+d(g_{2}(m), y)-2M.\]
Both of the right hand sides of the above two inequalities are equal to $d(g(n), y)+m-2M$. Since the first inequality holds for any geodesic ray $g_{1}\in {\cal L}_{T}(y, b)$, we have
\begin{align*}
\alpha_{b}(g(n), y)&=\limsup_{m\rightarrow \infty}(d(g(n), G(y, b)_{m})-m)\\
                   &\geq d(g(n), y)-2M.
\end{align*}
Similarly, we have
\[\alpha_{c}(g(n), y)\geq d(g(n), y)-2M.\]
On the other hand, it follows that
\begin{align*}
\alpha_{a}(g(n), y)&=\limsup_{m\rightarrow \infty}(d(g(n), G(y, a)_{m})-m)\\
                   &\geq -d(g(n), y)
\end{align*}
because $d(y, y')=m$ for any $y'\in G(y, a)_{m}$. Thus, we have
\[F_{abc}^{y}(g(n))\geq d(g(n), y)-4M.\]
Similarly, we can find constants $M'$ and $M''$ satisfying
\begin{align*}
F_{abc}^{y}(g'(n))&\geq d(g'(n), y)-4M',\\  
F_{abc}^{y}(g''(n))&\geq d(g''(n), y)-4M'',
\end{align*}
for any geodesic rays $g'\in {\cal L}_{T}(z, b)$, $g''\in {\cal L}_{T}(z, c)$ and $n\in {\Bbb N}$. We may replace $\max \{ M, M', M''\}$ by $M$ for simplicity. Then we have
\[F_{abc}^{y}(v)\geq d(v, y)-4M\]
for any $v\in G(z, a)\cup G(z, b)\cup G(z, c)$.  

Define the constant
\[N=\max \{ d(z, G(a, b)), d(z, G(b, c)), d(z, G(c, a))\}.\]
In general, we have 
\[\vert F_{abc}^{y}(z_{1})-F_{abc}^{y}(z_{2})\vert \leq 3d(z_{1}, z_{2})+288\delta \]
for any $z_{1}, z_{2}\in X$ by Lemma \ref{hyp-basic4} (iv) and Lemma \ref{beta}.

Let $w\in X(a, b, c)$ be any point. First, we assume that $w\in G(x, b)$ with $x\in G(a, b)$. Let $h$ be a geodesic in ${\cal L}_{T}(a, b)$ passing through $x$. Then we see that there exists $w'\in h$ such that $d(w, w')\leq \delta_{0}$ by Lemma \ref{lem:bow3} (i). It follows from Lemma \ref{lem:bow3} (i) that there exists a point $z'\in h$ with $d(z', z)\leq N+\delta_{0}$. By $\delta'$-slimness of the triangles $zz'a$ or $zz'b$ (see Lemma \ref{hyp-basic1}), we have $w_{0}\in G(z, a)\cup G(z, b)$ such that $d(w', w_{0})\leq N+\delta_{0}+\delta'$ and thus, $d(w, w_{0})\leq N+2\delta_{0}+\delta'$. Then
\begin{align*}
F_{abc}^{y}(w) & \geq F_{abc}^{y}(w_{0})-3(N+2\delta_{0}+\delta')-288\delta \\
               & \geq d(w_{0}, y)-4M-3(N+2\delta_{0}+\delta')-288\delta \\
               & \geq d(w, y)-4M-4(N+2\delta_{0}+\delta')-288\delta.
\end{align*}
Similarly, for any point $w$ in the set 
\[\left( \bigcup_{x\in G(b, c)}G(x, c)\right) \cup \left( \bigcup_{x\in G(c, a)}G(x, a)\right),\]
we have 
\[F_{abc}^{y}(w)\geq d(w, y)-4M-4(N+\delta_{0}+\delta')-288\delta. \]
Hence, $M_{0}=4M+4(N+\delta_{0}+\delta')+288\delta$ is a desired constant.
\end{pf*}

\begin{defn}
We define 
\begin{align*}
MV_{abc}^{y}&=\min \{ F_{abc}^{y}(w): w\in X(a, b, c)\},\\
MS_{abc}^{y}&=\{ w\in X(a, b, c): F_{abc}^{y}(w)=MV_{abc}^{y}\} ,\\
MS_{abc}&=\bigcup_{y\in X}MS_{abc}^{y}.
\end{align*}
\end{defn}

It follows from Propositions \ref{bounded} and \ref{F-is-proper} that the set $MS_{abc}^{y}$ is non-empty and finite for any $y\in X$ and $(a, b, c)\in \delta X$. Therefore, $MS_{abc}$ is non-empty.

\begin{prop}
For any $(a, b, c)\in \delta X$, the set $MS_{abc}$ is finite.
\end{prop}

\begin{pf} 
This is almost the same as \cite[Lemma 6.5]{adams1}.

Fix $y\in X$. By Lemma \ref{beta-fund} (iii) and Lemma \ref{beta}, we have
\[\vert \alpha_{a}(x, y)+\alpha_{a}(y, z)-\alpha_{a}(x, z)\vert \leq 200\delta+24\delta=224\delta,\]
for any $x, z\in X$. We have similar inequalities for $\alpha_{b}$ and $\alpha_{c}$. Summing these three inequalities, we have
\[|F_{abc}^{y}(x)+F_{abc}^{z}(y)-F_{abc}^{z}(x)|\leq 672\delta.\]
Define a function $G_{z}$ on $X$ by the formula
\[G_{z}=F_{abc}^{z}-F_{abc}^{z}(y),\]
and denote
\[m_{z}=\min \{ G_{z}(w): w\in X(a, b, c)\} =MV_{abc}^{z}-F_{abc}^{z}(y).\]
Then we have $MS_{abc}^{z}=\{ w\in X(a, b, c): G_{z}(w)=m_{z}\}$ and
\[|F_{abc}^{y}(x)-G_{z}(x)|\leq 672\delta \]
for any $x, z\in X$. This implies
\[|MV_{abc}^{y}-m_{z}|\leq 672\delta.\]
Hence, we have
\begin{align*}
MS_{abc}^{z} & \subseteq \{ w\in X(a, b, c): G_{z}(w)\leq MV_{abc}^{y}+672\delta \} \\
             & \subseteq \{ w\in X(a, b, c): F_{abc}^{y}(w)\leq MV_{abc}^{y}+1344\delta \}
\end{align*}
for any $z\in X$. This shows that 
\[MS_{abc}\subseteq \{ w\in X(a, b, c): F_{abc}^{y}(w)\leq MV_{abc}^{y}+1344\delta \}. \] 
The right hand side is bounded by Proposition \ref{F-is-proper} and thus, finite by Proposition \ref{bounded}. 
\end{pf}

\begin{defn}
Let ${\cal F}$ denote the set of all non-empty finite subsets of $X$. Define $MS\colon \delta X\rightarrow {\cal F}$\index{$MS$} by
\[MS(a, b, c)=MS_{abc}.\]
We call the map the {\it MIN set map}\index{MIN set map}. 
\end{defn}

Suppose that $\delta X$ is a standard Borel space and that ${\cal F}$ is equipped with the discrete Borel structure. We show that the map $MS$ is measurable. It can be seen as follows: first, note that the set 
\[S=\{ (w, (a, b, c))\in X\times \delta X: w\in X(a, b, c)\}\]
is a Borel subset of $X\times \delta X$ by the remark after Lemma \ref{chi-is-measurable}. 

\begin{lem}
Suppose that $p\colon Y\rightarrow Z$ is a Borel map between two standard Borel spaces such that $p^{-1}(z)$ is countable for any $z\in Z$. Let $f$ be a real-valued Borel function on $Y$. Then the function 
\[Z\ni z\mapsto \inf_{y\in p^{-1}(z)}f(y)\]
is measurable. 
\end{lem}

\begin{pf}
It is enough to show that the function
\[Z\ni z\mapsto \sup_{y\in p^{-1}(z)}f(y)\]
is measurable by considering $-f$ instead of $f$. Moreover, we may assume that $f$ is the characteristic function on a Borel subset $A\subseteq Y$. It follows from Theorem \ref{thm-standard-borel-space} (iv) that the subset $p(A)$ is measurable. Since the above function is the characteristic function on $p(A)$, it completes the proof.     
\end{pf}

Applying the above lemma for the projection $S\rightarrow \delta X$ onto the second coordinate, we see that the function $MV^{y}$ on $\delta X$ is measurable for any $y\in X$ and that the set
\[\{ (w, (a, b, c))\in S: MV_{abc}^{y}=F_{abc}^{y}(w)\}\]
is measurable in $S$. Thus, the map $MS^{y}\colon \delta X\rightarrow {\cal F}$ is measurable for any $y\in X$ and so is the map $MS$.

If $X$ is (the 1-skeleton of) the curve complex $C$ of a surface $M$ with $\kappa(M)\geq 0$, then the map $MS$ is equivariant for the actions of the mapping class group $\Gamma(M)$ by definition.

We mainly use the following map in what follows. Let ${\cal F}'$ be the set of all non-empty finite subsets of (the vertex set of) $C$ whose diameters are more than or equal to 3. The set ${\cal F}'$ is an invariant subset of ${\cal F}$ for the action of $\Gamma(M)$. Define $MS'\colon \delta C\rightarrow {\cal F}'$\index{$MS'$} by
\[MS'(a, b, c)=\{ x\in X(a, b, c): d(x, MS(a, b, c))\leq 3\}.\]
Note that the right hand side is finite by Proposition \ref{bounded}. The map $MS'$ is also measurable and equivariant for the actions of $\Gamma(M)$.


\section{Preliminaries on discrete measured equivalence relations}

In this section, we shall recall fundamentals of discrete measured equivalence relations.

\begin{defn}
Let $(X, \mu)$ be a standard Borel space with a finite positive measure. We say that a Borel subset ${\cal R}\subseteq X\times X$ is a {\it (non-singular) discrete measured equivalence relation}\index{discrete!measured equivalence relation}\index{non-singular!discrete measured equivalence relation} on $(X, \mu)$ if the following conditions are satisfied:
\begin{enumerate}
\renewcommand{\labelenumi}{\rm(\roman{enumi})}
\item ${\cal R}$ defines an equivalence relation on $X$, that is, if we write $x\sim y$ when $(x, y)\in {\cal R}$, then the relation $\sim$ on $X$ is an equivalence relation on $X$;
\item For any $x\in X$, the ${\cal R}$-equivalence class ${\cal R}x$\index{$R x$@${\cal R}x$} of $x$ is at most countable;
\item For any Borel subset $Y\subseteq X$ with $\mu(Y)=0$, we have $\mu({\cal R}Y)=0$, where 
\[{\cal R}A=\bigcup_{x\in A}{\cal R}x\]\index{$R A$@${\cal R}A$}
for a Borel subset $A\subseteq X$. The set ${\cal R}A$ is called the ${\cal R}$-{\it saturation}\index{R-saturation@${\cal R}$-!saturation} of $A$. 
\end{enumerate}
\end{defn}

Remark that the set ${\cal R}A$ is a Borel subset of $X$ for any Borel subset $A$ (use Theorem \ref{thm-standard-borel-space} (iv)). A Borel subset of $X$ with the form ${\cal R}A$ for some Borel subset $A$ is said to be {\it ${\cal R}$-invariant}\index{R-invariant@${\cal R}$-!invariant}. If any ${\cal R}$-invariant Borel subset of $X$ is either null or conull, then the relation ${\cal R}$ is said to be {\it ergodic}\index{ergodic relation}.

Let $[{\cal R}]$\index{$R ]$@$[{\cal R}]$} denote the set of all Borel automorphisms $g$ on $X$ such that $(g(x), x)\in {\cal R}$ for a.e.\ $x\in X$ and $[[{\cal R}]]$\index{$R ]] $@$[[{\cal R}]]$} denote the set of all Borel isomorphisms $g$ between two Borel subsets $A, B\subseteq X$ such that $(g(x), x)\in {\cal R}$ for a.e.\ $x\in A$.

If there exists an equivalent $\sigma$-finite measure $\nu$ to $\mu$ satisfying $\nu(A)=\nu(B)$ for any element $g\colon A\rightarrow B$ in $[[{\cal R}]]$, then we say that ${\cal R}$ is {\it measure-preserving}\index{measure-preserving} or that the measure $\nu$ is {\it invariant}\index{invariant!measure} for ${\cal R}$.

We define a measure $\tilde{\mu}$ on ${\cal R}$ by the formula
\[\tilde{\mu}(C)=\int_{X}\sum_{y\in {\cal R}x}\chi_{C}(x, y)d\mu (x)\]
for any Borel subset $C$ of ${\cal R}$, where $\chi_{C}$ is the characteristic function on $C$.

\begin{ex}\label{ex-thm-fm}
Let $(X, \mu)$ be a standard Borel space with a finite positive measure. Suppose that a discrete group $G$ has a {\it non-singular}\index{non-singular!action} Borel action on $(X, \mu)$, that is, the action preserves the class of the measure $\mu$. Then the subset  
\[{\cal R}=\{ (gx, x)\in X\times X: x\in X,\ g\in G\}\]
defines a discrete measured equivalence relation. Conversely, any discrete measured equivalence relation can be obtained in this way for some discrete group $G$ and its action \cite[Theorem 1]{fm}.
\end{ex}

\begin{defn}
Let ${\cal R}$ and ${\cal S}$ be discrete measured equivalence relations on standard Borel spaces $(X, \mu)$ and $(Y, \nu)$ with finite positive measures, respectively. We say that two relations ${\cal R}$ and ${\cal S}$ are {\it isomorphic}\index{isomorphic} if there exist conull Borel subsets $X'\subseteq X$, $Y'\subseteq Y$ and a Borel isomorphism $f\colon X'\rightarrow Y'$ such that $f_{*}\mu$ and $\nu$ are equivalent and $(f(x), f(y))\in {\cal S}$ if and only if $(x, y)\in {\cal R}$ for any $x, y\in X'$.
\end{defn}

\begin{defn}\label{defn-restriction}
Let ${\cal R}$ be a discrete measured equivalence relation on $(X, \mu)$. For a Borel subset $Y$ with positive measure, the discrete measured equivalence relation ${\cal R}\cap (Y\times Y)$ on $(Y, \mu|_{Y})$ is called the {\it restriction}\index{restriction} of ${\cal R}$ to $(Y, \mu|_{Y})$ and is denoted by $({\cal R})_{Y}$\index{$R Y$@$({\cal R})_{Y}$}, where $\mu|_{Y}$ is the restriction of $\mu$ to $Y$.
\end{defn}

\begin{defn}
Let ${\cal R}$ and ${\cal S}$ be discrete measured equivalence relations on standard Borel spaces $(X, \mu)$ and $(Y, \nu)$ with finite positive measures, respectively. We say that two relations ${\cal R}$ and ${\cal S}$ are {\it weakly isomorphic}\index{weakly!isomorphic} if there exist Borel subsets $A\subseteq X$, $B\subseteq Y$ with positive measure such that ${\cal R}A=X$, ${\cal S}B=Y$ and two restrictions $({\cal R})_{A}$ to $(A, \mu|_{A})$ and $({\cal S})_{B}$ to $(B, \nu|_{B})$ are isomorphic.
\end{defn}

Given two non-singular Borel actions of discrete groups $G_{1}$, $G_{2}$ on standard Borel spaces $(X_{1}, \mu_{1})$, $(X_{2}, \mu_{2})$ with $\sigma$-finite measures, respectively, we can define (weakly) orbit equivalence between them as in Definition \ref{defn-oe}. If we denote by ${\cal R}_{1}$ and ${\cal R}_{2}$ the relations generated by the actions of $G_{1}$ and $G_{2}$, respectively, then ${\cal R}_{1}$ and ${\cal R}_{2}$ are (weakly) isomorphic if and only if the two actions are (weakly) orbit equivalent.

\begin{defn}
Let ${\cal R}$ be a discrete measured equivalence relation on a standard Borel space $(X, \mu)$ with a finite positive measure.
\begin{enumerate}
\item[(i)] The relation ${\cal R}$ is said to be {\it finite}\index{finite relation} if each equivalence class of ${\cal R}$ is at most finite.
\item[(ii)] We say that ${\cal R}$ is {\it recurrent}\index{recurrent relation} if for any Borel subset $Y\subseteq X$ with positive measure and a.e.\  $x\in Y$, the intersection ${\cal R}x\cap Y$ is infinite.
\item[(iii)] A Borel subset $F\subseteq X$ is called a {\it fundamental domain}\index{fundamental domain} for the relation ${\cal R}$ if the intersection ${\cal R}x\cap F$ consists of exactly one point for a.e.\ $x\in X$.
\end{enumerate}
\end{defn}

\begin{prop}[\ci{Lemmas 2.3, 2.12}{adams2}]
Suppose that ${\cal R}$ is a discrete measured equivalence relation on a standard Borel space $(X, \mu)$ with a finite positive measure.
\begin{enumerate}
\item[(i)] If the relation ${\cal R}$ is finite, then it has a fundamental domain.
\item[(ii)] There exists an essentially unique ${\cal R}$-invariant Borel subset $Y\subseteq X$ such that $({\cal R})_{Y}$ has a fundamental domain and $({\cal R})_{X\setminus Y}$ is recurrent.
\end{enumerate}
\end{prop}

\begin{defn}
Let ${\cal R}$ be a discrete measured equivalence relation on $(X, \mu)$. We say that the relation ${\cal R}$ is of {\it type I}\index{type!I} if it has a fundamental domain and is of {\it type II}\index{type!II} if it is recurrent and measure-preserving.

Suppose that ${\cal R}$ is of type II. If there exists an equivalent finite positive measure to $\mu$ which is invariant for ${\cal R}$, then it is said to be of {\it type} ${\it II}_{1}$\index{type II1@type!${\rm II}_{1}$}. 
\end{defn}

\begin{defn}
Let ${\cal R}$ be a discrete measured equivalence relation on $(X, \mu)$ and $G$ be a standard Borel group. A Borel map  
\[\rho \colon {\cal R}\rightarrow G\]
is called a {\it cocycle}\index{cocycle} if there exists a conull Borel subset $X'$ of $X$ such that for a.e.\ $x\in X$, the cocycle identity
\[\rho(x, y)\rho(y, z)=\rho(x, z)\]
is satisfied for any $y, z\in {\cal R}x\cap X'$. We can show that if $\rho \colon {\cal R}\rightarrow G$ is a cocycle and two relations ${\cal R}$ and ${\cal R}'$ are isomorphic, then we can define a cocycle from ${\cal R}'$ into $G$. 
\end{defn}

Remark that we can find an ${\cal R}$-invariant conull Borel subset $X'$ of $X$ such that the above cocycle identity is satisfied for any $x\in X'$ and any $y, z\in {\cal R}x$ as follows: it follows from the fact in Example \ref{ex-thm-fm} that there exists a discrete group $H$ and its non-singular Borel action on $(X, \mu)$ such that
\[{\cal R}=\{ (hx, x)\in X\times X: x\in X,\ h\in H\}.\]
For $(h_{1}, h_{2})\in H\times H$, put
\[X_{(h_{1}, h_{2})}=\{ x\in X: \rho(h_{1}h_{2}x, h_{2}x)\rho(h_{2}x, x)=\rho(h_{1}h_{2}x, x)\}.\]
Then $X_{(h_{1}, h_{2})}$ is conull. If we define two Borel subsets
\[X''=\bigcap_{(h_{1}, h_{2})\in H\times H}X_{(h_{1}, h_{2})},\ \ X'=\bigcap_{h\in H}hX'',\]
then $X'$ is a desired subset. 

In what follows, although we often need to take an ${\cal R}$-invariant Borel subset of $X$ in this way, we do not mention it.

\begin{ex}
Let ${\cal R}$ be a discrete measured equivalence relation on $(X, \mu)$ given by an essentially free, non-singular Borel action of a discrete group $G$, where the word ``{\it essentially free}''\index{essentially free action} means that the Borel subset 
\[X_{g}=\{ x\in X: gx=x\}\]
is null for any $g\in G$. Then $X'=X\setminus (\bigcup_{g\in G}X_{g})$ is a $G$-invariant conull Borel subset on which $G$ acts freely. Thus, the formula
\[\rho(gx, x)=g\]
for $x\in X'$ and $g\in G$ defines a Borel cocycle $\rho \colon {\cal R}\rightarrow G$.

Note that by Furman's celebrated work \cite[Theorem D]{furman2}, there exist very concrete examples of an equivalence relation of type ${\rm II}_{1}$ which can not be (weakly) isomorphic to a relation given by an essentially free, measure-preserving Borel action of any discrete group $G$. 
\end{ex}

\begin{defn}
Let ${\cal R}$ be a discrete measured equivalence relation on $(X, \mu)$ and $G$ be a standard Borel group. Suppose that we have a Borel cocycle $\rho \colon {\cal R}\rightarrow G$ and a Borel space $S$ with a Borel $G$-action. Let ${\cal S}$ be a subrelation of ${\cal R}$. Then a Borel map $\varphi \colon X\rightarrow S$ is said to be {\it $\rho$-invariant}\index{rho-invariant Borel map@$\rho$-invariant!Borel map} for ${\cal S}$ if it satisfies
\[\rho(x, y)\varphi(y)=\varphi(x)\]
for a.e.\ $(x, y)\in {\cal S}$.  
\end{defn}

The following lemma is useful to extend a $\rho$-invariant Borel map as in the proof of Lemma \ref{inv-lem}:

\begin{lem}[\ci{Lemma 2.7}{adams2}]\label{lem-adams-map}
Let ${\cal R}$ be a discrete measured equivalence relation on $(X, \mu)$. Let $Y\subseteq X$ be a Borel subset with ${\cal R}Y=X$. Then there exists a Borel map $f\colon X\rightarrow Y$ satisfying the following conditions:
\begin{enumerate}
\item[(i)] $(f(x), x)\in {\cal R}$ for a.e.\ $x\in X$; 
\item[(ii)] $f(x)=x$ for any $x\in Y$.
\end{enumerate}
\end{lem}

\begin{pf}
As mentioned in Example \ref{ex-thm-fm}, we have a countable subgroup $G$ of $[{\cal R}]$ generating the relation ${\cal R}$. Fix an order on $G$ such that the identity in $G$ is minimal. For $x\in X$, let $g_{x}=\min \{ g\in G: g(x)\in Y\}$ and define $f(x)=g_{x}(x)$, which is a desired map. 
\end{pf}

\begin{lem}\label{inv-lem}
Let ${\cal R}$ be a discrete measured equivalence relation on $(X, \mu)$ and $G$ be a standard Borel group. Suppose that we have a Borel cocycle $\rho \colon {\cal R}\rightarrow G$ and a Borel space $S$ with a Borel $G$-action.
\begin{enumerate}
\renewcommand{\labelenumi}{\rm(\roman{enumi})}
\item Suppose that we have a sequence $\{ A_{n}\}_{n\in {\Bbb N}}$ of Borel subsets of $X$ such that $({\cal R})_{A_{n}}$ has a $\rho$-invariant Borel map $\varphi_{n}\colon A_{n}\rightarrow S$ for each $n\in {\Bbb N}$. Then $({\cal R})_{A}$ also has a $\rho$-invariant Borel map, where $A=\bigcup A_{n}$;
\item Let ${\cal M}$ be the set of all Borel subsets $A$ of $X$ with positive measure such that $({\cal R})_{A}$ has a $\rho$-invariant Borel map $A\rightarrow S$. Suppose that ${\cal M}$ is non-empty. Define 
\[m=\sup_{A\in {\cal M}}\mu(A).\]
Then there exists an ${\cal R}$-invariant Borel subset $B\in {\cal M}$ such that $m=\mu(B)$.  
\end{enumerate}
\end{lem}

\begin{pf}
For the assertion (i), we may assume that each $A_{n}$ is ${\cal R}$-invariant because $\varphi_{n}$ can be extendable to ${\cal R}A_{n}$, using a Borel map in Lemma \ref{lem-adams-map} as follows: define a $\rho$-invariant Borel map $\varphi_{n}'\colon {\cal R}A_{n}\rightarrow S$ by the formula  
\[\varphi_{n}'(x)=\rho(x, f(x))\varphi_{n}(f(x)),\]
where $f\colon {\cal R}A_{n}\rightarrow A_{n}$ is a Borel map such that
\begin{enumerate}
\item[(a)] $(f(x), x)\in {\cal R}$ for a.e.\ $x\in {\cal R}A_{n}$;
\item[(b)] $f(x)=x$ for any $x\in A_{n}$.
\end{enumerate}

Note that we have 
\[A= \bigsqcup_{n\in {\Bbb N}}\left( A_{n}\setminus \left( \bigcup_{m=1}^{n-1}A_{m}\right) \right).\]
Let us denote 
\[B_{1}=A_{1},\ \ B_{n}=A_{n}\setminus \left( \bigcup_{m=1}^{n-1}A_{m}\right),\ n\geq 2\]
and define $\varphi(x)=\varphi_{n}(x)$ for $x\in B_{n}$ and $n\in {\Bbb N}$. Then $\varphi \colon A\rightarrow S$ is $\rho$-invariant for $({\cal R})_{A}$. 

For the assertion (ii), there exists a sequence $\{ B_{n}\}$ of Borel subsets in ${\cal M}$ with $\mu(B_{n})\rightarrow m$. It follows from (i) that $B= \bigcup {\cal S}B_{n}$ is a desired subset.
\end{pf}

Next, we shall give a definition of amenability for discrete measured equivalence relations. Let ${\cal R}$ be a discrete measured equivalence relation on a standard Borel space $(X, \mu)$ with a finite positive measure. Suppose that $E$ is a separable Banach space. We denote the group of isometric automorphisms of $E$ by ${\rm Iso}(E)$\index{$Iso E$@${\rm Iso}(E)$} given the strong operator topology and denote the closed unit ball of the dual $E^{*}$ of $E$ by $E_{1}^{*}$\index{$E 1 *$@$E_{1}^{*}$} given the weak* topology. Assume that we have a Borel cocycle $\alpha \colon {\cal R}\rightarrow {\rm Iso}(E)$. Remark that ${\rm Iso}(E)$ is a separable metrizable group and the induced Borel structure is standard \cite[Lemma 1.1]{zim1}.

Moreover, we assume that we have a compact convex subset $A_{x}\subseteq E_{1}^{*}$ for $x\in X$ such that $\{ (x, a)\in X\times E^{*}_{1}: x\in X, \ a\in A_{x}\}$ is a Borel subset of $X\times E^{*}_{1}$ and for a.e.\ $(x, y)\in {\cal R}$, we have 
\[\alpha^{*}(x, y)A_{y}=A_{x},\]
where $\alpha^{*}(x, y)$ means the adjoint cocycle $\alpha^{*}(x, y)=(\alpha(x, y)^{-1})^{*}$. We denote by $F(X, \{ A_{x}\})$ the set of all measurable maps $\varphi \colon X\rightarrow E_{1}^{*}$ such that $\varphi(x)\in A_{x}$ for a.e.\ $x\in X$. We call $F(X, \{ A_{x}\})$ an {\it affine ${\cal R}$-space}\index{affine!${\cal R}$-space}.

\begin{defn}[\ci{Definition 3.1}{zim0}]
Let ${\cal R}$ be a discrete measured equivalence relation on a standard Borel space $(X, \mu)$ with a finite positive measure. The relation ${\cal R}$ is said to be {\it amenable}\index{amenable!relation} if any affine ${\cal R}$-space $F(X, \{ A_{x}\})$ has a fixed point, that is, there exists an element $\varphi \in F(X, \{ A_{x}\})$ such that
\[\rho(x, y)\varphi(y)=\varphi(x)\]
for a.e.\ $(x, y)\in {\cal R}$.
\end{defn}

This definition is equivalent to that in \cite{cfw} stated in terms of invariant means. One can find the proof of the equivalence of these two definitions in \cite[Theorem 4.2.7]{ana} (see also Theorem \ref{amenable-and-fixed}). We collect some fundamental results on amenability.

\begin{prop}\label{prop-amenable-fund}
Let ${\cal R}$ be a discrete measured equivalence relation on a standard Borel space $(X, \mu)$ with a finite positive measure. 
\begin{enumerate}
\item[(i)] If ${\cal R}$ is of type I, then it is amenable.
\item[(ii)] Suppose that ${\cal R}$ is given by an essentially free, non-singular Borel action of a discrete group $G$. If $G$ is amenable, then ${\cal R}$ is amenable. Conversely, if ${\cal R}$ is amenable and of type ${\it II}_{1}$, then $G$ is amenable.
\item[(iii)] The restriction $({\cal R})_{Y}$ to a Borel subset $Y\subseteq X$ with positive measure is amenable if and only if the relation $({\cal R})_{{\cal R}Y}$ is amenable.
\item[(iv)] Any subrelation of an amenable equivalence relation is amenable.
\end{enumerate} 
\end{prop}

The assertions (i) and (iii) follow by definition and we refer the reader to \cite[Proposition 4.3.3]{zim2} for (ii). The assertion (iv) comes from the fact that any amenable equivalence relation can be expressed as the union of an increasing sequence of finite subrelations and conversely, such an equivalence relation is amenable \cite{cfw}.

\begin{lem}[\ci{Lemma 2.16}{adams2}]\label{lem-ame-rec-exist}
Let ${\cal R}$ be a discrete measured equivalence relation on $(X, \mu)$. If ${\cal R}$ is recurrent, then there exists a recurrent amenable subrelation of ${\cal R}$.
\end{lem}


\section{Reducible elements in the mapping class group}\label{section-reducible-elements}

We shall recall some results from a group-theoretic aspect of the mapping class group and its action on ${\cal PMF}$ developed in \cite{ivanov1}. Let $M$ be a compact orientable surface of type $(g, p)$ with $\kappa(M)\geq 0$. Let $\Gamma(M)$ be the mapping class group of $M$. 

\begin{defn}
We shall call a diffeomorphism $F\colon M\rightarrow M$ {\it pure}\index{pure!diffeomorphism} if for some closed one-dimensional submanifold $L$ of $M$ (which may be empty), the following conditions are satisfied:
\begin{enumerate}
\item[(i)] All components of $L$ are non-trivial (i.e., each of them does not deform on $M$ to a point) and non-peripheral and are pairwise non-isotopic; 
\item[(ii)] The diffeomorphism $F$ is fixed (i.e., is the identity) on $L$ and it does not rearrange the components of $M\setminus L$ and it induces on each component of $M_{L}$ a diffeomorphism isotopic to either a pseudo-Anosov or the identity diffeomorphism. 
\end{enumerate}
Here, $M_{L}$ is the resulting surface of cutting $M$ along $L$. An element of $\Gamma(M)$ which contains a pure diffeomorphism is said to be {\it pure}\index{pure!element}. A subgroup of $\Gamma(M)$ consisting only pure elements is said to be {\it pure}\index{pure!subgroup}.
\end{defn}

Note that each pure element other than the trivial one has infinite order (see Proposition \ref{lem-blm}).

\begin{thm}[\ci{Corollaries 1.5, 1.8, 3.6}{ivanov1}]\label{pure-important}
\begin{enumerate}
\renewcommand{\labelenumi}{\rm(\roman{enumi})}
\item For any integer $m\geq 3$, the kernel $\Gamma(M; m)$\index{$\ C M m $@$\Gamma(M;m)$} of the natural $\Gamma(M)$-action on the first homology group $H_{1}(M; {\Bbb Z}/{m\Bbb Z})$ is torsion-free and consists of pure elements.
\item If $g\in \Gamma(M)$ is pure and $F\in {\cal PMF}$, then $g^{n}F=F$ for some $n\in {\Bbb Z}\setminus \{ 0\}$ implies $gF=F$.
\end{enumerate}
\end{thm}
Therefore, there exists a torsion-free subgroup of $\Gamma(M)$ with finite index consisting of pure elements.

\begin{cor}\label{permutation}
If $g\in \Gamma(M)$ is pure and $A$ is a finite subset of $S(M)$, then $gA=A$ implies $g\alpha =\alpha$ for each $\alpha \in V(C)$ such that there exists $\sigma \in A$ with $\alpha \in \sigma$.
\end{cor}

Recall that $S(M)$ is the set of all simplices of the curve complex $C=C(M)$ of $M$ if $\kappa(M)>0$ and $S(M)=V(C(M))$ if $\kappa(M)=0$. For $\sigma \in S(M)$, the union of some disjoint circles representing the vertices of $\sigma$ is called a {\it realization} of $\sigma$\index{realization}.

\begin{thm}
Suppose that $g\in \Gamma(M)$ and $\sigma \in S(M)$ satisfies $g\sigma =\sigma$. Let $L$ be a realization of $\sigma$. Then there exists a diffeomorphism $\varphi$ of $M$ in the class $g$ such that $\varphi(L)=L$.   
\end{thm}

This is a well-known result. We recommend the reader to see \cite[Theorem 5.2]{kos} for the proof.

In this case, $\varphi$ induces a diffeomorphism $\varphi_{L}$ on the cut surface $M_{L}$. Then $\varphi_{L}$ defines an isotopy class in the mapping class group $\Gamma(M_{L})$, which is defined similarly by the group of all isotopy classes of orientation-preserving self-diffeomorphisms on $M_{L}$. This induces a well-defined homomorphism $p_{\sigma}$ from the subgroup $\Gamma(M)_{\sigma}=\{ g\in \Gamma(M): g\sigma =\sigma \}$ into $\Gamma(M_{L})$.

\begin{prop}[\ci{Lemma 2.1 (1)}{blm}, \ci{Corollary 4.1.B}{ivanov2}]\label{lem-blm}
If one uses the above notation, then the kernel of the homomorphism $p_{\sigma}\colon \Gamma(M)_{\sigma}\rightarrow \Gamma(M_{L})$ for $\sigma \in S(M)$ is the free abelian group generated by the Dehn twists of all curves in $\sigma$ (see below for the definition of Dehn twists).
\end{prop}

\begin{thm}[\ci{Theorem 1.2}{ivanov1}]\label{comp-leave}
Let $g\in \Gamma(M; m)$ for some $m\geq 3$ and suppose that $\sigma \in S(M)$ satisfies $g\sigma =\sigma$. Let $L$ be a realization of $\sigma$ and $\varphi$ be a diffeomorphism on $M$ in the isotopy class $g$ such that $\varphi(L)=L$. Then $\varphi$ leaves all components of $L \cup \partial M$ and $M_{L}$.
\end{thm}

\begin{rem}\label{rem-comp-leave}
In particular, the image of the subgroup $\{ g\in \Gamma(M; m): g\sigma =\sigma \}$ for $m\geq 3$ under the homomorphism $p_{\sigma}$ is contained in the subgroup $\prod_{i}\Gamma(Q_{i})$ of $\Gamma(M_{L})$, where $\{ Q_{i}\}_{i}$ is the set of all components of $M_{L}$.
\end{rem}

The following corollary follows from Lemma \ref{lem-pants-fund}: 

\begin{cor}\label{cor-pants-stabilizer}
If $\sigma \in S(M)$ satisfies $|\sigma |=3g+p-3$ and we have $g\in \Gamma(M; m)$ for some $m\geq 3$ with $g\sigma =\sigma$, then a realization $L$ of $\sigma$ decomposes $M$ into pairs of pants and $g$ is in the free abelian group generated by the Dehn twists of all curves in $\sigma$.
\end{cor}

The Dehn twist about an element of $V(C)$ is an example of pure elements. We shall recall the definition of the Dehn twist around a simple closed curve on $M$.

We start with a description of a standard twist diffeomorphism of the annulus. Let $A$ be the annulus in ${\Bbb R}^{2}$ given by the inequality $1\leq r\leq 2$ in the standard polar coordinates $(r, \theta)$ in ${\Bbb R}^{2}$. Its boundary $\partial A$ consists of two components $\partial_{1}A$, $\partial_{2}A$ defined by the equations $r=1$, $r=2$, respectively. Let us fix a smooth function $\varphi \colon {\Bbb R}\rightarrow {\Bbb R}$ such that $\varphi(x)=0$ for $x\leq 1$, $\varphi(x)=2\pi$ for $x\geq 2$, $0\leq \varphi \leq 2\pi$ and the derivative $\varphi'$ is $\geq 0$. Let us define the diffeomorphism $T\colon A\rightarrow A$ by the formula
\[T(r, \theta)=(r, \theta+\varphi(r)),\]
which is called the {\it standard twist diffeomorphism}\index{standard!twist diffeomorphism} of the annulus $A$.

Let $e\colon A\rightarrow M$ be an orientation-preserving embedding of $A$ into a surface $M$. Take a diffeomorphism $e\circ T\circ e^{-1}\colon e(A)\rightarrow e(A)$ and extend it to a diffeomorphism $T_{e}\colon M\rightarrow M$ by letting it be the identity map on the complement of $e(A)$. The isotopy class of $T_{e}$ depends only on the isotopy class of the (unoriented) simple closed curve $C=e(a)$ on $M$ as far as the embedding $e$ is orientation-preserving, where $a=\{(r, \theta)\in A:r=3/2\}$. By a slight abuse of the language, we call this diffeomorphism a {\it twist diffeomorphism}\index{twist diffeomorphism} about the curve $C$ and call the isotopy class of the diffeomorphism the {\it Dehn twist\index{Dehn twist}} about the curve $C$. Clearly, this isotopy class depends only on the isotopy class of the curve $C$. Moreover, if the curve $C$ is in the trivial isotopy class or isotopic to one of the components in $\partial M$, then the Dehn twist about it is trivial. The reader should be referred to \cite[Section 4.1]{ivanov2} for more details.

Let $\sigma =\{ \alpha_{1}, \ldots ,\alpha_{m}\}$ be an element in $S(M)$ and $t_{i}$ be the isotopy class of the Dehn twist about the curve $\alpha_{i}$. Ivanov \cite{ivanov1} uses the following lemma to study a dynamical behavior of the action of a Dehn twist on ${\cal PMF}$.

\begin{lem}[\ci{Lemma 4.2}{ivanov1}]\label{lem-dehn-criterion}
Let $n_{1}, \ldots, n_{m}\in {\Bbb Z}$ and $t=t_{1}^{n_{1}}\cdots t_{m}^{n_{m}}\in \Gamma(M)$ be the product of the Dehn twists about the curves in $\sigma \in S(M)$. Then for all $\beta, \gamma \in {\cal MF}$, we have 
\begin{align*}
\sum_{i=1}^{m}(|n_{i}|-2)i(\gamma, \alpha_{i})i(\alpha_{i}, \beta)-i(\gamma, \beta)&\leq i(t(\gamma), \beta)\\
&\leq \sum_{i=1}^{m}|n_{i}|i(\gamma, \alpha_{i})i(\alpha_{i}, \beta)+i(\gamma, \beta).
\end{align*}
\end{lem}

\begin{cor}\label{cor-intersecting-not-invariant}
Let $\alpha, \beta \in V(C)$ with $i(\alpha, \beta)\neq 0$. We denote by $t_{\alpha}\in \Gamma(M)$ the Dehn twist about $\alpha$. Then $t_{\alpha}^{n}\beta \neq \beta$ for any $n\in {\Bbb Z}$ with $|n|\geq 3$.  
\end{cor}

\begin{cor}\label{trivial-stabilizer}
Let $\sigma$ be an element with $|\sigma |=3g+p-3$ and suppose that we have an element $\beta_{\alpha}\in V(C)$ for each $\alpha \in \sigma$ such that $i(\alpha, \beta_{\alpha})\neq 0$. Then the subgroup 
\[\{ g\in \Gamma(M): g\alpha =\alpha,\ g\beta_{\alpha}=\beta_{\alpha} {\rm \ for \ all \ } \alpha \in \sigma \} \]
of $\Gamma(M)$ is finite.
\end{cor}

\begin{pf}
Let $\Gamma$ be the above subgroup. Note that $M$ is decomposed into pairs of pants by a realization of $\sigma$ and the mapping class group of a pair of pants is finite (see Lemma \ref{lem-pants-fund}). By Proposition \ref{lem-blm}, the free abelian group generated by the Dehn twists around all curves in $\sigma$ is a subgroup of finite index of the group
\[\{ g\in \Gamma(M): g\alpha =\alpha {\rm \ for \ all \ } \alpha \in \sigma \}. \]
Assume that $\sigma =\{ \alpha_{1}, \ldots ,\alpha_{m}\} \in S(M)$ and $t=t_{1}^{n_{1}}\cdots t_{m}^{n_{m}} \in \Gamma$, where $t_{i}$ is the Dehn twist around $\alpha_{i}$ and $n_{i}\in {\Bbb Z}$. For each $\alpha_{j} \in \sigma$ and $k\in {\Bbb N}$, apply Lemma \ref{lem-dehn-criterion} to $t^{k}=t_{1}^{n_{1}k}\cdots t_{m}^{n_{m}k}$ with $\beta =\gamma =\beta_{\alpha_{j}}$. We have
\[\sum_{i=1}^{m}(|n_{i}k|-2)i(\beta_{\alpha_{j}}, \alpha_{i})^{2}\leq i(t^{k}(\beta_{\alpha_{j}}), \beta_{\alpha_{j}})=0\]
since $t\in \Gamma$ fixes $\beta_{\alpha_{j}}$. This inequality deduces a contradiction for a sufficiently large $k$ unless $n_{j}=0$. This means that $\Gamma$ is finite.
\end{pf}


\section[Subrelations of the two types]{Subrelations of the two types: irreducible and amenable ones and reducible ones}

In this section, we begin to study a discrete measured equivalence relation generated by an essentially free, non-singular Borel action of the mapping class group on a standard Borel space with a finite positive measure. We will show that its subrelations can be classified by their invariant Borel maps for the induced cocycle into the space of probability measures on ${\cal PMF}$. When the subrelations have invariant Borel maps, they can be classified into the two types, irreducible and amenable ones and reducible ones (or mixed ones), by the property of invariant Borel maps.

\subsection{Irreducible and amenable subrelations}

We consider the following assumption:

\begin{assumption}\label{assumption-star}\index{$(\ast $@$(\ast)$}
We call the following assumption $(\ast)$: let $\Gamma$ be a subgroup of $\Gamma(M; m)$, where $M$ is a surface with $\kappa(M)\geq 0$ and $m\geq 3$ is an integer. Let $(Y, \mu)$ be a standard Borel space with a finite positive measure and ${\cal S}$ be a recurrent discrete measured equivalence relation on $(Y, \mu)$. Suppose that we have a Borel cocycle
\[\rho \colon {\cal S}\rightarrow \Gamma \]
with finite kernel. It means that the kernel of $\rho$
\[\{ (x, y)\in {\cal S}: \rho(x, y)=e\} \]
is a finite subrelation of ${\cal S}$. 
\end{assumption}

Let $C$ be (the vertex set of) the curve complex of a surface $M$ with $\kappa(M)\geq 0$ and $\Gamma(M)$ be its mapping class group.

\begin{lem}\label{finite}
Let ${\cal F}'$ be the set of all non-empty finite subsets of $C$ whose diameters are more than or equal to $3$. Let $M({\cal F}')$ be the space of probability measures on ${\cal F}'$ with the $\ell^{1}$-norm topology. Then there exists a $\Gamma(M)$-equivariant Borel map $G_{0}\colon M({\cal F}')\rightarrow {\cal F}'$.
\end{lem}

\begin{pf}
The construction appears in the proof of \cite[Proposition 5.1]{adams2}. 

For each $\eta \in M({\cal F}')$, define 
\[G_{0}(\eta)=\bigcup \{ F\in {\cal F}': \eta(F)=\max_{{\cal F}'}\eta \}, \]
which is a desired map. Note that the $\ell^{1}$-norm topology and the pointwise convergence topology on the space $M({\cal F}')$ coincide.
\end{pf}

\begin{rem}\label{rem-borel-str-measurability}
In this remark, we introduce a Borel structure on the set $M(\partial C)$ (resp. $M(\delta C)$) of all probability measures on $\partial C$ (resp. $\delta C$) and prove measurability of the induced map $MS'_{*}\colon M(\delta C)\rightarrow M({\cal F}')$.

Recall that we have the $\Gamma(M)$-equivariant map $\pi \colon {\cal MIN}\rightarrow \partial C$ (see Chapter \ref{chapter:amenable-action}, Section \ref{boundary-of-curve-complex}). Let $M({\cal MIN})$ (resp. $M(({\cal MIN})^{3})$) denote the subset of $M({\cal PMF})$ (resp. $M(({\cal PMF})^{3})$) consisting of all probability measures whose supports are contained in ${\cal MIN}$ (resp. $({\cal MIN})^{3}$), which is a Borel subset of $M({\cal PMF})$ (resp. $M(({\cal PMF})^{3})$) by Corollary \ref{cor-third-measurable}. 

Let $M(\partial C)$ (resp. $M((\partial C)^{3})$) denote the space of all probability measures on $\partial C$ (resp. $(\partial C)^{3}$) with the Borel structure induced by the map $\pi_{*}\colon M({\cal MIN})\rightarrow M(\partial C)$ (resp. $\pi_{*}\colon M({\cal MIN}^{3})\rightarrow M((\partial C)^{3})$). Remark that $\pi_{*}$ is surjective since there exists a Borel section $S$ of $\pi$ such that the restriction of $\pi$ to $S$ is Borel isomorphism (see Proposition \ref{boundary-standard}). Let $p\colon \partial C\rightarrow S$ be the inverse Borel map. Then we get the Borel map $p \circ \pi \colon {\cal MIN}\rightarrow S$. Extending it to any Borel map from ${\cal PMF}$ into ${\cal PMF}$ and applying Proposition \ref{prop-first-measurable}, we can see that the induced map
\[p_{*}\circ \pi_{*}\colon M({\cal MIN})\rightarrow M(S)\]
is measurable. It follows that the map $p_{*}\colon M(\partial C)\rightarrow M(S)$ is measurable. Since the inverse of $p_{*}$ is the restriction of $\pi_{*}$, we see that $p_{*}$ is a Borel isomorphism. In particular, the Borel space $M(\partial C)$ is standard. Similarly, we can see that the Borel spaces $M(S^{3})$ and $M((\partial C)^{3})$ are Borel isomorphic under the map $\pi_{*}$.

Denote by $M(\delta C)$ the subset of $M((\partial C)^{3})$ consisting of all probability measures on $(\partial C)^{3}$ whose supports are contained in $\delta C$. We can see that $M(\delta C)$ is a Borel subset of $M((\partial C)^{3})$ by the above Borel isomorphism between $M((\partial C)^{3})$ and $M(S^{3})$ and by Corollary \ref{cor-third-measurable}. 

Composing two maps 
\begin{align*}
MS_{*}'&\colon M(\delta C)\rightarrow M({\cal F}'),\\  
G_{0}&\colon M({\cal F}')\rightarrow {\cal F}'
\end{align*}
constructed in Subsection \ref{MIN-set-map} and Lemma \ref{finite}, we get a $\Gamma(M)$-equivariant map
\[G_{1}\colon M(\delta C)\rightarrow {\cal F}'.\]
We can see that the map $MS_{*}'\circ \pi_{*}\colon M(\pi^{-1}(\delta C))\rightarrow M({\cal F}')$ is measurable by extending $MS'\circ \pi \colon \pi^{-1}(\delta C)\rightarrow {\cal F}'$ to a Borel map from ${\cal PMF}$ and applying Proposition \ref{prop-first-measurable}. This means that the map $MS'_{*}$ is measurable and so is $G_{1}$.
\end{rem}

\begin{prop}\label{invariant-amenable}
With the assumption $(\ast)$, suppose that there exists a $\rho$-invariant Borel map $\varphi \colon Y\rightarrow M(\partial C)$. Then the cardinality of the support of $\varphi(x)\in M(\partial C)$ is at most two for a.e.\ $x\in Y$. 
\end{prop}

\begin{pf}
Assume that the statement does not hold. Taking the restriction of $\varphi$ to an ${\cal S}$-invariant Borel subset of $Y$, we may assume that 
\[|{\rm supp}(\varphi(x))|\geq 3\]
for any $x\in Y$, where ${\rm supp}(\nu)$\index{$supp \ m$@${\rm supp}(\nu)$} denotes the support of a measure $\nu$. Then we would have a $\rho$-invariant Borel map
\[Y\rightarrow M(\delta C)\]
associating $x\in Y$ to the normalization of $(\varphi(x)\times \varphi(x)\times \varphi(x))\vert_{\delta C}$. Composing this map with $G_{1}\colon M(\delta C)\rightarrow {\cal F}'$, we get a $\rho$-invariant Borel map
\[\varphi' \colon Y\rightarrow {\cal F}'.\]
Let $F\in {\cal F}'$ be a subset of $C$ with $\mu((\varphi')^{-1}(F))>0$. Since $F$ has two points whose distance is more than or equal to $3$, the stabilizer of $F$ is finite by Lemma \ref{fill-curves-fix}. Since for a.e.\ $(x, y)\in ({\cal S})_{(\varphi')^{-1}(F)}$, the element $\rho(x, y)\in \Gamma$ fixes any point in $F$ by Corollary \ref{permutation}, the element $\rho(x, y)\in \Gamma$ is trivial by the fact that $\Gamma$ is torsion-free. This contradicts the recurrence of ${\cal S}$ because the kernel of $\rho$ is finite.   
\end{pf}

We denote 
\[\partial_{2}C=(\partial C\times \partial C)/\sim,\]\index{$\ z 2 C$@$\partial_{2}C$}
where the symbol $\sim$ means the equivalence relation given by the coordinate interchanging action of the symmetric group of two letters. Each element in $\partial_{2}C$ can be naturally identified with an atomic measure on $\partial C$ such that each atom has measure $1$ or $1/2$. Under this identification, we see that $\partial_{2}C$ is Borel isomorphic to a Borel subset of $M(\partial C)$, using a Borel isomorphism $M(S)\simeq M(\partial C)$ and Corollary \ref{cor-third-measurable}, where $S$ denotes a Borel section of $\pi \colon {\cal MIN}\rightarrow \partial C$ in Proposition \ref{boundary-standard}.

\begin{prop}\label{maximal}
Assume $(\ast)$.
\begin{enumerate}
\renewcommand{\labelenumi}{\rm(\roman{enumi})}
\item Suppose that there exists a $\rho$-invariant Borel map $Y\rightarrow M({\cal MIN})$. Then there exists a $\rho$-invariant Borel map $\varphi_{0}\colon Y\rightarrow M({\cal MIN})$ such that for any $\rho$-invariant Borel map $\varphi\colon Y\rightarrow M({\cal MIN})$, we have
\[{\rm supp}(\pi_{*}\varphi(x))\subseteq {\rm supp}(\pi_{*}\varphi_{0}(x))\]
for a.e.\ $x\in Y$.
\item Suppose that there exists a $\rho$-invariant Borel map $Y\rightarrow M(\partial C)$. Then there exists a $\rho$-invariant Borel map $\varphi_{0}\colon Y\rightarrow \partial_{2}C \subseteq M(\partial C)$ such that for any $\rho$-invariant Borel map $\varphi \colon Y\rightarrow M(\partial C)$, we have
\[{\rm supp}(\varphi(x))\subseteq {\rm supp}(\varphi_{0}(x))\]
for a.e.\ $x\in Y$.
\end{enumerate}
\end{prop}

\begin{pf}
The proof is essentially the same as that for \cite[Lemma 3.2]{adams2}.

Let $I$ denote the set of all $\rho$-invariant Borel maps $\varphi \colon Y\rightarrow M({\cal MIN})$. For each $\varphi \in I$, put
\begin{align*}
S(\varphi)&=\{ x\in Y: \vert {\rm supp}(\pi_{*}\varphi(x))\vert =2\}, \\
K&=\sup \{ \mu(S(\varphi)): \varphi \in I\}. 
\end{align*}
We claim that if $\{ \psi_{n}\}_{n\in {\Bbb N}}$ is a sequence of elements in $I$, then there exists an element $\psi \in I$ such that $S(\psi)=\bigcup S(\psi_{n})$. Let $S_{n}=S(\psi_{n})$ and $S_{n}'=S_{n}\setminus (S_{1}\cup \cdots \cup S_{n-1})$. Define $\psi \colon Y\rightarrow M({\cal MIN})$ by $\psi(x)=\psi_{n}(x)$ if $x\in S_{n}'$ and $\psi(x)=\psi_{0}(x)$ if $x\in Y\setminus (\bigcup S_{n}')$. Then $\psi \in I$ and $S(\psi)=\bigcup S(\psi_{n})$. 

It follows that there exists an element $\varphi_{0}\in I$ such that $K=\mu(S(\varphi_{0}))$. We will show that $\varphi_{0}$ is a desired map.

Let $\varphi \in I$. Remark that $S(\varphi)\setminus S(\varphi_{0})$ is a null set. It follows from Proposition \ref{invariant-amenable} that the cardinality of the support of the measure $\pi_{*}((\varphi(x)+\varphi_{0}(x))/2)$ is at most two for a.e.\ $x\in Y$, and we conclude that ${\rm supp}(\pi_{*}\varphi(x))\subseteq {\rm supp}(\pi_{*}\varphi_{0}(x))$ for a.e.\ $x\in S(\varphi_{0})$. 

Assume that there exists a Borel subset $A\subseteq Y\setminus S(\varphi_{0})$ with positive measure such that ${\rm supp}(\pi_{*}\varphi(x))\subsetneq {\rm supp}(\pi_{*}\varphi_{0}(x))$. We may assume that $A$ is contained in $Y\setminus S(\varphi)$ and is invariant for ${\cal S}$. Define $\varphi'\in I$ by $\varphi'(x)=(\varphi(x)+\varphi_{0}(x))/2$ for $x\in A$ and $\varphi'(x)=\varphi_{0}(x)$ for $x\in Y\setminus A$. Then $S(\varphi')=S(\varphi_{0})\cup A$ and $\mu(S(\varphi'))>K$, which is a contradiction.     

The assertion (ii) can be verified similarly by using Proposition \ref{invariant-amenable}.
\end{pf}

\begin{cor}\label{maximal-cor}
Assume $(\ast)$. 
\begin{enumerate}
\renewcommand{\labelenumi}{\rm(\roman{enumi})}
\item Suppose that there exists a $\rho$-invariant Borel map $Y\rightarrow M({\cal MIN})$. Then for any Borel subset $Y'$ of $Y$ with positive measure and any $\rho$-invariant map $\varphi \colon Y'\rightarrow M({\cal MIN})$ for $({\cal S})_{Y'}$, we have
\[{\rm supp}(\pi_{*}\varphi(x))\subseteq {\rm supp}(\pi_{*}\varphi_{0}(x))\]
for a.e.\ $x\in Y'$.
\item Suppose that there exists a $\rho$-invariant Borel map $Y\rightarrow M(\partial C)$. Then for any Borel subset $Y'$ of $Y$ with positive measure and any $\rho$-invariant map $\varphi \colon Y'\rightarrow M(\partial C)$ for $({\cal S})_{Y'}$, we have
\[{\rm supp}(\varphi(x))\subseteq {\rm supp}(\varphi_{0}(x))\]
for a.e.\ $x\in Y'$.
\end{enumerate}
\end{cor}

\begin{pf}
We show only the assertion (i) because the proof of (ii) is similar.

The map $\varphi \colon Y'\rightarrow M({\cal MIN})$ can be extended to the $\rho$-invariant Borel map $\varphi'\colon {\cal S}Y'\rightarrow M({\cal MIN})$ defined by  
\[x\mapsto \rho(x, f(x))\varphi(f(x)),\]
where $f\colon {\cal S}Y'\rightarrow Y'$ is a Borel map such that
\begin{enumerate}
\item[(a)] $(f(x), x)\in {\cal S}$ for a.e.\ $x\in {\cal S}Y'$;
\item[(b)] $f(x)=x$ for any $x\in Y'$
\end{enumerate}
(see Lemma \ref{lem-adams-map}). Define a map $\varphi''\colon Y\rightarrow M({\cal MIN})$ by
\begin{align*}
{\cal S}Y'\ni x&\mapsto \varphi'(x),\\
Y\setminus {\cal S}Y'\ni x&\mapsto \varphi_{0}(x).
\end{align*}
Then $\varphi''$ is $\rho$-invariant and we have 
\[{\rm supp}(\pi_{*}\varphi''(x))\subseteq {\rm supp}(\pi_{*}\varphi_{0}(x))\]
for a.e.\ $x\in Y$ by Proposition \ref{maximal}. Since $\varphi''(x)=\varphi(x)$ for $x\in Y'$, it completes the proof. 
\end{pf}

\begin{rem}
Under the assumption $(\ast)$, it follows from Proposition \ref{maximal} that there exists a unique $\rho$-invariant Borel map $\varphi_{\cal S}\colon Y\rightarrow \partial_{2}C$ satisfying the maximal property if there exists a $\rho$-invariant Borel map $Y\rightarrow M(\partial C)$. Corollary \ref{maximal-cor} means that $\varphi_{({\cal S})_{Y'}}$ is the restriction of $\varphi_{\cal S}$ to $Y'$ for any Borel subset $Y'$ of $Y$ with positive measure.
\end{rem}

\begin{lem}\label{amenable-lemma}
Let $\nu$ be a quasi-invariant probability measure on $\partial_{2}C$ for the $\Gamma(M)$-action. Then the $\Gamma(M)$-action on $(\partial_{2}C, \nu)$ is amenable.
\end{lem}

\begin{pf}
This proof is the same as that of \cite[Lemma 5.2]{adams2}. 

Let $F\subseteq \partial C\times \partial C$ be a fundamental domain for the action of the symmetric group of two letters on $\partial C\times \partial C$. Let $\sigma$ be the non-trivial element in the symmetric group. Define $F'=\sigma F$. Identifying $\partial_{2}C$ with $F$ (resp. $F'$) by way of the restriction of the projection $\partial C\times \partial C\rightarrow \partial_{2}C$, we identify the measure $\nu$ on $\partial_{2}C$ with a measure $\nu_{F}$ on $F$ (resp. $\nu_{F'}$ on $F'$). Let $\nu'=(\nu_{F}+\nu_{F'})/2$.     

Then the action of $\Gamma(M)$ on $(\partial C\times \partial C, \nu')$ is non-singular. Further, the natural map from $\partial C\times \partial C$ into $\partial_{2}C$ carries $\nu'$ to $\nu$. Since this map is $\Gamma(M)$-equivariant and has finite fibers, it suffices to show that the action of $\Gamma(M)$ on $(\partial C\times \partial C, \nu')$ is amenable.

Let $p\colon \partial C\times \partial C \rightarrow \partial C$ denote the projection onto the first coordinate. Since the map $p$ is $\Gamma(M)$-equivariant, by Theorem \ref{aeg-thm}, it suffices to show that the action of $\Gamma(M)$ on $(\partial C, p_{*}\nu')$ is amenable, which is proved in Theorem \ref{amenable-action-cc}. 
\end{pf}

Next, we consider the following assumption:

\begin{assumption}\index{$(\ast p $@$(\ast)'$}
We call the following assumption $(\ast)'$: let $(X, \mu)$ be a standard Borel space with a finite positive measure. Let $\Gamma$ be a subgroup of $\Gamma(M; m)$, where $M$ is a surface with $\kappa(M)\geq 0$ and $m\geq 3$ is an integer. Suppose that ${\cal R}$ is a discrete measured equivalence relation on $(X, \mu)$ and that we have a Borel cocycle
\[\rho \colon {\cal R}\rightarrow \Gamma\]
with finite kernel. Let $F$ be a fundamental domain of the finite equivalence relation $\ker \rho$. We assume that there exists an essentially free, non-singular Borel action of $\Gamma$ on $(F, \mu|_{F})$ generating the relation $({\cal R})_{F}$ and whose induced cocycle 
\[({\cal R})_{F}\rightarrow \Gamma, \ \ (gx, x)\mapsto g, \ g\in \Gamma, \ x\in F\]
is equal to the restriction of $\rho$ to $({\cal R})_{F}$. Let ${\cal S}$ be a recurrent subrelation of $({\cal R})_{Y}$ on $(Y, \mu|_{Y})$, where $Y$ is a Borel subset of $X$ with positive measure.
\end{assumption}

Remark that the recurrent relation ${\cal S}$ on $(Y, \mu|_{Y})$ in the above assumption satisfies $(\ast)$.

\begin{prop}\label{irreducible-amenable}
Assume  $(\ast)'$. 
\begin{enumerate}
\renewcommand{\labelenumi}{\rm(\roman{enumi})}
\item Suppose that there exists a $\rho$-invariant Borel map $\varphi \colon Y\rightarrow M({\cal MIN})$ for ${\cal S}$. Then the relation ${\cal S}$ is amenable.
\item Suppose that there exists a $\rho$-invariant Borel map $\varphi'\colon Y\rightarrow M(\partial C)$ for ${\cal S}$. Then the relation ${\cal S}$ is amenable.
\end{enumerate}
\end{prop}

\begin{pf}
We may assume that $F=X$ by Proposition \ref{prop-amenable-fund} (iii). The assertion (i) follows from (ii) by using $\pi \colon {\cal MIN}\rightarrow \partial C$. Hence, we show the assertion (ii).

Since ${\cal S}$ is recurrent, ${\rm supp}(\varphi'(x))$ consists of at most two points for a.e.\ $x\in Y$ by Proposition \ref{invariant-amenable}. It induces a $\rho$-invariant Borel map 
\[\varphi'\colon Y\rightarrow \partial_{2}C\]
for ${\cal S}$. We identify ${\cal S}$ with the relation on $(X, \mu)$ defined by the union
\[\{ (x, x)\in (X\setminus Y)\times (X\setminus Y)\} \cup {\cal S}\]
and extend $\varphi'$ to the map on $X$ defined by 
\[\varphi'(x)=a_{0}\]
for $x\in X\setminus Y$, where $a_{0}\in \partial_{2}C$ is a fixed point. Then the extended map $\varphi'$ is also $\rho$-invariant for ${\cal S}$. 

Consider the {\it Mackey range}\index{Mackey range} $S$ of the cocycle $\rho \colon {\cal S}\rightarrow \Gamma$ (see \cite[p.76]{zim2}). The space $S$ is a Borel $\Gamma$-space with a quasi-invariant positive measure and $L^{\infty}(S)$ can be identified with the weak* closed subalgebra 
\[\{ f\in L^{\infty}(X\times \Gamma): f(x, g)=f(y, \rho(y, x)g) \ {\rm for \ a.e.\ }(x, y)\in {\cal S} \ {\rm and} \ g\in \Gamma \}\]
of $L^{\infty}(X\times \Gamma)$ as $\Gamma$-modules, where the $\Gamma$-action on $X\times \Gamma$ is given by
\[g\cdot (x, g_{1})=(x, g_{1}g^{-1})\]
for $x\in Y$ and $g, g_{1}\in \Gamma$. The equivalence relation ${\cal R}_{1}$ on $X\times \Gamma$ defined by the formula
\[(x, g)\sim (y, \rho(y, x)g)\]
for $(x, y)\in {\cal R}$ and $g\in \Gamma$ is an equivalence relation of type I with the fundamental domain $\{ (x, e)\in X\times \Gamma :x\in X\}$. Recall that the relation ${\cal R}$ is given by an essentially free Borel action of $\Gamma$. The space $S$ is identified with the quotient space of $X\times \Gamma$ by the equivalence relation ${\cal R}_{2}$ defined by
\[(x, g)\sim (y, \rho(y, x)g)\]
for $(x, y)\in {\cal S}$ and $g\in \Gamma$ because ${\cal R}_{2}$ is a subrelation of the equivalence relation ${\cal R}_{1}$ and thus, of type I. The class of the measure on $S$ is the same as one on a fundamental domain $F_{2}$ for ${\cal R}_{2}$ when $S$ is identified with $F_{2}$. Denote the projection of $(x, g)\in X\times \Gamma$ by $[x, g]\in S$. Then ${\cal S}$ can be identified with a Borel subgroupoid 
\[{\cal G}=\{ ([x, e], \rho(x, y))\in S\rtimes \Gamma : (x, y)\in {\cal S}\}\]
of $S\rtimes \Gamma$ by the following lemmas (see Appendix \ref{general-amenable} for the definition of the groupoid $S\rtimes \Gamma$).

\begin{lem}
The subset ${\cal G}$ is a Borel subgroupoid of $S\rtimes \Gamma$ on the unit space
\[{\cal G}^{(0)}=\{ [x, e]\in S: x\in X\} \]
with the range and source maps and the inverse defined by
\[r([x, e], \rho(x, y))=[x, e], \ \ s([x, e], \rho(x, y))=\rho(x, y)^{-1}[x, e]=[y, e],\]
\[([x, e], \rho(x, y))^{-1}=([y, e], \rho(y, x)),\]
respectively, and the composition defined by
\[([x, e], \rho(x, y))\cdot ([y, e], \rho(y, z))=([x, e], \rho(x, z))\]
for $(x, y), (y, z)\in {\cal S}$.
\end{lem}

\begin{pf}
It is easy to show that ${\cal G}$ is a groupoid with the above operations. We must show that ${\cal G}$ is Borel. It follows from the equality
\[{\cal G}=\bigcup_{g\in \Gamma}\{ ([x, e], g)\in S\times \Gamma : x\in p_{1}(\rho^{-1}(g)\cap {\cal S})\} \] 
that ${\cal G}$ is a Borel subset of $S\times \Gamma$, where $p_{1}\colon {\cal S}\rightarrow X$ is the projection onto the first coordinate. Remark that $p_{1}$ and the projection $X\times \Gamma \rightarrow S$ has the property that the image of any Borel subset is a Borel subset (see Theorem \ref{thm-standard-borel-space} (iv)). It is also easy to see that the groupoid operations defined as above are measurable.
\end{pf}

\begin{lem}
The map $f\colon {\cal S}\rightarrow {\cal G}$ defined by the formula
\[f(x, y)=([x, e], \rho(x, y))\]
for $(x, y)\in {\cal S}$ is a Borel isomorphism such that the image of $\mu$ by $f$ and the measure on $S$ are equivalent.
\end{lem}

\begin{pf}
It is easy to prove that $f$ is a groupoid isomorphism. It is enough to show that $f^{-1}$ is a Borel map (see Theorem \ref{thm-standard-borel-space} (ii)). For a Borel subset $B\subseteq {\cal S}$, we see that
\[f(B)=\bigcup_{g\in \Gamma}\{ ([x, e], g)\in S\times \Gamma : x\in p_{1}(\rho^{-1}(g)\cap B)\}.\]
The right hand side is a Borel subset. The claim on the measure is clear since the projection $X\times \Gamma \rightarrow S$ carries a null set to a null set and the inverse image of a null set by the projection is also a null set.
\end{pf}

Return to the proof of Proposition \ref{irreducible-amenable}. Using the $\rho$-invariant Borel map $\varphi'\colon X\rightarrow \partial_{2}C$ for ${\cal S}$, we construct a Borel map $\varphi''\colon S\rightarrow \partial_{2}C$ by the formula
\[[x, g]\mapsto g^{-1}\varphi'(x).\]
Then $\varphi''$ is well-defined and $\Gamma$-equivariant. By Lemma \ref{amenable-lemma} and Theorem \ref{aeg-thm}, the groupoid $S\rtimes \Gamma$ is amenable. Since the groupoid ${\cal G}$ is principal, that is, the group 
\[{\cal G}_{w}^{w}=\{ \gamma \in {\cal G}: r(\gamma)=s(\gamma)=w\}\]
is trivial for a.e.\ $w\in S$, it is a subrelation of the induced equivalence relation on $S$
\[\{ (r(\gamma), s(\gamma))\in S\times S: \gamma \in S\rtimes \Gamma \}, \]
which is a quotient of the groupoid $S\rtimes \Gamma$ and is amenable because $S\rtimes \Gamma$ is amenable. Thus, the subrelation ${\cal S}$ is also amenable.  
\end{pf}

Finally, we note that a hyperbolic group also satisfies the same property as that of the mapping class group stated above. We can show the following proposition, using Adams' argument \cite{adams2} and applying the same argument as in the proof of Proposition \ref{irreducible-amenable}.

\begin{assumption}\index{$(\ast p h $@$(\ast)_{\rm h}'$}
We call the following assumption $(\ast)_{\rm h}'$: let $(X, \mu)$ be a standard Borel space with a finite positive measure. Let $\Gamma$ be a subgroup of a hyperbolic group $\Gamma_{0}$. Suppose that ${\cal R}$ is a discrete measured equivalence relation on $(X, \mu)$ and that we have a Borel cocycle
\[\rho \colon {\cal R}\rightarrow \Gamma\]
with finite kernel. Let $F$ be a fundamental domain of the finite equivalence relation $\ker \rho$. We assume that there exists an essentially free, non-singular Borel action of $\Gamma$ on $(F, \mu|_{F})$ whose induced cocycle 
\[({\cal R})_{F}\rightarrow \Gamma, \ \ (gx, x)\mapsto g, \ g\in \Gamma, \ x\in F\]
is equal to the restriction of $\rho$ to $({\cal R})_{F}$. Let ${\cal S}$ be a recurrent subrelation of $({\cal R})_{Y}$, where $Y$ is a Borel subset of $X$ with positive measure.
\end{assumption}

The following proposition is a slight generalization of \cite[Lemma 5.3]{adams2}. We use it in the following sections.

\begin{prop}\label{irreducible-amenable-hyp}
With the assumption $(\ast)_{\rm h}'$, suppose that there exists a $\rho$-invariant Borel map $\varphi \colon Y\rightarrow M(\partial \Gamma_{0})$ for ${\cal S}$. Then the relation ${\cal S}$ is amenable.
\end{prop}


\subsection{Reducible subrelations}

Let $M$ be a surface with $\kappa(M)\geq 0$. We denote $\Phi ={\cal PMF}\setminus {\cal MIN}$\index{$\ U=PMF MIN$@$\Phi ={\cal PMF}\setminus {\cal MIN}$} and define $S(M)$ to be the set of all simplices of the curve complex $C$ if $\kappa(M)>0$, and $S(M)=V(C(M))$ if $\kappa(M)=0$. First, we will construct a $\Gamma(M)$-equivariant Borel map $H\colon \Phi \rightarrow S(M)$. The following theorem says that every non-minimal measured foliation has at least one simple closed curve along the foliation with a certain additional property.

\begin{thm}[\ci{Corollary 2.15}{ivanov1}]\label{thm-H-essential}
For any $F\in \Phi$, there exists $\alpha \in V(C)$ satisfying the following: if $G\in {\cal PMF}$, then $i(G, \alpha)\neq 0$ implies $i(G, F)\neq 0$.
\end{thm}

\begin{defn}\label{construction-of-H}
For $F\in \Phi$, the subset of $V(C)$ consisting of the curves $\alpha$ satisfying the condition of Theorem \ref{thm-H-essential} is called the {\it canonical reduction system}\index{canonical reduction system!for a measured foliation} for $F$. Remark that it is an element in $S(M)$ because any element $\alpha$ of the canonical reduction system for $F$ satisfies $i(\alpha, F)=0$. Let $H \colon \Phi \rightarrow S(M)$\index{$H \ U S M$@$H\colon \Phi \rightarrow S(M)$} be the map associating the canonical reduction system. The $\Gamma(M)$-equivariance of $H$ is clear by definition.
\end{defn}

\begin{lem}
The map $H$ is measurable.
\end{lem}

\begin{pf}
For $\alpha \in V(C)$, if we denote 
\[A=\{ F\in {\cal PMF}: i(F, \alpha)\neq 0\},\]
then we see that
\begin{align*}
B=&\{ F\in  \Phi : \alpha {\rm \ is \ contained \ in \ the \ canonical \ reduction \ system \ for \ }F\} \\
=& \{ F\in \Phi : i(F, \cdot)\neq 0 \ {\rm on \ the \ open \ subset} \ A\ {\rm of}\ {\cal PMF}\}.
\end{align*}
Let us denote 
\[C=\{ (F, G)\in ({\cal PMF})^{2}: i(F, G)=0\}.\]
Then  
\[\Phi \setminus B=p(({\cal PMF}\times A)\cap C)\cap \Phi,\]
where $p\colon ({\cal PMF})^{2}\rightarrow {\cal PMF}$ is the projection onto the first coordinate. Since the intersection $({\cal PMF}\times A)\cap C$ can be expressed as the union of an increasing sequence of compact subsets in $({\cal PMF})^{2}$, we see that $\Phi \setminus B$ is a Borel subset.  
\end{pf}

\begin{lem}\label{partition}
With the assumption $(\ast)$, suppose that there exists a $\rho$-invariant Borel map
\[\varphi \colon Y\rightarrow M({\cal PMF})\]
such that $\varphi(y)(\Phi)=1$ for a.e.\ $y\in Y$. Then there exists a countable Borel partition $Y=\bigsqcup_{\sigma \in S(M)}Y_{\sigma}$ such that
\[\rho(({\cal S})_{Y_{\sigma}})\subseteq \Gamma_{\sigma},\]
where 
\[\Gamma_{\sigma}=\{ g\in \Gamma : g\sigma= \sigma \} =\{ g\in \Gamma : g\alpha =\alpha {\rm \ for \ any \ } \alpha \in \sigma \}\]
(see Corollary \ref{permutation}). 
\end{lem}

\begin{pf}
The Borel map $H\colon \Phi \rightarrow S(M)$ induces the Borel map $H_{*}\colon M(\Phi)\rightarrow M(S(M))$, where $M(\Phi)$ denotes the Borel subset of $M({\cal PMF})$ consisting of all probability measures on ${\cal PMF}$ whose supports are contained in $\Phi$ (see Proposition \ref{prop-first-measurable} and Corollary \ref{cor-third-measurable}). Then
\[H_{*}\varphi \colon Y\rightarrow M(S(M))\]
is a $\rho$-invariant Borel map. As in Lemma \ref{finite}, there exists a $\Gamma$-equivariant Borel map from $M(S(M))$ into ${\cal F}(S(M))$, the set of all non-empty finite subsets of $S(M)$, defined by
\[M(S(M))\ni \eta \mapsto \{ \sigma \in S(M): \eta(\sigma)=\max_{S(M)}\eta \} \in {\cal F}(S(M)).\]
Thus, we have a $\rho$-invariant Borel map
\[\psi \colon Y\rightarrow {\cal F}(S(M)).\]
Take $A\in {\cal F}(S(M))$ such that $\psi^{-1}(A)$ has positive measure. For a.e.\ $(x, y)\in {\cal S}$ with $x, y\in \psi^{-1}(A)$, we have $\rho(x, y)A=A$ by the $\rho$-invariance of $\psi$. By Corollary \ref{permutation}, the element $\rho(x, y)\in \Gamma$ fixes any curve in $A$. This completes the proof. 
\end{pf}

\begin{lem}\label{reducible}
With the assumption $(\ast)$, suppose that the cocycle $\rho \colon {\cal S}\rightarrow \Gamma$ is $\Gamma_{\sigma}$-valued for some $\sigma \in S(M)$, that is, $\rho(x, y)\in \Gamma_{\sigma}$ for a.e.\ $(x, y)\in {\cal S}$. 
\begin{enumerate}
\item[(i)] Any $\rho$-invariant Borel map $\psi \colon Y\rightarrow M({\cal PMF})$ satisfies $(\psi(y))(\Phi)=1$ for a.e.\ $y\in Y$.  
\item[(ii)] There exist no $\rho$-invariant Borel maps $Y\rightarrow \partial_{2}C$.
\end{enumerate}
\end{lem}

\begin{pf}
First, we define a $\Gamma(M)$-equivariant Borel map
\[G_{2}\colon C\times \partial_{2}C\rightarrow {\cal F}', \]
where the $\Gamma(M)$-action on $C\times \partial_{2}C$ is given by the formula
\[g\cdot (x, t)=(gx, gt)\]
for $g\in \Gamma(M)$, $x\in C$ and $t\in \partial_{2}C$. Recall that ${\cal F}'$ is the set of all non-empty finite subsets of $C$ whose diameters are more than or equal to 3. Define
\[G_{2}(x, t)=\{ x\}\cup \left( G(x, t)\cap \{ y\in C: d(y, x)=3\} \right)\]
 for $x\in C$ and $t\in \partial_{2}C$, where each element in $\partial_{2}C$ is identified with a subset of $\partial C$ whose cardinality is 1 or 2 and the symbol $G(x, t)$ means the union $\bigcup_{a\in t}G(x, a)$. It is easy to see that $G_{2}$ is a Borel map by Corollary \ref{measurable-cor}.

We show the assertion (i). Assume that we have a $\rho$-invariant Borel map
\[\psi \colon Y\rightarrow M({\cal PMF})\] 
such that there exists a Borel subset $Y'$ of $Y$ with positive measure satisfying $\psi(y)(\Phi)<1$ for any $y\in Y'$. We may assume that $Y'$ is an ${\cal S}$-invariant Borel subset and $\psi(y)({\cal MIN})=1$ for any $y\in Y'$. Then the map $\psi$ induces a $\rho$-invariant Borel map
\[\psi' \colon Y'\rightarrow \partial_{2}C\]  
for $({\cal S})_{Y'}$ by Proposition \ref{invariant-amenable}. Let $\alpha$ be any curve in $\sigma$, which is fixed by any element in $\Gamma_{\sigma}$ (see Corollary \ref{permutation}). Since $\rho$ is $\Gamma_{\sigma}$-valued, the Borel map $\psi''\colon Y'\rightarrow C\times \partial_{2}C$ defined by the formula
\[\psi''(y)=(\alpha, \psi'(y))\]
is $\rho$-invariant for $({\cal S})_{Y'}$. Composing this map with $G_{2}$, we obtain a $\rho$-invariant Borel map $Y'\rightarrow {\cal F}'$. This contradicts the recurrence of $({\cal S})_{Y'}$ as in the proof of Proposition \ref{invariant-amenable}.

The assertion (ii) has been verified in the above argument. 
\end{pf}

\begin{thm}\label{alternative}
With the assumption $(\ast)$, suppose that we admit a $\rho$-invariant Borel map $\varphi \colon Y\rightarrow M({\cal PMF})$ such that $\varphi(y)(\Phi)=1$ for a.e.\ $y\in Y$. If $Y'$ is a Borel subset of $Y$ with positive measure and $\psi \colon Y'\rightarrow M({\cal PMF})$ is a $\rho$-invariant Borel map for $({\cal S})_{Y'}$, then $\psi(y)(\Phi)=1$ for a.e.\ $y\in Y'$.
\end{thm}

\begin{pf}
Use Lemmas \ref{partition} and \ref{reducible}. 
\end{pf}

\begin{defn}
With the assumption $(\ast)$, the relation ${\cal S}$ is said to be {\it reducible}\index{reducible!subrelation} if it has a $\rho$-invariant Borel map $\varphi \colon Y\rightarrow M({\cal PMF})$ such that $\varphi(y)(\Phi)=1$ for a.e.\ $y\in Y$. 
\end{defn}

Conversely, the following corollary can be shown immediately:

\begin{cor}\label{alternative-cor}
With the assumption $(\ast)$, suppose that there exists a $\rho$-invariant Borel map $\varphi \colon Y\rightarrow M({\cal PMF})$ such that $\varphi(y)({\cal MIN})=1$ for a.e.\ $y\in Y$. If $Y'$ is a Borel subset of $Y$ with positive measure and $\psi \colon Y'\rightarrow M({\cal PMF})$ is a $\rho$-invariant Borel map for $({\cal S})_{Y'}$, then $\psi(y)({\cal MIN})=1$ for a.e.\ $y\in Y'$. 
\end{cor}

\begin{defn}
With the assumption $(\ast)$, we say that the subrelation ${\cal S}$ is {\it irreducible and amenable}\index{irreducible and amenable!subrelation} if it has a $\rho$-invariant Borel map $\varphi \colon Y\rightarrow M({\cal PMF})$ such that $\varphi(y)({\cal MIN})=1$ for a.e.\ $y\in Y$. 
\end{defn}

\begin{cor}\label{cor-irreducible-amenable-cor}
With the assumption $(\ast)'$, suppose that there exists a $\rho$-invariant Borel map $\varphi \colon Y\rightarrow \partial_{2}C$ for ${\cal S}$. Then ${\cal S}$ is irreducible and amenable.
\end{cor}

\begin{pf}
It follows from Proposition \ref{irreducible-amenable} (ii) that ${\cal S}$ is amenable. Thus, ${\cal S}$ has a $\rho$-invariant Borel map $\psi \colon Y\rightarrow M({\cal PMF})$ for ${\cal S}$. If $\psi(y)(\Phi)>0$ for $y$ in a Borel subset $Y'$ with positive measure, then there would exist a $\rho$-invariant Borel map $Y'\rightarrow S(M)$ for $({\cal S})_{Y'}$. This contradicts Lemma \ref{reducible} (ii).
\end{pf}

\begin{rem}\label{rem-disjoint-sum-irreducible-amenable-reducible}
Under the assumption $(\ast)$, if the relation ${\cal S}$ is amenable and any restriction of ${\cal S}$ to a Borel subset of $Y$ with positive measure is not reducible, then it is irreducible and amenable by definition. Under the assumption $(\ast)'$, the converse also holds by Proposition \ref{irreducible-amenable} and Corollary \ref{alternative-cor}.

Using Lemma \ref{inv-lem}, we see that under the assumption $(\ast)$, the relation ${\cal S}$ has an essentially unique Borel partition $Y=Y_{1}\sqcup Y_{2}\sqcup Y_{3}$ of $Y$ consisting of three ${\cal S}$-invariant Borel subsets $Y_{1}$, $Y_{2}$ and $Y_{3}$ satisfying the following conditions: the restriction $({\cal S})_{Y_{1}}$ is irreducible and amenable and the restriction $({\cal S})_{Y_{2}}$ is reducible. For any Borel subset $Y_{3}'$ of $Y_{3}$ with positive measure, the restriction $({\cal S})_{Y_{3}'}$ has no invariant Borel maps into $M({\cal PMF})$.  
\end{rem}


\section{Canonical reduction systems for reducible subrelations}

Let $M$ be a surface with $\kappa(M)\geq 0$. In \cite{ivanov1}, the canonical reduction system was introduced for any pure subgroup $\Gamma$ of the mapping class group $\Gamma(M)$ as follows: 

\begin{defn}[\ci{7.2}{ivanov1}]\label{def-crs-of-group}
Let $M$ be a surface with $\kappa(M)\geq 0$ and $\Gamma$ be a pure subgroup of the mapping class group $\Gamma(M)$. An isotopy class $\alpha \in V(C)$ is called an {\it essential reduction class}\index{essential!reduction class} for $\Gamma$ if the following two conditions are satisfied:
\begin{enumerate}
\renewcommand{\labelenumi}{\rm(\roman{enumi})}
\item $g\alpha =\alpha$ for any $g\in \Gamma$;
\item If $\beta \in V(C)$ and $i(\alpha, \beta)\neq 0$, then $g\beta \neq \beta$ for some $g\in \Gamma$.
\end{enumerate}
We call the set of all essential reduction classes for $\Gamma$ the {\it canonical reduction system}\index{canonical reduction system!for a pure subgroup} for $\Gamma$ and write it by $\sigma(\Gamma)$\index{$\ r \ C$@$\sigma(\Gamma)$}. Remark that the canonical reduction system is an element of $S(M)$. 
\end{defn}

Let $m$ be an integer and $\Gamma(M; m)$ be the kernel of the natural action of $\Gamma(M)$ on $H_{1}(M, {\Bbb Z}/m{\Bbb Z})$. Recall that $\Gamma(M; m)$ with $m\geq 3$ is a pure torsion-free subgroup of $\Gamma(M)$ with finite index (see Theorem \ref{pure-important}). 

Let $G$ be any subgroup of $\Gamma(M)$. If $\Gamma'$ and $\Gamma''$ are two pure normal subgroups of $G$ with finite index, then $\sigma(\Gamma')=\sigma(\Gamma'')$ (see \cite[7.4]{ivanov1}). Moreover, there always exists such a subgroup. If $\Gamma$ is the intersection $G\cap \Gamma(M; m)$ for some $m\geq 3$, then $\Gamma$ is a pure normal subgroup of $G$ with finite index. Thus, we can define the {\it canonical reduction system}\index{canonical reduction system!for a subgroup} for $G$ by the set $\sigma(\Gamma)$. 

If $G$ is finite, then $\sigma(G)$ is empty since $\sigma(\{ e\} )$ is empty.

\begin{thm}[\ci{Corollary 7.17}{ivanov1}]\label{thm-crs-non-empty}
An infinite subgroup $\Gamma$ of $\Gamma(M)$ is reducible if and only if its canonical reduction system is non-empty.
\end{thm}

Recall that a subgroup $\Gamma$ of $\Gamma(M)$ is said to be reducible if there exists $\sigma \in S(M)$ such that $g\sigma =\sigma$ for any $g\in \Gamma$ (see Theorem \ref{subgroup-classification}).

\begin{rem}\label{rem-main-lem2}
Canonical reduction systems are used to prove the following assertion (\cite[Corollary 7.13]{ivanov1}): suppose that $\Gamma_{1}$ is a subgroup of $\Gamma(M;m)$ with $m\geq 3$ and $\Gamma_{2}$ is an infinite normal subgroup of $\Gamma_{1}$. If $\Gamma_{2}$ is reducible, then 
\[g\sigma(\Gamma_{2})=\sigma(g\Gamma_{2}g^{-1})=\sigma(\Gamma_{2})\]
for any $g\in \Gamma_{1}$. Thus, $\Gamma_{1}$ is also reducible. Moreover, we can show $\sigma(\Gamma_{2})\subseteq \sigma(\Gamma_{1})$ by Corollary \ref{permutation}.
\end{rem}

Let $\Gamma$ be an infinite reducible subgroup of $\Gamma(M;m)$ with $m\geq 3$. Let $Q$ be a component of the resulting surface by cutting $M$ along a realization of $\sigma(\Gamma)$. Then we have the homomorphism $p\colon \Gamma \rightarrow \Gamma(Q)$ by Remark \ref{rem-comp-leave}.

\begin{thm}[\ci{Corollary 7.18}{ivanov1}]\label{thm-ivanov-last-hope}
The image $p(\Gamma)$ either is trivial or contains a pseudo-Anosov element and is not reducible.
\end{thm}

We will construct similar canonical reduction systems for reducible subrelations.  

\begin{defn}
Under the assumption $(\ast)$, for $\alpha \in V(C)$ and a Borel subset $A$ of $Y$ with positive measure, we say that the pair $(\alpha, A)$ is $\rho$-{\it invariant} (for ${\cal S}$)\index{rho-invariant pair@$\rho$-invariant!pair} if there exists some countable Borel partition $\{ A_{j}\}$ of $A$ such that the constant map $A\rightarrow \{ \alpha \}$ is  $\rho$-invariant for $({\cal S})_{A_{j}}$ for each $j$, that is,
\[\rho(x, y)\alpha =\alpha\]
for a.e.\ $(x, y)\in ({\cal S})_{A_{j}}$. 

We say that a $\rho$-invariant pair $(\alpha, A)$ is {\it essential}\index{essential rho-invariant pair@essential!$\rho$-invariant pair} if $B$ is a Borel subset of $A$ with positive measure and a curve $\beta \in V(C)$ satisfies $i(\alpha, \beta)\neq 0$, then the pair $(\beta, B)$ is not $\rho$-invariant.
\end{defn}

\begin{lem}\label{reducible-lem}
With the assumption $(\ast)$, the relation ${\cal S}$ is reducible if and only if for any Borel subset $A\subseteq Y$ with positive measure, there exists a $\rho$-invariant pair $(\alpha, B)$ with $B\subseteq A$. 
\end{lem}

\begin{pf}
It follows from Lemma \ref{partition} that if ${\cal S}$ is reducible, then for any Borel subset $A\subseteq Y$ with positive measure, there exists a $\rho$-invariant pair $(\alpha, B)$ with $B\subseteq A$. It follows from Theorem \ref{alternative} and Remark \ref{rem-disjoint-sum-irreducible-amenable-reducible} that the converse also holds. 
\end{pf}

The following lemma is easy to prove:

\begin{lem}\label{ess-making}
Assume $(\ast)$.
\begin{enumerate}
\renewcommand{\labelenumi}{\rm(\roman{enumi})}
\item If $(\alpha, A)$ is $\rho$-invariant, then so is $(\alpha, B)$ for any Borel subset $B$ of $A$ with positive measure. The same statement is true for essential $\rho$-invariant pairs.
\item Let $J$ be a countable index set and $A_{j}$ be a Borel subset of $Y$ with positive measure indexed by $j\in J$. If we have $\alpha \in V(C)$ such that $(\alpha, A_{j})$ is $\rho$-invariant for each $j\in J$, then the pair $(\alpha, \bigcup_{j\in J}A_{j})$ is also $\rho$-invariant. The same statement is true for essential $\rho$-invariant pairs.
\end{enumerate}
\end{lem}

In the next theorem, we show that the existence of a $\rho$-invariant pair implies the existence of an essential one. This assertion can be regarded as an analogue of Theorem \ref{thm-crs-non-empty}. For the proof, we shall recall some notions (see \cite[7.8]{ivanov1}).

Let $L$ be a submanifold of the surface $M$ which is a realization of some element in $S(M)$. Let $Q$ be a component of $M_{L}$, where $M_{L}$ denotes the resulting surface of cutting $M$ along $L$. Let $p_{L}\colon M_{L}\rightarrow M$ denote the canonical map. For $\delta \in V(C(M))$, we define a finite subset $r(\delta, Q)$\index{$r \ d Q$@$r(\delta, Q)$} of the set $V(C(Q))$ of all non-trivial isotopy classes of non-peripheral simple closed curves on $Q$. 

Let $\delta \in V(C(M))$ and represent the isotopy class $\delta$ by a circle $D$ that intersects each of the components of $L$ in the least possible number of points. Put $D_{L}=p_{L}^{-1}(D)$. The manifold $D_{L}$ consists of some intervals or it is a circle (if $D\cap L=\emptyset$). 

If either $D_{L}\cap Q=\emptyset$ or $D_{L}$ is a circle which lies in $Q$ and is peripheral for $Q$, then put $r(\delta, Q)=\emptyset$. 
If $D_{L}$ is a non-peripheral circle lying in $Q$, put $r(\delta, Q)=\{\delta \}$.

In the remaining cases, the intersection $D_{L}\cap Q$ consists of some intervals. For each such interval $I$, consider a regular neighborhood in $Q$ of the union of the interval $I$ and those components of $\partial Q$ on which the ends of $I$ lie (e.g., see \cite{hudson} for the definition of regular neighborhoods). Let $N_{I}$ denote the regular neighborhood. Then $N_{I}$ is a disk with two holes. Let $r'(\delta, Q)$ be the set of isotopy classes of components of the manifolds $\partial N_{I}\setminus \partial Q$, where $I$ runs through the set of all components of $D_{L}\cap Q$. Define $r(\delta, Q)$ as the resulting set of discarding from $r'(\delta, Q)$ the isotopy classes of trivial or peripheral circles of $Q$. We will regard $r(\delta, Q)$ as a subset of $V(C(M))$ using the embedding $V(C(Q))\hookrightarrow V(C(M))$. It is clear that this definition depends only on $\delta$ and the isotopy class of $Q$.

Let $F\colon M\rightarrow M$ be a diffeomorphism such that $F(L)=L$ and the induced diffeomorphism $M_{L}\rightarrow M_{L}$ takes $Q$ to $Q$. If $f\in \Gamma(M)$ denotes the isotopy class of $F$, then we have the equality
\[f(r(\delta, Q))=r(f\cdot \delta, Q)\]
by definition.

\begin{lem}[\ci{Lemma 7.9}{ivanov1}]\label{pants}
Let $L$ and $Q$ be the same as above and $\delta$ be an element in $V(C(M))$. If $r(\delta, Q)=\emptyset$, then one of the following three cases occurs:
\begin{enumerate}
\renewcommand{\labelenumi}{\rm(\roman{enumi})}
\item there is a circle in the class $\delta$ that does not intersect $Q$;
\item $\delta$ is the isotopy class of one of the components of $L$;
\item $Q$ is a disk with two holes.
\end{enumerate} 
\end{lem}

The following lemma will be used in the next chapter:

\begin{lem}[\ci{Lemma 7.10}{ivanov1}]\label{lem-pants-r}
Let $L$ and $Q$ be the same as above and $\delta$ be an element in $V(C(M))$. If $\alpha$ is an element in $V(C(Q))$ such that $i(\alpha, \delta)\neq 0$, then there exists an element $\delta'\in r(\delta, Q)$ such that $i(\alpha, \delta')\neq 0$.
\end{lem}

Using these notions, we prove the following main result in this section:

\begin{thm}\label{thm-ess-pair}
With the assumption $(\ast)$, if there exists a $\rho$-invariant pair, then there exists an essential $\rho$-invariant pair.
\end{thm}

\begin{pf}
This is essentially due to the proof of the group case \cite[Theorem 7.11]{ivanov1}. 

We assume that ${\cal S}$ has no essential $\rho$-invariant pairs, although ${\cal S}$ has a $\rho$-invariant one. We will deduce a contradiction.

For $\sigma \in S(M)$, define the subrelation ${\cal S}_{\sigma}=\rho^{-1}(\Gamma_{\sigma})$ of ${\cal S}$, where 
\[\Gamma_{\sigma}=\{ g\in \Gamma : g\sigma =\sigma \}=\{ g\in \Gamma : g\alpha =\alpha {\rm \ for \ any \ } \alpha \in \sigma \}. \]
Let ${\cal M}$ be the set of all Borel subsets of $Y$ with positive measure. Put
\[n_{0}=\sup \{ |\sigma|: \sigma \in S(M),\ E\in {\cal M},\ \rho(({\cal S})_{E})\subseteq \Gamma_{\sigma}\}.\]
Remark that since ${\cal S}$ has a $\rho$-invariant pair, there exist $\alpha \in V(C)$ and $E\in {\cal M}$ such that $\rho(({\cal S})_{E})\subseteq \Gamma_{{\alpha}}$. Let $\sigma \in S(M)$ and $Y_{\sigma}\in {\cal M}$ be elements such that $|\sigma |=n_{0}$ and $\rho(({\cal S})_{Y_{\sigma}})\subseteq \Gamma_{\sigma}$. It follows that $(\alpha, Y_{\sigma})$ is $\rho$-invariant for each $\alpha \in \sigma$. By our assumption, it is not essential. 

Let $L$ be a one-dimensional closed submanifold of $M$ representing $\sigma$. We will show that $M$ is decomposed into pairs of pants (i.e., disks with two holes) if it is cut along $L$. 

Let $Q$ be one component of $M_{L}$. Let $p_{L}\colon M_{L}\rightarrow M$ be the canonical map. Then $p_{L}(Q)\cap L$ is non-empty and consists of components of $L$. Let $\alpha$ be the isotopy class of a component of $p_{L}(Q)\cap L$. Then $\alpha$ is an element in $\sigma$. Since $(\alpha, Y_{\sigma})$ is not an essential $\rho$-invariant pair, there exist $\beta \in V(C(M))$, a Borel subset $A\subseteq Y_{\sigma}$ with positive measure and its countable Borel partition $A=\bigsqcup A_{j}$ such that $i(\alpha, \beta)\neq 0$ and $\beta$ is $\rho$-invariant for the relation $\bigsqcup ({\cal S})_{A_{j}}$. 

Consider the finite subset $r(\beta, Q)\subseteq V(C(Q))$. Then it follows from Theorem \ref{comp-leave} that $\rho(x, y)Q=Q$ for a.e.\ $(x, y)\in \bigsqcup ({\cal S})_{A_{j}}$ and that
\[\rho(x, y)r(\beta, Q)=r(\rho(x, y)\beta, Q)=r(\beta, Q)\]
for a.e.\ $(x, y)\in \bigsqcup ({\cal S})_{A_{j}}$. Since $\rho$ is valued in a pure subgroup and $r(\beta, Q)$ is a finite set, we have 
\[\rho(x, y)\beta'=\beta'\]
for a.e.\ $(x, y)\in \bigsqcup({\cal S})_{A_{j}}$ and any $\beta' \in r(\beta, Q)$ by Corollary \ref{permutation}. If $r(\beta, Q)$ were non-empty, then $(\beta', A)$ would be $\rho$-invariant for any $\beta'\in r(\beta, Q)$. This contradicts the maximality of $\vert \sigma \vert$ since $\beta'$ is represented by a non-trivial and non-peripheral curve on $Q$. Hence, we see that $r(\beta, Q)=\emptyset$. Since $i(\alpha, \beta)\neq 0$, the component $Q$ is a pair of pants by Lemma \ref{pants}. The claim follows.

Since $M$ is decomposed into pairs of pants by $L$, the subgroup $\Gamma_{\sigma}$ is contained in the free abelian subgroup $D_{\sigma}$ generated by the Dehn twists of all curves in $\sigma$ (see Corollary \ref{cor-pants-stabilizer}). We will show that $(\alpha, A')$ is an essential $\rho$-invariant pair for some $\alpha \in \sigma$ and Borel subset $A'$ of $Y_{\sigma}$ with positive measure. 

Assume that the claim is not true. 

\begin{lem}\label{essential-deny}
For each $\alpha \in \sigma$, there exist a countable Borel partition
\[Y_{\sigma}=\bigsqcup_{j}B_{j}^{(\alpha)}\]
(up to null sets) and $\beta_{j}^{(\alpha)}\in V(C)$ such that $i(\alpha, \beta_{j}^{(\alpha)})\neq 0$ and the pair $(\beta_{j}^{(\alpha)}, B_{j}^{(\alpha)})$ is $\rho$-invariant. 
\end{lem}

\begin{pf}
Let ${\cal M}'$ be the set of all Borel subsets $B$ of $Y_{\sigma}$ with positive measure such that there exist a countable Borel partition
\[B=\bigsqcup_{j}B_{j}^{(\alpha)}\]
and $\beta_{j}^{(\alpha)}\in V(C)$ such that $i(\alpha, \beta_{j}^{(\alpha)})\neq 0$ and the pair $(\beta_{j}^{(\alpha)}, B_{j}^{(\alpha)})$ is $\rho$-invariant. Define the number
\[m=\sup_{B\in {\cal M}'}\mu(B).\]
Remark that ${\cal M}'$ is non-empty and if $\{ B_{n}\}_{n\in {\Bbb N}}$ is a sequence of elements in ${\cal M}'$, then the union $\bigcup B_{n}$ is also in ${\cal M}'$. We can show that there exists a Borel subset $B\in {\cal M}'$ with $\mu(B)=m$ by this fact. Moreover, we can deduce $\mu(B)=\mu(Y_{\sigma})$ from the assumption that $(\alpha, A')$ is not essential $\rho$-invariant for each $\alpha \in \sigma$ and any Borel subset $A'$ of $Y_{\sigma}$ with positive measure.
\end{pf}   

It follows from this lemma that there exist a Borel subset $B$ of $Y_{\sigma}$ with positive measure and $\beta^{(\alpha)}\in V(C)$ for $\alpha \in \sigma$ such that $i(\alpha, \beta^{(\alpha)})\neq 0$ and $\beta^{(\alpha)}$ is $\rho$-invariant for $({\cal S})_{B}$. Therefore, we have
\[\rho(({\cal S})_{B})\subseteq D_{\sigma}\cap \left( \bigcap_{\alpha \in \sigma}\Gamma_{\beta^{(\alpha)}}\right).\]
The right hand side is trivial by Corollary \ref{trivial-stabilizer}. This contradicts the recurrence of ${\cal S}$. 
\end{pf}

Under the assumption $(\ast)$, suppose that  ${\cal S}$ is reducible. For $\alpha \in V(C)$, let ${\cal M}_{\alpha}$ be the set of all Borel subsets $A$ of $Y$ such that $(\alpha, A)$ is an essential $\rho$-invariant pair and denote
\[m_{\alpha}=\sup_{A\in {\cal M}_{\alpha}}\mu(A).\]
It follows from Theorem \ref{thm-ess-pair} that ${\cal M}_{\alpha}$ is non-empty for some $\alpha \in V(C)$ and from Lemma \ref{ess-making} (ii) that there exists a Borel subset $X_{\alpha}$ for any $\alpha \in V(C)$ such that if $(\alpha, A)$ is any essential $\rho$-invariant pair, then we have $\mu(A\setminus X_{\alpha})=0$ and if $X_{\alpha}$ has positive measure, then $(\alpha, X_{\alpha})$ is an essential $\rho$-invariant pair. Moreover, the equality
\[Y=\bigcup_{\alpha \in V(C)}X_{\alpha}\]
up to null sets follows from Lemma \ref{reducible-lem} and Theorem \ref{thm-ess-pair} since ${\cal S}$ is reducible.

If $\alpha_{1}, \alpha_{2}\in V(C)$ satisfy that $X_{\alpha_{1}}\cap X_{\alpha_{2}}$ has positive measure, then $i(\alpha_{1}, \alpha_{2})=0$ by definition. Thus, we can define a Borel map $\varphi_{\rho}\colon Y\rightarrow S(M)$ by the formula
\[\varphi_{\rho}(x)=\{ \alpha \in V(C): x\in X_{\alpha}\} \]
up to null sets. 

\begin{defn}\label{crs}
In the above notation, we call the Borel map $\varphi_{\rho}\colon Y\rightarrow S(M)$ the {\it canonical reduction system} (or {\it CRS}) for the cocycle $\rho \colon {\cal S}\rightarrow \Gamma$. \index{canonical reduction system!for a reducible subrelation}\index{CRS} As far as there is no possibility of confusion, we call it the CRS for ${\cal S}$ simply.
\end{defn}

\begin{lem}\label{crs-sublem}
With the assumption $(\ast)$, let $(\alpha, B)$ be a $\rho$-invariant pair and $g\in [[{\cal S}]]$ such that $B\subseteq {\rm dom}(g)$ and 
\[\rho(gx, x)=\gamma \in \Gamma \]
for a.e.\ $x\in B$. Then the pair $(\gamma \alpha, gB)$ is $\rho$-invariant. If $(\alpha, B)$ is an essential $\rho$-invariant pair, then so is $(\gamma \alpha, gB)$.
\end{lem}

\begin{pf}
The $\rho$-invariance of $(\gamma \alpha, gB)$ comes from
\begin{align*}
\rho(gx, gy)\gamma \alpha &=\rho(gx, x)\rho(x, y)\rho(y, gy)\gamma \alpha \\
                              &=\gamma \rho(x, y)\alpha \\
                              &=\gamma \alpha
\end{align*}
for a.e.\ $(x, y)\in \bigsqcup ({\cal S})_{B_{j}}$, where $B=\bigsqcup B_{j}$ is a countable Borel partition such that $\alpha$ is $\rho$-invariant for $\bigsqcup ({\cal S})_{B_{j}}$. 

We assume that $(\alpha, B)$ is an essential $\rho$-invariant pair. Let $\beta \in V(C)$ with $i(\gamma \alpha, \beta)\neq 0$ and $B'$ be a Borel subset of $gB$ with positive measure. We must show that $(\beta, B')$ is not $\rho$-invariant. 

If $(\beta, B')$ were $\rho$-invariant, then we could show as above that $(\gamma^{-1}\beta, g^{-1}B')$ is also $\rho$-invariant. On the other hand, since $i(\alpha, \gamma^{-1}\beta)\neq 0$, $g^{-1}B'\subseteq B$ and $(\alpha, B)$ is essential, the pair $(\gamma^{-1}\beta, g^{-1}B')$ is not $\rho$-invariant, which is a contradiction.
\end{pf}

\begin{lem}\label{crs-fund}
With the assumption $(\ast)$, suppose that ${\cal S}$ is reducible and let $\varphi_{\rho}\colon Y\rightarrow S(M)$ be the CRS for ${\cal S}$.
\begin{enumerate}
\renewcommand{\labelenumi}{\rm(\roman{enumi})}
\item The map $\varphi_{\rho}$ is maximal in the sense that if $(\alpha, A)$ is an essential $\rho$-invariant pair, then 
\[\mu(A\setminus \varphi_{\rho}^{-1}(\{ \sigma \in S(M):\alpha \in \sigma \} ))=0.\] 
\item The map $\varphi_{\rho}$ is $\rho$-invariant for ${\cal S}$.
\item[(iii)] Let $A$ be a Borel subset of $Y$ with positive measure. Then the CRS for the restriction $({\cal S})_{A}$ is equal to the restriction to $A$ of $\varphi_{\rho}\colon Y\rightarrow S(M)$. 
\end{enumerate}
\end{lem}

\begin{pf}
Since we have the equality
\[\varphi_{\rho}^{-1}(\{ \sigma \in S(M): \alpha \in \sigma \})=X_{\alpha}\]
with the notation right before Definition \ref{crs}, the assertion (i) has already been verified. We prove the assertion (ii). Fix $\alpha \in V(C)$. Let $g\in [[{\cal S}]]$ be a partial Borel isomorphism such that the intersection $X_{\alpha}\cap {\rm dom}(g)$ has positive measure. Let $\gamma$ be any element in $\Gamma$ such that the Borel subset
\[A_{\gamma}=\{ x\in X_{\alpha}\cap {\rm dom}(g): \rho(gx, x)=\gamma \}\]
has positive measure. It is enough to show $\gamma \alpha \in \varphi_{\rho}(gx)$ for a.e.\ $x\in A_{\gamma}$. However, since we have $gA_{\gamma}\subseteq X_{\gamma \alpha}$ (up to null sets) by Lemma \ref{crs-sublem} and the above assertion (i), the assertion (ii) holds.

For the assertion (iii), note the following remarks, which follow by definition:

\begin{enumerate}
\item[(a)] For any Borel subsets $B\subseteq A$ with positive measure and $\alpha \in V(C)$, the pair $(\alpha, B)$ is a (essential) $\rho$-invariant pair for ${\cal S}$ if and only if it is a (essential) $\rho$-invariant pair for $({\cal S})_{A}$.
\item[(b)] For $\alpha \in V(C)$ and a Borel subset $B$ of $Y$ with positive measure, the pair $(\alpha, B)$ is an essential $\rho$-invariant pair for ${\cal S}$ if and only if $\alpha$ is an isotopy class in $\varphi_{\rho}(y)\in S(M)$ for a.e.\ $y\in B$.
\end{enumerate}

We denote by $\varphi_{\rho, A}\colon A\rightarrow S(M)$ the CRS for the restriction $({\cal S})_{A}$. If $\alpha \in V(C)$ and $E$ is a Borel subset of $A$ with positive measure such that $\alpha \in \varphi_{\rho}(x)$ for $x\in E$, then $(\alpha, E)$ is an essential $\rho$-invariant pair for ${\cal S}$ by the remark (b) and thus, for $({\cal S})_{A}$ by the remark (a). This means $\alpha \in \varphi_{\rho, A}(x)$ for a.e.\ $x\in E$ by the remark (b). 

Conversely, if $\alpha \in V(C)$ and $E$ is a Borel subset of $A$ with positive measure such that $\alpha \in \varphi_{\rho, A}(x)$ for $x\in E$, then $(\alpha, E)$ is an essential $\rho$-invariant pair for $({\cal S})_{A}$ by the remark (b) and thus, for ${\cal S}$ by the remark (a). This means that $\alpha \in \varphi_{\rho}(x)$ for a.e.\ $x\in E$ by the remark (b). 
\end{pf}


\section[Indecomposability of equivalence relations]{Indecomposability of equivalence relations generated by the mapping class group}

In this section, we give necessary conditions for an equivalence relation generated by an essentially free, non-singular Borel action of the mapping class group to be decomposed into a product of two recurrent relations. For an equivalence relation generated by a hyperbolic group, Adams \cite{adams2} showed that the amenability is a necessary condition for the decomposability (see Theorem \ref{adams2-main}). Moreover, it is also a sufficient condition \cite{cfw}. We will show an analogous result for equivalence relations generated by the mapping class group.

\subsection{Generalities for normality}\label{subsection-normal}

In this subsection, we recall fundamentals of normal subrelations developed by Feldman-Sutherland-Zimmer \cite{fsz}. 

Let ${\cal R}$ be a non-singular discrete measured equivalence relation on a standard Borel space $(X, \mu)$ with a finite positive measure. Let ${\cal S}$ be its subrelation. Note that ${\cal S}$ induces an equivalence relation $\sim$ on each ${\cal R}$-class ${\cal R}x$, $x\in X$ defined as follows: for $y, z\in {\cal R}x$, we write $y\sim z$ if $(y, z)\in {\cal S}$. Let $I(x)$ denote the cardinality of the quotient ${\cal R}x/\sim$. Then the function $I\colon x\mapsto I(x) \in {\Bbb N}\cup \{ \infty \}$ is measurable and ${\cal R}$-invariant \cite[Lemma 1.1]{fsz}. It is called the {\it index function}\index{index function} for the pair $({\cal R}, {\cal S})$. 

We denote $X_{a}=I^{-1}(a)$ for $a\in {\Bbb N}\cup \{ \infty \}$. Then for each $a\in {\Bbb N}\cup \{ \infty \}$, there exist Borel maps $\phi^{a}_{j}\colon X_{a}\rightarrow X_{a}$ such that for a.e.\ $x\in X_{a}$, the family $\{ {\cal S}\phi^{a}_{j}(x): j\in {\Bbb N},\ 0\leq j< a\}$ is a partition of ${\cal R}x$ \cite[Lemma 1.1]{fsz}. An element in the family $\{ \phi_{j}^{a}\}$ is called a {\it choice function}\index{choice function}.

In \cite{fsz}, the subrelation ${\cal S}$ is said to be normal in ${\cal R}$ if we can find a family $\{ \phi_{j}^{a}\}$ of choice functions such that if $(x, y)\in ({\cal S})_{X_{a}}$, then $(\phi_{j}^{a}(x), \phi_{j}^{a}(y))\in ({\cal S})_{X_{a}}$. Several equivalent conditions for the normality are known \cite[Theorem 2.2]{fsz}. However, for our purpose, we adopt a slightly different definition of normal subrelations as follows:

\begin{defn}
Let ${\cal R}$ be a non-singular discrete measured equivalence relation on a standard Borel space $(X, \mu)$ with a finite positive measure. Let ${\cal S}$ be its subrelation. We denote by ${\rm End}_{{\cal R}}({\cal S})'$\index{$E n d R S p $@${\rm End}_{{\cal R}}({\cal S})'$} the set of all Borel maps $g$ from a Borel subset ${\rm dom}(g)$ of $X$ to $X$ such that 
\begin{enumerate}
\renewcommand{\labelenumi}{\rm(\roman{enumi})}
\item $(g(x), x)\in {\cal R}$ for a.e.\ $x\in {\rm dom}(g)$.
\item $(g(x), g(y))\in {\cal S}$ if and only if $(x, y)\in {\cal S}$ for a.e.\ $x, y\in {\rm dom}(g)$
\end{enumerate}
Let us denote by $[[{\cal R}]]_{\cal S}$\index{$R ]] S$@$[[{\cal R}]]_{\cal S}$} the set of all partial Borel isomorphisms in ${\rm End}_{\cal R}({\cal S})'$.
\end{defn}

Remark that any endomorphism in ${\rm End}_{\cal R}({\cal S})'$ has at most countable fibers. 

\begin{defn}
Let ${\cal R}$ be a non-singular discrete measured equivalence relation on a standard Borel space $(X, \mu)$ with a finite positive measure. Its subrelation ${\cal S}$ is said to be {\it normal}\index{normal subrelation} in ${\cal R}$ if we can find a countable family $\{ \phi_{n}\}$ of ${\rm End}_{\cal R}({\cal S})'$ such that for a.e.\ $(x, y)\in {\cal R}$, there exists $n$ such that $x\in {\rm dom}(\phi_{n})$ and $(\phi_{n}(x), y)\in {\cal S}$. We call such $\{ \phi_{n}\}$ a {\it family of normal choice functions}\index{family of normal choice functions} for the pair $({\cal R}, {\cal S})$.
\end{defn}

We shall give fundamental results about normal subrelations (in our sense). For two Borel maps $f$, $g$ from Borel subsets ${\rm dom}(f)$, ${\rm dom}(g)$ of $X$ to $X$, respectively, their composition is defined by the restricted map
\[g\circ f\colon {\rm dom}(f)\cap f^{-1}({\rm dom}(g))\rightarrow X.\] 
The following lemma is easy to prove:

\begin{lem}\label{generate-normal}
Let ${\cal R}$ be a non-singular discrete measured equivalence relation on a standard Borel space $(X, \mu)$ with a finite positive measure. Let ${\cal S}$ be its subrelation.
\begin{enumerate}
\item[(i)] If $g\in [[{\cal R}]]_{\cal S}$, then $g^{-1}\in [[{\cal R}]]_{\cal S}$.
\item[(ii)] If $g\in {\rm End}_{\cal R}({\cal S})'$, then the restriction of $g$ to any Borel subset of ${\rm dom}(g)$ also belongs to ${\rm End}_{\cal R}({\cal S})'$.
\item[(iii)] The composition of two elements in ${\rm End}_{\cal R}({\cal S})'$ also belongs to ${\rm End}_{\cal R}({\cal S})'$. In particular, if $g\in {\rm End}_{\cal R}({\cal S})'$ and $f$ is a Borel map from a Borel subset ${\rm dom}(f)$ of $X$ into $X$ such that $(f(x), x)\in {\cal S}$ for a.e.\ $x\in {\rm dom}(f)$, then the compositions $g\circ f$, $f\circ g$ are also in ${\rm End}_{\cal R}({\cal S})'$.
\end{enumerate}
\end{lem}

\begin{lem}\label{lem-normal-subgroup}
Let $G$ be a discrete group and $H$ be its normal subgroup. If a pair $({\cal R}, {\cal S})$ of a discrete measured equivalence relation and its subrelation on $(X, \mu)$ is defined by an essentially free, non-singular Borel action of $G$ and its restriction to $H$, then the subrelation ${\cal S}$ is normal in ${\cal R}$.
\end{lem}

\begin{pf}
Let $\{ g_{n}\}$ be a set of all representatives of the set $G/H$ of left cosets. We show that $\{ g_{n}\}$ is a family of normal choice functions for the pair $({\cal R}, {\cal S})$. We may assume that the action of $G$ is free.

If $(x, y)\in {\cal S}$, then we have a unique $h\in H$ such that $y=hx$. Then 
\[(g_{n}x, g_{n}y)=(g_{n}x, g_{n}hx)=(g_{n}x, h'g_{n}x)\in {\cal S}\] 
for some $h'\in H$. Conversely, if $(g_{n}x, g_{n}y)\in {\cal S}$ for $x, y\in X$, then we have unique $h, h'\in H$ such that $g_{n}y=hg_{n}x=g_{n}h'x$. It follows that $y=h'x$, which means $(x, y)\in {\cal S}$. Thus, we have verified $g_{n}\in [[{\cal R}]]_{\cal S}$ for each $n$.

For any $(x, y)\in {\cal R}$, we can find a unique $g\in G$ such that $y=gx$. Moreover, we can find unique $n$ and $h\in H$ such that $g=hg_{n}$ since the subgroup $H$ is normal in $G$. Then $(g_{n}x, y)\in {\cal S}$.
\end{pf}

\begin{lem}\label{normal-restriction}
Let ${\cal R}$ be a non-singular discrete measured equivalence relation on $(X, \mu)$ and ${\cal S}$ be its normal subrelation. Then for any Borel subset $A$ of $X$ with positive measure, the restriction $({\cal S})_{A}$ is a normal subrelation of $({\cal R})_{A}$. 
\end{lem}

\begin{pf}
Let $B={\cal S}A$ be the ${\cal S}$-saturation of $A$. Then there exists a Borel map $f\colon B\rightarrow A$ such that
\begin{enumerate}
\renewcommand{\labelenumi}{\rm(\roman{enumi})}
\item $(f(x), x)\in {\cal S}$ for a.e.\ $x\in B$;
\item $f(x)=x$ for any $x\in A$
\end{enumerate}
by Lemma \ref{lem-adams-map}. Let $\{ \phi_{n}\}\subseteq {\rm End}_{{\cal R}}({\cal S})'$ be a family of normal choice functions for the pair $({\cal R}, {\cal S})$. Let us define 
\[D_{n}=\{ x\in A\cap {\rm dom}(\phi_{n}): \phi_{n}(x)\in B \}\]
and define a Borel map $\phi_{n}'\colon D_{n}\rightarrow A$ by the formula
\[\phi_{n}'(x)=f\circ \phi_{n}(x)\]
for $x\in D_{n}$. Then $\phi_{n}'\in {\rm End}_{\cal R}({\cal S})'$ by Lemma \ref{generate-normal} (iii).

For $(x, y)\in ({\cal R})_{A}$, there exists $n$ such that $x\in {\rm dom}(\phi_{n})$ and $(\phi_{n}(x), y)\in {\cal S}$. It follows from $y\in A$ that $\phi_{n}(x)\in {\cal S}A=B$, that is, $x\in D_{n}$. Thus, $(\phi_{n}'(x), y)\in ({\cal S})_{A}$, which implies that $({\cal S})_{A}$ is normal in $({\cal R})_{A}$.
\end{pf}

\begin{lem}\label{generate-generate-normal}
Let $G$ be a countable subset of ${\rm End}_{\cal R}({\cal S})'$ and ${\cal S}_{G}$ be the equivalence relation generated by ${\cal S}$ and $G$. Then ${\cal S}$ is normal in ${\cal S}_{G}$.
\end{lem}

\begin{pf}
We may assume that the domain of any element in $G$ is ${\cal S}$-invariant, using the map constructed in Lemma \ref{lem-adams-map}. Thanks to Theorem \ref{thm-standard-borel-space} (iii), each $g\in G$ can be expressed as a disjoint sum of elements in $[[{\cal R}]]_{\cal S}$, that is, there exist $g_{n}\in [[{\cal R}]]_{\cal S}$, $n\in {\Bbb N}$ such that ${\rm dom}(g_{n})$ is disjoint each other for $n\in {\Bbb N}$ and 
\[{\rm dom}(g)=\bigsqcup_{n\in {\Bbb N}}{\rm dom}(g_{n}).\] 
Moreover, the restriction of $g$ to ${\rm dom}(g_{n})$ is equal to $g_{n}$. We extend the inverse map $g_{n}^{-1}\in [[{\cal R}]]_{\cal S}$ to a map $f_{n}^{g}\in {\rm End}_{\cal R}({\cal S})'$ with domain ${\cal S}$-invariant, using Lemma \ref{lem-adams-map}. The map $f_{n}^{g}$ satisfies 
\[f_{n}^{g}(x)=g_{n}^{-1}(x)\]
for $x\in {\rm ran}(g_{n})$. Let us denote $G'=\{ g, f_{n}^{g}\}_{g\in G}\subseteq {\rm End}_{\cal R}({\cal S})'$.

Next, we prove the following claim: for a.e.\ $(x, y)\in {\cal S}_{G}$, we can find a finite sequence $x=x_{0}, x_{1}, \ldots, x_{l}=y$ such that $(x_{2i}, x_{2i+1})\in {\cal S}$ and $(x_{2i-1}, x_{2i})$ is in the graph 
\[\{ (x, g(x))\in X\times X: x\in {\rm dom}(g)\}\]
of some element $g\in G'$. By definition, for $(x, y)\in {\cal S}_{G}$, we can find a finite sequence $x=x_{0}, x_{1}, \ldots, x_{l}=y$ such that $(x_{2i}, x_{2i+1})\in {\cal S}$ and either $(x_{2i-1}, x_{2i})$ or $(x_{2i}, x_{2i-1})$ is in the graph of some element $g\in G$. If $(x_{2i}, x_{2i-1})$ is in the graph of $g\in G$, then there exists a unique $n\in {\Bbb N}$ such that $x_{2i}\in {\rm dom}(g_{n})$ and $x_{2i-1}=g_{n}(x_{2i})$. It follows from the definition of the map $f_{n}^{g}$ that $(x_{2i-1}, x_{2i})$ is in the graph of $f_{n}^{g}$, which shows the claim.      

Therefore, replacing $G$ by $G'$, we may assume that any element in $G$ has an ${\cal S}$-invariant domain and that for a.e.\ $(x, y)\in {\cal S}_{G}$, we can find a finite sequence $x=x_{0}, x_{1}, \ldots, x_{l}=y$ such that $(x_{2i}, x_{2i+1})\in {\cal S}$ and $(x_{2i-1}, x_{2i})$ is in the graph of some element in $G$.

Let $\Omega$ be the countable set of all words consisting of letters in $G$. Each word in $\Omega$ naturally defines an element in ${\rm End}_{\cal R}({\cal S})'$ by Lemma \ref{generate-normal} (iii). We prove that the set $\Omega$ defines a family of normal choice functions for the pair $({\cal S}_{G}, {\cal S})$. 

For a.e.\ $(x, y)\in {\cal S}_{G}$, we can find a finite sequence $x=x_{0}, \ldots, x_{l}=y$ such that $(x_{2i}, x_{2i+1})\in {\cal S}$ and $(x_{2i-1}, x_{2i})$ is in the graph of some element $g_{i}\in G$. For the proof that $\Omega$ is a family of normal choice functions, it suffices to show by induction on $i$ that we have
\[g_{i}\circ g_{i-1}\circ \cdots \circ g_{1}(x_{1})\in {\rm dom}(g_{i+1}),\]
\[(x_{2i+1}, g_{i}\circ g_{i-1}\circ \cdots \circ g_{1}(x_{1}))\in {\cal S}\]  because each $g_{i}\circ \cdots \circ g_{1}$ defines an element in $\Omega$. When $i=1$, since $g_{1}(x_{1})=x_{2}$, $(x_{2}, x_{3})\in {\cal S}$ and $x_{3}\in {\rm dom}(g_{2})$, we see that $g_{1}(x_{1})\in {\cal S}({\rm dom}(g_{2}))={\rm dom}(g_{2})$. Suppose that it is true in the case of $i$. Since $g_{i+1}(x_{2i+1})=x_{2i+2}$, $(x_{2i+2}, x_{2i+3}) \in {\cal S}$ and $(x_{2i+1}, g_{i}\circ \cdots \circ g_{1}(x_{1}))\in {\cal S}$, we have $(x_{2i+3}, g_{i+1}\circ \cdots \circ g_{1}(x_{1}))\in {\cal S}$. It follows from $x_{2i+3}\in {\rm dom}(g_{i+2})$ that $g_{i+1}\circ \cdots \circ g_{1}(x_{1})\in {\rm dom}(g_{i+2})$.       
\end{pf}


\subsection{Main lemmas}  

In this subsection, we prove two central lemmas (Lemmas \ref{main-lem1} and \ref{main-lem2}) for the proof of the main theorems in this chapter. In both the lemmas, we use the maximal properties for invariant Borel maps for irreducible and amenable subrelations and reducible subrelations developed before. Lemma \ref{main-lem2} can be seen as an analogue of the claim in Remark \ref{rem-main-lem2}.  

\begin{assumption}\index{$(\ast p p $@$(\ast)''$}
We call the following assumption $(\ast)''$: let $\Gamma$ be a subgroup of $\Gamma(M; m)$, where $M$ is a surface with $\kappa(M)\geq 0$ and $m\geq 3$ is an integer. Let $(X, \mu)$ be a standard Borel space with a finite positive measure and ${\cal R}$ be a non-singular discrete measured equivalence relation on $(X, \mu)$. Suppose that we have a Borel cocycle
\[\rho \colon {\cal R}\rightarrow \Gamma \]
with finite kernel. Let ${\cal S}$ be a recurrent subrelation of ${\cal R}$. 
\end{assumption}

\begin{lem}\label{main-lem1}
With the assumption $(\ast)''$, suppose that ${\cal S}$ has a $\rho$-invariant Borel map $X\rightarrow \partial_{2}C$ for ${\cal S}$. Let $\varphi_{0}\colon X\rightarrow \partial_{2}C$ be the ``maximal'' $\rho$-invariant Borel map for ${\cal S}$ in Proposition \ref{maximal} (ii). Let $G$ be a countable subset of $[[{\cal R}]]_{\cal S}$. Then the map $\varphi_{0}$ is $\rho$-invariant for the relation ${\cal S}_{G}$ generated by ${\cal S}$ and $G$.
\end{lem}

\begin{pf}
The proof is similar to that for \cite[Lemma 3.4]{adams2}. We denote $G^{-1}=\{ g^{-1}\in [[{\cal R}]]_{\cal S}: g\in G\}$. For $g\in G\cup G^{-1}$, define a Borel map $\psi_{g}\colon {\rm dom}(g)\rightarrow \partial_{2}C$ by the formula
\[\psi_{g}(x)=\rho(x, gx)\varphi_{0}(gx)\]
for $x\in {\rm dom}(g)$. Then for $(x, y)\in ({\cal S})_{{\rm dom}(g)}$, we have 
\begin{align*}
\rho(x, y)\psi_{g}(y) &= \rho(x, y)\rho(y, gy)\varphi_{0}(gy)\\
                      &= \rho(x, gx)\rho(gx, gy)\varphi_{0}(gy)\\
                      &= \rho(x, gx)\varphi_{0}(gx)=\psi_{g}(x).
\end{align*}
Thus, the map $\psi_{g}$ is $\rho$-invariant for $({\cal S})_{{\rm dom}(g)}$. By Corollary \ref{maximal-cor} (ii), we have 
\[{\rm supp}(\psi_{g}(x))\subseteq {\rm supp}(\varphi_{0}(x))\]
for a.e.\ $x\in {\rm dom}(g)$. 

Similarly, the Borel map ${\rm ran}(g)\ni y\mapsto \rho(y, g^{-1}y)\varphi_{0}(g^{-1}y)$ is $\rho$-invariant for $({\cal S})_{{\rm ran}(g)}$. Thus, we have 
\[{\rm supp}(\rho(y, g^{-1}y)\varphi_{0}(g^{-1}y))\subseteq {\rm supp}(\varphi_{0}(y))\]
for a.e.\ $y\in {\rm ran}(g)$. These two inclusions show that $\psi_{g}(x)=\varphi_{0}(x)$ for a.e.\ $x\in {\rm dom}(g)$, which means that $\varphi_{0}$ is $\rho$-invariant for the relation ${\cal S}_{G}$. 
\end{pf}

\begin{lem}\label{main-sublem}
With the assumption $(\ast)''$, let $(\alpha, B)$ be a $\rho$-invariant pair for ${\cal S}$ and $g\in [[{\cal R}]]_{\cal S}$ such that $B\subseteq {\rm dom}(g)$ and 
\[\rho(gx, x)=\gamma \in \Gamma \]
for a.e.\ $x\in B$. Then $(\gamma \alpha, gB)$ is a $\rho$-invariant pair for ${\cal S}$. If $(\alpha, B)$ is an essential $\rho$-invariant pair, then so is $(\gamma \alpha, gB)$.
\end{lem}

\begin{pf}
This is almost the same as the proof of Lemma \ref{crs-sublem}. We have 
\begin{align*}
\rho(gx, gy)\gamma \alpha &=\rho(gx, x)\rho(x, y)\rho(y, gy)\gamma \alpha \\
                              &=\gamma \rho(x, y)\alpha \\
                              &=\gamma \alpha
\end{align*}
for $(x, y)\in \bigsqcup ({\cal S})_{B_{j}}$, where $B=\bigsqcup B_{j}$ is a countable Borel partition such that $\alpha$ is $\rho$-invariant for $\bigsqcup ({\cal S})_{B_{j}}$. Since $(gx, gy)\in {\cal S}$ for $x, y\in {\rm dom}(g)$ implies $(x, y)\in {\cal S}$, the above equality means that $\gamma \alpha$ is invariant for $\bigsqcup ({\cal S})_{gB_{j}}$, that is, the pair $(\gamma \alpha, gB)$ is $\rho$-invariant for ${\cal S}$.

We assume that $(\alpha, B)$ is an essential $\rho$-invariant pair. Let $\beta \in V(C)$ with $i(\gamma \alpha, \beta)\neq 0$ and $B'$ be a Borel subset of $gB$ with positive measure. We must show that $(\beta, B')$ is not $\rho$-invariant for ${\cal S}$. 

If $(\beta, B')$ is $\rho$-invariant for ${\cal S}$, then we can show as above that $(\gamma^{-1}\beta, g^{-1}B')$ is also $\rho$-invariant for ${\cal S}$. On the other hand, since $i(\alpha, \gamma^{-1}\beta)\neq 0$, $g^{-1}B'\subseteq B$ and $(\alpha, B)$ is essential, the pair $(\gamma^{-1}\beta, g^{-1}B')$ is not $\rho$-invariant for ${\cal S}$, which is a contradiction.
\end{pf}

\begin{lem}\label{main-lem2}
With the assumption $(\ast)''$, suppose that ${\cal S}$ is reducible and let $\varphi_{\cal S}\colon X\rightarrow S(M)$ be the CRS for ${\cal S}$. Let $G$ be a countable subset of $[[{\cal R}]]_{\cal S}$. Then the map $\varphi_{{\cal S}}$ is $\rho$-invariant for the relation ${\cal S}_{G}$ generated by ${\cal S}$ and $G$. In particular, ${\cal S}_{G}$ is reducible.
\end{lem}
Moreover, if $\varphi_{\cal S}$ is constant, then any curve in the value of $\varphi_{\cal S}$ is invariant for ${\cal S}_{G}$ by the pureness of $\Gamma$ (see Corollary \ref{permutation}).

\begin{pf}
For $g\in G\cup G^{-1}$, define a Borel map $\psi_{g}\colon {\rm dom}(g)\rightarrow S(M)$ by 
\[\psi_{g}(x)=\rho(x, gx)\varphi_{{\cal S}}(gx).\]
For $(x, y)\in ({\cal S})_{{\rm dom}(g)}$, we have 
\begin{align*}
\rho(x, y)\psi_{g}(y) &= \rho(x, y)\rho(y, gy)\varphi_{{\cal S}}(gy)\\
                      &= \rho(x, gx)\rho(gx, gy)\varphi_{{\cal S}}(gy)\\
                      &= \rho(x, gx)\varphi_{{\cal S}}(gx)=\psi_{g}(x)
\end{align*}
since $\varphi_{\cal S}$ is $\rho$-invariant for ${\cal S}$ by Lemma \ref{crs-fund} (ii). Thus, $\psi_{g}$ is $\rho$-invariant for $({\cal S})_{{\rm dom}(g)}$. 

There exists a countable Borel partition 
\[{\rm dom}(g)=\bigsqcup_{\gamma \in \Gamma}A_{\gamma}\]
such that $\rho(x, gx)=\gamma$ for $x\in A_{\gamma}$. Let $\sigma \in S(M)$ be an element such that the measure of the set $\psi_{g}^{-1}(\sigma)$ is positive. Let $\alpha$ be any element in $\sigma$ and $\gamma \in \Gamma$ be an element with $\mu(\psi_{g}^{-1}(\sigma)\cap A_{\gamma})>0$. Then we will show that $(\alpha, \psi_{g}^{-1}(\sigma)\cap A_{\gamma})$ is an essential $\rho$-invariant pair for ${\cal S}$. 

The $\rho$-invariance of $(\alpha, \psi_{g}^{-1}(\sigma)\cap A_{\gamma})$ for ${\cal S}$ follows easily from the fact that $\psi_{g}$ is $\rho$-invariant for $({\cal S})_{{\rm dom}(g)}$ and Corollary \ref{permutation}. By definition, we have
\begin{align*}
\psi^{-1}_{g}(\sigma) &=\{ x\in {\rm dom}(g): \psi_{g}(x)=\sigma \} \\
                      &=\{ x\in {\rm dom}(g): \rho(x, gx)\varphi_{{\cal S}}(gx)=\sigma \},
\end{align*}
which implies 
\[\psi_{g}^{-1}(\sigma)\cap A_{\gamma}\subseteq \{ x\in {\rm dom}(g): \varphi_{{\cal S}}(gx)=\gamma^{-1}\sigma \} =g^{-1}(\varphi^{-1}_{{\cal S}}(\gamma^{-1}\sigma)\cap {\rm ran}(g)).\]
It follows from the definition of $\varphi_{{\cal S}}$ and the above inclusion that $(\gamma^{-1}\alpha, g(\psi_{g}^{-1}(\sigma)\cap A_{\gamma}))$ is an essential $\rho$-invariant pair for ${\cal S}$. Hence, $(\alpha, \psi_{g}^{-1}(\sigma)\cap A_{\gamma})$ is an essential $\rho$-invariant pair for ${\cal S}$ by Lemma \ref{main-sublem}. 

By the maximality of the Borel set
\[X_{\alpha}=\varphi_{\cal S}^{-1}(\{ \sigma \in S(M): \alpha \in \sigma \}) \]
(see Lemma \ref{crs-fund} (i)), we have $\psi_{g}^{-1}(\sigma)\cap A_{\gamma}\subseteq X_{\alpha}$ (up to null sets). Since $\gamma \in \Gamma$ is an arbitrary element such that $\psi_{g}^{-1}(\sigma)\cap A_{\gamma}$ has positive measure and ${\rm dom}(g)=\bigsqcup_{\gamma \in \Gamma}A_{\gamma}$, we have $\psi_{g}^{-1}(\sigma)\subseteq X_{\alpha}\cap {\rm dom}(g)$. Since $\alpha \in \sigma$ is arbitrary, we see that 
\[\psi_{g}^{-1}(\sigma)\subseteq {\rm dom}(g)\cap \left( \bigcap_{\alpha \in \sigma}X_{\alpha}\right).\] 

If we had $\alpha' \in V(C)\setminus \sigma$ such that $\psi_{g}^{-1}(\sigma)\cap X_{\alpha'}$ has positive measure, then $i(\alpha, \alpha')=0$ for any $\alpha \in \sigma$ since $\psi_{g}^{-1}(\sigma)\subseteq \bigcap_{\alpha \in \sigma}X_{\alpha}$. Moreover, there would exist $\gamma \in \Gamma$ such that $E=\psi^{-1}_{g}(\sigma)\cap X_{\alpha'}\cap A_{\gamma}$ has positive measure. Then $(\alpha', E)$ is an essential $\rho$-invariant pair for ${\cal S}$. For a.e.\ $x\in E$, we have 
\[\sigma =\psi_{g}(x)=\gamma \varphi_{{\cal S}}(gx)\]
and thus, $\varphi_{{\cal S}}(gx)=\gamma^{-1}\sigma$. On the other hand, by Lemma \ref{main-sublem}, $(\gamma^{-1}\alpha', gE)$ is also an essential $\rho$-invariant pair for ${\cal S}$. This contradicts $\varphi_{{\cal S}}(gx)=\gamma^{-1}\sigma$ for a.e.\ $x\in E$ and $\gamma^{-1}\alpha'\notin \gamma^{-1}\sigma$.

We have shown that for any $\alpha'\in V(C)\setminus \sigma$, the Borel subset $\psi_{g}^{-1}(\sigma)\cap X_{\alpha'}$ is null. Hence, we see that $\psi_{g}^{-1}(\sigma)=\varphi_{\cal S}^{-1}(\sigma)\cap {\rm dom}(g)$ for any $\sigma \in S(M)$ (up to null sets), that is, $\psi_{g}=\varphi_{{\cal S}}$ a.e.\ on ${\rm dom}(g)$. This completes the proof.   
\end{pf}


\subsection{Necessary conditions for the existence of amenable normal subrelations}

First, we show that equivalence relations generated by an essentially free, measure-preserving action of the mapping class group can be neither amenable nor reducible.

\begin{assumption}\index{$(\diamond $@$(\diamond)$}
We call the following assumption $(\diamond)$: let $(X, \mu)$ be a standard Borel space with a finite positive measure. Let $\Gamma$ be a sufficiently large subgroup of $\Gamma(M; m)$, where $M$ is a surface with $\kappa(M)\geq 0$ and $m\geq 3$ is an integer. Suppose that we have an essentially free, measure-preserving Borel action of $\Gamma$ on $(X, \mu)$. Let ${\cal R}$ be the induced equivalence relation on $(X, \mu)$.
\end{assumption}

Recall that a subgroup of the mapping class group is said to be sufficiently large if it has an independent pair of pseudo-Anosov elements (see Theorem \ref{subgroup-classification}).

\begin{prop}\label{typeII}
With the assumption $(\diamond)$, the restriction $({\cal R})_{A}$ is neither amenable nor reducible for any Borel subset $A$ of $X$ with positive measure.
\end{prop}

\begin{pf}
Since $\Gamma$ is not amenable, neither is $({\cal R})_{A}$. Let $\rho \colon {\cal R}\rightarrow \Gamma$ be the induced cocycle. Consider an amenable subrelation ${\cal S}$ generated by a pseudo-Anosov element $g$ in $\Gamma$. Since ${\cal R}$ is of type ${\rm II}_{1}$, the relation ${\cal S}$ is recurrent. Let $\{ F_{+}, F_{-}\}\subseteq {\cal MIN}$ be the set of the fixed points of $g$. We define the Borel map $\varphi \colon X\rightarrow M({\cal PMF})$ by setting $\varphi(x)$ to be the measure such that its support is $\{ F_{+}, F_{-}\}$ and each atom has measure $1/2$ for any $x\in X$. Then $\varphi$ is $\rho$-invariant for ${\cal S}$. If $({\cal R})_{A}$ were reducible, then there would exist a $\rho$-invariant Borel map $\psi \colon A\rightarrow S(M)$ for $({\cal R})_{A}$. The map $\psi$ is also $\rho$-invariant for $({\cal S})_{A}$. This contradicts Theorem \ref{alternative}. Remark that $S(M)$ can be identified with a subset of ${\cal PMF}\setminus {\cal MIN}$ (see Chapter \ref{chapter:amenable-action}, Section \ref{mcg}) and thus, each element in $S(M)$ can be naturally regarded as an atomic measure on ${\cal PMF}\setminus {\cal MIN}$. 
\end{pf}

Under the assumption $(\ast)$, we say that ${\cal S}$ is a {\it disjoint sum of an irreducible and amenable relation and a reducible relation}\index{disjoint sum of an irreducible and amenable relation and a reducible relation} if there exists a Borel partition $Y=Y_{1}\sqcup Y_{2}$ such that $({\cal S})_{Y_{1}}$ is irreducible and amenable and $({\cal S})_{Y_{2}}$ is reducible (where either $Y_{1}$ or $Y_{2}$ may be null). In this case, remark that such a Borel partition is unique up to null sets and each $Y_{i}$ is ${\cal S}$-invariant (see Remark \ref{rem-disjoint-sum-irreducible-amenable-reducible}).

\begin{thm}\label{main-normal}
With the assumption $(\ast)'$, suppose that we have a subrelation ${\cal T}$ of $({\cal R})_{Y}$ containing ${\cal S}$ as a subrelation. Moreover, assume that ${\cal S}$ is normal in ${\cal T}$ and is amenable. Then ${\cal T}$ is a disjoint sum of an irreducible and amenable relation and a reducible relation.
\end{thm}

\begin{pf}
Since ${\cal S}$ is amenable, there exists a $\rho$-invariant Borel map $\varphi \colon Y\rightarrow M({\cal PMF})$ for ${\cal S}$. By Theorem \ref{alternative} and Corollary \ref{alternative-cor}, the Borel subsets
\begin{align*}
Y_{1}&=\{ x\in Y: \varphi(x)({\cal MIN})=1\},\\
Y_{2}&=\{ x\in Y: \varphi(x)({\cal PMF}\setminus {\cal MIN})=1\} 
\end{align*}
satisfy $Y=Y_{1}\sqcup Y_{2}$ and each $Y_{i}$ is ${\cal S}$-invariant. 

If $\mu(Y_{1})>0$, then $({\cal S})_{Y_{1}}$ is normal in $({\cal T})_{Y_{1}}$. Note that any Borel map in ${\rm End}_{({\cal T})_{Y{1}}}(({\cal S})_{Y_{1}})'$ can be expressed as a countable disjoint sum of elements in  $[[({\cal T})_{Y_{1}}]]_{({\cal S})_{Y_{1}}}$ as in the proof of Lemma \ref{generate-generate-normal}. Let $G$ be a countable subset of $[[({\cal T})_{Y_{1}}]]_{({\cal S})_{Y_{1}}}$ constructed as above from a family of normal choice functions for the pair $(({\cal T})_{Y_{1}}, ({\cal S})_{Y_{1}})$ and apply Lemma \ref{main-lem1}. Since $({\cal T})_{Y_{1}}$ is generated by $({\cal S})_{Y_{1}}$ and the choice functions, the relation $({\cal T})_{Y_{1}}$ has a $\rho$-invariant Borel map $Y_{1}\rightarrow \partial_{2}C$. It follows from Corollary \ref{cor-irreducible-amenable-cor} that $({\cal T})_{Y_{1}}$ is irreducible and amenable. 

Similarly, if $\mu(Y_{2})>0$, then we have a $\rho$-invariant Borel map $Y_{2}\rightarrow S(M)$ for $({\cal T})_{Y_{2}}$, using Lemma \ref{main-lem2}.   
\end{pf}

The next corollary follows from Proposition \ref{typeII} and Theorem \ref{main-normal}. It is an analogue of \cite[Theorem 4.8]{fsz}. 

\begin{cor}\label{cor-ame-normal}
With the assumption $(\diamond)$, suppose that we have a discrete measured equivalence relation ${\cal R}'$ weakly isomorphic to ${\cal R}$. Then ${\cal R}'$ admits no recurrent amenable normal subrelations.
\end{cor}

Using these results, we show that the mapping class group of a surface $M$ with $\kappa(M)\geq 0$ can not be measure equivalent to groups of some type. We shall recall the notion of measure equivalence among discrete groups. We recommend the reader to see \cite{furman1}, \cite{furman2} and \cite{gab2} for details. 

\begin{defn}[\ci{0.5.E}{gromov2}]\label{defn-me}
Two discrete groups $\Gamma$ and $\Lambda$ are said to be {\it measure equivalent}\index{measure equivalent} if there exist commuting, measure-preserving, essentially free actions of $\Gamma$ and $\Lambda$ on some standard Borel space $(\Omega, m)$ with a $\sigma$-finite non-zero positive measure such that the action of each of the groups $\Gamma$ and $\Lambda$ admits a fundamental domain with finite measure.
\end{defn}

In fact, this defines an equivalence relation among discrete groups \cite[Section 2]{furman1}. Let $G$ be a locally compact second countable group and $\Gamma_{1}$ and $\Gamma_{2}$ be two lattices (i.e., discrete subgroup with cofinite volume) in $G$. Then $G$ is unimodular \cite[Lemma 2.32]{sauer} and the left and right multiplicative actions of $\Gamma_{1}$ and $\Gamma_{2}$ on $G$ with the Haar measure, respectively, show that $\Gamma_{1}$ and $\Gamma_{2}$ are measure equivalent \cite[Example 1.2]{furman1}. Also, a discrete group is measure equivalent to its subgroup of finite index and to its quotient group by a finite normal subgroup. We have presented several notable results on measure equivalence in Chapter \ref{introduction}. The reference \cite{gab-survey} is a quite detailed survey in which various recent progress about measure equivalence are treated. 

\begin{rem}\label{rem-qi-me}
The notion of measure equivalence can be regarded as a measure-theoretical analogue of a coarse geometric notion of quasi-isometry among finitely generated groups by the following observation of Gromov:

\begin{thm}[\ci{${\rm 0.2.C_{2}'}$}{gromov2}, \ci{Theorem 2.14}{sauer}]
Let $\Gamma$ and $\Lambda$ be two finitely generated groups. Then they are quasi-isometric if and only if they are {\rm topologically equivalent}\index{topologically equivalent!finitely generated groups} in the following sense: there exist commuting, continuous actions of $\Gamma$ and $\Lambda$ on some locally compact space such that the action of each of the groups is properly discontinuous and cocompact.
\end{thm}

Note that two cocompact lattices in the same locally compact group are topologically equivalent. However, in general, two lattices in the same locally compact group are not topologically equivalent. As mentioned in Chapter \ref{introduction}, the quasi-isometry classes of irreducible lattices in semisimple Lie groups can be completely understood through combination of several papers by several authors (see \cite{farb} for details). 
\end{rem}

The notion of measure equivalence can be expressed in terms of discrete measured equivalence relations as follows:  

\begin{thm}[\cite{furman2}, \ci{Theorem 2.3}{gab2}]\label{me-eq-oe}
Let $\Gamma$ and $\Lambda$ be two infinite discrete groups. Then $\Gamma$ and $\Lambda$ are measure equivalent if and only if they admit essentially free, measure-preserving ergodic actions on standard probability spaces which generate weakly isomorphic relations.   
\end{thm}

The next corollary follows from Lemmas \ref{lem-normal-subgroup}, \ref{normal-restriction} and Corollary \ref{cor-ame-normal}:

\begin{cor}\label{cor-me-ame-normal}
Let $\Gamma$ be a sufficiently large subgroup of $\Gamma(M; m)$, where $M$ is a surface with $\kappa(M)\geq 0$ and $m\geq 3$ is an integer. Then $\Gamma$ can not be measure equivalent to any group which has an infinite amenable normal subgroup.
\end{cor}

Finally, we note that the above assertions are satisfied also for hyperbolic groups instead of the mapping class group in the following form. We state only the analogue of Corollaries \ref{cor-ame-normal} and \ref{cor-me-ame-normal}, which can be shown by using Proposition \ref{irreducible-amenable-hyp} and the same idea for the proof of the above assertions. Needless to say, these are essentially due to Adams \cite{adams2} (see also Remark \ref{rem-ms-ham}).

\begin{thm}\label{prop-hyp-main}
With the assumption $(\ast)_{\rm h}'$, suppose that we have a subrelation ${\cal T}$ of $({\cal R})_{Y}$ containing ${\cal S}$ as a subrelation. Moreover, assume that ${\cal S}$ is normal in ${\cal T}$ and is amenable. Then the relation ${\cal T}$ is also amenable.  
\end{thm}

\begin{thm}\label{thm-hyp-main}
Let $\Gamma$ be a non-elementary subgroup of a hyperbolic group. 
\begin{enumerate}
\item[(i)] Suppose that we have an essentially free, measure-preserving Borel action of $\Gamma$ on a standard Borel space $(X, \mu)$ with a finite positive measure. Let ${\cal R}$ be the induced equivalence relation on $(X, \mu)$ and suppose that we have a discrete measured equivalence relation ${\cal R}'$ weakly isomorphic to ${\cal R}$. Then ${\cal R}'$ admits no non-amenable subrelation ${\cal S}$ which has a recurrent amenable normal subrelations of ${\cal S}$.
\item[(ii)] The group $\Gamma$ can not be measure equivalent to any group containing a non-amenable subgroup $G$ which has an infinite amenable normal subgroup of $G$.
\end{enumerate}
\end{thm}


\subsection{Indecomposability}\label{subsection-indec}

Given a discrete measured equivalence relation ${\cal R}_{i}$ on a standard Borel space $(X_{i}, \mu_{i})$ with a finite positive measure for $i=1, 2$, we denote by ${\cal R}_{1}\times {\cal R}_{2}$\index{$R 1 a R 2$@${\cal R}_{1}\times {\cal R}_{2}$} the relation
\[\{ ((x_{1}, x_{2}), (y_{1}, y_{2}))\in (X_{1}\times X_{2})^{2}: (x_{1}, y_{1})\in {\cal R}_{1}, \ (x_{2}, y_{2})\in {\cal R}_{2}\} \]
on the product space $(X_{1}\times X_{2}, \mu_{1}\times \mu_{2})$. Note that if ${\cal S}_{i}$ is a normal subrelation of ${\cal R}_{i}$ for $i=1, 2$, then both ${\cal R}_{1}\times {\cal S}_{2}$ and ${\cal S}_{1}\times {\cal R}_{2}$ are normal in ${\cal R}_{1}\times {\cal R}_{2}$.

\begin{thm}\label{main-thm}
With the assumption $(\ast)'$, let ${\cal R}_{i}$ be a non-singular recurrent discrete measured equivalence relation on a standard Borel space $(X_{i}, \mu_{i})$ with a finite positive measure for $i=1, 2$. Suppose that we have a Borel subset $Z$ of $X_{1}\times X_{2}$ with positive measure and that $({\cal R}_{1}\times {\cal R}_{2})_{Z}$ and ${\cal S}$ are weakly isomorphic. Then ${\cal S}$ is a disjoint sum of an irreducible and amenable relation and a reducible relation.
\end{thm}

\begin{pf}
In this proof, for a subrelation ${\cal R}_{1}'$ of ${\cal R}_{1}$, the symbol ${\cal R}_{1}'$ means the relation ${\cal R}_{1}'\times {\cal D}_{X_{2}}$ on $(X_{1}\times X_{2}, \mu_{1}\times \mu_{2})$, where ${\cal D}_{X_{2}}$ denotes the trivial relation on $X_{2}$. Similarly, for a subrelation ${\cal R}_{2}'$ of ${\cal R}_{2}$, the symbol ${\cal R}_{2}'$ means the relation ${\cal D}_{X_{1}}\times {\cal R}_{2}'$ on $(X_{1}\times X_{2}, \mu_{1}\times \mu_{2})$. 

We may suppose that there exists a Borel isomorphism $f\colon Z \rightarrow Y$ which induces the isomorphism between the two relations $({\cal R}_{1}\times {\cal R}_{2})_{Z}$ and ${\cal S}$. We denote
\[f({\cal R}')=\{ (f(x), f(y))\in {\cal S}: (x, y)\in {\cal R}'\} \]
for a subrelation ${\cal R}'$ of $({\cal R}_{1}\times {\cal R}_{2})_{Z}$. Then $f({\cal R}')$ is a subrelation of ${\cal S}$.

Let ${\cal T}$ be a recurrent amenable subrelation of ${\cal R}_{1}$ (see Lemma \ref{lem-ame-rec-exist}). By the amenability, there exists a $\rho$-invariant Borel map
\[\varphi \colon Y\rightarrow M({\cal PMF})\]
for the relation $f(({\cal T})_{Z})$. By Theorem \ref{alternative} and Corollary \ref{alternative-cor}, if we define  
\begin{align*}
Z_{1}&=\{ x\in Z: \varphi(f(x))({\cal MIN})=1\},\\   
Z_{2}&=\{ x\in Z: \varphi(f(x))({\cal PMF}\setminus {\cal MIN})=1\},
\end{align*}
then $Z=Z_{1}\sqcup Z_{2}$ and both $Z_{1}$ and $Z_{2}$ are $({\cal T})_{Z}$-invariant.  

If $(\mu_{1}\times \mu_{2})(Z_{1})>0$, since $f(({\cal T})_{Z_{1}})$ is normal in $f(({\cal T}\times {\cal R}_{2})_{Z_{1}})$, then $f(({\cal T}\times {\cal R}_{2})_{Z_{1}})$ is irreducible and amenable by Corollary \ref{cor-irreducible-amenable-cor} and Lemma \ref{main-lem1}. In particular, $f(({\cal R}_{2})_{Z_{1}})$ is irreducible and amenable. In a similar way, we can show that $f(({\cal R}_{1}\times {\cal R}_{2})_{Z_{1}})$ is irreducible and amenable since ${\cal R}_{2}$ is normal in ${\cal R}_{1}\times {\cal R}_{2}$. 

If $(\mu_{1}\times \mu_{2})(Z_{2})>0$, then $f(({\cal R}_{1}\times {\cal R}_{2})_{Z_{2}})$ is reducible similarly. 
\end{pf}

The next corollary follows from Proposition \ref{typeII} and Theorem \ref{main-thm}. It is an analogue of Theorems \ref{zim3-main} and \ref{adams2-main}.

\begin{cor}
Assume $(\diamond)$. Let ${\cal R}_{i}$ be a non-singular recurrent discrete measured equivalence relation on a standard Borel space $(X_{i}, \mu_{i})$ with finite positive measure for $i=1, 2$. Then there exist no Borel subset $Z$ of $X_{1}\times X_{2}$ with positive measure such that ${\cal R}$ and $({\cal R}_{1}\times {\cal R}_{2})_{Z}$ are weakly isomorphic. 
\end{cor}

We can prove the following corollary with the same idea for the proofs of Theorems \ref{main-normal} and \ref{main-thm}:

\begin{cor}\label{cor-me-product}
Let $\Gamma$ be a sufficiently large subgroup of $\Gamma(M)$, where $M$ is a surface with $\kappa(M)\geq 0$. Then $\Gamma$ is not measure equivalent to any group of the form $G_{1}\times G_{2}$ such that each $G_{i}$ is infinite and either $G_{1}$ or $G_{2}$ has an infinite amenable subgroup.
\end{cor}

\begin{pf}
We give only a sketch of the proof. Suppose that the corollary is not true and that we have a group $G_{1}\times G_{2}$ measure equivalent to $\Gamma$ such that $G_{1}$ contains an infinite amenable subgroup $G$. We may assume that $\Gamma$ is a subgroup of $\Gamma(M;m)$ with $m\geq 3$. Then there would exist weakly isomorphic ergodic relations ${\cal R}_{i}$ of type ${\rm II}_{1}$ on $(X_{i}, \mu_{i})$ for $i=1, 2$ such that ${\cal R}_{1}$ and ${\cal R}_{2}$ are generated by essentially free, measure-preserving actions of $\Gamma$ and $G_{1}\times G_{2}$, respectively. Let $f\colon A_{1}\rightarrow A_{2}$, $A_{i}\subseteq X_{i}$ be a Borel isomorphism which induces the weak isomorphism between ${\cal R}_{1}$ and ${\cal R}_{2}$. We denote by $\rho \colon {\cal R}_{1}\rightarrow \Gamma$ the induced cocycle. 

Let us denote by ${\cal S}$ the amenable recurrent subrelation of ${\cal R}_{2}$ generated by the action of the subgroup $G$. The subrelation $f^{-1}(({\cal S})_{A_{2}})$ has the $\rho$-invariant Borel map from $A_{1}$ into $\partial_{2}C$ or $S(M)$ with the maximal property. Since $G$ is a normal subgroup of the group $G\times G_{2}$, using the same technique in the proof of Theorem \ref{main-normal}, we can show the existence of a $\rho$-invariant Borel map from $A_{1}$ into $\partial_{2}C$ or $S(M)$ for the subrelation $f^{-1}(({\cal S}_{1})_{A_{2}})$ of $({\cal R}_{1})_{A_{1}}$, where ${\cal S}_{1}$ denotes the relation generated by the subgroup $G\times G_{2}$. 

Let ${\cal S}_{2}$ denote the subrelation of ${\cal R}_{2}$ generated by the action of $G_{2}$. Since there exists a $\rho$-invariant Borel map from $A_{1}$ into $\partial_{2}C$ or $S(M)$ for the relation $f^{-1}(({\cal S}_{2})_{A_{2}})\subseteq f^{-1}(({\cal S}_{1})_{A_{2}})$, we have the $\rho$-invariant Borel map for $f^{-1}(({\cal S}_{2})_{A_{2}})$ with the maximal property. Similarly, it follows from the normality of $G_{2}$ in $G_{1}\times G_{2}$ that we can find a $\rho$-invariant Borel map from $A_{1}$ into $\partial_{2}C$ or $S(M)$ for the relation $({\cal R}_{1})_{A_{1}}$. Thus, we see that $({\cal R}_{1})_{A_{1}}$ is a disjoint sum of an irreducible and amenable relation and a reducible relation. This contradicts Proposition \ref{typeII}.     
\end{pf}


\section{Comparison with hyperbolic groups}

In this last section, we remark a certain difference between the mapping class group and hyperbolic groups.

\begin{thm}\label{thm-hyp-mcg-not-me}
The mapping class group of a non-exceptional surface is not measure equivalent to any non-elementary subgroup of a hyperbolic group.
\end{thm}

\begin{pf}
We apply Theorem \ref{thm-hyp-main} (ii). Let $M$ be a non-exceptional surface. Choose $\alpha, \beta_{1}, \beta_{2}\in V(C)$ such that $i(\alpha, \beta_{1})=i(\alpha, \beta_{2})=0$ and $i(\beta_{1}, \beta_{2})\neq 0$. We can choose such curves because $M$ is non-exceptional. Then the subgroup generated by the Dehn twists about $\beta_{1}$ and $\beta_{2}$ is non-amenable by the remark right before Theorem \ref{tits-alternative-mcg}. Moreover, each element in the subgroup commutes with the Dehn twist about $\alpha$. Hence, the mapping class group of $M$ has a non-amenable subgroup $G$ which contains an infinite amenable normal subgroup of $G$.  
\end{pf}

\begin{rem}
This difference between the mapping class group $\Gamma(M)$ for a non-exceptional surface $M$ and a hyperbolic group comes from the difference of properties of the centralizers (or normalizers) for their elements as mentioned in the beginning of this chapter. 

It also affects a difference between their group von Neumann algebras. A von Neumann algebra is said to be {\it diffuse}\index{diffuse} if it does not contain a minimal projection. We say that a von Neumann algebra ${\cal M}$ is said to be {\it solid}\index{solid} if for any diffuse subalgebra ${\cal A}$, its relative commutant
\[\{ x\in {\cal M}: ax=xa {\rm \ for \ any \ } a\in {\cal A}\}\]
is injective (see \cite{ozawa2}). The notion of injectivity for a von Neumann algebra corresponds to the notion of amenability for a group. It follows from the existence of non-amenable centralizers of infinite order elements in $\Gamma(M)$ (see the proof of Theorem \ref{thm-hyp-mcg-not-me}) that the group von Neumann algebra of $\Gamma(M)$ is not solid. On the other hand, Ozawa \cite{ozawa2} verified the solidity of the group von Neumann algebras associated with hyperbolic groups.

The solidity of the group von Neumann algebras for hyperbolic groups gives the following celebrated result: a factor ${\cal M}$ (i.e., a von Neumann algebra whose center consists of scalars) is said to be {\it prime}\index{prime} if ${\cal M}$ is isomorphic to a factor of the form ${\cal M}_{1}\otimes {\cal M}_{2}$, then either ${\cal M}_{1}$ or ${\cal M}_{2}$ is finite-dimensional. A discrete group $G$ is said to be {\it ICC}\index{ICC} if the conjugacy class of any non-trivial element in $G$ is infinite. Note that the group von Neumann algebra of $G$ is a factor if and only if $G$ is ICC (see \cite[Chapter V, Proposition 7.9]{tak1}). Since any solid factor is prime, the group factor of a non-elementary hyperbolic ICC group is prime. This result can be viewed as a von Neumann algebra analogue of Theorem \ref{adams2-main} due to Adams. Although Ozawa's proof is completely different from Adams' one, the amenability of the action on the boundary plays a crucial role in both of the proofs.

Hence, it is natural to ask whether the group factors associated with various ICC subgroups of the mapping class group are prime or not (see Remark \ref{rem-ms-ham} for the existence of ICC subgroups of the mapping class group).  
\end{rem}

\begin{rem}
Non-elementary hyperbolicity for finitely generated groups is not invariant under measure equivalence. For example, the free product ${\Bbb Z}^{2}\ast {\Bbb Z}$ is measure equivalent to the free group of rank $2$ \cite[Section 2.2]{gab-survey}. 
\end{rem}

\begin{rem}\label{rem-ms-ham}
Monod-Shalom \cite[Section 7]{ms} introduced the class ${\cal C}_{\rm reg}$ consisting of all discrete groups $G$ for which the second bounded cohomology group $H_{\rm b}^{2}(G, \ell^{2}(G))$ is non-trivial. They showed that being in the class ${\cal C}_{\rm reg}$ is an invariant under measure equivalence \cite[Theorem 1.18]{ms}. 

This class contains all non-elementary subgroups of a hyperbolic group \cite[Theorem 3]{mms}. Moreover, any torsion-free group in the class ${\cal C}_{\rm reg}$ is ICC \cite[Proposition 7.11]{ms}. On the other hand, this class does not contain any group of the following form \cite[Proposition 7.10]{ms}:
\begin{enumerate}
\item[(i)] the direct product of two infinite groups.
\item[(ii)] a group containing an infinite amenable normal subgroup.
\end{enumerate}

Hamenst\"adt \cite[Corollary B]{ham2} shows that any sufficiently large subgroup in the mapping class group of a surface $M$ with $\kappa(M)\geq 0$ is in the class ${\cal C}_{\rm reg}$. Therefore, from the above Monod-Shalom's result, any torsion-free sufficiently large subgroup is ICC. Moreover, we can deduce the following results which contain Corollaries \ref{cor-me-ame-normal} and \ref{cor-me-product}: any sufficiently large subgroup of $\Gamma(M)$ is measure equivalent neither to a discrete group which has an infinite amenable normal subgroup nor to the direct product of two infinite groups.

\end{rem}




\chapter[Classification in terms of measure equivalence I]{Classification of the mapping class groups in terms of measure equivalence I}\label{chapter-best}

In this chapter, we give some classification result about the mapping class groups from the viewpoint of measure equivalence. For this, we need a more detailed investigation on subrelations of equivalence relations generated by the mapping class group. 

In Section \ref{section-reducible-revisited}, we study reducible subrelations again. We classify subsurfaces of a surface cut by all curves in the canonical reduction system of a reducible subrelation into three types, trivial ones, irreducible and amenable ones and irreducible and non-amenable ones. 

In the case of a reducible subgroup $G$ contained in $\Gamma(M;m)$ for some $m\geq 3$, we have a decomposition $M_{L}=\bigsqcup_{i}Q_{i}$ of $M$ by cutting $M$ along a realization $L$ of $\sigma(G)$. Moreover, we obtain natural homomorphisms 
\[p\colon G\rightarrow \prod_{i}\Gamma(Q_{i}),\ \ p_{j}\colon G\rightarrow \Gamma(Q_{j}),\]
where $p_{j}$ is the composition of $p$ and the projection $\prod_{i}\Gamma(Q_{i})\rightarrow \Gamma(Q_{j})$ (see Remark \ref{rem-comp-leave}). We classify subsurfaces $Q_{j}$ into the following three types:
\begin{enumerate}
\item[(i)] if $p_{j}(G)$ is trivial, then $Q_{j}$ is said to be {\it trivial}\index{trivial subsurface};
\item[(ii)] if $p_{j}(G)$ is amenable and contains a pseudo-Anosov element in $\Gamma(Q_{j})$, then $Q_{j}$ is said to be {\it irreducible and amenable}\index{irreducible and amenable!subsurface};
\item[(iii)] if $p_{j}(G)$ contains an independent pair of pseudo-Anosov elements in $\Gamma(Q_{j})$, then $Q_{j}$ is said to be {\it irreducible and non-amenable}\index{irreducible and non-amenable!subsurface}.
\end{enumerate}
Remark that the image $p_{j}(G)$ either trivial or contains a pseudo-Anosov element in $\Gamma(Q_{j})$ (see Theorem \ref{thm-ivanov-last-hope}). The classification of subsurfaces for a reducible subrelation is an analogue of the above classification. However, in the case of relations, there seem no analogues of the image $p_{j}(G)$. Thus, we need a different approach for the definition of subsurfaces of the three types.

In Section \ref{section-irreducible-amenable-subsurface}, we show that for a reducible subrelation ${\cal S}$ on $(X, \mu)$, the action of ${\cal S}$ on an irreducible and amenable subsurface satisfies similar properties to the action of an irreducible and amenable subrelation appearing in Chapter \ref{chapter-indec}. We prove that if $Q$ is an irreducible and amenable subsurface and if the action of ${\cal S}$ on the space of probability measures on ${\cal PMF}(Q)$ has an invariant Borel map $\varphi \colon X\rightarrow M({\cal PMF}(Q))$, then the support of $\varphi(x)$ is contained in ${\cal MIN}(Q)$ and the support of $\pi_{*}\varphi(x)$ consists of at most two points for a.e.\ $x\in X$, where $\pi_{*}$ is the induced map by $\pi \colon {\cal MIN}(Q)\rightarrow \partial C(Q)$.

In Section \ref{section-amenable-action-general}, we prove that non-amenable reducible subrelations always have irreducible and non-amenable subsurfaces.

Finally, in Section \ref{section-classification}, we give a partial result of Theorem \ref{thm-int-classification} (see Corollary \ref{cor-classification1}). The following theorem is a group-theoretic result corresponding to the key observation for the proof (see Appendix \ref{app-group} for the proof): 

\begin{thm}\label{main-group-version-pre}
Let $M$ be a compact orientable surface of type $(g, p)$ satisfying $\kappa(M)\geq 0$. Let $\Gamma(M)$ be the mapping class group of $M$. If $G_{1}\times \cdots \times G_{n}$ is a subgroup of $\Gamma(M)$ with each $G_{i}$ non-amenable, then 
\[n\leq g+\left[ \frac{g+p-2}{2}\right].\]
Moreover, there exists a subgroup of the form $G_{1}\times \cdots \times G_{n}$ with each $G_{i}$ non-amenable and 
\[n=g+\left[ \frac{g+p-2}{2}\right], \]
where for $a\in {\Bbb R}$, we denote by $[a]$ the maximal integer less than or equal to $a$. 
\end{thm}

The number $g+[(g+p-2)/2]$ is equal to the maximal number of components of $M_{L}$ which is not a pair of pants among realizations $L$ of all elements in $S(M)\cup \{ \emptyset \}$. 

If $G_{1}\times \cdots \times G_{n}$ is a subgroup of $\Gamma(M;m)$ for some $m\geq 3$ such that each $G_{i}$ is non-amenable and $n\geq 2$, then it is easy to see that the subgroup $G_{1}\times \cdots \times G_{n}$ is reducible and that any two curves in the union $\sigma =\bigcup_{i}\sigma(G_{i})$ do not intersect. The key ingredient for the proof of Theorem \ref{main-group-version-pre} is the following two facts: 
\begin{enumerate}
\item[(a)] Since $G_{i}$ is non-amenable, each $G_{i}$ has an irreducible and non-amenable subsurface, which is also a component of $M_{L}$, where $L$ is a realization of $\sigma$. 
\item[(b)] For each $i\neq j$, irreducible and non-amenable subsurfaces for $G_{i}$ and $G_{j}$ are disjoint each other. 
\end{enumerate}
Hence, the first part of Theorem \ref{main-group-version-pre} follows since a pair of pants can not be an irreducible and non-amenable subsurface by definition. It is easy to find a subgroup $G_{1}\times \cdots \times G_{n}$ as in the second part of Theorem \ref{main-group-version-pre}. We give a proof of Corollary \ref{cor-classification1} along this idea for the proof of Theorem \ref{main-group-version-pre}. (In Section \ref{section-amenable-action-general}, we show a corresponding property to the fact (a) in terms of relations.)

In general, relations generated by an essentially free action of the direct product of two groups can not be expressed as a product ${\cal R}_{1}\times {\cal R}_{2}$ of relations of the form treated in Chapter \ref{chapter-indec}, Section \ref{subsection-indec}. Therefore, in the case of relations, there seem no notions corresponding to the direct product of groups. 

Instead of considering direct products, we use the notion of normal subrelations. For two groups $H_{1}$ and $H_{2}$, each of them is normal in the direct product $H_{1}\times H_{2}$. Correspondingly, consider the situation where two relations ${\cal S}_{1}$ and ${\cal S}_{2}$ on the same standard Borel space satisfy that each ${\cal S}_{i}$ is normal in ${\cal S}_{1}\vee {\cal S}_{2}$. We regard this situation as an analogue of the direct product of two groups. It goes without saying that when a group $H$ has two subgroups $H_{1}$, $H_{2}$ such that $H_{1}$ and $H_{2}$ generate $H$ and each $H_{i}$ is normal in $H$, the group $H$ is not necessarily the direct product of $H_{1}$ and $H_{2}$. Therefore, for the proof of Corollary \ref{cor-classification1}, we need more technical and artificial conditions in the above situation (see Definitions \ref{defn-dagger} and \ref{defn-dagger-n}).

Needless to say, the proof of Theorem \ref{main-group-version-pre} is much easier than the proof for the case of measure equivalence in this chapter. Hence, we recommend the reader who wants to understand the outline of this chapter to see Appendix \ref{app-group}, first.


\section{Reducible subrelations, revisited}\label{section-reducible-revisited}

\begin{defn}\label{defn-geometric-things}
Let $M$ be a compact orientable surface with $\kappa(M)\geq 0$. 
\begin{enumerate}
\item[(i)] A {\it subsurface}\index{subsurface} of $M$ is (an isotopy class of) the following compact submanifold $Q$ of $M$: there exist $\sigma \in S(M)\cup \{ \emptyset \}$, its realization $L$ and an open tubular neighborhood $N$ of $L$ such that $Q$ is a component of $M\setminus N$. We have the natural embedding $Q\rightarrow M$. In what follows, for simplicity, we often denote by $M_{\sigma}$ the resulting surface $M_{L}$ by cutting $M$ along $L$. We identify a component of $M_{\sigma}$ with the corresponding subsurface of $M$.

\item[(ii)] Let $Q$ be a component of $M_{L}$ with a realization $L$ on $M$ of some $\sigma \in S(M)\cup \{ \emptyset \}$. Let $p_{L}\colon M_{L}\rightarrow M$ be the natural map. A component of the set $p_{L}(Q)\cap L$ is called a {\it relative boundary component}\index{relative boundary component} of $Q$. We denote by $\partial_{M}Q$\index{$\ z M Q $@$\partial_{M}Q$} the set of all relative boundary components of $Q$. 

\item[(iii)] We denote by $D=D(M)$\index{$D=D M$@$D=D(M)$} the set of all isotopy classes of subsurfaces in $M$, on which the mapping class group $\Gamma(M)$ acts naturally. Let ${\cal F}_{0}(D)$\index{$F 0 D $@${\cal F}_{0}(D)$} denote the set of all finite subsets $F$ of $D$ (including the empty set) such that if $Q_{1}, Q_{2}\in F$ and $Q_{1}\neq Q_{2}$, then $Q_{1}$ and $Q_{2}$ can be realized disjointly on $M$. 

\item[(iv)] We say that $F\in {\cal F}_{0}(D)$ {\it fills}\index{fill} $M$ when the following condition is satisfied: there exist $\sigma \in S(M)\cup \{ \emptyset \}$, its realization $L$ and an open tubular neighborhood $N$ such that $F$ is the set of all components of $M\setminus N$. 
\end{enumerate}
\end{defn}

\begin{defn}\label{defn-fill-surface}
Suppose that $\{ \alpha_{1}, \ldots, \alpha_{n}\}$ is a finite set consisting of at least two elements in $V(C(M))$ and satisfies the following condition: for any $k, l\in \{1, \ldots, n\}$, there exists a sequence $m_{1}, \ldots, m_{r}\in \{ 1, \ldots, n\}$ such that $m_{1}=k$, $m_{r}=l$ and $i(\alpha_{m_{j}}, \alpha_{m_{j+1}})\neq 0$ for any $j$. We define the {\it subsurface filled by the curves}\index{subsurface!filled by $\{ \alpha_{k}\}$} $\{ \alpha_{k}\}$ as the following subsurface $Q$: if $G$ is a subgroup of $\Gamma(M)$ generated by the Dehn twists of all curves in $\{ \alpha_{k}\}$, then let $Q$ be a component of $M_{\sigma(G)}$ containing the curves $\alpha_{k}$. 
\end{defn}

First, note that any curve $\alpha_{k}$ is disjoint from any curve in $\sigma(G)$ by Corollary \ref{cor-intersecting-not-invariant}. Secondly, remark that for any realization $L_{k}$ of $\alpha_{k}$, the union $\bigcup_{k=1}^{n}L_{k}$ is connected. We call a finite subset $\{ \alpha_{k}\}$ of $V(C(M))$ satisfying the above condition a {\it connected chain}\index{connected chain}. It follows from this remark that the connected component $Q$ in Definition \ref{defn-fill-surface} is uniquely determined. We write $Q$ as $S(\{ \alpha_{k}\})$\index{$S \a k$@$S(\{ \alpha_{k}\})$}.

\begin{lem}\label{lem-fill-subsurface}
Let $\{ \alpha_{k}\}$ be a connected chain in $V(C(M))$ and $Q=S(\{ \alpha_{k}\})$.
\begin{enumerate}
\item[(i)] If $\beta$ is an element in $V(C(Q))$, then there exists $k$ such that $i(\alpha_{k}, \beta)\neq 0$.
\item[(ii)] If $g\in \Gamma(M)$ is a pure torsion-free element fixing (the set or any curves in) $\{ \alpha_{k}\}$, then $g$ fixes $Q$ and any element in $V(C(Q))$. 
\end{enumerate} 
\end{lem}

\begin{pf}
Let $G$ be the intersection of $\Gamma(M;m)$ and the subgroup generated by the Dehn twists of all curves $\alpha_{k}$. Then we have the natural homomorphism $p\colon G\rightarrow \Gamma(Q)$ (see Remark \ref{rem-comp-leave}). If $\beta \in V(C(Q))$ satisfied $i(\alpha_{k}, \beta)=0$ for any $k$, then the image $p(G)$ would be an infinite reducible subgroup fixing the curve $\beta \in V(C(Q))$. This contradicts Theorem \ref{thm-ivanov-last-hope} and proves the assertion (i).   

If $g\in \Gamma(M)$ is a pure torsion-free element fixing any curves in $\{ \alpha_{k}\}$, then it is clear that $g$ fixes $Q$ by definition. Let $G'$ be the intersection of $\Gamma(M;m)$ and the subgroup generated by $g$ and the Dehn twists about all curves in $\{ \alpha_{k}\}$. We have a homomorphism 
\[p\colon G'\rightarrow \Gamma(Q)\]
which is an extension of $p\colon G\rightarrow \Gamma(Q)$ defined above. Since $G$ is a normal subgroup of $G'$, we see that $\sigma(G)\subseteq \sigma(G')$ by Remark \ref{rem-main-lem2}. It follows from the assertion (i) and Corollary \ref{cor-intersecting-not-invariant} that there exist no curves in $\sigma(G')$ which is also an element in $V(C(Q))$. Thus, $Q$ is also a component of $M_{\sigma(G')}$. It follows from Theorem \ref{thm-ivanov-last-hope} that $p(G')$ is finite and from Theorem \ref{pure-important} (ii) that $g$ fixes any element in $V(C(Q))$. 
\end{pf}

\begin{rem}\label{rem-kappa-increase}
Let $\{ \alpha_{k}\}$ be a connected chain in $V(C(M))$ and $Q=S(\{ \alpha_{k}\})$. We suppose that there exists a relative boundary component $\beta$ of $Q$. (If there exist no such components, then $Q=M$.) 

For any $\gamma \in V(C(M))$ with $i(\gamma, \beta)\neq 0$, we consider the finite subset $r(\gamma, Q)$ of $V(C(Q))$ (see the comment right before Lemma \ref{pants}). Since $Q$ is not a pair of pants, the set $r(\gamma, Q)$ is non-empty by Lemma \ref{pants}. It follows from the definition of the set $r(\gamma, Q)$ that for each $k$, if $i(\gamma, \alpha_{k})=0$, then $i(\gamma', \alpha_{k})=0$ for any $\gamma'\in r(\gamma, Q)$. By Lemma \ref{lem-fill-subsurface} (i), there exists $k$ such that $i(\gamma, \alpha_{k})\neq 0$.

Thus, we can consider the subsurface $R$ filled by the curves $\{ \alpha_{k}\} \cup \{ \gamma \}$. Let $G_{Q}$, $G_{R}$ be the groups generated by the Dehn twists of all curves in $\{ \alpha_{k}\}$, $\{ \alpha_{k}\}\cup \{ \gamma \}$, respectively. 

\begin{lem}
\begin{enumerate}
\item[(i)] Any two curves in the union $\sigma(G_{Q})\cup \sigma(G_{R})$ do not intersect. Thus, $\sigma(G_{Q})\cup \sigma(G_{R})$ is an element in $S(M)$.
\item[(ii)] The subsurface $R$ contains $Q$. In other words, when we realize $\sigma(G_{Q})\cup \sigma(G_{R})$ on $M$ and identify each component of the realization with the corresponding element in $\sigma(G_{R})\cup \sigma(G_{Q})$, the corresponding open component to $R$ of $M \setminus \sigma(G_{R})$ contains the corresponding open component to $Q$ of $M\setminus \sigma(G_{Q})$.
\item[(iii)] The relative boundary component $\beta$ of $Q$ is an element in $V(C(R))$.
\end{enumerate}
\end{lem}

\begin{pf}
For any $\alpha \in \sigma(G_{Q})$ and $\alpha'\in \sigma(G_{R})$, we have $i(\alpha, \alpha')=0$ because if $i(\alpha, \alpha')\neq 0$, then this would contradict the essentiality of $\alpha$ for $G_{Q}$. 

If $R$ did not contain $Q$, then there would exist $\alpha'\in \sigma(G_{R})$ which is in $V(C(Q))$. Since $\alpha'$ is invariant for the subgroup $G_{R}$ which contains $G_{Q}$, this contradicts Corollary \ref{cor-intersecting-not-invariant} and Lemma \ref{lem-fill-subsurface} (i). It completes the proof of the assertion (ii).

It follows from $\beta \in \sigma(G_{Q})$ that $i(\beta, \alpha')=0$ for any $\alpha'\in \sigma(G_{R})$. Since $i(\beta, \gamma)\neq 0$ and $\gamma \in V(C(R))$, we see that $\beta$ is not an element in $\sigma(G_{R})$ and that $\beta$ can be realized on $R$. 
\end{pf}

Since the subsurface $R$ contains $\beta$ as a non-trivial, non-peripheral circle, we see that $\kappa(Q)<\kappa(R)$. Recall that the number $\kappa(M)+1=3g+p-3$ is equal to the maximal number of non-trivial isotopy classes of non-peripheral simple closed curves on $M$ which can be realized disjointly on $M$.    
\end{rem}

In what follows, the symbol ${\cal S}$ means a reducible relation on a standard Borel space $(Y, \mu)$ with a finite positive measure satisfying the assumption $(\ast)$. Let us denote by $\varphi \colon Y\rightarrow S(M)$ the CRS for ${\cal S}$ defined in Definition \ref{crs}.

\begin{rem}\label{rem-curve-comp-leave}
Let $\sigma \in S(M)$ be an element with $\mu(\varphi^{-1}(\sigma))>0$. Set $Z=\varphi^{-1}(\sigma)$. Recall that the Borel map $\varphi$ is $\rho$-invariant for ${\cal S}$ by Lemma \ref{crs-fund} (ii). It follows that $\rho(x, y)\sigma =\sigma$ for a.e.\ $(x, y)\in ({\cal S})_{Z}$. Moreover, it follows from Theorem \ref{comp-leave} that $\rho(x, y)$ fixes any curve in $\sigma$ and any component of $M_{\sigma}$. 
\end{rem}

\begin{lem}\label{lem:trivial}
Assume that $\varphi$ is constant and let $\sigma \in S(M)$ be its value. Let $Q$ be a component of $M_{\sigma}$. Suppose that there exist $\alpha \in V(C(M))$ with $r(\alpha, Q)\neq \emptyset$ and $A\subseteq Y$ with positive measure such that $(\alpha, A)$ is $\rho$-invariant. Then there exists a countable Borel partition $A=\bigsqcup A_{n}$ such that 
\[\rho(x, y)\beta =\beta\]
for any $\beta \in V(C(Q))$ and a.e.\ $(x, y)\in ({\cal S})_{A_{n}}$. In particular, the pair $(\beta, A)$ is $\rho$-invariant for any curve $\beta$ in $V(C(Q))$.
\end{lem}

\begin{pf}
Since $r(\alpha, Q)\neq \emptyset$, the subsurface $Q$ is not a pair of pants (see Lemma \ref{pants}). 

By $\rho$-invariance of $(\alpha, A)$, there exists a countable Borel partition $A=\bigsqcup A_{n}$ such that $\alpha$ is invariant for $({\cal S})_{A_{n}}$ for all $n$. Since
\[\rho(x, y)r(\alpha, Q)=r(\rho(x, y)\alpha, Q)=r(\alpha, Q)\]
for a.e.\ $(x, y)\in ({\cal S})_{A_{n}}$ by Remark \ref{rem-curve-comp-leave} and $r(\alpha, Q)$ is finite, we see that any $\alpha'\in r(\alpha, Q)$ is also invariant for $({\cal S})_{A_{n}}$. Thus, by replacing $\alpha$ by $\alpha'$, we may assume that $\alpha$ is a curve in $V(C(Q))$. 

It follows from $\alpha \notin \sigma$ that $(\alpha, A')$ is not essential for any Borel subset $A'$ of $A$ with positive measure. Hence, there exist a countable Borel partition $A=\bigsqcup B_{n}$ and $\alpha_{n}'\in V(C(M))$ such that $i(\alpha_{n}', \alpha)\neq 0$ and $(\alpha_{n}', B_{n})$ is $\rho$-invariant (see the proof of Lemma \ref{essential-deny}). It follows from $i(\alpha_{n}', \alpha)\neq 0$ that $\alpha_{n}'$ is not a relative boundary component of $Q$. Since $Q$ is not a pair of pants, we have $r(\alpha_{n}', Q)\neq \emptyset$ by Lemma \ref{pants}. By Lemma \ref{lem-pants-r}, there exists $\alpha_{n}\in r(\alpha_{n}', Q)$ such that $i(\alpha_{n}, \alpha)\neq 0$. As above, we can show that $(\alpha_{n}, B_{n})$ is $\rho$-invariant. 

For each $n$, if $\partial_{Q}(S(\alpha, \alpha_{n}))$ is empty, then $S(\alpha, \alpha_{n})=Q$. Moreover, by Lemma \ref{lem-fill-subsurface} (ii), we see that elements in $\Gamma$ fixing $\alpha$ and $\alpha_{n}$ fixes also any curve in $V(C(Q))$. This completes the proof.

Next, we assume that $\partial_{Q}(S(\alpha, \alpha_{n}))$ is non-empty. Take a boundary component $\beta_{n}$ in $\partial_{Q}(S(\alpha, \alpha_{n}))$. Since $(\beta_{n}, B_{n})$ is $\rho$-invariant by Lemma \ref{lem-fill-subsurface} (ii), but not essential, there exist a countable Borel partition $B_{n}=\bigsqcup_{m}B_{nm}$ and $\alpha_{nm}\in V(C(M))$ such that $i(\beta_{n}, \alpha_{nm})\neq 0$ and $(\alpha_{nm}, B_{nm})$ is $\rho$-invariant. As above, we can assume that $\alpha_{nm}$ is a curve in $V(C(Q))$. Since $i(\beta_{n}, \alpha_{nm})\neq 0$, we can construct the subsurface $S(\alpha, \alpha_{n}, \alpha_{nm})$ of $Q$ and for each $n$, we have 
\[\kappa(S(\alpha, \alpha_{n}))<\kappa(S(\alpha, \alpha_{n}, \alpha_{nm}))\leq \kappa(Q)\] 
by Remark \ref{rem-kappa-increase}.

Continuing this process, we obtain a countable Borel partition $A=\bigsqcup C_{n}$ and families $\{ \alpha_{k}^{n}\}_{k=1}^{l_{n}}$ consisting of finite simple closed curves in $V(C(Q))$ with $l_{n}\leq \kappa(Q)$ such that we can construct the subsurface $S(\alpha, \{ \alpha_{k}^{n}\}_{k})$ equal to $Q$ and that the pairs $(\alpha, C_{n})$ and $(\alpha_{k}^{n}, C_{n})$ are $\rho$-invariant. Since elements in $\Gamma$ which fix $\alpha$ and $\alpha_{k}^{n}$ for all $k=1, \ldots, l_{n}$ fix also any curve in $V(C(Q))$, it completes the proof.  
\end{pf}

\begin{prop}\label{maximal-alpha}
With the assumption $(\ast)$, for any $\alpha \in V(C(M))$, there exists an essentially unique Borel subset $Y_{\alpha}\subseteq Y$ such that for any Borel subset $A\subseteq Y\setminus Y_{\alpha}$ with positive measure, the pair $(\alpha, A)$ is not $\rho$-invariant and moreover, if $Y_{\alpha}$ has positive measure, then $(\alpha, Y_{\alpha})$ is $\rho$-invariant.
\end{prop}

\begin{pf}
Let ${\cal P}_{\alpha}$ be the family of all Borel subsets $A\subseteq Y$ such that $\mu(A)>0$ and $(\alpha, A)$ is $\rho$-invariant. It satisfies that if $A_{n}\in {\cal P}_{\alpha}$, then $\bigcup A_{n}\in {\cal P}_{\alpha}$ by Lemma \ref{ess-making} (ii). If ${\cal P}_{\alpha}$ is empty, let us define $Y_{\alpha}=\emptyset$. Otherwise, let $Y_{\alpha}$ be a Borel subset satisfying 
\[\mu(Y_{\alpha})=\sup \{ \mu(A): A\in {\cal P}_{\alpha}\}.\]
This $Y_{\alpha}$ satisfies the desired property. 
\end{pf}

\begin{rem}\label{rem-invariant-domain}
It follows from Lemma \ref{crs-fund} (iii), Lemma \ref{lem:trivial} and Proposition \ref{maximal-alpha} that if $Q$ is a component of $M_{\sigma}$ for some $\sigma \in S(M)$, then 
\[Y_{\alpha}\cap \varphi^{-1}(\sigma)=Y_{\beta}\cap \varphi^{-1}(\sigma)\]
for any $\alpha, \beta \in V(C(Q))$ up to null sets, where $\varphi$ is the CRS for ${\cal S}$.
\end{rem}

\begin{thm}\label{alternative-trivial-irreducible}
With the assumption $(\ast)$, suppose that ${\cal S}$ is reducible and let $\varphi \colon Y\rightarrow S(M)$ be the CRS for ${\cal S}$. Then there exist essentially unique two $\rho$-invariant Borel maps
\begin{align*}
\varphi_{t}&\colon Y\rightarrow {\cal F}_{0}(D),\\
\varphi_{i}&\colon Y\rightarrow {\cal F}_{0}(D)
\end{align*}
such that
\begin{enumerate}
\renewcommand{\labelenumi}{\rm(\roman{enumi})}
\item any element in $\varphi_{t}(y)\cup \varphi_{i}(y)$ is a component of $M_{\varphi(y)}$ for a.e.\ $y\in Y$;
\item the union $\varphi_{t}(y)\cup \varphi_{i}(y)$ fills $M$ and $\varphi_{t}(y)\cap \varphi_{i}(y)=\emptyset$ for a.e.\ $y\in Y$;
\item if $Q$ is in $F\in {\cal F}_{0}(D)$ with $\mu(\varphi_{t}^{-1}(F))>0$, then either $Q$ is a pair of pants or the pair $(\alpha, \varphi_{t}^{-1}(F))$ is $\rho$-invariant for any $\alpha \in V(C(Q))$;
\item if $Q$ is in $F\in {\cal F}_{0}(D)$ with $\mu(\varphi_{i}^{-1}(F))>0$, then $Q$ is not a pair of pants and $(\alpha, A)$ is not $\rho$-invariant for any $\alpha \in V(C(M))$ with $r(\alpha, Q)\neq \emptyset$ and any Borel subset $A\subseteq \varphi_{i}^{-1}(F)$ with positive measure.
\end{enumerate}
\end{thm} 

\begin{pf}
First, we define $\varphi_{t}$ on $\varphi^{-1}(\sigma)$ for $\sigma \in S(M)$. Let $Y_{\alpha}$ be the Borel subset in Proposition \ref{maximal-alpha} for $\alpha \in V(C(M))$. For each component $Q$ of $M_{\sigma}$ which is not a pair of pants, take $\alpha_{Q} \in V(C(Q))$. Define
\begin{align*}
\varphi_{t}(y)=\{ Q: y\in & Y_{\alpha_{Q}}\cap \varphi^{-1}(\sigma)\} \\
                          & \cup \{ {\rm a \ component \ of \ } M_{\sigma} {\rm \ which \ is \ a \ pair \ of \ pants}\}
\end{align*}
for $y\in \varphi^{-1}(\sigma)$. This definition is independent of the choice of $\alpha_{Q}$ by Remark \ref{rem-invariant-domain}. Define
\[\varphi_{i}(y)=\{ {\rm a \ component \ of \ } M_{\sigma}\} \setminus \varphi_{t}(y)\]
for $y\in \varphi^{-1}(\sigma)$. Since $\varphi^{-1}(\sigma)$ and $Y_{\alpha_{Q}}$ are Borel subsets of $Y$ for each $\sigma \in S(M)$ and component $Q$ of $M_{\sigma}$, both the maps $\varphi_{t}$ and $\varphi_{i}$ are measurable. 

It is easy to see the conditions (i), (ii) and (iii). We assume that the condition (iv) is not satisfied. Then there would exist $F\in {\cal F}_{0}(D)$, $Q\in F$, $A\subseteq \varphi^{-1}_{i}(F)$ with positive measure and $\alpha \in V(C(M))$ with $r(\alpha, Q)\neq \emptyset$ such that $(\alpha, A)$ is $\rho$-invariant. We may assume that all $\varphi$, $\varphi_{t}$ and $\varphi_{i}$ are constant on $A$ and thus, there exists $\sigma \in S(M)$ such that $\varphi(y)=\sigma$ for any $y\in A$ and $F$ is a set of components of $M_{\sigma}$. Moreover, we may assume $\alpha \in V(C(Q))$, using the same argument at the beginning of the proof of Lemma \ref{lem:trivial}. However, since it follows that $A\subseteq Y_{\alpha}$, we see that $Q\in \varphi_{t}(y)$ for $y\in A$ by definition, which is a contradiction. It proves the condition (iv).

Finally, we show the $\rho$-invariance of the maps $\varphi_{t}$ and $\varphi_{i}$. Let $g\in [[{\cal S}]]$ and 
\[A_{\gamma}=\{ x\in {\rm dom}(g): \rho(gx, x)=\gamma \}\]
for $\gamma \in \Gamma$. By Lemma \ref{crs-sublem}, if $(\alpha, B)$ with $B\subseteq A_{\gamma}$ is $\rho$-invariant, then so is $(\gamma \alpha, gB)$. Hence, $gB\subseteq Y_{\gamma \alpha}$ up to null sets. It follows that $g(Y_{\alpha}\cap A_{\gamma})\subseteq Y_{\gamma \alpha}$. By $\rho$-invariance of $\varphi$, we see that
\[g(Y_{\alpha_{Q}}\cap \varphi^{-1}(\sigma)\cap A_{\gamma})\subseteq Y_{\gamma \alpha_{Q}}\cap \varphi^{-1}(\gamma \sigma)=Y_{\alpha_{(\gamma Q)}}\cap \varphi^{-1}(\gamma \sigma),\] 
where the last equality comes from Remark \ref{rem-invariant-domain}. This implies the $\rho$-invariance of $\varphi_{t}$. It follows that $\varphi_{i}$ is also $\rho$-invariant.

The uniqueness of $\varphi_{t}$ and $\varphi_{i}$ is clear.
\end{pf}

\begin{defn}
We call $\varphi_{t}$ (resp. $\varphi_{i}$) the {\it system of trivial} (resp. {\it irreducible}) {\it subsurfaces}\index{system!of trivial subsurfaces}\index{system!of irreducible subsurfaces} for ${\cal S}$. In short, we call it the {\it T} (resp. {\it I}) {\it system}\index{T!system}\index{I!system} for ${\cal S}$. We often call an element in $\varphi_{t}(x)$ (resp. $\varphi_{i}(x)$) a {\it trivial} (resp. {\it irreducible}) {\it subsurface}\index{trivial subsurface}\index{irreducible subsurface} for $x\in Y$ and in short, a {\it T} (resp. {\it I}) {\it subsurface}\index{T!subsurface}\index{I!subsurface}.
\end{defn} 

When we identify a subsurface with a component of the cut surface along some curves, we call T and I subsurfaces {\it T} and {\it I components}\index{T!component}\index{I!component}, respectively.

\begin{lem}\label{restriction-of-invariants-first}
With the assumption $(\ast)$, suppose that ${\cal S}$ is reducible. Let us denote by $\varphi$, $\varphi_{t}$ and $\varphi_{i}$ the CRS, T and I systems for ${\cal S}$, respectively. Let $A\subseteq Y$ be a Borel subset with positive measure. Then the CRS, T and I systems for the relation $({\cal S})_{A}$ are the restrictions to $A$ of $\varphi$, $\varphi_{t}$ and $\varphi_{i}$, respectively.
\end{lem}

\begin{pf}
For the CRS, we have shown the lemma in Lemma \ref{crs-fund} (iii). 

Let $\varphi_{t, A}$ and $\varphi_{i, A}$ be the T and I systems for $({\cal S})_{A}$. It follows from Lemma \ref{crs-fund} (iii) and Theorem \ref{alternative-trivial-irreducible} (i), (ii) that $\varphi_{t, A}(x)\cup \varphi_{i, A}(x)=\varphi_{t}(x)\cup \varphi_{i}(x)$ for a.e.\ $x\in A$.

Let $E$ be a Borel subset of $A$ with positive measure such that all $\varphi$, $\varphi_{t}$, $\varphi_{i}$, $\varphi_{t, A}$ and $\varphi_{i, A}$ are constant on $E$ and denote their values on $E$ by $\Phi$, $\Phi_{t}$, $\Phi_{i}$, $\Phi_{t, A}$ and $\Phi_{i, A}$, respectively. Then we have
\[\Phi_{t, A}\cup \Phi_{i, A}=\Phi_{t}\cup \Phi_{i}.\]
If $Q\in \Phi_{t}$ and $\alpha \in V(C(Q))$, then $(\alpha, E)$ is $\rho$-invariant for ${\cal S}$ and thus, for $({\cal S})_{A}$. This means that $Q\in \Phi_{i, A}$ cannot happen and thus, $Q\in \Phi_{t, A}$. Conversely, if $Q\in \Phi_{t, A}$ and $\alpha \in V(C(Q))$, then $(\alpha, E)$ is $\rho$-invariant for $({\cal S})_{A}$ and thus, for ${\cal S}$. This means that $Q\in \Phi_{t}$. Hence, $\varphi_{t}=\varphi_{t, A}$ on $E$. Since $E$ is an arbitrary Borel subset of $A$ on which all the CRS's, T and I systems for ${\cal S}$ and $({\cal S})_{A}$ are constant, we see that $\varphi_{t}=\varphi_{t, A}$ on $A$. It follows that $\varphi_{i}=\varphi_{i, A}$ on $A$.    
\end{pf}

With the assumption $(\ast)$, if $\sigma \in S(M)$ is ${\cal S}$-invariant, then $\rho$ is a cocycle into the stabilizer $\Gamma_{\sigma}$ of $\sigma$ and $Q$ is ${\cal S}$-invariant for each component $Q$ of $M_{\sigma}$. Thus, we have a Borel cocycle ${\cal S}\rightarrow \Gamma(Q)$ by composing $\rho$ with the homomorphism from $\Gamma_{\sigma}$ into $\Gamma(Q)$ (see Remark \ref{rem-comp-leave}). We denote this cocycle by $\rho_{Q}$.

Let $S$ be a Borel $\Gamma(Q)$-space and $Y'$ be a Borel subset of $Y$ with positive measure. For simplicity, we say that a Borel map $\psi \colon Y'\rightarrow S$ is $\rho_{Q}$-invariant (for ${\cal S}$)\index{rho-invariant Borel map@$\rho$-invariant!Borel map} if it is $\rho_{Q}$-invariant for $({\cal S})_{Y'}$.

\begin{lem}\label{irreducible-partition}
With the assumption $(\ast)$, suppose that ${\cal S}$ is reducible. Assume that all $\varphi$, $\varphi_{t}$ and $\varphi_{i}$ are constant on a Borel subset $Y'$ of $Y$ with positive measure and denote their values by the same symbols. Let $Q\in \varphi_{i}$. 

\begin{enumerate}
\item[(a)] There exist essentially unique two ${\cal S}$-invariant Borel subsets $Y_{1}$, $Y_{2}$ of $Y'$ such that
\begin{enumerate}  
\item[(i)] $Y'=Y_{1}\sqcup Y_{2}$ up to null sets;
\item[(ii)] if $Y_{1}$ has positive measure, then there exist two $\rho_{Q}$-invariant Borel maps
\[\psi \colon Y_{11}\rightarrow M({\cal PMF}(Q)), \ \ \psi'\colon Y_{12}\rightarrow \partial_{2}C(Q)\]
such that $Y_{1}=Y_{11}\cup Y_{12}$ (either $Y_{11}$ or $Y_{12}$ may be null). Moreover, any of such a map $\psi$ is $M({\cal MIN}(Q))$-valued; 
\item[(iii)] if $Y_{2}$ has positive measure, then for any Borel subset $A\subseteq Y_{2}$ with positive measure, there exist neither $\rho_{Q}$-invariant Borel maps $A\rightarrow M({\cal PMF}(Q))$ nor $A\rightarrow \partial_{2}C(Q)$.
\end{enumerate}

\item[(b)] If $Y''$ is a Borel subset of $Y'$ with positive measure, then the above Borel partition corresponding to $Y''$ is $Y''=(Y''\cap Y_{1})\sqcup (Y''\cap Y_{2})$.
\end{enumerate}
\end{lem}

\begin{pf}
By Lemma \ref{inv-lem} (ii), we can find $Y_{1}$ and $Y_{2}$ satisfying all the conditions in (a) except for the latter assertion in the condition (ii). 

In the condition (ii), if the map $\psi$ were not $M({\cal MIN}(Q))$-valued, then we would have a Borel subset $A\subseteq Y_{1}$ with positive measure such that $\psi(x)({\cal MIN}(Q))<1$ for any $x\in A$. Recall the $\Gamma(Q)$-equivariant Borel maps
\[M({\cal PMF}(Q)\setminus {\cal MIN}(Q))\rightarrow M(S(Q)),\]
\[M(S(Q))\rightarrow {\cal F}(S(Q))\]
in the proof of Lemma \ref{partition}, where ${\cal F}(S(Q))$ is the set of all non-empty finite subsets of $S(Q)$. By restricting the measure $\psi(x)$ for $x\in A$ to the set ${\cal PMF}(Q)\setminus {\cal MIN}(Q)$ and normalizing, we obtain a $\rho_{Q}$-invariant Borel map
\[A\rightarrow M({\cal PMF}(Q)\setminus {\cal MIN}(Q)).\]
Composing these maps, we get a $\rho_{Q}$-invariant Borel map $f\colon A\rightarrow {\cal F}(S(Q))$. Let $F \in {\cal F}(S(Q))$ be an element such that $f^{-1}(F)$ has positive measure. Since $({\cal S})_{f^{-1}(F)}$ fixes any curve in elements in $F$, this contradicts $Q\in \varphi_{i}$ and the property (iv) of $\varphi_{i}$ in Theorem \ref{alternative-trivial-irreducible}. The essential uniqueness of $Y_{1}$ and $Y_{2}$ is clear.

The assertion (b) is easy to prove.
\end{pf}

\begin{rem}
If $Q\in D$ and $\gamma \in \Gamma(M)$, then we have the natural continuous maps
\begin{align*}
\gamma \colon {\cal MF}(Q)&\rightarrow {\cal MF}(\gamma Q),\\
\gamma \colon C(Q)&\rightarrow C(\gamma Q)
\end{align*}
defined by the formulas
\begin{align*}
(\gamma(t))(\alpha)&=t(\gamma^{-1}\alpha),\\
\gamma(x)&=\gamma \cdot x
\end{align*}
for $\alpha \in V(C(\gamma Q))$, $t\in {\cal MF}(Q)\subseteq {\Bbb R}_{\geq 0}^{V(C(Q))}$ and $x\in V(C(Q))\subseteq V(C(M))$. The map $\gamma \colon {\cal MF}(Q)\rightarrow {\cal MF}(\gamma Q)$ induces the map $\gamma \colon {\cal PMF}(Q)\rightarrow {\cal PMF}(\gamma Q)$. These maps induce the continuous maps 
\begin{align*}
\gamma \colon M({\cal PMF}(Q))&\rightarrow M({\cal PMF}(\gamma Q)),\\
\gamma \colon \partial_{2}C(Q)&\rightarrow \partial_{2}C(\gamma Q).
\end{align*}
\end{rem}

\begin{rem}\label{rem-mendo}
This remark will be used in the proof of the next theorem.

With the assumption $(\ast)$, suppose that ${\cal S}$ is reducible. Let us denote by $\varphi$, $\varphi_{t}$ and $\varphi_{i}$ the CRS, T and I systems for ${\cal S}$. Suppose that $Y'$ is a Borel subset of $Y$ with positive measure such that all $\varphi$, $\varphi_{t}$ and $\varphi_{i}$ are constant on $Y'$ and denote their values on $Y'$ by $\sigma$, $F_{1}$ and $F_{2}$, respectively. Let $Q\in F_{2}$ and $Y'=Y_{1}^{Q}\sqcup Y_{2}^{Q}$ be the Borel partition for $Q$ in Lemma \ref{irreducible-partition} (a). 

Let $g\in [[{\cal S}]]$ and 
\[A_{\gamma}=\{ x\in {\rm dom}(g): \rho(gx, x)=\gamma \}\]  
for $\gamma \in \Gamma$. Assume that $A_{\gamma}\cap Y'$ has positive measure for a fixed $\gamma \in \Gamma$. Then 
\[\varphi(gx)=\gamma \sigma,\ \varphi_{t}(gx)=\gamma F_{1},\ \varphi_{i}(gx)=\gamma F_{2}\]
for a.e.\ $x\in A_{\gamma}$ since $\varphi$, $\varphi_{t}$ and $\varphi_{i}$ are $\rho$-invariant. In particular, $\varphi$, $\varphi_{t}$ and $\varphi_{i}$ are constant on $gA_{\gamma}$. Moreover, $\gamma Q\in \gamma F_{2}$.

We show that $g(A_{\gamma}\cap Y')=g(A_{\gamma}\cap Y_{1}^{Q})\sqcup g(A_{\gamma}\cap Y_{2}^{Q})$ is equal to the Borel partition $g(A_{\gamma}\cap Y')=Y_{1}^{\gamma Q}\sqcup Y_{2}^{\gamma Q}$ for $\gamma Q$ in Lemma \ref{irreducible-partition} (a). Let  
\[\psi \colon Y_{11}^{Q}\rightarrow M({\cal PMF}(Q)), \ \ \psi'\colon Y_{12}^{Q}\rightarrow \partial_{2}C(Q)\]
be $\rho$-invariant Borel maps in Lemma \ref{irreducible-partition} (a) such that $Y_{1}^{Q}=Y_{11}^{Q}\cup Y_{12}^{Q}$ (either $Y_{11}^{Q}$ or $Y_{12}^{Q}$ may be null). Define the maps
\begin{align*}
\psi_{g}&\colon g(A_{\gamma}\cap Y_{11}^{Q})\rightarrow M({\cal PMF}(\gamma Q)),\\
\psi_{g}'&\colon g(A_{\gamma}\cap Y_{12}^{Q})\rightarrow \partial_{2}C(\gamma Q)
\end{align*}
by the formulas
\begin{align*}
\psi_{g}(gx)&=\gamma (\psi(x)),\\
\psi_{g}'(gy)&=\gamma (\psi'(y))
\end{align*}
for $x\in A_{\gamma}\cap Y_{11}^{Q}$ and $y\in A_{\gamma}\cap Y_{12}^{Q}$. Then for $(gx, gy)\in {\cal S}$ with $gx, gy\in g(A_{\gamma}\cap Y_{11}^{Q})$, we have
\begin{align*}
\rho(gx, gy)\psi_{g}(gy)&=\rho(gx, x)\rho(x, y)\rho(y, gy)\gamma (\psi(y))\\
                        &=\gamma \rho(x, y)\psi(y)\\
                        &=\gamma(\psi(x))=\psi_{g}(gx).
\end{align*}
This means that $\psi_{g}$ is $\rho$-invariant. Similarly, $\psi_{g}'$ is also $\rho$-invariant. It follows that $g(A_{\gamma}\cap Y_{1}^{Q})\subseteq Y_{1}^{\gamma Q}$. If $Z_{1}=Y_{1}^{\gamma Q}\setminus g(A_{\gamma}\cap Y_{1}^{Q})$ had positive measure, then $g^{-1}(Z_{1})\subseteq A_{\gamma}\cap Y_{2}^{Q}$. By the definition of $Y_{1}^{\gamma Q}$, there exist $\rho$-invariant Borel maps
\[\overline{\psi}\colon Z_{11}\rightarrow M({\cal PMF}(\gamma Q)), \ \ \overline{\psi}'\colon Z_{12}\rightarrow \partial_{2}C(\gamma Q)\]
with $Z_{1}=Z_{11}\cup Z_{12}$. We can show as above that both of the maps 
\begin{align*}
g^{-1}(Z_{11})\ni x&\mapsto \gamma^{-1}(\overline{\psi}(gx))\in M({\cal PMF}(Q)),\\
g^{-1}(Z_{12})\ni x&\mapsto \gamma^{-1}(\overline{\psi}'(gx))\in \partial_{2}C(Q)
\end{align*}
are $\rho$-invariant. This contradicts $g^{-1}(Z_{1})\subseteq A_{\gamma}\cap Y_{2}^{Q}$. Therefore, $g(A_{\gamma}\cap Y_{1}^{Q})=Y_{1}^{\gamma Q}$ and thus, $g(A_{\gamma}\cap Y_{2}^{Q})=Y_{2}^{\gamma Q}$.

\end{rem}

\begin{thm}\label{alternative-irreducible-amenable-non-amenable}
With the assumption $(\ast)$, suppose that ${\cal S}$ is reducible. Let $\varphi$ be the CRS for ${\cal S}$. Then there exist essentially unique two $\rho$-invariant Borel maps
\begin{align*}
\varphi_{ia}&\colon Y\rightarrow {\cal F}_{0}(D),\\
\varphi_{in}&\colon Y\rightarrow {\cal F}_{0}(D)
\end{align*}
such that
\begin{enumerate}
\renewcommand{\labelenumi}{\rm(\roman{enumi})} 
\item $\varphi_{i}(y)=\varphi_{ia}(y)\cup \varphi_{in}(y)$ and $\varphi_{ia}(y)\cap \varphi_{in}(y)=\emptyset$ for a.e.\ $y\in Y$;
\item if $Q$ is in $F\in {\cal F}_{0}(D)$ with $\mu(\varphi_{ia}^{-1}(F))>0$, then there exist two $\rho_{Q}$-invariant Borel maps
\[\psi \colon A_{1}\rightarrow M({\cal PMF}(Q)), \ \ \psi'\colon A_{2}\rightarrow \partial_{2}C(Q)\]
such that $\varphi_{ia}^{-1}(F)=A_{1}\cup A_{2}$ (either $A_{1}$ or $A_{2}$ may be null). Moreover, any of such a map $\psi$ is $M({\cal MIN}(Q))$-valued;
\item if $Q$ is in $F\in {\cal F}_{0}(D)$ with $\mu(\varphi_{in}^{-1}(F))>0$, then there exist neither $\rho_{Q}$-invariant Borel maps $A\rightarrow M({\cal PMF}(Q))$ nor $A\rightarrow \partial_{2}C(Q)$ for any Borel subset $A\subseteq \varphi_{in}^{-1}(F)$ with positive measure.
\end{enumerate}
\end{thm}

\begin{pf}
We define $\varphi_{ia}$ on each Borel subset $Y'=\varphi^{-1}(\sigma)\cap \varphi_{t}^{-1}(F_{1})\cap \varphi_{i}^{-1}(F_{2})$ with positive measure for $\sigma \in S(M)$ and $F_{1}, F_{2}\in {\cal F}_{0}(D)$. If $F_{2}=\emptyset$, then define $\varphi_{ia}=\emptyset$ on $Y'$. If $Q\in F_{2}$, then let $Y'=Y_{1}^{Q}\sqcup Y_{2}^{Q}$ be the Borel partition in Lemma \ref{irreducible-partition} (a). Let us define 
\[\varphi_{ia}(y)=\{ Q\in F_{2}: y\in Y_{1}^{Q}\}\]
for $y\in Y'$ and define
\[\varphi_{in}(y)=F_{2}\setminus \varphi_{ia}(y)\]
for $y\in Y'$. Since $Y_{1}^{Q}$ is a Borel subset of $Y'$, the map $\varphi_{ia}$ is measurable on $Y'$ and so is the map $\varphi_{in}$. Thus, they are measurable on $Y$. Moreover, these maps satisfy all the conditions (i), (ii) and (iii).

Next, we show the $\rho$-invariance for the map $\varphi_{ia}\colon Y\rightarrow {\cal F}_{0}(D)$. Suppose that $\sigma \in S(M)$ and $F_{1}, F_{2}\in {\cal F}_{0}(D)$ are elements such that the Borel subset $Y'=\varphi^{-1}(\sigma)\cap \varphi_{t}^{-1}(F_{1})\cap \varphi_{i}^{-1}(F_{2})$ has positive measure. Let $g\in [[{\cal S}]]$ and
\[A_{\gamma}=\{ x\in {\rm dom}(g): \rho(gx, x)=\gamma \}\]
for $\gamma \in \Gamma$. Assume that $A_{\gamma}\cap Y'$ has positive measure for a fixed $\gamma \in \Gamma$. Then 
\[\varphi(gx)=\gamma \sigma,\ \varphi_{t}(gx)=\gamma F_{1},\ \varphi_{i}(gx)=\gamma F_{2}\]
for a.e.\ $x\in A_{\gamma}\cap Y'$ and 
\[Y''=\varphi^{-1}(\gamma \sigma)\cap \varphi^{-1}_{t}(\gamma F_{1})\cap \varphi_{i}^{-1}(\gamma F_{2})\]
contains $g(A_{\gamma}\cap Y')$.

If $Q\in F_{2}$, then we have the Borel partition $Y'=Y_{1}^{Q}\sqcup Y_{2}^{Q}$ in Lemma \ref{irreducible-partition} (a). Since $\gamma Q\in \gamma F_{2}$, we have the Borel partition $Y''=Z_{1}^{\gamma Q}\sqcup Z_{2}^{\gamma Q}$ for $\gamma  Q$ in Lemma \ref{irreducible-partition} (a). It follows from Lemma \ref{irreducible-partition} (b) and Remark \ref{rem-mendo} that $g(A_{\gamma}\cap Y_{j}^{Q})\subseteq Z_{j}^{\gamma Q}$ for $j=1, 2$. This means that $\gamma Q\in \varphi_{ia}(gx)$ for a.e.\ $x\in A_{\gamma}\cap Y_{1}^{Q}$ and that $\gamma \varphi_{ia}(x)\subseteq \varphi_{ia}(gx)$ for a.e.\ $x\in A_{\gamma}$. Considering $g^{-1}$, we can deduce the equality $\gamma \varphi_{ia}(x)=\varphi_{ia}(gx)$ for a.e.\ $x\in A_{\gamma}$. Thus, we have shown the $\rho$-invariance of the map $\varphi_{ia}$. It follows that $\varphi_{in}$ is also $\rho$-invariant.

The uniqueness of $\varphi_{ia}$ and $\varphi_{in}$ can be verified easily.
\end{pf}

\begin{defn}
We call $\varphi_{ia}$ (resp. $\varphi_{in}$) the {\it system of irreducible and amenable} (resp. {\it irreducible and non-amenable}) {\it subsurfaces}\index{system! of irreducible and amenable subsurfaces}\index{system!of irreducible and non-amenable subsurfaces} for ${\cal S}$. In short, we call it the {\it IA} (resp. {\it IN}) {\it system}\index{IA!system}\index{IN!system} for ${\cal S}$. We often call an element in $\varphi_{ia}(x)$ (resp. $\varphi_{in}(x)$) an {\it irreducible and amenable} (resp. {\it irreducible and non-amenable}) {\it subsurface}\index{irreducible and amenable!subsurface}\index{irreducible and non-amenable!subsurface} for $x\in Y$ and in short, an {\it IA} (resp. {\it IN}) {\it subsurface}\index{IA!subsurface}\index{IN!subsurface}.
\end{defn}

When we identify a subsurface with a component of the cut surface along some curves, we call IA and IN subsurfaces {\it IA} and {\it IN components}\index{IA!component}\index{IN!component}, respectively.

\begin{lem}\label{restriction-of-invariants}
With the assumption $(\ast)$, suppose that ${\cal S}$ is reducible. Let us denote by $\varphi_{ia}$ and $\varphi_{in}$ the IA and IN systems for ${\cal S}$, respectively. Let $A\subseteq Y$ be a Borel subset with positive measure. Then the IA and IN systems for the relation $({\cal S})_{A}$ are the restrictions to $A$ of $\varphi_{ia}$ and $\varphi_{in}$, respectively.
\end{lem}

\begin{pf}
Let $\varphi_{ia, A}$ and $\varphi_{in, A}$  be the IA and IN systems for $({\cal S})_{A}$, respectively, and let $E$ be a Borel subset of $A$ with positive measure such that all the CRS, T, IA and IN systems for ${\cal S}$ and $({\cal S})_{A}$ are constant on $E$. We denote the values of the IA and IN systems on $E$ for ${\cal S}$ and $({\cal S})_{A}$ by $\Phi_{ia}$, $\Phi_{in}$, $\Phi_{ia, A}$ and $\Phi_{in, A}$, respectively. It follows from Lemma \ref{restriction-of-invariants-first} and Theorem \ref{alternative-irreducible-amenable-non-amenable} (i) that 
\[\Phi_{ia}\cup \Phi_{in}=\Phi_{ia, A}\cup \Phi_{in, A}.\]
If $Q\in \Phi_{ia}$, then there exist $\rho$-invariant Borel maps
\[\psi \colon E_{1}\rightarrow M({\cal PMF}(Q)), \ \ \psi'\colon E_{2}\rightarrow \partial_{2}C(Q)\] 
for ${\cal S}$ with $E=E_{1}\cup E_{2}$. These Borel maps are $\rho$-invariant also for $({\cal S})_{A}$. Thus, $Q\in \Phi_{in, A}$ cannot happen and thus, $Q\in \Phi_{ia, A}$. Hence, $\Phi_{ia}\subseteq \Phi_{ia, A}$.  Similarly, we can prove $\Phi_{ia, A}\subseteq \Phi_{ia}$. This means that $\varphi_{ia}=\varphi_{ia, A}$ on $E$. Since $E$ is an arbitrary Borel subset of $A$ on which all the CRS, T, IA and IN systems for ${\cal S}$ and $({\cal S})_{A}$ are constant, we see that $\varphi_{ia}=\varphi_{ia, A}$ on $A$. It follows that $\varphi_{in}=\varphi_{in, A}$ on $A$.
\end{pf}


\section{Irreducible and amenable subsurfaces}\label{section-irreducible-amenable-subsurface}

In this section, we show that for an IA subsurface $Q$ for ${\cal S}$, the cocycle $\rho_{Q}$ satisfies similar properties to ones for irreducible and amenable subrelations. Propositions \ref{invariant-amenable'}, \ref{maximal'} and Lemma \ref{main-lem1'} (i) correspond to Propositions \ref{invariant-amenable}, \ref{maximal} and Lemma \ref{main-lem1}, respectively.

First, assume $(\ast)$ and that ${\cal S}$ is reducible and let $\varphi$, $\varphi_{t}$, $\varphi_{ia}$ and $\varphi_{in}$ be the CRS, T, IA and IN systems for ${\cal S}$, respectively. Suppose that they are constant on $Y$ and denote the values of $\varphi$, $\varphi_{t}$, $\varphi_{ia}$ and $\varphi_{in}$ by the same symbols.

Assume that $\varphi_{ia}$ is non-empty. Fix a component $Q$ in $\varphi_{ia}$. Remark that $Q$ is not a pair of pants. Then there exists a Borel cocycle
\[\rho_{Q}\colon {\cal S}\rightarrow \Gamma(Q).\]
By Theorem \ref{alternative-irreducible-amenable-non-amenable}, there exist two $\rho_{Q}$-invariant Borel maps
\[\psi \colon A_{1}\rightarrow M({\cal PMF}(Q)), \ \ \psi'\colon A_{2}\rightarrow \partial_{2}C(Q)\]
such that $Y=A_{1}\cup A_{2}$ (either $A_{1}$ or $A_{2}$ may be null). Since $\psi$ is $M({\cal MIN}(Q))$-valued, we obtain the Borel map
\[\pi_{*}\psi \colon A_{1}\rightarrow M(\partial C(Q)),\]
composing $\psi$ with the induced Borel map $\pi_{*}$ by $\pi \colon {\cal MIN}(Q)\rightarrow \partial C(Q)$ in Chapter \ref{chapter:amenable-action}, Section \ref{boundary-of-curve-complex} (see also Remark \ref{rem-borel-str-measurability}). 

\begin{prop}\label{invariant-amenable'}
In the above assumption, the cardinality of the support of the measure $\pi_{*}\psi (y)$ on $\partial C(Q)$ is at most two for a.e.\ $y\in A_{1}$. Thus, the map $\pi_{*}\psi$ induces a $\rho_{Q}$-invariant Borel map $A_{1}\rightarrow \partial_{2}C(Q)$. 
\end{prop}

\begin{pf}
Assume that it is not true. By normalization, we may assume that the cardinality of the support of $\pi_{*}\psi (y)$ would be more than two for $y$ in some Borel subset $A_{1}'$ of $A_{1}$ with positive measure. As in the proof of Proposition \ref{invariant-amenable}, we get a $\rho_{Q}$-invariant Borel map $\psi''\colon A_{1}'\rightarrow {\cal F}'(C(Q))$, where ${\cal F}'(C(Q))$ is the set of all non-empty finite subsets of $V(C(Q))$ whose diameters are more than or equal to 3. 

Let $F\in {\cal F}'(C(Q))$ be an element with $\mu({\psi''}^{-1}(F))>0$. Take $\alpha \in F$. It follows from the finiteness of $F$ and from the pureness of $\Gamma$ that $\alpha$ is invariant for $({\cal S})_{{\psi''}^{-1}(F)}$. Thus, $(\alpha, {\psi''}^{-1}(F))$ is $\rho$-invariant. This contradicts the assumption that $Q$ is an I subsurface for ${\cal S}$.
\end{pf}

Thus, combining $\psi'\colon A_{2}\rightarrow \partial_{2}C(Q)$, we obtain a $\rho_{Q}$-invariant Borel map $Y\rightarrow \partial_{2}C(Q)$ by Lemma \ref{inv-lem} (i).

\begin{prop}\label{maximal'}
We can construct a $\rho_{Q}$-invariant Borel map $\psi_{0}\colon Y\rightarrow \partial_{2}C(Q)$ such that for any $\rho_{Q}$-invariant Borel map $\psi \colon Y\rightarrow \partial_{2}C(Q)$, we have 
\[{\rm supp}(\psi(y))\subseteq {\rm supp}(\psi_{0}(y))\]
for a.e.\ $y\in Y$.
\end{prop}

The method to construct this $\psi_{0}$ is the same as in Proposition \ref{maximal}.

Next, we assume $(\ast)''$ and that ${\cal S}$ is normal in ${\cal R}$ and ${\cal S}$ is reducible. Let $\varphi$, $\varphi_{t}$, $\varphi_{ia}$ and $\varphi_{in}$ be the CRS, T, IA and IN systems for ${\cal S}$, respectively. Suppose that they are constant on $X$ and denote their values by the same symbols. Suppose that $\varphi_{ia}$ is non-empty and let $Q$ be a component of $\varphi_{ia}$. Let $\psi_{0}\colon X\rightarrow \partial_{2}C(Q)$ be the $\rho_{Q}$-invariant Borel map for ${\cal S}$ in Proposition \ref{maximal'}. 

Note that ${\cal R}$ is reducible and that $Q$ is invariant also for ${\cal R}$ since the Borel map $\varphi$ is invariant for ${\cal R}$ by Lemma \ref{main-lem2} and $\rho$ is valued in the group $\Gamma(M; m)$ (see Theorem \ref{comp-leave}). Thus, we have a Borel cocycle
\[\rho_{Q}\colon {\cal R}\rightarrow \Gamma(Q)\]
whose restriction to ${\cal S}$ is the same as the cocycle already considered.

\begin{lem}\label{main-lem1'}
With the above assumption,
\begin{enumerate}
\item[(i)] the Borel map $\psi_{0}\colon X\rightarrow \partial_{2}C(Q)$ is $\rho_{Q}$-invariant for ${\cal R}$. 
\item[(ii)] if we denote by $\psi \colon X\rightarrow S(M)$ the CRS for ${\cal R}$, then $\varphi \subseteq \psi(x)$ for a.e.\ $x\in X$. 
\item[(iii)] if we denote by $\psi_{ia}\colon X\rightarrow {\cal F}_{0}(D)$ the IA system for ${\cal R}$, then $Q\in \psi_{ia}(x)$ for a.e.\ $x\in X$. 
\end{enumerate}
\end{lem}

\begin{pf}
The proof of the assertion (i) is the same as in Lemma \ref{main-lem1}.

Since the CRS $\varphi$ for ${\cal S}$ is invariant for ${\cal R}$, if $\alpha$ is any curve in $\varphi$, then the pair $(\alpha, X)$ is $\rho$-invariant for ${\cal R}$. Thus, $(\alpha, X)$ is an essential $\rho$-invariant pair for ${\cal R}$ because ${\cal R}$ contains ${\cal S}$, which proves the assertion (ii). 

The assertion (iii) follows from (i), (ii) and the fact that for a.e.\ $x\in X$, any element $\alpha \in V(C(Q))$ is not an element in $\psi(x)$ because $Q$ is an I subsurface for ${\cal S}$. 
\end{pf}

Finally, we give a simplification of Theorems \ref{alternative-trivial-irreducible} and \ref{alternative-irreducible-amenable-non-amenable}.

\begin{cor}\label{cor-simple}
With the assumption $(\ast)$, suppose that ${\cal S}$ is reducible. Let $\varphi$, $\varphi_{t}$, $\varphi_{ia}$ and $\varphi_{in}$ be the CRS, T, IA and IN systems for ${\cal S}$. Let $\sigma \in S(M)$, $F_{1}, F_{2}, F_{3}\in {\cal F}_{0}(D)$ be elements such that the Borel subset
\[A=\varphi^{-1}(\sigma)\cap \varphi_{t}^{-1}(F_{1})\cap \varphi_{ia}^{-1}(F_{2})\cap \varphi_{in}^{-1}(F_{3})\]
has positive measure. Then for each component $Q$ of $M_{\sigma}$, we see that
\begin{enumerate}
\item[(i)] $Q\in F_{1}$ if and only if either $Q$ is a pair of pants or $Q$ is not a pair of pants and there exists a $\rho_{Q}$-invariant Borel map $A\rightarrow S(Q)$.
\item[(ii)] $Q\in F_{2}$ if and only if there exists a $\rho_{Q}$-invariant Borel map $A\rightarrow \partial_{2}C(Q)$ and there exist no $\rho_{Q}$-invariant Borel maps $B\rightarrow S(Q)$ for any Borel subset $B$ of $A$ with positive measure.
\item[(iii)] $Q\in F_{3}$ if and only if there exist neither $\rho_{Q}$-invariant Borel maps $B\rightarrow S(Q)$ nor $B\rightarrow \partial_{2}C(Q)$ for any Borel subset $B$ of $A$ with positive measure.
\end{enumerate}   
\end{cor} 

\begin{pf}
We show the assertion (i). If $Q\in F_{1}$ and $Q$ is not a pair of pants, then it follows from Theorem \ref{alternative-trivial-irreducible} (iii) that the pair $(\alpha, A)$ is $\rho$-invariant for any $\alpha \in V(C(Q))$. Let $\alpha \in V(C(Q))$. It follows that there exists a countable Borel partition $A=\bigsqcup A_{n}$ such that 
\[\rho(x, y)\alpha =\alpha\]
for a.e.\ $(x, y)\in ({\cal S})_{A_{n}}$. Thus, the Borel map
\[f_{n}\colon A_{n}\rightarrow S(Q)\]
defined by $f_{n}(x)=\{ \alpha \}$ for any $x\in A_{n}$ is $\rho$-invariant. It follows from Lemma \ref{inv-lem} (i) that there exists a $\rho$-invariant Borel map $A\rightarrow S(Q)$.

Conversely, if there exists a $\rho$-invariant Borel map $f\colon A\rightarrow S(Q)$, then let $\sigma \in S(Q)$ be an element with $\mu(f^{-1}(\sigma))>0$ and let $\alpha \in \sigma$. Since we have
\[\rho(x, y)\alpha =\alpha\]
for a.e.\ $(x, y)\in ({\cal S})_{f^{-1}(\sigma)}$, the pair $(\alpha, f^{-1}(\sigma))$ is $\rho$-invariant. If $Q$ were in $F_{2}\cup F_{3}$, it would contradict Theorem \ref{alternative-trivial-irreducible} (iv).   

The assertions (ii) and (iii) are easy to prove.
\end{pf}


\section{Amenable, reducible subrelations}\label{section-amenable-action-general}

Under the assumption $(\ast)$, if ${\cal S}$ is reducible and amenable, then the irreducible and non-amenable system $\varphi_{in}$ for ${\cal S}$ must be empty by the definition of amenable relations. In this section, we verify the converse under the assumption $(\ast)'$, similarly to Proposition \ref{irreducible-amenable}.

\subsection{Amenability for the action of a reducible subgroup}

In this subsection, we show an analogue of Theorem \ref{amenable-action-cc}.
   
Let $M$ be a compact orientable surface with $\kappa(M)\geq 0$. Let $\sigma \in S(M)$ and $\Gamma_{\sigma}$ be the $\Gamma(M;m)$-stabilizer of $\sigma$, where $m\geq 3$ is an integer. By the pureness of $\Gamma(M;m)$, any element in $\Gamma_{\sigma}$ fixes any curve in $\sigma$. Let $\{ Q_{i}\}_{i}$ be the set of all components of $M_{\sigma}$ except for pairs of pants. We have a natural homomorphism from $\Gamma_{\sigma}$ into $\prod_{i}\Gamma(Q_{i})$, which induces the action of $\Gamma_{\sigma}$ on $\prod_{i}C(Q_{i})$ and on $\prod_{i}\partial C(Q_{i})$. Let us denote $X_{\sigma}=\prod_{i}V(C(Q_{i}))$ and $\partial X_{\sigma}=\prod_{i}\partial C(Q_{i})$ for simplicity.

\begin{thm}
The action of $\Gamma_{\sigma}$ on $(\partial X_{\sigma}, \mu)$ is amenable for any quasi-invariant probability measure $\mu$ on $\partial X_{\sigma}$.
\end{thm}

\begin{pf}
In this proof, we identify $C(Q_{i})$ with the vertex set of the complex $C(Q_{i})$. In the proof of Theorem \ref{amenable-action-cc-non-excep} (see also Remark \ref{rem-excep-funct-const}), we have constructed the Borel functions
\[f_{n}^{(i)}\colon \partial C(Q_{i})\times C(Q_{i})\rightarrow [0, 1]\]
for $n\in {\Bbb N}$ such that
\begin{enumerate}
\renewcommand{\labelenumi}{\rm(\roman{enumi})}
\item for any $a\in \partial C(Q_{i})$, we have
\[\sum_{x\in C(Q_{i})}f_{n}^{(i)}(a, x)=1;\]
\item for any $g\in \Gamma(Q_{i})$, we have
\[\lim_{n\rightarrow \infty}\sup_{a\in \partial C(Q_{i})}\sum_{x\in C(Q_{i})}|f_{n}^{(i)}(a, x)-f_{n}^{(i)}(g^{-1}a, g^{-1}x)|=0.\]
\end{enumerate}
Using these functions, we define the Borel function $F_{n}\colon \partial X_{\sigma} \times X_{\sigma}\rightarrow [0, 1]$ by the formula
\[F_{n}((a_{i}), (x_{i}))=\prod_{i}f_{n}^{(i)}(a_{i}, x_{i})\]
for $(a_{i})\in \partial X_{\sigma}$ and $(x_{i})\in X_{\sigma}$. Then we see that
\[\sum_{(x_{i})\in X_{\sigma}}F_{n}((a_{i}), (x_{i}))=1\]
for any $(a_{i})\in \partial X_{\sigma}$ and that 
\[\lim_{n\rightarrow \infty}\sup_{a\in \partial X_{\sigma}}\sum_{x\in X_{\sigma}}|F_{n}(a, x)-F_{n}(g^{-1}a, g^{-1}x)|=0\]
for any $g\in \prod_{i}\Gamma(Q_{i})$. This implies
\[\lim_{n\rightarrow \infty}\sum_{x\in X_{\sigma}}\int_{\partial X_{\sigma}}h(a)\varphi(a, x)(F_{n}(g^{-1}a, g^{-1}x)-F_{n}(a, x))d\mu(a)=0\]
for any $h\in L^{1}(\partial X_{\sigma})$, $\varphi \in L^{\infty}(\partial X_{\sigma}\times X_{\sigma})$ and $g\in \prod_{i}\Gamma(Q_{i})$. It follows that the sequence $\{ F_{n}\}$ is an approximate weakly invariant mean for the $(\partial X_{\sigma}\rtimes \prod_{i}\Gamma(Q_{i}))$-space $\partial X_{\sigma}\times X_{\sigma}$ (see Definition \ref{app-inv-mean}).

Using the functions in the proof of Lemma \ref{lem-amenable-action-stabilizer}, we can prove that for any $x\in X_{\sigma}$, the action of the group 
\[\{ g\in \prod_{i}\Gamma(Q_{i}): gx=x\}\]
on $(\partial X_{\sigma}, \mu)$ is amenable.

By these facts, we can show that the action of $\prod_{i}\Gamma(Q_{i})$ on $(\partial X_{\sigma}, \mu)$ is amenable in the same way as in the proof of Theorem \ref{amenable-action-cc-non-excep}.

On the other hand, the kernel of the natural homomorphism from $\Gamma_{\sigma}$ into $\prod_{i}\Gamma(Q_{i})$ is contained in the free abelian subgroup generated by the Dehn twists around all curves in $\sigma$ (see Proposition \ref{lem-blm}). In particular, it is amenable. It follows from Proposition \ref{prop-ame-cri-other} that the action of $\Gamma_{\sigma}$ on $(\partial X_{\sigma}, \mu)$ is also amenable.
\end{pf}

\begin{cor}\label{amenable-action-cc2-reducible}
The action of $\Gamma_{\sigma}$ on $(\prod_{i}\partial_{2}C(Q_{i}), \mu)$ is amenable for any quasi-invariant probability measure $\mu$ on $\prod_{i}\partial_{2}C(Q_{i})$.
\end{cor}

\begin{pf}
The idea for the proof is the same as that for Lemma \ref{amenable-lemma}.

Consider the coordinate interchanging action of the symmetric group $G_{i}$ of two letters on $\partial C(Q_{i})\times \partial C(Q_{i})$. Define $G=\prod_{i}G_{i}$ which acts on $\prod_{i}(\partial C(Q_{i}))^{2}$. Let $F_{i}\subseteq (\partial C(Q_{i}))^{2}$ be a fundamental domain for the action of $G_{i}$ and denote $F=\prod_{i}F_{i}$. By way of the restriction of the projection 
\[q\colon \prod_{i}(\partial C(Q_{i}))^{2}\rightarrow \prod_{i}\partial_{2}C(Q_{i}),\]
we identify $F$ with $\prod_{i}\partial_{2}C(Q_{i})$ and identify the measure $\mu_{F}$ on $F$ with the measure $\mu$. Let us define the measure $\mu'$ on $\prod_{i}(\partial C(Q_{i}))^{2}$ by the formula
\[\mu'=2^{-K}\sum_{g\in G}g_{*}\mu_{F},\]
where $K$ is the number of components of $M_{\sigma}$ which are not pairs of pants. Then $q_{*}\mu'=\mu$. 

The rest of the proof is the same as in Lemma \ref{amenable-lemma}.
\end{pf}


\subsection{Amenability as a relation and non-existence of irreducible and non-amenable subsurfaces}

\begin{thm}\label{irr-non-ame-emp-ame}
With the assumption $(\ast)'$, suppose that ${\cal S}$ is reducible and let $\varphi$, $\varphi_{t}$, $\varphi_{ia}$ and $\varphi_{in}$ be the CRS, T, IA and IN systems for ${\cal S}$, respectively. If $\varphi_{in}$ is empty a.e.\ on $Y$, then ${\cal S}$ is amenable.
\end{thm}

\begin{pf}
We may suppose that $F=X$, that is, the relation ${\cal S}$ is a recurrent subrelation on $Y(\subseteq X)$ of the relation ${\cal R}$ generated by an essentially free, non-singular Borel action of $\Gamma$ on $(X, \mu)$. Replacing any Borel subset of $Y$ on which all $\varphi$, $\varphi_{t}$, $\varphi_{ia}$ and $\varphi_{in}$ are constant by $Y$, we may assume that all $\varphi$, $\varphi_{t}$, $\varphi_{ia}$ and $\varphi_{in}$ are constant on $Y$ (see Lemmas \ref{restriction-of-invariants-first} and \ref{restriction-of-invariants}). Let $\sigma \in S(M)$ be the value of $\varphi$. Denote the values of $\varphi_{t}$, $\varphi_{ia}$ and $\varphi_{in}$ by the same symbols. Remark that $\rho$ is $\Gamma_{\sigma}$-valued, where $\Gamma_{\sigma}$ is the $\Gamma$-stabilizer of $\sigma$ and that each subsurface $Q\in \varphi_{t}\cup \varphi_{ia}\cup \varphi_{in}$ is preserved by any element in $\Gamma_{\sigma}$. Thus, $\rho$ induces a Borel cocycle
\[\rho_{Q}\colon {\cal S}\rightarrow \Gamma_{\sigma}\rightarrow \Gamma(Q).\]

If $Q$ is in $\varphi_{t}$ and not a pair of pants, then it follows from Lemma \ref{lem:trivial} that there exists a countable Borel partition $Y=\bigsqcup Y_{n}$ such that for any $n\in {\Bbb N}$, the restriction of $\rho_{Q}$ to $({\cal S})_{Y_{n}}$ is a map into the subgroup
\[Z=\{ g\in \Gamma(Q): g\alpha =\alpha {\rm \ for \ any \ } \alpha \in V(C(Q))\}. \] 
Note that $g\in Z$ satisfies $ga=a$ for any $a\in \partial C(Q)$. 

It follows that there exists a countable Borel partition $Y=\bigsqcup Y_{n}$ and a $\rho_{Q}$-invariant Borel map $Y_{n}\rightarrow \partial_{2}C(Q)$ for each $n$ and any subsurface $Q\in \varphi_{t}$ which is not a pair of pants. 

Moreover, for any subsurface $R\in \varphi_{ia}$, there exists a $\rho_{R}$-invariant Borel map $Y\rightarrow \partial_{2}C(R)$. Thus, there exists a $\rho_{\sigma}$-invariant Borel map 
\[\psi_{n}\colon Y_{n}\rightarrow \prod_{i}\partial_{2}C(Q_{i}),\]
where $i$ runs through the index set of the set $\{ Q_{i}\}$ of all components of $M_{\sigma}$ which are not pairs of pants and $\rho_{\sigma}\colon {\cal S}\rightarrow \prod_{i}\Gamma(Q_{i})$ is the induced cocycle $\prod_{i}\rho_{Q_{i}}$. Therefore, we obtain a $\rho_{\sigma}$-invariant Borel map 
\[\psi \colon Y\rightarrow \prod_{i}\partial_{2}C(Q_{i})\]
by Lemma \ref{inv-lem} (i). 

We can show amenability of ${\cal S}$ by using $\psi$ and Corollary \ref{amenable-action-cc2-reducible} and by applying the idea for the proof of Proposition \ref{irreducible-amenable}. We give only a sketch of the proof.
We identify ${\cal S}$ with the relation on $X$ defined by the union
\[\{(x, x)\in (X\setminus Y)\times (X\setminus Y)\} \cup {\cal S}.\]

The space $X\times \Gamma_{\sigma}$ has the action of $\Gamma_{\sigma}$ which is defined by the formula
\[g\cdot (x, g_{1})=(x, g_{1}g^{-1})\]
for $x\in X$ and $g, g_{1}\in \Gamma_{\sigma}$. Then the Mackey range $S$ of the cocycle $\rho :{\cal S}\rightarrow \Gamma_{\sigma}$ can be identified with the quotient space of $X\times \Gamma_{\sigma}$ by the equivalence relation
\[(x, g)\sim (y, \rho(y, x)g)\]
for $(x, y)\in {\cal S}$, $g\in \Gamma_{\sigma}$. Note that this is a subrelation of the equivalence relation defined by the formula
\[(x, g)\sim (y, \rho(y, x)g)\]
for $(x, y)\in \rho^{-1}(\Gamma_{\sigma})$ and $g\in \Gamma_{\sigma}$, which has a fundamental domain $\{ (x, e)\in X\times \Gamma_{\sigma}: x\in X\}$. Moreover, the Mackey range $S$ has the $\Gamma_{\sigma}$-action induced from the $\Gamma_{\sigma}$-action on $X\times \Gamma_{\sigma}$ defined above.

We can show that ${\cal S}$ can be identified with a Borel subgroupoid of the groupoid $S\rtimes \Gamma_{\sigma}$ by way of the map ${\cal S}\ni (x, y)\mapsto ([x, e], \rho(x, y))\in S\rtimes \Gamma_{\sigma}$, where $[x, g]\in S$ denotes the projection of $(x, g)\in X\times \Gamma_{\sigma}$. 

Using the Borel map $\psi \colon Y\rightarrow \prod_{i}\partial_{2}C(Q_{i})$ constructed above, we define a $\rho_{\sigma}$-invariant Borel map
\[\psi'\colon X\rightarrow \prod_{i}\partial_{2}C(Q_{i})\]
for the extended ${\cal S}$ on $(X, \mu)$ by $\psi'=\psi$ on $Y$ and $\psi'=a_{0}$ on $X\setminus Y$, where $a_{0}\in \prod_{i}\partial_{2}C(Q_{i})$ is a fixed point. We can construct a $\Gamma_{\sigma}$-equivariant Borel map $\psi''\colon S\rightarrow \prod_{i}\partial_{2}C(Q_{i})$ by $[x, g]\mapsto g^{-1}\psi'(x)$. By Corollary \ref{amenable-action-cc2-reducible} and Theorem \ref{aeg-thm}, the Borel groupoid $S\rtimes \Gamma_{\sigma}$ is amenable, which implies that ${\cal S}$ is amenable.    
\end{pf}

The above argument shows also the following proposition:

\begin{prop}\label{prop-inv-partial-amenable}
With the assumption $(\ast)'$, suppose that ${\cal S}$ is reducible and that there exists $\sigma \in S(M)$ such that $\rho(x, y)$ fixes $\sigma$ for a.e.\ $(x, y)\in {\cal S}$. Let $\{ Q_{i}\}_{i}$ be the set of all components of $M_{\sigma}$ which are not pairs of pants and $\rho_{\sigma}\colon {\cal S}\rightarrow \prod_{i}\Gamma(Q_{i})$ be the induced cocycle $\prod_{i}\rho_{Q_{i}}$. Moreover, we assume that there exists a $\rho_{\sigma}$-invariant Borel map 
\[\psi \colon Y\rightarrow \prod_{i}\partial_{2}C(Q_{i}).\]
Then the relation ${\cal S}$ is amenable.
\end{prop}


\section{Classification}\label{section-classification}

In this section, we give some classification result for the mapping class groups in terms of measure equivalence.

\begin{assumption}\label{assumption-diamond'}
We call the following assumption $(\diamond)'$\index{$(\diamond p $@$(\diamond)'$}: let $(X, \mu)$ be a standard Borel space with a finite positive measure. Let $\Gamma$ be a subgroup of $\Gamma(M; m)$, where $M$ is a surface with $\kappa(M)\geq 0$ and $m\geq 3$ is an integer. Suppose that ${\cal R}$ is a discrete measured equivalence relation on $(X, \mu)$ and that we have a Borel cocycle
\[\rho \colon {\cal R}\rightarrow \Gamma\]
with finite kernel. Let $F$ be a fundamental domain of the finite equivalence relation $\ker \rho$. We assume that there exists an essentially free, non-singular Borel action of $\Gamma$ on $(F, \mu|_{F})$ generating the relation $({\cal R})_{F}$ and whose induced cocycle 
\[({\cal R})_{F}\rightarrow \Gamma, \ \ (gx, x)\mapsto g, \ g\in \Gamma, \ x\in F\]
is equal to the restriction of $\rho$ to $({\cal R})_{F}$.
\end{assumption}

This assumption is almost the same as the assumption $(\ast)'$ in which we have a recurrent subrelation ${\cal S}$ of $({\cal R})_{Y}$ on a Borel subset $Y$ of $X$ with positive measure.

\begin{defn}\label{defn-dagger}
Let ${\cal R}$ be a discrete measured equivalence relation on a standard Borel space $(X, \mu)$ with a finite positive measure. We say that a family 
\[\{ Y, \{ {\cal S}_{1}^{1}, {\cal S}_{1}^{2}\}, \{ {\cal S}_{2}^{1}, {\cal S}_{2}^{2}\} \}\]
satisfies $(\dagger)$\index{$(\dagger $@$(\dagger)$} if $Y$ is a Borel subset of $X$ with positive measure and ${\cal S}_{i}^{j}$ $(i, j=1, 2)$ are amenable recurrent subrelations of the relation $({\cal R})_{Y}$ on $Y$ satisfying the following conditions: let us denote ${\cal R}_{i, A}=({\cal S}^{1}_{i})_{A}\vee({\cal S}^{2}_{i})_{A}$ for a Borel subset $A\subseteq Y$ with positive measure and $i=1, 2$. Then 
\begin{enumerate}
\renewcommand{\labelenumi}{\rm(\roman{enumi})}
\item the relation ${\cal R}_{i, A}$ is non-amenable for any Borel subset $A\subseteq Y$ with positive measure and $i=1, 2$;
\item for any Borel subsets $A\subseteq B$ of $Y$ with positive measure, the subrelation $({\cal S}_{1}^{j})_{A}$ is normal in $({\cal S}_{1}^{j})_{A}\vee ({\cal R}_{2, B})_{A}$ and the subrelation $({\cal S}_{2}^{j})_{A}$ is normal in $({\cal S}_{2}^{j})_{A}\vee ({\cal R}_{1, B})_{A}$ for $j=1, 2$;
\item the relation ${\cal R}_{i, A}$ is normal in ${\cal R}_{1, A}\vee {\cal R}_{2, A}$ for any Borel subset $A$ of $Y$ with positive measure and $i=1, 2$.
\end{enumerate}
\end{defn}

\begin{rem}
Note that if a family $\{ Y, \{ {\cal S}_{1}^{1}, {\cal S}_{1}^{2}\}, \{ {\cal S}_{2}^{1}, {\cal S}_{2}^{2}\} \}$ satisfies $(\dagger)$, then the family $\{ A, \{ ({\cal S}_{1}^{1})_{A}, ({\cal S}_{1}^{2})_{A}\}, \{ ({\cal S}_{2}^{1})_{A}, ({\cal S}_{2}^{2})_{A}\} \}$ also satisfies $(\dagger)$ for any Borel subset $A$ of $Y$ with positive measure. 

Suppose that ${\cal R}$ (resp. ${\cal R}'$) is a discrete measured equivalence relation on $(X, \mu)$ (resp. $(X', \mu')$) and there exists a Borel isomorphism $f\colon X\rightarrow X'$ inducing the isomorphism between ${\cal R}$ and ${\cal R}'$. If a family $\{ Y, \{ {\cal S}_{1}^{1}, {\cal S}_{1}^{2}\}, \{ {\cal S}_{2}^{1}, {\cal S}_{2}^{2}\} \}$ of a Borel subset $Y$ of $X$ and subrelations of $({\cal R})_{Y}$ satisfies $(\dagger)$, then the family $\{ f(Y), \{ f({\cal S}_{1}^{1}), f({\cal S}_{1}^{2})\}, \{ f({\cal S}_{2}^{1}), f({\cal S}_{2}^{2})\} \}$ also satisfies $(\dagger)$.  
\end{rem}

\begin{rem}\label{rem-usankusai}
With the assumption $(\diamond)'$, think of amenable recurrent subrelations ${\cal S}^{1}$, ${\cal S}^{2}$ of the relation $({\cal R})_{Y}$ on $(Y, \mu)$ satisfying the following conditions: let us denote ${\cal R}_{A}=({\cal S}^{1})_{A}\vee ({\cal S}^{2})_{A}$ for a Borel subset $A$ of $Y$ with positive measure. Then
\begin{enumerate}
\item[(i)] the relation ${\cal R}_{A}$ is non-amenable for any Borel subset $A$ of $Y$ with positive measure;
\item[(ii)] the relation $({\cal S}^{j})_{A}$ is normal in ${\cal R}_{A}$ for any Borel subset $A$ of $Y$ with positive measure and $j=1, 2$.
\end{enumerate}
If we put ${\cal S}_{i}^{j}={\cal S}^{j}$ for $i, j=1, 2$, then the family $\{ Y, \{ {\cal S}_{1}^{1}, {\cal S}_{1}^{2}\}, \{ {\cal S}_{2}^{1}, {\cal S}_{2}^{2}\} \}$ would satisfy $(\dagger)$. However, by the following argument, we see that there exist no such subrelations ${\cal S}^{1}$, ${\cal S}^{2}$ (see Remark \ref{rem-usankusai-right}). 
\end{rem}

The next lemma shows a fundamental fact about a family satisfying $(\dagger)$ under the assumption $(\diamond)'$.

\begin{lem}\label{lem-lem-lem}
Assume $(\diamond)'$ and that we have a family $\{ Y, \{ {\cal S}_{1}^{1}, {\cal S}_{1}^{2}\}, \{ {\cal S}_{2}^{1}, {\cal S}_{2}^{2}\} \}$ of a Borel subset $Y$ of $X$ and subrelations of $({\cal R})_{Y}$ satisfying $(\dagger)$. Then the relation ${\cal R}_{1, A}\vee {\cal R}_{2, A}$ is reducible for any Borel subset $A$ of $Y$ with positive measure. Moreover, if $\varphi^{(1)}$ and $\varphi^{(2)}$ are the CRS's for ${\cal R}_{1, A}$ and ${\cal R}_{2, A}$, respectively, then 
\begin{enumerate}
\item[(i)] both the Borel maps $\varphi^{(1)}$, $\varphi^{(2)}\colon Y\rightarrow S(M)$ are $\rho$-invariant for ${\cal R}_{1, A}\vee {\cal R}_{2, A}$.
\item[(ii)] the union $\varphi^{(1)}(x)\cup \varphi^{(2)}(x)$ is an element in $S(M)$ for a.e.\ $x\in A$.
\end{enumerate}
\end{lem}

\begin{pf}
Since $({\cal S}_{1}^{1})_{A}$ is normal in $({\cal S}_{1}^{1})_{A}\vee {\cal R}_{2, A}$ and $({\cal S}_{1}^{1})_{A}$ is amenable and recurrent, it follows from Theorem \ref{main-normal} that $({\cal S}_{1}^{1})_{A}\vee {\cal R}_{2, A}$ is a disjoint sum of an irreducible and amenable relation and a reducible relation. Since $(({\cal S}_{1}^{1})_{A}\vee {\cal R}_{2, A})_{A'}$ is non-amenable for any Borel subset $A'$ of $A$ with positive measure, the relation $({\cal S}_{1}^{1})_{A}\vee {\cal R}_{2, A}$ is reducible. In particular, ${\cal R}_{2, A}$ is reducible. Since the CRS $\varphi^{(2)}$ for ${\cal R}_{2, A}$ is a $\rho$-invariant Borel map for ${\cal R}_{1, A}\vee {\cal R}_{2, A}$ by Lemma \ref{main-lem2}, we see that ${\cal R}_{1, A}\vee {\cal R}_{2, A}$ is reducible.

The assertion (i) for $\varphi^{(2)}$ has already shown. The proof for $\varphi^{(1)}$ is the same. If the assertion (ii) were not true, then we would have a Borel subset $B$ of $A$ with positive measure, $\sigma_{1}, \sigma_{2}\in S(M)$ and $\alpha_{1}\in \sigma_{1}$, $\alpha_{2}\in \sigma_{2}$ such that $\varphi^{(1)}(x)=\sigma_{1}$, $\varphi^{(2)}(x)=\sigma_{2}$ for any $x\in B$ and $i(\alpha_{1}, \alpha_{2})\neq 0$. By definition, $(\alpha_{1}, B)$ is an essential $\rho$-invariant pair for ${\cal R}_{1, A}$. On the other hand, it follows from the assertion (i) that the pair $(\alpha_{2}, B)$ is $\rho$-invariant for ${\cal R}_{1, A}$. This contradicts the essentiality of the pair $(\alpha_{1}, B)$.
\end{pf}

The following lemma will be used in the proof of Theorem \ref{main-main-thm}.

\begin{lem}\label{in-contained-in}
With the assumption $(\ast)''$, suppose that ${\cal R}$ is reducible. Let $\varphi$, $\varphi_{in}$ and $\psi$, $\psi_{in}$ be the CRS's, IN systems for ${\cal S}$ and ${\cal R}$, respectively. Then for a.e.\ $x\in X$, the union $\varphi(x)\cup \psi(x)$ is an element in $S(M)$ and each subsurface in $\varphi_{in}(x)$ is contained in some subsurface in $\psi_{in}(x)$.
\end{lem}

It is easy to see that the union $\varphi(x)\cup \psi(x)$ is an element in $S(M)$ for a.e.\  $x\in X$. Otherwise, it would contradict the essentiality of elements in $\varphi(x)$ for ${\cal S}$. The last conclusion in the lemma means that when we realize the union $\varphi(x)\cup \psi(x)\in S(M)$ on $M$ and identify the isotopy class of elements in $\varphi(x)\cup \psi(x)$ with its realization, each IN component of $M\setminus \varphi(x)$ for ${\cal S}$ is contained in some IN component of $M\setminus \psi(x)$ for ${\cal R}$.

\begin{pf}
By Lemmas \ref{restriction-of-invariants-first} and \ref{restriction-of-invariants}, we may assume that all the CRS's, T, IA and IN systems for ${\cal R}$ and ${\cal S}$ are constant. Let $\sigma, \tau \in S(M)$ be the values of the CRS's for ${\cal R}$, ${\cal S}$, respectively. Realize $\sigma \cup \tau \in S(M)$ disjointly on $M$ and fix it. We identify the curves in the realization with the isotopy classes which contain them. In what follows, we often identify a subsurface with the corresponding open component in the complement of some of these fixed curves on $M$.

Let $Q$ be an IN subsurface for ${\cal S}$. Note that $Q$ does not contain any curve in $\sigma$. Otherwise, it would contradict Corollary \ref{cor-simple} since any curve in $\sigma$ is invariant for ${\cal R}$. It follows that there exists a unique subsurface $R$ of $M\setminus \sigma$ containing $Q$. It is clear that $R$ is not in the T system for ${\cal R}$. 

Assume that $R$ is an IA subsurface for ${\cal R}$. If $Q=R$, then $Q$ would not be an IN subsurface for ${\cal S}$, which is a contradiction. Thus, we have $Q\subsetneq R$ and there exists a component $\alpha$ of $\partial_{R}Q$. The curve $\alpha$ is invariant for ${\cal S}$ because $\alpha$ is in the CRS for ${\cal S}$. Since $R$ is an IA subsurface for ${\cal R}$, there exists a $\rho$-invariant Borel map $Y\rightarrow \partial_{2}C(R)$ for ${\cal R}$. By restricting the domain of the $\Gamma(R)$-equivariant Borel map
\[G_{2}\colon V(C(R))\times \partial_{2}C(R)\rightarrow {\cal F}'(C(R))\]
in the proof of Lemma \ref{reducible} to the set $\{ \alpha \}\times \partial_{2}C(R)$, we can construct a Borel map
\[\partial_{2}C(R)\rightarrow {\cal F}'(C(R))\]
which is equivariant for the action of the $\Gamma$-stabilizer of $\alpha$ and $R$. Recall that ${\cal F}'(C(R))$ is the set of all non-empty finite subsets of $V(C(R))$ whose diameters are more than or equal to 3. Composing the above two maps, we obtain a $\rho$-invariant Borel map 
\[f\colon Y\rightarrow {\cal F}'(C(R))\]
for ${\cal S}$. Let $F\in {\cal F}'(C(R))$ be a subset with $\mu(f^{-1}(F))>0$. If $(x, y)\in ({\cal S})_{f^{-1}(F)}$, then $\rho(x, y)$ fixes $F$ (see Remark \ref{rem-curve-comp-leave}) and thus, fixes any curve in $F$. Since the diameter of $F$ is more than or equal to 3, the quotient of $\rho(x, y)$ in $\Gamma(R)$ has finite order (see Lemma \ref{fill-curves-fix}). Thus, for any curve $\beta$ in $V(C(R))$, there exists $n\in {\Bbb N}$ such that $\rho(x, y)^{n}\beta =\beta$, which implies $\rho(x, y)\beta =\beta$ by the pureness of $\Gamma$ (see Theorem \ref{pure-important} (ii)). Therefore, for any curve $\beta$ in $V(C(R))$, the pair $(\beta, f^{-1}(F))$ is $\rho$-invariant for ${\cal S}$. It contradicts the assumption that $Q$ is an IN subsurface for ${\cal S}$. 

Hence, $R$ is an IN subsurface for ${\cal R}$. 
\end{pf}

The next theorem is the most important observation for the classification in this chapter:

\begin{thm}\label{main-main-thm}
Assume $(\diamond)'$ and that we have a family $\{ Y, \{ {\cal S}_{1}^{1}, {\cal S}_{1}^{2}\}, \{ {\cal S}_{2}^{1}, {\cal S}_{2}^{2}\} \}$ of a Borel subset $Y$ of $X$ and subrelations of $({\cal R})_{Y}$ satisfying $(\dagger)$. Then there exists a Borel subset 
$A\subseteq Y$ with positive measure satisfying the following conditions: let 
$\varphi^{(i)}$ be the CRS for ${\cal R}_{i, A}$ for $i=1, 2$ and realize the union $\varphi^{(1)}(x)\cup \varphi^{(2)}(x)\in S(M)$ on $M$ (see Lemma \ref{lem-lem-lem}).
\begin{enumerate}
\item[(i)] For $(i, j)=(1, 2), (2, 1)$, each IN component of $M\setminus \varphi^{(i)}(x)$ does not contain any curve in $\varphi^{(1)}(x)\cup \varphi^{(2)}(x)$ for a.e.\ $x\in A$. Thus, each IN component of $M\setminus \varphi^{(i)}(x)$ is contained in some component of $M\setminus \varphi^{(j)}(x)$. 
\item[(ii)] Any IN component of $M\setminus \varphi^{(1)}(x)$ is disjoint from any IN component of $M\setminus \varphi^{(2)}(x)$ for a.e. $x\in A$.
\end{enumerate}
\end{thm}

Note that the condition (i) is automatic since for any Borel subset $A$ of $Y$ with positive measure and $(i, j)=(1, 2), (2, 1)$, the CRS for ${\cal R}_{j, A}$ is $\rho$-invariant for ${\cal R}_{i, A}$ by Lemma \ref{main-lem2}. Hence, we need to find a Borel subset $A$ with positive measure satisfying the condition (ii).

For the proof, we prepare some results of a family $\{ Y, \{ {\cal S}_{1}^{1}, {\cal S}_{1}^{2}\}, \{ {\cal S}_{2}^{1}, {\cal S}_{2}^{2}\} \}$ of a Borel subset $Y$ of $X$ and subrelations of $({\cal R})_{Y}$ satisfying $(\dagger)$ under the assumption $(\diamond)'$. As in the proof of Lemma \ref{in-contained-in}, we often identify a subsurface of $M$ with an open component in the complement of some fixed curves in $M$. 

We can find a Borel subset $Y'$ of $Y$ with positive measure such that all the CRS, T, IA and IN systems for $({\cal R}_{i, Y})_{Y'}$ are constant for $i=1, 2$. For simplicity, let us denote ${\cal R}_{i}=({\cal R}_{i, Y})_{Y'}$ and denote the CRS and IN system for ${\cal R}_{i}$ by $\psi^{(i)}$ and $\psi_{in}^{(i)}$, respectively. Since they are constant, we write their values by the same symbols. Remark that $\psi^{(1)}\cup \psi^{(2)}\in S(M)$ by Lemma \ref{lem-lem-lem} (ii).

Realize all curves in $\psi^{(1)}\cup \psi^{(2)}\in S(M)$ disjointly on $M$ and fix them. Then we can realize each subsurface in $\psi_{in}^{(i)}$ as an open component of the set $M\setminus \psi^{(i)}$. 

Let $Q_{1}\in \psi_{in}^{(1)}$. Since any curve $\alpha \in \psi^{(2)}$ is in the CRS for $({\cal R}_{1, Y}\vee {\cal R}_{2, Y})_{Y'}$ by Lemma \ref{main-lem1'} (ii) and the pair $(\alpha, Y')$ is $\rho$-invariant for $({\cal R}_{1, Y})_{Y'}={\cal R}_{1}$, we see that there exist no curves $\alpha$ in $\psi^{(2)}$ with $\alpha \in V(C(Q_{1}))$. 

Therefore, there exists a unique component $Q_{2}$ of $M\setminus \psi^{(2)}$ containing $Q_{1}$. If $Q_{2}$ is an IN subsurface for ${\cal R}_{2}$, then we can prove $Q_{1}=Q_{2}$ by repeating the above argument for $Q_{2}$. Thus, we have shown the following lemma:

\begin{lem}\label{boundary-equal-lem}
If $Q_{i}\in \psi_{in}^{(i)}$ is realized as an open component of $M\setminus \psi^{(i)}$ for $i=1, 2$, then either $Q_{1}\cap Q_{2}=\emptyset$ or $Q_{1}=Q_{2}$.
\end{lem}

Using this lemma, we give a proof of Theorem \ref{main-main-thm}.

\begin{pf*}{{\sc Proof of Theorem \ref{main-main-thm}.}}
We use the same notation as above. Moreover, let $B$ be a Borel subset of $Y'$ with positive measure such that all the CRS, T, IA and IN systems for $({\cal S}_{i}^{j})_{B}$ for $i, j=1, 2$ are constant.

We assume that the intersection
\[\left( \bigcup \psi_{in}^{(1)}\right)\cap \left(\bigcup \psi_{in}^{(2)}\right)\subseteq D(M)\]
is non-empty and let $Q$ be a subsurface in the intersection. (If the intersection is empty, then let $\pi^{(1)}$ and $\pi^{(2)}$ be the IN systems for ${\cal R}_{1, Y'}$ and ${\cal R}_{2, Y'}$, respectively. By applying Lemma \ref{in-contained-in} for the relation ${\cal R}_{i}$ and its subrelation ${\cal R}_{i, Y'}$, we see that for a.e.\ $x\in Y'$, any subsurface in $\pi^{(i)}(x)$ is contained in some subsurface in $\psi_{in}^{(i)}(x)$. Since the intersection 
\[\left( \bigcup \psi_{in}^{(1)}\right)\cap \left(\bigcup \psi_{in}^{(2)}\right)\]
is empty, the IN systems $\pi^{(1)}$, $\pi^{(2)}$ satisfy the condition (ii) by Lemma \ref{boundary-equal-lem}, which means the completion of the proof of Theorem \ref{main-main-thm}.)

Note that for any curve $\alpha$ in the CRS for $({\cal S}_{i}^{j})_{B}$, the pair $(\alpha, B)$ is $\rho$-invariant for $({\cal R}_{k})_{B}$ for each $j=1, 2$, $(i, k)=(1, 2), (2, 1)$ by Lemma \ref{main-lem2} because $({\cal S}_{i}^{j})_{B}$ is normal in $({\cal S}_{i}^{j})_{B}\vee ({\cal R}_{k})_{B}$. In particular, any curve in the CRS for $({\cal S}_{i}^{j})_{B}$ does not intersect any curve in $\psi^{(k)}$ for each $j=1, 2$, $(i, k)=(1, 2), (2, 1)$.

By the above property of the essential reduction classes for $({\cal S}_{1}^{1})_{B}$, we have the following three possible cases (i), (iia) and (iib): 
\begin{enumerate}
\item[(i)] there exists an essential reduction class $\alpha_{1}$ for $({\cal S}_{1}^{1})_{B}$ represented by a curve in $V(C(Q))$;
\item[(ii)] there exist no such curves. 
\end{enumerate}
In the case (ii), cut $M$ along a realization of the CRS for $({\cal S}_{1}^{1})_{B}$. Since there exist no essential reduction classes for $({\cal S}_{1}^{1})_{B}$ in $V(C(Q))$, we have a component $R_{1}$ of the cut surface containing $Q$. Since the relation $({\cal S}_{1}^{1})_{B}$ is amenable, we have the following two possibilities:  
\begin{enumerate}
\item[(iia)] $R_{1}$ is a T subsurface for $({\cal S}_{1}^{1})_{B}$;
\item[(iib)] $R_{1}$ is an IA subsurface for $({\cal S}_{1}^{1})_{B}$.
\end{enumerate}
If $Q\subsetneq R_{1}$, then any component of $\partial_{R_{1}}Q$ is invariant for $({\cal S}_{1}^{1})_{B}$. Thus, $R_{1}$ is a T subsurface for $({\cal S}_{1}^{1})_{B}$. It is the case (iia). Thus, in the case (iib), we have the equality $R_{1}=Q$.

Similarly, we can consider the three cases (i)', (iia)' and (iib)' for $({\cal S}_{1}^{2})_{B}$. 

In the case (i), by Lemma \ref{main-lem2} and the assumption that $({\cal S}_{1}^{1})_{B}$ is normal in $({\cal S}_{1}^{1})_{B}\vee ({\cal R}_{2})_{B}$, we see that $\alpha_{1}$ is invariant for $({\cal R}_{2})_{B}$. Thus, $Q$ is a T subsurface for $({\cal R}_{2})_{B}$, which is a contradiction.

In the case (iib), there exists a $\rho$-invariant Borel map $B\rightarrow \partial_{2}C(Q)$ for $({\cal S}_{1}^{1})_{B}$. By Lemma \ref{main-lem1'} (i) and the assumption that $({\cal S}_{1}^{1})_{B}$ is normal in $({\cal S}_{1}^{1})_{B}\vee ({\cal R}_{2})_{B}$, there exists a $\rho$-invariant Borel map $B\rightarrow \partial_{2}C(Q)$ for $({\cal S}_{1}^{1})_{B}\vee ({\cal R}_{2})_{B}$. Hence, $Q$ is not an IN subsurface for $({\cal R}_{2})_{B}$, which is a contradiction.

Similarly, we can deduce contradictions from both the cases (i)' and (iib)'. Thus, $Q$ must be contained in T subsurfaces both for $({\cal S}_{1}^{1})_{B}$ and $({\cal S}_{1}^{2})_{B}$. 

For each component $Q$ of the intersection $(\bigcup \psi_{in}^{(1)}) \cap (\bigcup \psi_{in}^{(2)})$, let $\{ \alpha_{Q}, \beta_{Q}\}$ be a pair of curves in $V(C(Q))$ filling $Q$ (see the comment right before Lemma \ref{fill-curves-fix}). Since $Q$ is contained in T subsurfaces both for $({\cal S}_{1}^{1})_{B}$ and $({\cal S}_{1}^{2})_{B}$, there exists a Borel subset $A\subseteq B$ with positive measure such that for any component $Q$ of the intersection $(\bigcup \psi_{in}^{(1)})\cap (\bigcup \psi_{in}^{(2)})$, both $\alpha_{Q}$ and $\beta_{Q}$ are invariant both for $({\cal S}_{1}^{1})_{A}$ and $({\cal S}_{1}^{2})_{A}$. Thus, $\alpha_{Q}$ and $\beta_{Q}$ are invariant for ${\cal R}_{1, A}=({\cal S}_{1}^{1})_{A}\vee ({\cal S}_{1}^{2})_{A}$.

We show that this $A$ is a desired subset. If $A$ were not a desired subset, then by Lemma \ref{boundary-equal-lem}, there would exist a Borel subset $A'$ of $A$ with positive measure satisfying the following conditions:
\begin{enumerate}
\item[(a)] all the CRS's, T, IA and IN systems for ${\cal R}_{i, A}$ are constant on $A'$ for $i=1, 2$;
\item[(b)] if $\varphi_{in}^{(i)}$ denotes the value of the IN system for ${\cal R}_{i, A}$ on $A'$ for $i=1, 2$, then the intersection 
\[\left(\bigcup \varphi_{in}^{(1)}\right)\cap \left(\bigcup \varphi_{in}^{(2)}\right)\subseteq D(M)\]
is non-empty.
\end{enumerate}

Let $Q'$ be a component of the intersection $(\bigcup \varphi_{in}^{(1)})\cap (\bigcup \varphi_{in}^{(2)})$. Since $({\cal R}_{i, A})_{A'}$ is a subrelation of $({\cal R}_{i})_{A'}$, it follows from Lemma \ref{in-contained-in} that there exists $Q_{i}\in \psi_{in}^{(i)}$ such that $Q'\subseteq Q_{i}$ for each $i=1, 2$. Since $Q_{1}$ and $Q_{2}$ intersect, we see that $R=Q_{1}=Q_{2}$ by Lemma \ref{boundary-equal-lem}. 

On the other hand, by the choice of $A$, there exists a filling pair $\{ \alpha_{R}, \beta_{R}\}$ of curves in $R$ such that both curves in the pair are invariant for ${\cal R}_{1, A}$. Since $Q'$ is not a pair of pants, either $r(\alpha_{R}, Q')$ or $r(\beta_{R}, Q')$ is non-empty by Lemma \ref{pants}. It follows from Lemma \ref{lem:trivial} that for any curve $\gamma$ in $V(C(Q'))$, the pair $(\gamma, A')$ is $\rho$-invariant for $({\cal R}_{1, A})_{A'}$, which contradicts the condition that $Q'$ is an IN subsurface for $({\cal R}_{1, A})_{A'}$.        
\end{pf*}

\begin{cor}\label{cor-main-main-thm}
If $A$ is the Borel subset in Theorem \ref{main-main-thm}, then the theorem holds also for any Borel subset of $A$ with positive measure.
\end{cor}

\begin{pf}
For any Borel subset $A'\subseteq A$ with positive measure, since ${\cal R}_{i, A'}$ is a subrelation of $({\cal R}_{i, A})_{A'}$, any IN subsurface for ${\cal R}_{i, A'}$ is contained in some IN subsurface for $({\cal R}_{i, A})_{A'}$ by Lemma \ref{in-contained-in}. It follows from Lemma \ref{restriction-of-invariants} and Theorem \ref{main-main-thm} that IN subsurfaces for $({\cal R}_{1, A})_{A'}$ and $({\cal R}_{2, A})_{A'}$ are disjoint. Thus, so do IN subsurfaces for ${\cal R}_{1, A'}$ and ${\cal R}_{2, A'}$.
\end{pf}

For extracting an invariant for equivalence relations generated by the mapping class groups, we need the following geometric lemma.

Let $M=M_{g, p}$ be a compact orientable surface of type $(g, p)$. Assume $\kappa(M)=3g+p-4\geq 0$. Let $\sigma$ be an element in $S(M)\cup \{ \emptyset \}$ and realize all the curves in $\sigma$ disjointly on $M$. We denote by $M_{\sigma}$ the resulting surface obtained by cutting $M$ along each curve in $\sigma$. If $\sigma$ consists of a single element $\alpha \in V(C(M))$, then we denote by $M_{\alpha}$\index{$M \ a $@$M_{\alpha}$} instead of $M_{\sigma}$ for simplicity. Let $n(\sigma)$\index{$n \ r$@$n(\sigma)$} be the number of components of $M_{\sigma}$ which are not pairs of pants. For a surface $M$ of type $(g, p)$, we write 
\[n(M)=g+\left[ \frac{g+p-2}{2}\right], \]\index{$n M$@$n(M)$}
where for $a\in {\Bbb R}$, we denote by $[a]$\index{$a for a R$@$[a]$ for $a\in {\Bbb R}$} the maximal integer less than or equal to $a$.

Recall that a curve $\alpha \in V(C(M))$ is said to be {\it separating}\index{separating curve} if the surface $M$ is cut along $\alpha$, then the resulting surface has two components. Otherwise, the curve $\alpha$ is said to be {\it non-separating}\index{non-separating curve}.

\begin{lem}\label{main-geometric-lem}
Let $M$ be a surface of type $(g, p)$ with $\kappa(M)\geq 0$. 
\begin{enumerate}
\item[(i)] If $\alpha \in V(C(M))$ is non-separating, then the resulting surface cut by $\alpha$ is of type $(g-1, p+2)$.
\item[(ii)] If $\alpha \in V(C(M))$ is separating, then let $M_{g_{1}, p_{1}}$, $M_{g_{2}, p_{2}}$ be the two components of the resulting surface cut by $\alpha$. Then $g=g_{1}+g_{2}$ and $p+2=p_{1}+p_{2}$.
\item[(iii)] We have $n(\sigma)\leq n(M)$ for any $\sigma \in S(M)\cup \{ \emptyset \}$. Moreover, there exists $\sigma \in S(M)\cup \{ \emptyset \}$ such that $n(\sigma)=n(M)$.
\end{enumerate}
\end{lem}

\begin{pf}
If $\alpha \in V(C(M))$ and we cut $M$ along $\alpha$, then the following two cases could occur: the resulting surface $M_{\alpha}$ is connected or not. Remark that in general, the number of components in a pants decomposition of a surface $M_{g', p'}$ is equal to $-\chi(M_{g', p'})=2g'+p'-2$, where $\chi(M_{g', p'})$ denotes the Euler characteristic of the surface $M_{g', p'}$ (see \cite[Sections 2.2, 2.3 and 2.4]{ivanov2}). By comparing the number of components in a pants decomposition, we see that if $M_{\alpha}$ is connected, then it is of type $(g-1, p+2)$. If $M_{\alpha}$ is not connected, then it has two components of types $(g_{1}, p_{1})$ and $(g_{2}, p_{2})$, respectively, and we have $(p_{1}-1)+(p_{2}-1)=p$. It follows from
\[(2g_{1}+p_{1}-2)+(2g_{2}+p_{2}-2)=2g+p-2\]
that $g_{1}+g_{2}=g$. We have shown the assertions (i) and (ii).

We prove the assertion (iii). We try to find the cases where maximal value of $n(\sigma)$ is attained. Take $\sigma \in S(M)$. In general, if a surface $M'$ is decomposed into pairs of pants by one simple closed curve, then $M'$ is either of type $(0, 4)$ or $(1, 1)$. If there exists a component $Q$ of $M_{\sigma}$ which is neither of types $(0, 4)$, $(1, 1)$ nor $(0, 3)$ and if $\alpha$ is in $V(C(Q))$, then $n(\sigma \cup \{ \alpha \})\geq n(\sigma)$. Hence, we may assume that any component of $M_{\sigma}$ is one of types $(0, 4)$, $(1, 1)$ and $(0, 3)$. 

Let $n_{1}$, $n_{2}$ and $n_{3}$ be the numbers of components of $M_{\sigma}$ of types $(0, 4)$, $(1, 1)$ and $(0, 3)$, respectively. By definition, $n(\sigma)=n_{1}+n_{2}$. We obtain two pairs of pants from $M_{0, 4}$ and one pair of pants from $M_{1, 1}$ by pants decompositions. By comparing the numbers of components in a pants decomposition of $M$, we see that 
\[2n_{1}+n_{2}+n_{3}=2g+p-2.\]

It follows from the condition of the genus that $n_{2}\leq g$. Thus, it is necessary for the maximality of $n(\sigma)$ that $n_{2}$ is as large and $n_{3}$ is as small as possible. If we put $n_{2}=g$, then $2n_{1}+n_{3}=g+p-2$. If $g+p$ is even, then put $n_{3}=0$, otherwise, put $n_{3}=1$. Then $n(\sigma)$ is equal to $g+[(g+p-2)/2]=n(M)$.

In fact, this value can be attained by the following construction of $\sigma$. Assume $M\neq M_{1, 1}, M_{2, 0}$. First, we can choose $g$ curves which separate $M_{1, 1}$ from $M$ and are non-isotopic each other. The remaining surface is of type $(0, g+p)$. This surface can be decomposed into $(g+p)/2-1$ components of type $(0, 4)$ if $g+p$ is even and into $(g+p-3)/2$ components of type $(0, 4)$ and one component of type $(0, 3)$ if $g+p$ is odd. Finally, if $M=M_{1, 1}$, then it is clear that $n(\sigma)=1$ and if $M=M_{2, 0}$, then we can decompose $M$ into two components of type $(1, 1)$ and so $n(\sigma)=2$.    
\end{pf}

\begin{defn}\label{defn-dagger-n}
Let ${\cal R}$ be a discrete measured equivalence relation on a standard Borel space $(X, \mu)$ with a finite positive measure. For $n\in {\Bbb N}$, we say that a family
\[\{ Y, \{ {\cal S}_{1}^{1}, {\cal S}_{1}^{2}\}, \ldots, \{ {\cal S}_{n}^{1}, {\cal S}_{n}^{2}\} \}\] 
satisfies $(\dagger)_{n}$\index{$(\dagger n $@$(\dagger)_{n}$} if $Y$ is a Borel subset of $X$ with positive measure and $\{ {\cal S}_{k}^{1}, {\cal S}_{k}^{2}\}$ $(k=1, \ldots, n)$ is a pair of subrelations of $({\cal R})_{Y}$ on $Y$ such that for any distinct pair $k, l\in \{1, \ldots, n\}$, the family 
\[\{ Y, \{ {\cal S}_{k}^{1}, {\cal S}_{k}^{2}\}, \{ {\cal S}_{l}^{1}, {\cal S}_{l}^{2}\} \}\] 
satisfies $(\dagger)$.   
\end{defn}

\begin{rem}\label{rem-dagger-n}
Note that if a family 
\[\{ Y, \{ {\cal S}_{1}^{1}, {\cal S}_{1}^{2}\}, \ldots, \{ {\cal S}_{n}^{1}, {\cal S}_{n}^{2}\}\}\]
satisfies $(\dagger)_{n}$, then the family 
\[\{ A, \{ ({\cal S}_{1}^{1})_{A}, ({\cal S}_{1}^{2})_{A}\}, \ldots, \{ ({\cal S}_{n}^{1})_{A}, ({\cal S}_{n}^{2})_{A}\}\}\]
also satisfies $(\dagger)_{n}$ for any Borel subset $A$ of $Y$ with positive measure. 

Suppose that ${\cal R}$ (resp. ${\cal R}'$) is a discrete measured equivalence relation on $(X, \mu)$ (resp. $(X', \mu')$) and there exists a Borel isomorphism $f\colon X\rightarrow X'$ inducing the isomorphism between ${\cal R}$ and ${\cal R}'$. If a family 
\[\{ Y, \{ {\cal S}_{1}^{1}, {\cal S}_{1}^{2}\}, \ldots, \{ {\cal S}_{n}^{1}, {\cal S}_{n}^{2}\}\}\]
of a Borel subset $Y$ of $X$ and subrelations of $({\cal R})_{Y}$ satisfies $(\dagger)_{n}$, then the family 
\[\{ f(Y), \{ f({\cal S}_{1}^{1}), f({\cal S}_{1}^{2})\}, \ldots, \{ f({\cal S}_{n}^{1}), f({\cal S}_{n}^{2}\})\}\]
also satisfies $(\dagger)_{n}$. 

\end{rem}

\begin{cor}\label{main-important-cor}
Assume $(\diamond)'$ and that we have a family 
\[\{ Y, \{ {\cal S}_{1}^{1}, {\cal S}_{1}^{2}\}, \ldots, \{ {\cal S}_{n}^{1}, {\cal S}_{n}^{2}\}\}\]
of a Borel subset $Y$ of $X$ and subrelations of $({\cal R})_{Y}$ satisfying $(\dagger)_{n}$. 
\begin{enumerate}
\item[(i)]There exists a Borel subset $A$ of $Y$ with positive measure such that subsurfaces in the IN systems for ${\cal R}_{k, A}$ and ${\cal R}_{l, A}$ are disjoint for any distinct $k, l\in \{ 1, \ldots, n\}$.  
\item[(ii)] We have $n\leq n(M)$.
\end{enumerate}
\end{cor}

Here, we say that subsurfaces in the IN systems for ${\cal R}_{k, A}$ and ${\cal R}_{l, A}$ are disjoint if $Q_{k}$ and $Q_{l}$ are disjoint for any $Q_{k}\in \psi^{(k)}(x)$, $Q_{l}\in \psi^{(l)}(x)$ and a.e.\ $x\in A$, where $\psi^{(k)}$ and $\psi^{(l)}$ are the IN systems for ${\cal R}_{k, A}$ and ${\cal R}_{l, A}$, respectively. Remark that the assertion (i) holds also for any Borel subset of $A$ with positive measure similarly as in Corollary \ref{cor-main-main-thm}.

\begin{pf}
We may assume $n\geq 2$ because $n(M)=g+[(g+p-2)/2]\geq 1$ when $\kappa(M)\geq 0$. For a Borel subset $A\subseteq Y$ with positive measure and $k\in \{ 1, \ldots, n\}$, we put 
\[{\cal R}_{k, A}=({\cal S}_{k}^{1})_{A}\vee ({\cal S}_{k}^{2})_{A}.\]
By Theorem \ref{main-main-thm}, choose a Borel subset $A_{12}\subseteq X$ with positive measure such that subsurfaces in the IN systems for ${\cal R}_{1, A_{12}}$ and ${\cal R}_{2, A_{12}}$ are disjoint. Applying Theorem \ref{main-main-thm} to the family 
\[\{ A_{12}, \{ ({\cal S}_{1}^{1})_{A_{12}}, ({\cal S}_{1}^{2})_{A_{12}}\}, \{ ({\cal S}_{3}^{1})_{A_{12}}, ({\cal S}_{3}^{2})_{A_{12}}\} \},\]
we can choose a Borel subset $A_{13}\subseteq A_{12}$ with positive measure such that subsurfaces in the IN systems for ${\cal R}_{1, A_{13}}$ and ${\cal R}_{3, A_{13}}$ are disjoint. Remark that subsurfaces in the IN systems for ${\cal R}_{1, A_{13}}$ and ${\cal R}_{2, A_{13}}$ are also disjoint by Corollary \ref{cor-main-main-thm}.

By repeating this process for any pair $k, l\in \{ 1, \ldots, n\}$ with $k<l$, there exists a Borel subset $A\subseteq Y$ with positive measure such that subsurfaces in the IN systems for ${\cal R}_{k, A}$ and ${\cal R}_{l, A}$ are disjoint for any $k, l\in \{ 1, \ldots, n\}$ with $k\neq l$. This means that for a.e.\ $x\in A$ and $k\neq l$, when we denote by $\varphi^{(k)}$ the CRS for ${\cal R}_{k, A}$ and realize the union
\[\varphi(x)=\bigcup_{k=1}^{n}\varphi^{(k)}(x)\in S(M)\]
on $M$, any IN component of $M\setminus \varphi^{(k)}(x)$ for ${\cal R}_{k, A}$ is contained in some component of $M\setminus \varphi^{(l)}(x)$ and is disjoint from any IN component of $M\setminus \varphi^{(l)}(x)$ for ${\cal R}_{l, A}$. Hence, any IN component of $M\setminus \varphi^{(k)}(x)$ for ${\cal R}_{k, A}$ is also a component of $M\setminus \varphi(x)$ for a.e.\ $x\in A$ and any $k$. Since ${\cal R}_{k, A}$ is non-amenable, the IN system for ${\cal R}_{k, A}$ is non-empty for any $k$ by Theorem \ref{irr-non-ame-emp-ame}. The corollary can be verified by Lemma \ref{main-geometric-lem} (iii). 
\end{pf}

\begin{rem}\label{rem-usankusai-right}
We can verify the claim in Remark \ref{rem-usankusai}, using Corollary \ref{main-important-cor} as follows: with the assumption $(\diamond)'$, if we had amenable recurrent subrelations ${\cal S}^{1}$, ${\cal S}^{2}$ of ${\cal S}$ satisfying the conditions in Remark \ref{rem-usankusai}, then for any $n\in {\Bbb N}$, the family
\[\{ Y, \{ {\cal S}_{1}^{1}, {\cal S}_{1}^{2}\},\ldots, \{ {\cal S}_{n}^{1}, {\cal S}_{n}^{2}\} \}\]
defined by ${\cal S}_{i}^{j}={\cal S}^{j}$ for each $i=1, \ldots, n$ and $j=1, 2$ would satisfy $(\dagger)_{n}$. This contradicts Corollary \ref{main-important-cor} (ii). 

On the other hand, using Proposition \ref{prop-inv-partial-amenable}, we can show directly the following assertion: with the assumption $(\ast)'$, let ${\cal S}^{1}$ and ${\cal S}^{2}$ be amenable recurrent subrelations of ${\cal S}$ satisfying only the condition
\begin{enumerate}
\item[(ii)] the relation $({\cal S}^{j})_{A}$ is normal in ${\cal R}_{A}$ for any Borel subset $A$ of $Y$ with positive measure for each $j=1, 2$.
\end{enumerate}
Then the relation ${\cal R}_{A}=({\cal S}^{1})_{A}\vee ({\cal S}^{2})_{A}$ is amenable for some Borel subset $A$ of $Y$ with positive measure.
\end{rem}

Finally, we prove that in the assumption $(\diamond)'$, if $\Gamma =\Gamma(M;m)$, $F=X$ and the action of $\Gamma$ on $X$ is measure-preserving, then the relation ${\cal R}$ has a family of subrelations satisfying $(\dagger)_{n(M)}$.

\begin{lem}\label{lem-main-free-groups-1}
Let $G_{1}$ and $G_{2}$ be two infinite cyclic groups. Suppose that we have an essentially free, measure-preserving Borel action of the free product $G_{1}*G_{2}$ on a standard Borel space $(X, \mu)$ with a finite positive measure. Let us denote by ${\cal S}_{i}$ the equivalence relation on $(X, \mu)$ generated by the action of $G_{i}$ for $i=1, 2$. Then for any Borel subset $A$ of $X$ with positive measure, the relation $({\cal S}_{1})_{A}\vee ({\cal S}_{2})_{A}$ is non-amenable.
\end{lem}

\begin{pf}
This can be shown easily by using cost of a discrete measure-preserving equivalence relation as follows (see Appendix \ref{app-cost}): the relation $({\cal S}_{1})_{A}\vee ({\cal S}_{2})_{A}$ is the free product of the two subrelations $({\cal S}_{1})_{A}$ and $({\cal S}_{2})_{A}$ (see Definition \ref{defn-free-product-relation}). Since $({\cal S}_{1})_{A}$ and $({\cal S}_{2})_{A}$ are amenable and recurrent, both of their costs are equal to $\mu(A)$. Thus, we have 
\[(1/\mu(A)){\cal C}_{\mu}(({\cal S}_{1})_{A}\vee ({\cal S}_{2})_{A})=2\]
by Theorem \ref{thm-cost-free-product}. This implies that $({\cal S}_{1})_{A}\vee ({\cal S}_{2})_{A}$ is non-amenable by Theorem \ref{thm-lev-amenable}.
\end{pf}

\begin{lem}\label{lem-main-free-groups-2}
Let $G_{1}$ be a discrete group and $G_{2}$ and $G_{3}$ be two subgroups of a discrete group $G$ such that $G_{2}$ and $G_{3}$ generate $G$. Suppose that we have an essentially free, non-singular Borel action of the group $G_{1}\times G$ on a standard Borel space $(X, \mu)$ with a finite positive measure. Let us denote by ${\cal R}_{j}$ the equivalence relation on $(X, \mu)$ generated by the action of $G_{j}$ for $j=1, 2, 3$. Then for any pair of Borel subsets $A\subseteq B$ of $X$ with positive measure, the subrelation $({\cal R}_{1})_{A}$ is normal in $({\cal R}_{1})_{A}\vee (({\cal R}_{2})_{B}\vee ({\cal R}_{3})_{B})_{A}$.
\end{lem}

\begin{pf}
For simplicity, we denote ${\cal T}=({\cal R}_{2})_{B}\vee ({\cal R}_{3})_{B}$. Let $\{ g_{n}\}_{n\in {\Bbb N}}$ and $\{ h_{m}\}_{m\in {\Bbb N}}$ be enumerations of all elements in $G_{2}$ and $G_{3}$, respectively. We define Borel maps 
\begin{align*}
g_{n}'&\colon B\cap g_{n}^{-1}(B)\rightarrow B\cap g_{n}(B),\\
h_{m}'&\colon B\cap h_{m}^{-1}(B)\rightarrow B\cap h_{m}(B)
\end{align*}
by the restrictions of $g_{n}$, $h_{m}$, respectively. Then for a.e.\ $(x, y)\in {\cal T}$, we can find a word $w$ made of letters $g_{n}'$, $h_{m}'$ such that $x\in {\rm dom}(w)$ and $y=w(x)$. 

For any word $w$ made of letters $g_{n}'$, $h_{m}'$, we define a Borel map 
\[\omega_{w}\colon A\cap w^{-1}(A\cap {\rm ran}(w))\rightarrow w(A\cap {\rm dom}(w))\cap A\]
by the restriction of $w$. Then for a.e.\ $(x, y)\in ({\cal T})_{A}$, there exists $w$ such that $x\in {\rm dom}(\omega_{w})$ and $y=\omega_{w}(x)$.

We show that $\omega_{w}$ is in $[[{\cal T'}]]_{({\cal R}_{1})_{A}}$, where ${\cal T}'=({\cal R}_{1})_{A}\vee ({\cal T})_{A}$. It is clear that the Borel map $\omega_{w}$ corresponding to a word $w$ of letters $g_{n}'$, $h_{m}'$ satisfies that $(\omega_{w}(x), x)\in {\cal T}'$ for a.e.\ $x\in {\rm dom}(\omega_{w})$. For $x\in {\rm dom}(\omega_{w})$ and $f\in G_{1}$ with $fx\in {\rm dom}(\omega_{w})$, we have
\[\omega_{w}(fx)=gfx=fgx=f\omega_{w}(x),\]
where $g\in G$ is the corresponding element to the word $w$. This means that $(\omega_{w}(fx), \omega_{w}(x))\in ({\cal R}_{1})_{A}$. Thus, if $(x, y)\in ({\cal R}_{1})_{A}$ and $x, y\in {\rm dom}(\omega_{w})$, then we have $(\omega_{w}(x), \omega_{w}(y))\in ({\cal R}_{1})_{A}$. Since the inverse of the partial Borel isomorphism $\omega_{w}$ is equal to $\omega_{w^{-1}}$ for the word $w^{-1}$, we see that for a.e.\ $x, y\in {\rm dom}(\omega_{w})$, if $(\omega_{w}(x), \omega_{w}(y))\in ({\cal R}_{1})_{A}$, then $(x, y)\in ({\cal R}_{1})_{A}$. Therefore, we have shown that any $\omega_{w}$ is in $[[{\cal T'}]]_{({\cal R}_{1})_{A}}$.

For a.e.\ $(x, y)\in {\cal T}'$, we can find a finite sequence $x=x_{0}, x_{1}, \ldots, x_{l}=y$ such that $(x_{2j-1}, x_{2j})\in ({\cal R}_{1})_{A}$ and $(x_{2j}, x_{2j+1})\in ({\cal T})_{A}$ for each $j$. Moreover, we can find words $w_{j}$ of letters $g_{n}'$, $h_{m}'$ such that $x_{2j}\in {\rm dom}(\omega_{w_{j}})$ and $x_{2j+1}=\omega_{w_{j}}(x_{2j})$ for each $j$. It means that the relation $({\cal R}_{1})_{A}\vee ({\cal T})_{A}$ is generated by $({\cal R}_{1})_{A}$ and the countable set of the Borel maps $\omega_{w}$ for all words $w$. It follows from Lemma \ref{generate-generate-normal} that the subrelation $({\cal R}_{1})_{A}$ is normal in the relation ${\cal T}'=({\cal R}_{1})_{A}\vee ({\cal T})_{A}$. 
\end{pf}

\begin{lem}\label{lem-main-free-groups-3}
For $i=1, 2$, let $G_{i}^{1}$ and $G_{i}^{2}$ be two subgroups of a discrete group $G_{i}$ such that $G_{i}^{1}$ and $G_{i}^{2}$ generate $G_{i}$. Suppose that we have an essentially free, non-singular Borel action of the group $G_{1}\times G_{2}$ on a standard Borel space $(X, \mu)$ with a finite positive measure. Let us denote by ${\cal R}_{i}^{j}$ the equivalence relation on $(X, \mu)$ generated by the action of $G_{i}^{j}$ for $i, j=1, 2$. Let us denote ${\cal R}_{i, A}=({\cal R}_{i}^{1})_{A}\vee ({\cal R}_{i}^{2})_{A}$ for a Borel subset $A$ of $X$ with positive measure and $i=1, 2$. Then the relation ${\cal R}_{i, A}$ is normal in ${\cal R}_{1, A}\vee {\cal R}_{2, A}$ for each $i=1, 2$.
\end{lem}

\begin{pf}
By symmetry, it suffices to show that ${\cal R}_{1, A}$ is normal in ${\cal R}_{1, A}\vee {\cal R}_{2, A}$. Let $\{ g_{n}\}_{n\in {\Bbb N}}$ and $\{ h_{m}\}_{m\in {\Bbb N}}$ be enumerations of all elements in $G_{2}^{1}$ and $G_{2}^{2}$, respectively. As in the proof of Lemma \ref{lem-main-free-groups-2}, we define Borel maps
\begin{align*}
g_{n}'&\colon A\cap g_{n}^{-1}(A)\rightarrow A\cap g_{n}(A),\\
h_{m}'&\colon A\cap h_{m}^{-1}(A)\rightarrow A\cap h_{m}(A)
\end{align*}
by the restrictions of $g_{n}$, $h_{m}$, respectively. Then for a.e.\ $(x, y)\in ({\cal R}_{1}^{1})_{A}$ with $x, y\in {\rm dom}(g_{n}')$, we can find a unique $g\in G_{1}^{1}$ such that $y=gx$. Then we have
\[g_{n}'(y)=g_{n}gx=gg_{n}x=gg_{n}'(x).\]
It follows that $(g_{n}'(x), g_{n}'(y))\in ({\cal R}_{1}^{1})_{A}$. Similarly, we can show that for a.e.\ $x, y\in {\rm dom}(g_{n}')$, we have $(g_{n}'(x), g_{n}'(y))\in ({\cal R}_{1}^{2})_{A}$ if $(x, y)\in ({\cal R}_{1}^{2})_{A}$. It follows that for a.e.\ $x, y\in {\rm dom}(g_{n}')$ with $(x, y)\in {\cal R}_{1, A}$, we have $(g_{n}'(x), g_{n}'(y))\in {\cal R}_{1, A}$. 

Hence, the Borel map $g_{n}'$ belongs to $[[{\cal R}_{1, A}\vee {\cal R}_{2, A}]]_{{\cal R}_{1, A}}$ because there exists $n'\in {\Bbb N}$ such that $(g_{n}')^{-1}=g_{n'}'$. Similarly, we can prove that each $h_{m}'$ is also in $[[{\cal R}_{1, A}\vee {\cal R}_{2, A}]]_{{\cal R}_{1, A}}$. Since $\{ g_{n}'\}$, $\{ h_{m}'\}$ and ${\cal R}_{1, A}$ generate the relation ${\cal R}_{1, A}\vee {\cal R}_{2, A}$, the lemma follows from Lemma \ref{generate-generate-normal}. 
\end{pf}

\begin{cor}\label{cor-main-main-free-groups}
Let $M$ be a surface with $\kappa(M)\geq 0$. Suppose that an infinite subgroup $\Gamma$ of the mapping class group $\Gamma(M)$ is measure equivalent to a group $G$ containing a subgroup isomorphic to the direct product of $n$ free groups with rank $2$. Then $n\leq n(M)$.
\end{cor}

\begin{pf}
We may assume that $\Gamma$ is a subgroup of $\Gamma(M;m)$ for some integer $m\geq 3$. We suppose that there exists a group $G$ measure equivalent to $\Gamma$ which contains a subgroup $(G_{1}^{1}\ast G_{1}^{2})\times \cdots \times (G_{n}^{1}\ast G_{n}^{2})$ with each $G_{i}^{j}$ isomorphic to ${\Bbb Z}$. Then we have two equivalence relations ${\cal R}_{1}$ on $(X_{1}, \mu_{1})$ and ${\cal R}_{2}$ on $(X_{2}, \mu_{2})$ of type ${\rm II}_{1}$ satisfying the following conditions: the relation ${\cal R}_{1}$ (resp. ${\cal R}_{2}$) is generated by an essentially free, measure-preserving Borel action of $\Gamma$ on $(X_{1}, \mu_{1})$ (resp. $G$ on $(X_{2}, \mu_{2})$). Moreover, ${\cal R}_{1}$ and ${\cal R}_{2}$ are weakly isomorphic.

Let $f\colon A_{1}\rightarrow A_{2}$, $A_{i}\subseteq X_{i}$ be a partial Borel isomorphism inducing the weak isomorphism between ${\cal R}_{1}$ and ${\cal R}_{2}$. Let us denote by ${\cal S}_{i}^{j}$ the subrelation of ${\cal R}_{2}$ generated by the action of the subgroup $G_{i}^{j}$. It follows from Lemmas \ref{lem-main-free-groups-1}, \ref{lem-main-free-groups-2} and \ref{lem-main-free-groups-3} that the inverse image by $f$ of the family 
\[\{ A_{2}, \{ ({\cal S}_{1}^{1})_{A_{2}}, ({\cal S}_{1}^{2})_{A_{2}}\},\ldots, \{ ({\cal S}_{n}^{1})_{A_{2}}, ({\cal S}_{n}^{2})_{A_{2}}\} \}\]
of the set $A_{2}$ and subrelations of $({\cal R}_{2})_{A_{2}}$ satisfies $(\dagger)_{n}$. The desired inequality follows from Corollary \ref{main-important-cor} (ii). 
\end{pf}

\begin{cor}\label{cor-classification1}
Suppose that $M^{1}$ and $M^{2}$ are compact orientable surfaces with $\kappa(M^{1}), \kappa(M^{2})\geq 0$ and that the mapping class groups $\Gamma(M^{1})$ and $\Gamma(M^{2})$ are measure equivalent. Then $n(M^{1})=n(M^{2})$.
\end{cor}

\begin{pf}
Let $\sigma \in S(M)\cup \{ \emptyset \}$ be a disjoint family of curves on $M^{1}$ such that $n(\sigma)$ is maximal, that is, $n(\sigma)=n(M^{1})$ (see Lemma \ref{main-geometric-lem} (iii)). For each component $Q$ of $(M^{1})_{\sigma}$ of type $(0, 4)$ or $(1, 1)$, choose two intersecting curves on $Q$ and consider their Dehn twists. They generate a group containing a free subgroup of rank 2 (see Theorems \ref{tits-alternative-mcg} and \ref{commuting-dehn}). It follows that $\Gamma(M^{1})$ has $n(\sigma)=n(M^{1})$ free subgroups of rank $2$ commuting each other and that $n(M^{1})\leq n(M^{2})$ by Corollary \ref{cor-main-main-free-groups}. By symmetry, we obtain the equality $n(M^{1})=n(M^{2})$.
\end{pf}

The proof of Corollary \ref{cor-classification1} shows that the mapping class group $\Gamma(M_{g, p})$ contains a subgroup isomorphic to the direct product of $g+[(g+p-2)/2]$ free groups of rank $2$. Thus, Corollary \ref{cor-main-main-free-groups} shows the following corollary which has already shown by Gaboriau \cite[Corollaire 5.9]{gab-l2}. For $n\in {\Bbb N}\setminus \{ 0\}$, let us denote by ${\Bbb F}_{n}$ the free group of rank $n$.

\begin{cor}\label{cor-gab-free-groups}
Let $n_{1}, \ldots, n_{k}, m_{1}, \ldots, m_{l}$ be natural numbers more than or equal to $2$. If the direct product ${\Bbb F}_{n_{1}}\times \cdots \times {\Bbb F}_{n_{k}}$ is measure equivalent to a discrete group containing the direct product ${\Bbb F}_{m_{1}}\times \cdots \times {\Bbb F}_{m_{l}}$, then $k\geq l$.
\end{cor}


\chapter[Classification in terms of measure equivalence II]{Classification of the mapping class groups in terms of measure equivalence II}\label{chapter-best2}

In this chapter, we improve the classification of the mapping class groups in terms of measure equivalence in Corollary \ref{cor-classification1} in the following form (see Corollary \ref{cor-kappa-invariant}, Corollary \ref{cor-kappa-odd-class} (ii) and Corollary \ref{cor-kappa-even-class-mcg}): let us denote
\begin{equation*}
g_{0}(M)= \begin{cases}
           2 & {\rm if} \ g\leq 2, \\
           g & {\rm if} \ g>2
          \end{cases}
\end{equation*}
for a compact orientable surface $M$ of type $(g, p)$. 

\begin{thm}\label{thm-main-class-final}
Suppose that $M^{1}$ and $M^{2}$ are compact orientable surfaces with $\kappa(M^{1}), \kappa(M^{2})\geq 0$ and that the mapping class groups $\Gamma(M^{1})$ and $\Gamma(M^{2})$ are measure equivalent. Then $\kappa(M^{1})=\kappa(M^{2})$ and $g_{0}(M^{1})=g_{0}(M^{2})$.
\end{thm}

Given a compact orientable surface $M$ of type $(g, p)$ and a relation ${\cal R}$ generated by an essentially free, measure-preserving action of the mapping class group $\Gamma(M)$ on $(X, \mu)$, we study a family 
\[{\cal F}=\{ Y, \{ {\cal S}_{1}^{1}, {\cal S}_{1}^{2}\},\ldots, \{ {\cal S}_{n(M)}^{1}, {\cal S}_{n(M)}^{2}\} \}\]
of a Borel subset $Y$ of $X$ with positive measure and subrelations of $({\cal R})_{Y}$ satisfying $(\dagger)_{n(M)}$ (see Definition \ref{defn-dagger-n}), where $n(M)$ is equal to $g+[(g+p-2)/2]$ and is the maximal number with property described in Corollary \ref{main-important-cor} (ii). Roughly speaking, it follows from the proof of Corollary \ref{main-important-cor} that for such a family ${\cal F}$, there exists $\sigma \in S(M)$ with $n(\sigma)=n(M)$ such that each pair $\{ {\cal S}_{i}^{1}, {\cal S}_{i}^{2}\}$ corresponds to a unique component $Q_{i}$ of $M_{\sigma}$, which is an IN component for the relation ${\cal R}_{i, A}=({\cal S}_{i}^{1})_{A}\vee ({\cal S}_{i}^{2})_{A}$ for some Borel subset $A$ of $Y$ with positive measure. Recall that $n(\sigma)$ is the number of components of $M_{\sigma}$ which are not pairs of pants. In the construction of a family satisfying $(\dagger)_{n(M)}$ in Corollary \ref{cor-main-main-free-groups}, we may regard ${\cal S}_{i}^{1}$ and ${\cal S}_{i}^{2}$ as subrelations generated by Dehn twists about two intersecting curves in $V(C(Q_{i}))$. This correspondence enables us to use topological properties of the surface for the classification of the mapping class groups in terms of measure equivalence. In this chapter, we push this approach further.
 
In Section \ref{sec-geometric-lemmas}, we study $\sigma \in S(M)$ with $n(\sigma)=n(M)$. When the surface $M$ is cut along all curves in $\sigma$, the possibility of the types and numbers of the components in $M_{\sigma}$ is investigated. This possibility is fairly different between the cases where $\kappa(M)$ is odd and even. We also give another difference between them as follows (see Lemmas \ref{lem-even-nn} and \ref{lem-odd-nn}): assume that $\kappa(M)$ is even. If $\sigma \in S(M)$ and $\alpha \in V(C(M))\setminus \sigma$ satisfy $n(\sigma)=n(M)$ and $i(\alpha, \beta )=0$ for any $\beta \in \sigma$, then $n(\tau)\leq n(M)-1$ for any $\tau \in S(M)$ with $\alpha \in \tau$. On the other hand, if $\kappa(M)$ is odd, then for any $\alpha \in V(C(M))$, there exists $\tau \in S(M)$ such that $n(\tau)=n(M)$ and $\alpha \in \tau$.   

In Section \ref{sec-families-subrelations}, using geometric lemmas in Section \ref{sec-geometric-lemmas}, we find certain properties which any family ${\cal F}$ satisfying $(\dagger)_{n(M)}$ possesses. As mentioned above, one can think that subrelations in such a family ${\cal F}$ correspond to some curves in $V(C(M))$. We shall state one of the above geometric properties about curves in terms of relations as follows (see Theorem \ref{thm-even-n-1}): suppose that $\kappa(M)$ is even and that we have a family 
\[{\cal F}=\{ Y, \{ {\cal S}_{1}^{1}, {\cal S}_{1}^{2}\},\ldots, \{ {\cal S}_{n(M)}^{1}, {\cal S}_{n(M)}^{2}\} \}\]
of a Borel subset $Y$ of $X$ and subrelations of $({\cal R})_{Y}$ satisfying $(\dagger)_{n(M)}$. Then for any $i\in \{ 1, \ldots, n(M)\}$, there exists $s\in \{ 1, 2\}$ satisfying the following condition: assume that we have a family
\[\{ Y, \{ {\cal T}_{1}^{1}, {\cal T}_{1}^{2}\},\ldots, \{ {\cal T}_{n}^{1}, {\cal T}_{n}^{2}\} \}\]
of the Borel set $Y$ and subrelations of $({\cal R})_{Y}$ satisfies $(\dagger)_{n}$ and a certain condition involving the subrelations ${\cal S}_{i}^{s}$ (that is, ${\cal S}_{i}^{s}$ and ${\cal T}_{j}^{t}$ are ``commutative'' for any $j\in \{ 1, \ldots, n\}$ and $t\in \{ 1, 2\}$). Then $n\leq n(M)-1$. We give the definition of the ``commutativity'' in Definition \ref{defn-dagger-n-wrt}, which is formulated in terms of normal subrelations as well as that of $(\dagger)$ in Definition \ref{defn-dagger}.

On the other hand, another geometric property about curves shows that if $\kappa(M)$ is odd, then there exists a family 
\[{\cal F}=\{ X, \{ {\cal S}_{1}^{1}, {\cal S}_{1}^{2}\},\ldots, \{ {\cal S}_{n(M)}^{1}, {\cal S}_{n(M)}^{2}\} \}\] 
satisfying $(\dagger)_{n(M)}$ and the following condition: for any $i\in \{ 1, \ldots, n(M)\}$ and $s\in \{ 1, 2\}$, there exists a family
\[\{ X, \{ {\cal T}_{1}^{1}, {\cal T}_{1}^{2}\}, \ldots, \{ {\cal T}_{n(M)}^{1}, {\cal T}_{n(M)}^{2}\} \}\]
satisfying $(\dagger)_{n(M)}$ and the above stated condition with respect to the relation ${\cal S}_{i}^{s}$ (see Theorem \ref{thm-kappa-invariant-odd}). Remark that $n(M)=[(\kappa(M)-2)/2]$. Thus, in order to deduce the equality $\kappa(M^{1})=\kappa(M^{2})$ in Theorem \ref{thm-main-class-final}, it suffices to show that for two surfaces $M^{1}$ and $M^{2}$, if $n(M^{1})=n(M^{2})$ and $\kappa(M^{1})$ is even and $\kappa(M^{2})$ is odd, then their mapping class groups are not measure equivalent. The above remarkable difference about families satisfying $(\dagger)_{n(M)}$ between the cases where $\kappa(M)$ is even and odd enables us to prove the equality $\kappa(M^{1})=\kappa(M^{2})$ in Theorem \ref{thm-main-class-final}. This will be done in Section \ref{sec-app1}.

In Section \ref{sec-families-subrelations}, we give other two statements (Theorems \ref{thm-odd-n-1} and \ref{thm-even-n-2}) about families ${\cal F}$ satisfying $(\dagger)_{n(M)}$. They will be used in Sections \ref{sec-app2} and \ref{sec-app3}, where we deduce the equality $g_{0}(M^{1})=g_{0}(M^{2})$ in Theorem \ref{thm-main-class-final}.

In Sections \ref{sec-app2} and \ref{sec-app3}, we investigate the cases where $\kappa(M)$ is odd and even, respectively. The arguments in these sections are combinatorial. Our fundamental attitude for these sections is as follows: if we are going to prove that the mapping class groups $\Gamma(M^{1})$ and $\Gamma(M^{2})$ of two distinct surfaces $M^{1}$ and $M^{2}$ with $\kappa(M^{1})=\kappa(M^{2})$ are not measure equivalent, then first, we show that the relation ${\cal R}_{1}$ generated by the action of $\Gamma(M^{1})$ satisfies some property given in Section \ref{sec-families-subrelations}. On the other hand, we prove that the relation ${\cal R}_{2}$ generated by the action of $\Gamma(M^{2})$ has a family which does not satisfy the property. The concrete construction of such a family for the relation ${\cal R}_{2}$ is the main task of Sections \ref{sec-app2} and \ref{sec-app3}.

\section{Geometric lemmas}\label{sec-geometric-lemmas}

For a compact orientable surface $M$ of type $(g, p)$, we write $n(M)=g+[(g+p-2)/2]$. In this section, we study $\sigma \in S(M)$ with $n(\sigma)=n(M)$, where $n(\sigma)$ is the number of components of $M_{\sigma}$ which are not pairs of pants.

\begin{lem}\label{lem-even-n}
Suppose that $\kappa(M)\geq 0$ is even and we have $\sigma \in S(M)$ with $n(\sigma)=n(M)=g+(g+p-2)/2$. 
\begin{enumerate}
\item[(i)] The numbers of components of $M_{\sigma}$ of types $(1, 1)$ and $(0, 4)$ are $g$ and $(g+p-2)/2$, respectively. 
\item[(ii)] All curves in $\sigma$ are separating ones.
\end{enumerate}
\end{lem}

\begin{pf}
In general, if $Q$ is a surface with $\kappa(Q)\geq 0$ and $n(Q)=1$, then it is one of types $(0, 4)$, $(0, 5)$, $(1, 1)$ and $(1, 2)$. If $M_{\sigma}$ had a component which is neither of the above types nor of type $(0, 3)$, then the component could be decomposed into at least two components neither of which are pairs of pants. This contradicts $n(\sigma)=n(M)$ and Lemma \ref{main-geometric-lem} (iii). Thus, any component of $M_{\sigma}$ is one of the above five types.

Let $n_{0, 3}$, $n_{0, 4}$, $n_{0, 5}$, $n_{1, 1}$ and $n_{1, 2}$ be the numbers of components of $M_{\sigma}$ of types $(0, 3)$, $(0, 4)$, $(0, 5)$, $(1, 1)$ and $(1, 2)$, respectively. By comparing the number of components in a pants decomposition, we have 
\begin{equation}\label{eq-pants}
n_{0, 3}+2n_{0, 4}+3n_{0, 5}+n_{1, 1}+2n_{1, 2}=2g+p-2.
\end{equation}
It follows from the definition of $n(\sigma)$ that 
\begin{equation}\label{eq-n}
n_{0, 4}+n_{0, 5}+n_{1, 1}+n_{1, 2}=g+(g+p-2)/2.
\end{equation}   
By the condition of the genus, we have 
\begin{equation}\label{eq-genus}
n_{1, 1}+n_{1, 2}\leq g
\end{equation}
(see also Lemma \ref{main-geometric-lem} (i), (ii)). Calculating $(\ref{eq-pants})-(\ref{eq-n})\times 2$, we see that $n_{0, 3}+n_{0, 5}-n_{1, 1}=-g$. It follows from (\ref{eq-genus}) that $n_{0, 3}+n_{0, 5}+g=n_{1, 1}\leq g-n_{1, 2}\leq g$. Thus, we have $n_{0, 3}+n_{0, 5}\leq 0$, which means $n_{0, 3}=n_{0, 5}=0$. It follows that $n_{1, 1}=g$, $n_{1, 2}=0$ and $n_{0, 4}=(g+p-2)/2$, which shows the assertion (i).

If $\alpha \in \sigma$ were non-separating, then $M_{\alpha}$ would be of type $(g-1, p+2)$ by Lemma \ref{main-geometric-lem} (i). Then it is impossible for $M_{\sigma}$ to have $g$ components of type $(1, 1)$, which contradicts the assertion (i).    
\end{pf}

Next, we consider the case where $\kappa(M)$ is odd.

\begin{lem}\label{lem-odd-n}
Suppose that $\kappa(M)\geq 0$ is odd and we have $\sigma \in S(M)$ with $n(\sigma)=n(M)=g+(g+p-3)/2$. 
\begin{enumerate}
\item[(i)] Any component of $M_{\sigma}$ is one of types $(0, 3)$, $(0, 4)$, $(0, 5)$, $(1, 1)$ and $(1, 2)$.
\item[(ii)] Let $n_{0, 3}$, $n_{0, 4}$, $n_{0, 5}$, $n_{1, 1}$ and $n_{1, 2}$ be the numbers of components of $M_{\sigma}$ of types $(0, 3)$, $(0, 4)$, $(0, 5)$, $(1, 1)$ and $(1, 2)$, respectively. Then we have the following four possibilities of $(n_{0, 3}, n_{0, 4}, n_{0, 5}, n_{1, 1}, n_{1, 2})$: 
\begin{enumerate}
\item[(a)] $(1, (g+p-3)/2, 0, g, 0)$;
\item[(b)] $(0, (g+p-5)/2, 1, g, 0)$ (if $g+p\geq 5$);
\item[(c)] $(0, (g+p-1)/2, 0, g-1, 0)$ (if $g\geq 1$);
\item[(d)] $(0, (g+p-3)/2, 0, g-1, 1)$ (if $g\geq 1$).
\end{enumerate}
\end{enumerate}
\end{lem}

\begin{pf}
The proof of the assertion (i) is the same as the beginning of the proof of Lemma \ref{lem-even-n}. Moreover, as in the proof of Lemma \ref{lem-even-n}, we obtain the following two equalities and one inequality:
\begin{equation}\label{eq-pants-2}
n_{0, 3}+2n_{0, 4}+3n_{0, 5}+n_{1, 1}+2n_{1, 2} = 2g+p-2, 
\end{equation}
\begin{equation}\label{eq-n-2}
n_{0, 4}+n_{0, 5}+n_{1, 1}+n_{1, 2} = g+(g+p-3)/2,
\end{equation}
\begin{equation}\label{eq-genus-2}
n_{1, 1}+n_{1, 2} \leq g.
\end{equation}
Calculating $(\ref{eq-pants-2})-(\ref{eq-n-2})\times 2$, we see that $n_{0, 3}+n_{0, 5}-n_{1, 1}=-g+1$. It follows from (\ref{eq-genus-2}) that $n_{0, 3}+n_{0, 5}=n_{1, 1}-g+1\leq -n_{1, 2}+1\leq 1$.

If $n_{0, 3}=1$, then $n_{0, 5}=0$ and the case (a) occurs. If $n_{0, 5}=1$, then $n_{0, 3}=0$ and the case (b) occurs. If $n_{0, 3}=n_{0, 5}=0$, then we have $n_{1, 1}=g-1$, $n_{0, 4}+n_{1, 2}=(g+p-1)/2$ and $n_{1, 2}\leq 1$. It follows that either the case (c) or (d) occurs.
\end{pf}

\begin{rem}\label{rem-lem-odd-n}
In fact, for each possibility $(n_{0, 3}, n_{0, 4}, n_{0, 5}, n_{1, 1}, n_{1, 2})$ in Lemma \ref{lem-odd-n} (ii), we can construct $\sigma \in S(M)$ such that the numbers of components of $M_{\sigma}$ for each type are $(n_{0, 3}, n_{0, 4}, n_{0, 5}, n_{1, 1}, n_{1, 2})$. 
\end{rem}

The following two lemmas give a remarkable difference between the cases where $\kappa(M)$ is even and odd:

\begin{lem}\label{lem-even-nn}
Suppose that $\kappa(M)\geq 0$ is even and we have $\sigma \in S(M)$ with $n(\sigma)=n(M)$. Let $\alpha \in V(C(M))\setminus \sigma$ be a curve satisfying $i(\alpha, \beta)=0$ for any $\beta \in \sigma$. Then $n(\tau)\leq n(M)-1$ for any $\tau \in S(M)$ with $\alpha \in \tau$.
\end{lem}

\begin{pf}
It follows from Lemma \ref{lem-even-n} (i) that $M_{\sigma}$ has $g$ components of type $(1, 1)$ and $(g+p-2)/2$ components of type $(0, 4)$. Assume the existence of $\tau \in S(M)$ such that $\alpha \in \tau$ and $n(\tau)=n(M)$. We will deduce a contradiction.

By Lemma \ref{lem-even-n} (ii), the curve $\alpha$ is separating. Since $\alpha$ does not intersect any curve in $\sigma$, it is a curve on some component of $M_{\sigma}$. If the component were of type $(1, 1)$, then $\alpha$ would be a non-separating curve on $M$, which is a contradiction. Thus, $\alpha$ is a curve on some component $Q$ of $M_{\sigma}$ of type $(0, 4)$. 

Since $\alpha$ is separating, we have a decomposition $M_{\alpha}=Q_{1}\sqcup Q_{2}$. Suppose that $Q_{i}$ is of type $(g_{i}, p_{i})$ for $i=1, 2$. Then we have two equalities $g_{1}+g_{2}=g$ and $p_{1}+p_{2}=p+2$. 

We show that both $g_{1}+p_{1}$ and $g_{2}+p_{2}$ are odd. For each $i=1, 2$, there exists $g_{i}$ curves in $\sigma$ each of which separates a component of type $(1, 1)$ from $Q_{i}$ and the remaining surface is $M_{0, g_{i}+p_{i}}$. (Remark that each $Q_{i}$ is not of type $(2, 0)$ since $p_{i}\geq 1$.) Moreover, the remaining surface is decomposed into a certain number of components of type $(0, 4)$ and one component of type $(0, 3)$ by the other curves in $\sigma$ since $\alpha$ cut $Q$ into two pairs of pants. It follows that $g_{i}+p_{i}-2$, the number of components in a pants decomposition of $M_{0, g_{i}+p_{i}}$ is odd and that so is $g_{i}+p_{i}$.

Therefore, 
\begin{align*}
n(Q_{1})+n(Q_{2})&=g_{1}+[(g_{1}+p_{1}-2)/2]+g_{2}+[(g_{2}+p_{2}-2)/2]\\
                 &=g_{1}+(g_{1}+p_{1}-3)/2+g_{2}+(g_{2}+p_{2}-3)/2 \\
                 &=g+(g+p-2)/2-1=n(M)-1.
\end{align*}
On the other hand, any $\tau'\in S(M)$ with $\alpha \in \tau'$ satisfies
\[n(\tau')\leq n(Q_{1})+n(Q_{2})\leq n(M)-1,\]
which contradicts the existence of $\tau$.   
\end{pf}

\begin{rem}\label{rem-lem-even-nn}
We can show the following assertion by the above proof: suppose that $\kappa(M)\geq 0$ is even and we have $\sigma \in S(M)$ with $n(\sigma)=n(M)$. Let $Q$ be a component of $M_{\sigma}$ of type $(0, 4)$ and $\alpha \in V(C(Q))$. Then $\alpha$ is a separating curve on $M$. Otherwise, $M_{\sigma \cup \{ \alpha \}}$ would have at most $g-1$ components of type $(1, 1)$, which is a contradiction. Moreover, if we cut $M$ along $\alpha$, then we have two components both of whose complexities are odd.
\end{rem}

\begin{lem}\label{lem-odd-nn}
Suppose that $\kappa(M)\geq 0$ is odd. Let $\alpha \in V(C(M))$. Then there exists $\sigma \in S(M)$ such that $\alpha \in \sigma$ and $n(\sigma)=n(M)$.
\end{lem}

\begin{pf}
First, assume that $\alpha$ is a separating curve. Then we have a decomposition $M_{\alpha}=Q_{1}\sqcup Q_{2}$. Suppose that $Q_{i}$ is of type $(g_{i}, p_{i})$ for $i=1, 2$. Then we have two equalities $g_{1}+g_{2}=g$ and $p_{1}+p_{2}=p+2$. Since exactly one of the numbers $g_{1}+p_{1}$ and $g_{2}+p_{2}$ is even, we have
\begin{align*}
n(Q_{1})+n(Q_{2}) &=g_{1}+[(g_{1}+p_{1}-2)/2]+g_{2}+[(g_{2}+p_{2}-2)/2]\\
                  &=g+(g+p-3)/2=n(M).
\end{align*}
Since there exists $\sigma_{i}\in S(Q_{i})$ such that $n(\sigma_{i})=n(Q_{i})$ for each $i=1, 2$, the union $\sigma =\{ \alpha \}\cup \sigma_{1}\cup \sigma_{2}\in S(M)$ is a desired element.

Next, assume that $\alpha$ is a non-separating curve. Then $M_{\alpha}$ is of type $(g-1, p+2)$ and 
\begin{align*}
n(M_{g-1, p+2}) &=g-1+[((g-1)+(p+2)-2)/2]\\
                &=g-1+[(g+p-1)/2]=g-1+(g+p-1)/2=n(M).
\end{align*}
Since there exists $\sigma'\in S(M_{\alpha})$ with $n(\sigma')=n(M_{\alpha})$, the union $\sigma =\{ \alpha \}\cup \sigma'\in S(M)$ is a desired element.
\end{pf}


\section{Families of subrelations satisfying the maximal condition}\label{sec-families-subrelations}

Let $M$ be a compact orientable surface of type $(g, p)$ satisfying $\kappa(M)\geq 0$. Assume $(\diamond)'$ (see Assumption \ref{assumption-diamond'}). Hence, we have a subgroup $\Gamma$ of $\Gamma(M; m)$ with $m\geq 3$ and a discrete measured equivalence relation ${\cal R}$ on a standard Borel space $(X, \mu)$ with a finite positive measure. Moreover, we have a Borel cocycle $\rho \colon {\cal R}\rightarrow \Gamma$ with finite kernel. There exists a fundamental domain $F\subseteq X$ of the finite relation $\ker \rho$ such that we have an essentially free, non-singular Borel action of $\Gamma$ on $(F, \mu|_{F})$ generating the relation $({\cal R})_{F}$ and whose induced cocycle is equal to the restriction of $\rho$ to $({\cal R})_{F}$.     

We assume that there exists a family 
\[{\cal F}=\{ Y, \{ {\cal S}_{1}^{1}, {\cal S}_{1}^{2}\},\ldots, \{ {\cal S}_{n(M)}^{1}, {\cal S}_{n(M)}^{2}\} \}\]
of a Borel subset $Y$ of $X$ and subrelations of $({\cal R})_{Y}$ satisfying $(\dagger)_{n(M)}$. Note that the number $n(M)$ has the maximal property described in Corollary \ref{main-important-cor} (ii). Let us denote ${\cal S}_{i, A}=({\cal S}_{i}^{1})_{A}\vee ({\cal S}_{i}^{2})_{A}$ for $i\in \{ 1, \ldots, n(M)\}$ and a Borel subset $A$ of $Y$ with positive measure.

It follows from Corollary \ref{main-important-cor} (i) that there exists a Borel subset $Y_{1}$ of $Y$ with positive measure such that the IN systems for ${\cal S}_{i, Y_{1}}$ and ${\cal S}_{j, Y_{1}}$ are disjoint for any distinct $i, j\in \{ 1, \ldots, n(M)\}$. Moreover, there exists a Borel subset $Y_{2}$ of $Y_{1}$ with positive measure such that all the CRS, T, IA and IN systems for the relations ${\cal S}_{i, Y_{1}}$ and ${\cal S}_{i}^{s}$ are constant on $Y_{2}$ for any $i\in \{ 1, \ldots, n(M)\}$ and $s\in \{ 1, 2\}$. 

Let us denote (the value on $Y_{2}$ of) the CRS for ${\cal S}_{i, Y_{1}}$ by $\varphi_{i}$. Then $\varphi =\bigcup_{i=1}^{n(M)}\varphi_{i} \in S(M)$ and each IN component for $({\cal S}_{i, Y_{1}})_{Y_{2}}$ is also a component of $M_{\varphi}$ and they are disjoint for $i$ (see the proof of Corollary \ref{main-important-cor}). 

It follows from Lemma \ref{main-geometric-lem} (iii) and the definition of $n(M)$ that for each $i$, the relation $({\cal S}_{i, Y_{1}})_{Y_{2}}$ has a unique IN component $Q_{i}$, and conversely that each component of $M_{\varphi}$ not of type $(0, 3)$ is an IN component $Q_{i}$ for some $({\cal S}_{i, Y_{1}})_{Y_{2}}$.

\begin{prop}\label{prop-posi-s}
Fix $i\in \{ 1, \ldots, n(M)\}$. Suppose that $Q_{i}$ is either of type $(0, 4)$ or $(1, 1)$. For each $s\in \{ 1, 2\}$, one of the following three cases occurs:
\begin{enumerate}
\item[(I)] There exists a curve $\beta_{i}^{s}$ in the CRS for $({\cal S}_{i}^{s})_{Y_{2}}$ which is also an element in $V(C(Q_{i}))$;
\item[(II)] The CRS for $({\cal S}_{i}^{s})_{Y_{2}}$ contains all curves in $\partial_{M}Q_{i}$ and $Q_{i}$ is an IA component for $({\cal S}_{i}^{s})_{Y_{2}}$;
\item[(III)] If $t\in \{ 1, 2\} \setminus \{ s\}$, then the CRS for $({\cal S}_{i}^{t})_{Y_{2}}$ contains all curves in $\partial_{M}Q_{i}$ and $Q_{i}$ is an IA component for $({\cal S}_{i}^{t})_{Y_{2}}$. 
\end{enumerate}
\end{prop}      

Recall that $\partial_{M}Q_{i}$ denotes the set of relatively boundary components (see Definition \ref{defn-geometric-things}). Remark that if $\kappa(M)$ is even, then the curve $\beta_{i}^{s}$ in the case (I) satisfies the hypothesis on $\alpha$ in Lemma \ref{lem-even-nn}. 

For the proof of this proposition, we need the following lemma:

\begin{lem}\label{lem-posi-s}
Fix $i\in \{ 1, \ldots, n(M)\}$. Suppose that $Q_{i}$ is either of type $(0, 4)$ or $(1, 1)$. For any Borel subset $A$ of $Y_{2}$ with positive measure, the IN system for ${\cal S}_{i, A}$ is constant on (a conull subset of) $A$ and its value consists of only the component $Q_{i}$.
\end{lem}

\begin{pf}
Let $A$ be a Borel subset of $Y_{2}$ with positive measure. Since ${\cal S}_{i, A}$ is non-amenable, the IN system for ${\cal S}_{i, A}$ is non-empty for a.e.\  point on $A$ by Theorem \ref{irr-non-ame-emp-ame}. It follows from Lemma \ref{in-contained-in} that for a.e.\ $x\in A$, each IN component of ${\cal S}_{i, A}$ for $x$ is contained in some IN component of $({\cal S}_{i, Y_{1}})_{A}$ for $x$ because ${\cal S}_{i, A}$ is a subrelation of $({\cal S}_{i, Y_{1}})_{A}$. On the other hand, the IN system for $({\cal S}_{i, Y_{1}})_{A}$ is constant on $A$ and its value consists of only the component $Q_{i}$. Since $Q_{i}$ is either of type $(0, 4)$ or $(1, 1)$, this shows the lemma.
\end{pf}

\begin{pf*}{{\sc Proof of Proposition \ref{prop-posi-s}}}
Let $s\in \{ 1, 2\}$ and denote by $\tau^{s}\in S(M)$ (the value on $Y_{2}$ of) the CRS for $({\cal S}_{i}^{s})_{Y_{2}}$. Note that any curve in $\tau^{s}$ does not intersect any curve in $\partial_{M}Q_{i}$. Otherwise, there would exist curves $\alpha \in \tau^{s}$ and $\beta \in \partial_{M}Q_{i}$ with $i(\alpha, \beta)\neq 0$. It contradicts the essentiality of the pair $(\alpha, Y_{2})$ for $({\cal S}_{i}^{s})_{Y_{2}}$ since $\beta$ is invariant for $({\cal S}_{i, Y_{1}})_{Y_{2}}$ and thus, for $({\cal S}_{i}^{s})_{Y_{2}}$. 

We have the following two possibilities:
\begin{enumerate}
\item[(i)] There exists a curve in $\tau^{s}$ which is also in $V(C(Q_{i}))$;
\item[(ii)] There exist no such curves.
\end{enumerate}
The case (i) is equal to the case (I). In the case (ii), there exists a unique component $R^{s}$ of $M_{\tau^{s}}$ with $Q_{i}\subseteq R^{s}$. Then we have the following possibilities:
\begin{enumerate}
\item[(iia)] $R^{s}$ is a T component for $({\cal S}_{i}^{s})_{Y_{2}}$;
\item[(iib)] $R^{s}$ is an IA component for $({\cal S}_{i}^{s})_{Y_{2}}$.
\end{enumerate}
In the case (iib), any curve in $\partial_{M}Q_{i}$ can not be in $V(C(R^{s}))$ since it is invariant for $({\cal S}_{i}^{s})_{Y_{2}}$. Thus, we see that $Q_{i}=R^{s}$ and the case (iib) is the case (II).

In the case (iia), there exists a Borel subset $Y_{2}'$ of $Y_{2}$ with positive measure such that 
\[\rho_{R^{s}}(x, y)\alpha =\alpha \]
for a.e.\ $(x, y)\in ({\cal S}_{i}^{s})_{Y_{2}'}$ and any $\alpha \in V(C(R^{s}))$ by Lemma \ref{lem:trivial}, where $\rho_{R^{s}}$ is the composition of the cocycle $\rho \colon {\cal R}\rightarrow \Gamma$ and the homomorphism $p_{R^{s}}\colon \Gamma_{\tau^{s}}\rightarrow \Gamma(R^{s})$. 

Let $t\in \{ 1, 2\} \setminus \{ s\}$ and denote by $\tau^{t}\in S(M)$ (the value on $Y_{2}$ of) the CRS for $({\cal S}_{i}^{t})_{Y_{2}}$. We can consider the cases (i)', (iia)' and (iib)' similarly for the relation $({\cal S}_{i}^{t})_{Y_{2}}$. We show that in the case (iia), only the case (iib)' can occur.

If the case (i)' occurred, then there would exist a curve $\beta_{1}\in \tau^{t}$ which is also in $V(C(Q_{i}))$. By the choice of $Y_{2}'$, the curve $\beta_{1}$ is invariant for $({\cal S}_{i}^{s})_{Y_{2}'}$. Thus, $\beta_{1}$ is invariant for the relation ${\cal S}_{i, Y_{2}'}=({\cal S}_{i}^{s})_{Y_{2}'}\vee ({\cal S}_{i}^{t})_{Y_{2}'}$. On the other hand, it follows from Lemma \ref{lem-posi-s} that the IN system for ${\cal S}_{i, Y_{2}'}$ is constant a.e.\ on $Y_{2}'$ and its value consists of only the component $Q_{i}$, which is a contradiction.

If the case (iia)' occurred, then there would exist a Borel subset $Y_{2}''$ of $Y_{2}'$ with positive measure such that
\[\rho_{R^{t}}(x, y)\alpha =\alpha \]
for a.e.\ $(x, y)\in ({\cal S}_{i}^{t})_{Y_{2}''}$ and any $\alpha \in V(C(R^{t}))$ by Lemma \ref{lem:trivial}. Let $\beta_{2}\in V(C(Q_{i}))$. By the choice of $Y_{2}''$, the curve $\beta_{2}$ is invariant for the relation ${\cal S}_{i, Y_{2}''}=({\cal S}_{i}^{s})_{Y_{2}''}\vee ({\cal S}_{i}^{t})_{Y_{2}''}$. We can deduce a contradiction as above.

Therefore, we have shown that in the case (iia), only the case (iib)' can happen. This case corresponds to the case (III).
\end{pf*}

We introduce the following new condition:

\begin{defn}\label{defn-dagger-n-wrt}
Let ${\cal R}$ be a discrete measured equivalence relation on a standard Borel space $(X, \mu)$ with a finite positive measure and $Y$ be a Borel subset of $X$ with positive measure. Let ${\cal S}_{1}, \ldots, {\cal S}_{m}$ be a finite family of recurrent subrelations on $Y$ of $({\cal R})_{Y}$. Given $n\in {\Bbb N}$ and a family
\[{\cal G}=\{ Y, \{ {\cal T}_{1}^{1}, {\cal T}_{1}^{2}\}, \ldots, \{ {\cal T}_{n}^{1}, {\cal T}_{n}^{2}\} \}\] 
of the Borel set $Y$ and subrelations of $({\cal R})_{Y}$ satisfying $(\dagger)_{n}$, we say that the family ${\cal G}$ satisfies $(\dagger)_{n}$ {\it with respect to} ${\cal S}_{1}, \ldots, {\cal S}_{m}$\index{$(\dagger n with respect to n $@$(\dagger)_{n}$ with respect to ${\cal S}_{1}, \ldots, {\cal S}_{m}$} if the following conditions are satisfied: 
\begin{enumerate}
\item[(i)] For any Borel subsets $A\subseteq B$ of $Y$ with positive measure, the relation $({\cal S}_{j})_{A}$ is normal in $({\cal S}_{j})_{A}\vee ({\cal T}_{i, B})_{A}$ for each $i\in \{ 1, \ldots, n\}$ and $j\in \{ 1, \ldots, m\}$.   
\item[(ii)] For any Borel subset $A$ of $Y$ with positive measure, the relation ${\cal T}_{i, A}$ is normal in $({\cal S}_{j})_{A}\vee {\cal T}_{i, A}$ for each $i\in \{ 1, \ldots, n\}$ and $j\in \{ 1,\ldots, m\}$.
\end{enumerate}
Here, we denote by ${\cal T}_{i, A}$ the relation $({\cal T}_{i}^{1})_{A}\vee ({\cal T}_{i}^{2})_{A}$ for $i\in \{ 1, \ldots, n\}$ and a Borel subset $A$ of $Y$. 
\end{defn}

\begin{rem}\label{rem-dagger-wrt-restriction}
As in Remark \ref{rem-dagger-n}, we give the following remarks: in the above definition, if a family 
\[\{ Y, \{ {\cal T}_{1}^{1}, {\cal T}_{1}^{2}\}, \ldots, \{ {\cal T}_{n}^{1}, {\cal T}_{n}^{2}\} \}\]
satisfies $(\dagger)_{n}$ with respect to ${\cal S}_{1}, \ldots, {\cal S}_{m}$, then for any Borel subset $A$ of $Y$ with positive measure, the family
\[\{ A, \{ ({\cal T}_{1}^{1})_{A}, ({\cal T}_{1}^{2})_{A}\}, \ldots, \{ ({\cal T}_{n}^{1})_{A}, ({\cal T}_{n}^{2})_{A}\} \}\]
satisfies $(\dagger)_{n}$ with respect to $({\cal S}_{1})_{A}, \ldots, ({\cal S}_{m})_{A}$.

Suppose that ${\cal R}$ (resp. ${\cal R}'$) is a discrete measured equivalence relation on $(X, \mu)$ (resp. $(X', \mu')$) and there exists a Borel isomorphism $f\colon X\rightarrow X'$ inducing the isomorphism between ${\cal R}$ and ${\cal R}'$. Let $Y$ be a Borel subset of $X$ with positive measure and ${\cal S}_{1},\ldots, {\cal S}_{m}$ be recurrent subrelations of $({\cal R})_{Y}$ on $Y$. If a family 
\[\{ Y, \{ {\cal T}_{1}^{1}, {\cal T}_{1}^{2}\}, \ldots, \{ {\cal T}_{n}^{1}, {\cal T}_{n}^{2}\}\}\]
of a Borel subset $Y$ of $X$ and subrelations of $({\cal R})_{Y}$ satisfies $(\dagger)_{n}$ with respect to ${\cal S}_{1}, \ldots, {\cal S}_{m}$, then the family 
\[\{ f(Y), \{ f({\cal T}_{1}^{1}), f({\cal T}_{1}^{2})\}, \ldots, \{ f({\cal T}_{n}^{1}), f({\cal T}_{n}^{2}\})\}\]
also satisfies $(\dagger)_{n}$ with respect to $f({\cal S}_{1}), \ldots, f({\cal S}_{m})$.

\end{rem}

Return to the situation before Proposition \ref{prop-posi-s}. We have a family
\[{\cal F}=\{ Y, \{ {\cal S}_{1}^{1}, {\cal S}_{1}^{2}\},\ldots, \{ {\cal S}_{n(M)}^{1}, {\cal S}_{n(M)}^{2}\} \}\]
of a Borel subset $Y$ of $X$ and subrelations of $({\cal R})_{Y}$ satisfying $(\dagger)_{n(M)}$. Moreover, we have chosen Borel subsets $Y_{2}\subseteq Y_{1}$ of $Y$ with positive measure. Recall that $Q_{i}$ denotes the unique IN component for $({\cal S}_{i, Y_{1}})_{Y_{2}}$ for $i\in \{ 1, \ldots, n(M)\}$. 

Suppose that distinct elements $i_{1},\ldots, i_{m}, j_{1},\ldots, j_{m'}\in \{ 1, 2, \ldots, n(M)\}$ satisfy that each of $Q_{i_{1}}, \ldots, Q_{i_{m}}, Q_{j_{1}},\ldots, Q_{j_{m'}}$ is either of type $(0, 4)$ or $(1, 1)$. The numbers $m$ and $m'$ may be zero. Assume that there exist $s_{1},\ldots, s_{m}, t_{1},\ldots, t_{m'}\in \{ 1, 2\}$ such that $({\cal S}_{i_{1}}^{s_{1}})_{Y_{2}}, \ldots, ({\cal S}_{i_{m}}^{s_{m}})_{Y_{2}}$ satisfy the condition (I) in Proposition \ref{prop-posi-s} and $({\cal S}_{j_{1}}^{t_{1}})_{Y_{2}}, \ldots, ({\cal S}_{j_{m'}}^{t_{m'}})_{Y_{2}}$ satisfy the condition (II) in the same proposition. Let $\beta_{i_{k}}^{s_{k}}\in V(C(Q_{i_{k}}))$ be a curve in the CRS for $({\cal S}_{i_{k}}^{s_{k}})_{Y_{2}}$.

Suppose that we have a family 
\[\{ Y, \{ {\cal T}_{1}^{1}, {\cal T}_{1}^{2}\}, \ldots, \{ {\cal T}_{n}^{1}, {\cal T}_{n}^{2}\} \}\] 
of the Borel set $Y$ and subrelations of $({\cal R})_{Y}$ satisfying $(\dagger)_{n}$ with respect to the relations ${\cal S}_{i_{1}}^{s_{1}}, \ldots, {\cal S}_{i_{m}}^{s_{m}}, {\cal S}_{j_{1}}^{t_{1}},\ldots, {\cal S}_{j_{m'}}^{t_{m'}}$.

Then it follows from Corollary \ref{main-important-cor} (i) that there exists a Borel subset $Y_{3}$ of $Y_{2}$ with positive measure such that for any distinct $u, v\in \{ 1, \ldots, n\}$, the IN systems for ${\cal T}_{u, Y_{3}}$ and ${\cal T}_{v, Y_{3}}$ are disjoint a.e.\ on $Y_{3}$. Moreover, there exists a Borel subset $Y_{4}$ of $Y_{3}$ with positive measure such that for each $u$, the CRS, T, IA and IN systems for ${\cal T}_{u, Y_{3}}$ are constant on $Y_{4}$. Let us denote the values on $Y_{4}$ of the CRS and IN system for ${\cal T}_{u, Y_{3}}$ by $\Phi^{(u)}$ and $\Psi^{(u)}$, respectively. Then the union $\Phi =\bigcup_{u=1}^{n}\Phi^{(u)}$ is an element in $S(M)$ and each component of $\Psi^{(u)}$ is also a component of $M_{\Phi}$ (see the proof of Corollary \ref{main-important-cor}). We denote 
\[\Phi'=\{ \beta_{i_{1}}^{s_{1}},\ldots, \beta_{i_{m}}^{s_{m}}\},\ \ \Phi''=\bigcup_{l=1}^{m'}\partial_{M}Q_{j_{l}}, \ \ \Psi=\bigcup_{u=1}^{n}\Psi^{(u)}.\]
Then $\Phi', \Phi''\in S(M)$ and $\Psi \in {\cal F}_{0}(D(M))$ since any components in $\Psi^{(u)}$ and $\Psi^{(v)}$ are disjoint for all distinct $u, v\in \{ 1, \ldots, n\}$.

\begin{thm}\label{thm-general-class}
In the above notation, 
\begin{enumerate}
\item[(i)] each curve in $\Phi$ neither intersects $\beta_{i_{k}}^{s_{k}}$ nor curves in $\partial_{M}Q_{j_{l}}$ for any $k\in \{ 1, \ldots, m\}$ and $l\in \{ 1, \ldots, m'\}$. In particular, if we put $\Phi_{0}=\Phi \cup \Phi'\cup \Phi''$, then $\Phi_{0}\in S(M)$.
\item[(ii)] we have $r(\beta_{i_{k}}^{s_{k}}, Q_{j_{l}})=\emptyset$ for any $k\in \{ 1, \ldots, m\}$ and $l\in \{ 1, \ldots, m'\}$.
\item[(iii)] we have $r(\beta_{i_{k}}^{s_{k}}, Q)=\emptyset$ for any $Q\in \Psi$ and $k\in \{ 1, \ldots, m\}$.
\item[(iv)] for any $u\in \{ 1, \ldots, n\}$ and $l\in \{ 1, \ldots, m'\}$, the component $Q_{j_{l}}$ is a component of $M_{\Phi_{0}}$ and is contained in some component of $M_{\Phi^{(u)}}$ which is not in $\Psi^{(u)}$.
\item[(v)] each $Q\in \Psi$ is a component of $M_{\Phi_{0}}$.
\item[(vi)] we have $n\leq n(\Phi_{0})-m'$.
\end{enumerate}
\end{thm}

\begin{pf}
Since $({\cal S}_{i_{k}}^{s_{k}})_{Y_{3}}$ is normal in $({\cal S}_{i_{k}}^{s_{k}})_{Y_{3}}\vee {\cal T}_{u, Y_{3}}$, the curve $\beta_{i_{k}}^{s_{k}}$ is invariant for ${\cal T}_{u, Y_{3}}$. Hence, we have $r(\beta_{i_{k}}^{s_{k}}, Q)=\emptyset$ for any component $Q$ in $\Psi^{(u)}$. Otherwise, $Q$ would be a T component for $({\cal T}_{u, Y_{3}})_{Y_{4}}$ by Lemma \ref{lem:trivial}, which contradicts $Q\in \Psi^{(u)}$. Moreover, since $\beta_{i_{k}}^{s_{k}}$ is invariant for $({\cal T}_{u, Y_{3}})_{Y_{4}}$, any curve in $\Phi^{(u)}$ does not intersect $\beta_{i_{k}}^{s_{k}}$ by the essentiality of the curves in $\Phi^{(u)}$ for $({\cal T}_{u, Y_{3}})_{Y_{4}}$. 

On the other hand, since $({\cal S}_{j_{l}}^{t_{l}})_{Y_{3}}$ is normal in $({\cal S}_{j_{l}}^{t_{l}})_{Y_{3}}\vee {\cal T}_{u, Y_{3}}$, we can show as above that any curve in $\Phi^{(u)}$ does not intersect any relative boundary component of $Q_{j_{l}}$ which is contained in the CRS for $({\cal S}_{j_{l}}^{t_{l}})_{Y_{3}}$. This proves the assertion (i).

We have shown the assertion (iii). The assertion (ii) is clear because $i_{k}, j_{l}\in \{ 1, \ldots, n(M)\}$ are distinct and $Q_{i_{k}}$ and $Q_{j_{l}}$ are disjoint. 

Fix $u\in \{ 1, \ldots, n\}$ and $l\in \{ 1, \ldots, m'\}$. Since ${\cal T}_{u, Y_{3}}$ is normal in ${\cal T}_{u, Y_{3}}\vee ({\cal S}_{j_{l}}^{t_{l}})_{Y_{3}}$, the CRS of ${\cal T}_{u, Y_{3}}$ is invariant for $({\cal S}_{j_{l}}^{t_{l}})_{Y_{3}}$. In particular, any curve in $\Phi^{(u)}$ is invariant for $({\cal S}_{j_{l}}^{t_{l}})_{Y_{4}}$. It follows that any curve in $\Phi^{(u)}$ cannot be an element in $V(C(Q_{j_{l}}))$ because $Q_{j_{l}}$ is an IA component for $({\cal S}_{j_{l}}^{t_{l}})_{Y_{4}}$ by the condition (II) in Proposition \ref{prop-posi-s}. Thus, there exists a unique component $R_{l}^{u}$ of $M_{\Phi^{(u)}}$ with $Q_{j_{l}}\subseteq R_{l}^{u}$. Since $u\in \{ 1,\ldots, n\}$ is arbitrary, it follows that $Q_{j_{l}}$ is contained in some component of $M_{\Phi}$ and is a component of $M_{\Phi_{0}}$ by the assertion (ii).

We show that $R_{l}^{u}$ is not a component in $\Psi^{(u)}$ for any $u\in \{ 1, \ldots, n\}$ and $l\in \{ 1, \ldots, m'\}$, which proves the assertion (iv). Suppose $Q_{j_{l}}=R^{u}_{l}$. Since $({\cal S}_{j_{l}}^{t_{l}})_{Y_{3}}$ is normal in $({\cal S}_{j_{l}}^{t_{l}})_{Y_{3}}\vee {\cal T}_{u, Y_{3}}$ and $Q_{j_{l}}$ is an IA component for $({\cal S}_{j_{l}}^{t_{l}})_{Y_{3}}$, it follows from Lemma \ref{main-lem1'} (iii) that $Q_{j_{l}}$ is also an IA component for $({\cal S}_{j_{l}}^{t_{l}})_{Y_{3}}\vee {\cal T}_{u, Y_{3}}$ a.e.\ on $Y_{3}$. Thus, $R^{u}_{l}$ is not a component in $\Psi^{(u)}$ by Corollary \ref{cor-simple}.

Next, assume $Q_{j_{l}}\subsetneq R_{l}^{u}$. Then there exists a curve in $\partial_{M}Q_{j_{l}}$ which is also an element in $V(C(R_{l}^{u}))$. The curve is invariant for $({\cal T}_{u, Y_{3}})_{Y_{4}}$ because $({\cal S}_{j_{l}}^{t_{l}})_{Y_{4}}$ is normal in $({\cal S}_{j_{l}}^{t_{l}})_{Y_{4}}\vee ({\cal T}_{u, Y_{3}})_{Y_{4}}$. Hence, $R^{u}_{l}$ is not a component in $\Psi^{(u)}$. This completes the proof of the assertion (iv).

The assertion (v) comes from the assertions (iii), (iv) because each $Q\in \Psi$ is a component of $M_{\Phi}$. The inequality (vi) comes from the assertions (iv), (v) and the fact that any components in $\Psi^{(u)}$ and $\Psi^{(v)}$ are disjoint for all distinct $u, v\in \{ 1, \ldots, n\}$. 
\end{pf}

We apply this theorem to various situations as follows:

\begin{thm}\label{thm-even-n-1}
With the assumption $(\diamond)'$, suppose that $\kappa(M)$ is even and a family
\[{\cal F}=\{ Y, \{ {\cal S}_{1}^{1}, {\cal S}_{1}^{2}\},\ldots, \{ {\cal S}_{n(M)}^{1}, {\cal S}_{n(M)}^{2}\} \}\]
of a Borel subset $Y$ of $X$ and subrelations of $({\cal R})_{Y}$ satisfies $(\dagger)_{n(M)}$. Then for any $i\in \{ 1, \ldots, n(M)\}$, there exists $s\in \{ 1, 2\}$ satisfying the following condition: if a family 
\[{\cal G}=\{ Y, \{ {\cal T}_{1}^{1}, {\cal T}_{1}^{2}\}, \ldots, \{ {\cal T}_{n}^{1}, {\cal T}_{n}^{2}\} \}\] 
of the Borel set $Y$ and subrelations of $({\cal R})_{Y}$ satisfies $(\dagger)_{n}$ with respect to the relation ${\cal S}_{i}^{s}$, then $n\leq n(M)-1$.
\end{thm}

\begin{pf}
We take Borel subsets $Y_{2}\subseteq Y_{1}$ of $Y$ with positive measure for the family ${\cal F}$ as in the beginning of this section. Fix $i\in \{ 1, \ldots, n(M)\}$. Remark that the IN component $Q_{i}$ for $({\cal S}_{i, Y_{1}})_{Y_{2}}$ is either of type $(0, 4)$ or $(1, 1)$ by Lemma \ref{lem-even-n} (i). Apply Proposition \ref{prop-posi-s} to the relations $({\cal S}_{i}^{1})_{Y_{2}}$ and $({\cal S}_{i}^{2})_{Y_{2}}$. 

Suppose that $({\cal S}_{i}^{1})_{Y_{2}}$ satisfies the condition (I). It follows from Theorem \ref{thm-general-class} (vi) that given any family ${\cal G}$ satisfying $(\dagger)_{n}$ with respect to ${\cal S}_{i}^{1}$, we see that $n\leq n(\Phi \cup \{ \beta_{i}^{1}\})$. It follows from Lemma \ref{lem-even-nn} that $n(\Phi \cup \{ \beta_{i}^{1}\})\leq n(M)-1$ and thus, $n\leq n(M)-1$.

Suppose that $({\cal S}_{i}^{1})_{Y_{2}}$ satisfies the condition (II). It follows from Theorem \ref{thm-general-class} (vi) that given any family ${\cal G}$ satisfying $(\dagger)_{n}$ with respect to ${\cal S}_{i}^{1}$, we see that $n\leq n(M)-1$. 

If $({\cal S}_{i}^{1})_{Y_{2}}$ satisfies the condition (III), then the relation $({\cal S}_{i}^{2})_{Y_{2}}$ satisfies the condition (II) and we can apply Theorem \ref{thm-general-class} (vi) similarly. 
\end{pf}

\begin{thm}\label{thm-odd-n-1}
With the assumption $(\diamond)'$, suppose that $\kappa(M)$ is odd and $g\geq 3$. Assume that we have a family
\[{\cal F}=\{ Y, \{ {\cal S}_{1}^{1}, {\cal S}_{1}^{2}\},\ldots, \{ {\cal S}_{n(M)}^{1}, {\cal S}_{n(M)}^{2}\} \}\]
of a Borel subset $Y$ of $X$ and subrelations of $({\cal R})_{Y}$ satisfying $(\dagger)_{n(M)}$. Then after permuting the indices $1, \ldots, n(M)$, we see that for any distinct $i, j\in \{ 1, \ldots, g-1\}$, there exist $s, t\in \{ 1, 2\}$ satisfying the following condition: if a family 
\[{\cal G}=\{ Y, \{ {\cal T}_{1}^{1}, {\cal T}_{1}^{2}\}, \ldots, \{ {\cal T}_{n}^{1}, {\cal T}_{n}^{2}\} \}\] 
of the Borel set $Y$ and subrelations of $({\cal R})_{Y}$ satisfies $(\dagger)_{n}$ with respect to the relations ${\cal S}_{i}^{s}$ and ${\cal S}_{j}^{t}$, then $n\leq n(M)-1$.
\end{thm}

\begin{pf}
We use the same notation as in the beginning of this section. Since the CRS, T, IA and IN systems for $({\cal S}_{i, Y_{1}})_{Y_{2}}$ are constant for any $i\in \{ 1, \ldots, n(M)\}$, we can rearrange the indices $1, \ldots, n(M)$ so that the unique IN component $Q_{i}$ for $({\cal S}_{i, Y_{1}})_{Y_{2}}$ is of type $(1, 1)$ for any $i\in \{ 1, \ldots, g-1\}$. Remark that for any $\sigma \in S(M)$ with $n(\sigma)=n(M)$, there exist at least $g-1$ components of $M_{\sigma}$ of type $(1, 1)$ by Lemma \ref{lem-odd-n} (ii).  

Let $i, j\in \{ 1, \ldots, g-1\}$ be any distinct elements. We apply Proposition \ref{prop-posi-s} to the four relations $({\cal S}_{i}^{1})_{Y_{2}}$, $({\cal S}_{i}^{2})_{Y_{2}}$, $({\cal S}_{j}^{1})_{Y_{2}}$ and $({\cal S}_{j}^{2})_{Y_{2}}$. 

Suppose that one of the four relations, say $({\cal S}_{i}^{1})_{Y_{2}}$, satisfies the condition (II). It follows from Theorem \ref{thm-general-class} (vi) that given any family ${\cal G}$ satisfying $(\dagger)_{n}$ with respect to ${\cal S}_{i}^{1}$ and ${\cal S}_{j}^{1}$, we see that $n\leq n(M)-1$.

Suppose that there exist $s, t\in \{ 1, 2\}$ such that both $({\cal S}_{i}^{s})_{Y_{2}}$ and $({\cal S}_{j}^{t})_{Y_{2}}$ satisfy the condition (I). It follows from Theorem \ref{thm-general-class} (vi) that given any family ${\cal G}$ satisfying $(\dagger)_{n}$ with respect to ${\cal S}_{i}^{s}$ and ${\cal S}_{j}^{t}$, we see that $n\leq n(\Phi \cup \{ \beta_{i}^{s}, \beta_{j}^{t}\})$. When we cut the surface $M$ along the curves $\beta_{i}^{s}$ and $\beta_{j}^{t}$, the genus of the resulting surface is $g-2$ since both $Q_{i}$ and $Q_{j}$ are of type $(1, 1)$. Hence, the equality $n(\Phi \cup \{ \beta_{i}^{s}, \beta_{j}^{t}\})=n(M)$ cannot occur by Lemma \ref{lem-odd-n} (ii). It means that $n(\Phi \cup \{ \beta_{i}^{s}, \beta_{j}^{t}\})\leq n(M)-1$ and thus, $n\leq n(M)-1$. 
\end{pf}

\begin{thm}\label{thm-even-n-2}
With the assumption $(\diamond)'$, suppose that $\kappa(M)$ is even and $g\geq 3$. Assume that we have a family
\[{\cal F}=\{ Y, \{ {\cal S}_{1}^{1}, {\cal S}_{1}^{2}\},\ldots, \{ {\cal S}_{n(M)}^{1}, {\cal S}_{n(M)}^{2}\} \}\]
of a Borel subset $Y$ of $X$ and subrelations of $({\cal R})_{Y}$ satisfying $(\dagger)_{n(M)}$. Then after permuting the indices $1, \ldots, n(M)$, we see that for any distinct $i, j, k\in \{ 1, \ldots, g\}$, there exist $s, t, u\in \{ 1, 2\}$ satisfying the following condition: if a family
\[{\cal G}=\{ Y, \{ {\cal T}_{1}^{1}, {\cal T}_{1}^{2}\}, \ldots, \{ {\cal T}_{n}^{1}, {\cal T}_{n}^{2}\} \}\] 
of the Borel set $Y$ and subrelations of $({\cal R})_{Y}$ satisfies $(\dagger)_{n}$ with respect to the relations ${\cal S}_{i}^{s}$, ${\cal S}_{j}^{t}$ and ${\cal S}_{k}^{u}$, then $n\leq n(M)-2$.
\end{thm}

\begin{pf}
We use the same notation as in the beginning of this section. As in the proof of Theorem \ref{thm-odd-n-1}, we rearrange the indices $1, \ldots, n(M)$ so that the IN component $Q_{i}$ for $({\cal S}_{i, Y_{1}})_{Y_{2}}$ is of type $(1, 1)$ for any $i\in \{ 1, \ldots, g\}$. Remark that for any $\sigma \in S(M)$ with $n(\sigma)=n(M)$, there exist $g$ components of $M_{\sigma}$ of type $(1, 1)$ by Lemma \ref{lem-even-n} (i). 

Let $i, j, k\in \{ 1, \ldots, g\}$ be any distinct elements. We apply Proposition \ref{prop-posi-s} to the six relations $({\cal S}_{i}^{s})_{Y_{2}}$, $({\cal S}_{j}^{s})_{Y_{2}}$, $({\cal S}_{k}^{s})_{Y_{2}}$ for $s\in \{ 1, 2\}$. 

Suppose that there exist two distinct elements in $\{ i, j, k\}$, say $i$ and $j$, and $s, t\in \{ 1, 2\}$ such that both $({\cal S}_{i}^{s})_{Y_{2}}$ and $({\cal S}_{j}^{t})_{Y_{2}}$ satisfy the condition (II). It follows from Theorem \ref{thm-general-class} (vi) that given any family ${\cal G}$ satisfying $(\dagger)_{n}$ with respect to ${\cal S}_{i}^{s}$, ${\cal S}_{j}^{t}$ and ${\cal S}_{k}^{1}$, we see that $n\leq n(M)-2$.

The remaining cases are the following two ones:
\begin{enumerate}
\item[(a)] We can rearrange the indices $i, j, k$ so that there exist $s, t, u\in \{ 1, 2\}$ such that both $({\cal S}_{i}^{s})_{Y_{2}}$ and $({\cal S}_{j}^{t})_{Y_{2}}$ satisfy the condition (I) and $({\cal S}_{k}^{u})_{Y_{2}}$ satisfies the condition (II);
\item[(b)] There exist $s, t, u\in \{ 1, 2\}$ such that all the relations $({\cal S}_{i}^{s})_{Y_{2}}$, $({\cal S}_{j}^{t})_{Y_{2}}$ and $({\cal S}_{k}^{u})_{Y_{2}}$ satisfy the condition (I).
\end{enumerate}
In the case (a), given any family ${\cal G}$ satisfying $(\dagger)_{n}$ with respect to ${\cal S}_{i}^{s}$, ${\cal S}_{j}^{t}$ and ${\cal S}_{k}^{u}$, we see that 
\[n\leq n(\Phi \cup \{ \beta_{i}^{s}, \beta_{j}^{t}\} \cup \partial_{M}Q_{k})-1\leq n(M)-2\]
by Lemma \ref{lem-even-nn} and Theorem \ref{thm-general-class} (vi).

In the case (b), given any family ${\cal G}$ satisfying $(\dagger)_{n}$ with respect to ${\cal S}_{i}^{s}$, ${\cal S}_{j}^{t}$ and ${\cal S}_{k}^{u}$, we see that 
\[n\leq n(\Phi \cup \{ \beta_{i}^{s}, \beta_{j}^{t}, \beta_{k}^{u}\})\leq n(M_{\{ \beta_{i}^{s}, \beta_{j}^{t}, \beta_{k}^{u}\}})\]
by Theorem \ref{thm-general-class} (vi), where for a disconnected surface $M'$, we denote by $n(M')$ the sum of $n(Q)$ for all components $Q$ of $M'$. If the surface $M$ is cut along the curve $\beta_{i}^{s}$, then the resulting surface $M_{\beta_{i}^{s}}$ is connected since $Q_{i}$ is of type $(1, 1)$. Moreover, the complexity $\kappa(M_{\beta_{i}^{s}})$ is odd. Thus,  
\[n(M_{\{ \beta_{i}^{s}, \beta_{j}^{t}, \beta_{k}^{u}\}})\leq n((M_{\beta_{i}^{s}})_{\{ \beta_{j}^{t}, \beta_{k}^{u}\}})\leq n(M_{\beta_{i}^{s}})-1\leq n(M)-2\] 
by Lemma \ref{lem-odd-n} (ii) and Lemma \ref{lem-even-nn}. Hence, we have $n\leq n(M)-2$.    
\end{pf}


\section[Application I]{Application I (Invariance of complexity under measure equivalence)}\label{sec-app1}

In this section, we apply Theorem \ref{thm-even-n-1} and show that complexity is invariant under measure equivalence among the mapping class groups. For the proof, we need the following lemmas. In this section, let us denote by $t(\alpha)\in \Gamma(M)$ the Dehn twist about $\alpha$ for $\alpha \in V(C(M))$.

\begin{lem}\label{lem-free-generate-curves}
Let $M$ be a surface with $\kappa(M)\geq 0$ and $m\in {\Bbb N}$ be a natural number more than or equal to $3$. Then there exists $N\in {\Bbb N}$ satisfying the following condition: let $Q$ be a subsurface of $M$ of type $(0, 4)$ or $(1, 1)$ and let $\alpha, \beta \in V(C(Q))$ be curves with $i(\alpha, \beta)\neq 0$. Then both $t(\alpha)^{N}$ and $t(\beta)^{N}$ are in $\Gamma(M; m)$ and they generate the free subgroup of rank $2$ of $\Gamma(M;m)$.
\end{lem}

Remark that for any $k\in {\Bbb N} \setminus \{ 0\}$, the number $Nk$ also satisfies the condition in this lemma.

\begin{pf}
There exists $N'\in {\Bbb N}$ such that $g^{N'}\in \Gamma(M; m)$ for all $g\in \Gamma(M)$. Since both $t(\alpha)$ and $t(\beta)$ preserve each curve in $\partial_{M}Q$ and preserve $Q$, we have a homomorphism
\[p\colon G\rightarrow \Gamma(Q), \]
where $G$ denotes the subgroup of $\Gamma(M)$ generated by $t(\alpha)$ and $t(\beta)$. Note that $\Gamma(Q)$ has a subgroup of finite index isomorphic to the free group of rank $2$. Remark that in general, if $H_{1}$ is a group and $H_{2}$ is a subgroup of $H_{1}$ of finite index, then the kernel of the action of $H_{1}$ on $H_{1}/H_{2}$ is a normal subgroup of $H_{1}$ of finite index contained in $H_{2}$. Fix a free normal subgroup $G_{1}$ (resp. $G_{2}$) of $\Gamma(M_{0, 4})$ (resp. $\Gamma(M_{1, 1})$) of finite index. Then there exists $N''\in {\Bbb N}$ such that $g^{N''}\in G_{1}$ (resp. $G_{2}$) for any element $g\in \Gamma(M_{0, 4})$ (resp. $\Gamma(M_{1, 1})$). It follows that $p(t(\alpha)^{N''})$ and $p(t(\beta)^{N''})$ generate a free subgroup of rank $2$ in $\Gamma(Q)$ since $\alpha$ and $\beta$ satisfy $i(\alpha, \beta)\neq 0$ (see Theorem \ref{commuting-dehn}). Thus, $t(\alpha)^{N''}$ and $t(\beta)^{N''}$ generate a free subgroup of rank $2$ in $\Gamma(M)$. Hence, $N=N'N''$ is a desired number.
\end{pf}

\begin{lem}\label{lem-family-good}
Let $N\in {\Bbb N}$ be the natural number in Lemma \ref{lem-free-generate-curves} and let $\sigma \in S(M)$ be an element such that $M_{\sigma}$ has $n$ components $R_{1}, \ldots, R_{n}$ each of which is either of type $(0, 4)$ or $(1, 1)$. Let $\sigma'=\{ \alpha_{1},\ldots, \alpha_{p}\} \subseteq \sigma$ be a subset (which may be empty). Suppose that we have two intersecting curves $\beta_{k}^{1}, \beta_{k}^{2}\in V(C(R_{k}))$ for each $k\in \{ 1,\ldots, n\}$.

Assume $(\diamond)'$ and that $F=X$, $\Gamma =\Gamma(M;m)$ and the action of $\Gamma$ on $X$ is measure-preserving. We denote by ${\cal S}_{i}$ and ${\cal T}_{k}^{s}$ the subrelations of ${\cal R}$ generated by the actions of $t(\alpha_{i})^{N}$ and $t(\beta_{k}^{s})^{N}$ for $i\in \{ 1,\ldots, p\}$, $k\in \{ 1,\ldots, n\}$ and $s\in \{ 1, 2\}$, respectively.  
\begin{enumerate}
\item[(i)] The family 
\[{\cal G}=\{ X, \{ {\cal T}_{1}^{1}, {\cal T}_{1}^{2}\}, \ldots, \{ {\cal T}_{n}^{1}, {\cal T}_{n}^{2}\} \}\] 
satisfies $(\dagger)_{n}$. 
\item[(ii)] The family ${\cal G}$ satisfies $(\dagger)_{n}$ with respect to ${\cal S}_{1}, \ldots, {\cal S}_{p}$ if $\sigma'$ is non-empty.
\end{enumerate}
\end{lem}

\begin{pf}
Note that $t(\beta_{k}^{s})$ and $t(\beta_{k'}^{s'})$ are commutative in $\Gamma(M;m)$ for any distinct $k, k'\in \{ 1, \ldots, n_{0}\}$ and any $s, s'\in \{ 1, 2\}$ and that $t(\beta_{k}^{1})^{N}$ and $t(\beta_{k}^{2})^{N}$ generate the free subgroup of rank $2$ for any $k$. It follows from Lemmas \ref{lem-main-free-groups-1}, \ref{lem-main-free-groups-2} and \ref{lem-main-free-groups-3} that the family ${\cal G}$ satisfies $(\dagger)_{n}$. Moreover, since $t(\alpha_{i})$ and $t(\beta_{k}^{s})$ are commutative for any $i\in \{ 1,\ldots, p\}$, $k\in \{ 1,\ldots, n\}$ and $s\in \{ 1, 2\}$, it follows from Lemmas \ref{lem-main-free-groups-2} and \ref{lem-main-free-groups-3} that the family ${\cal G}$ satisfies $(\dagger)_{n}$ with respect to the relations ${\cal S}_{1}, \ldots, {\cal S}_{p}$. 
\end{pf}

\begin{thm}\label{thm-kappa-invariant-odd}
With the assumption $(\diamond)'$, suppose that $\kappa(M)$ is odd, $F=X$, $\Gamma =\Gamma(M; m)$ and the action of $\Gamma$ on $X$ is measure-preserving. Let us denote by $n_{0}=n(M)$. Then there exists a family
\[\{ X, \{ {\cal S}_{1}^{1}, {\cal S}_{1}^{2}\}, \ldots, \{ {\cal S}_{n_{0}}^{1}, {\cal S}_{n_{0}}^{2}\} \}\]
of the Borel set $X$ and subrelations of ${\cal R}$ satisfying $(\dagger)_{n_{0}}$ and the following condition: for any $i\in \{ 1, \ldots, n_{0}\}$ and $s\in \{ 1, 2\}$, there exists a family
\[\{ X, \{ {\cal T}_{1}^{1}, {\cal T}_{1}^{2}\}, \ldots, \{ {\cal T}_{n_{0}}^{1}, {\cal T}_{n_{0}}^{2}\} \}\]
satisfying $(\dagger)_{n_{0}}$ with respect to the relation ${\cal S}_{i}^{s}$.
\end{thm}

\begin{pf}
Let $N\in {\Bbb N}$ be the natural number in Lemma \ref{lem-free-generate-curves} and $\sigma \in S(M)$ be an element such that the numbers of components of $M_{\sigma}$ of types $(1, 1)$, $(0, 4)$ and $(0, 3)$ are $g$, $(g+p-3)/2$ and $1$, respectively (see Remark \ref{rem-lem-odd-n}). We denote by $Q_{1}, \ldots, Q_{n_{0}}$ the components of $M_{\sigma}$ either of type $(1, 1)$ or $(0, 4)$. Choose two intersecting curves $\beta_{i}^{1}, \beta_{i}^{2}\in V(C(Q_{i}))$ for each $i\in \{ 1, \ldots, n_{0}\}$. Let ${\cal S}_{i}^{s}$ be the subrelation of ${\cal R}$ generated by the action of $t(\beta_{i}^{s})^{N}$. Then the family
\[\{ X, \{ {\cal S}_{1}^{1}, {\cal S}_{1}^{2}\}, \ldots, \{ {\cal S}_{n_{0}}^{1}, {\cal S}_{n_{0}}^{2}\} \}\]
satisfies $(\dagger)_{n_{0}}$ by Lemma \ref{lem-family-good} (i). We show that this family is a desired one.

Take any $i\in \{ 1,\ldots, n(M)\}$ and $s\in \{ 1, 2\}$. It follows from Lemma \ref{lem-odd-nn} that there exists $\tau \in S(M)$ such that $\beta_{i}^{s}\in \tau$ and $n(\tau)=n(M)$. If $M_{\tau}$ has a component either of type $(0, 5)$ or $(1, 2)$, then add a certain curve $\alpha$ on the component to $\tau$ so that $M_{\tau \cup \{ \alpha \}}$ has components of only types $(0, 3)$, $(0, 4)$ and $(1, 1)$ (see Lemma \ref{lem-odd-n} (ii)). Hence, we may assume that $M_{\tau}$ has components of only types $(0, 3)$, $(0, 4)$ and $(1, 1)$. 

Let $R_{1}, \ldots, R_{n_{0}}$ be all components of $M_{\tau}$ either of type $(0, 4)$ or $(1, 1)$. As above, for $k\in \{ 1, \ldots, n_{0}\}$, choose two intersecting curves $\gamma_{k}^{t}\in V(C(R_{k}))$ for $t\in \{ 1, 2\}$. Let ${\cal T}_{k}^{t}$ be the subrelation of ${\cal R}$ generated by the action of $t(\gamma_{k}^{t})^{N}$. Then the family 
\[\{ X, \{ {\cal T}_{1}^{1}, {\cal T}_{1}^{2}\}, \ldots, \{ {\cal T}_{n_{0}}^{1}, {\cal T}_{n_{0}}^{2}\} \}\]
satisfies $(\dagger)_{n_{0}}$ with respect to the relation ${\cal S}_{i}^{s}$ by Lemma \ref{lem-family-good} (ii).  
\end{pf}

\begin{cor}\label{cor-kappa-invariant}
Let $M^{1}$ and $M^{2}$ be surfaces with $\kappa(M^{1}), \kappa(M^{2})\geq 0$. Suppose that the mapping class groups $\Gamma(M^{1})$ and $\Gamma(M^{2})$ are measure equivalent. Then $\kappa(M^{1})=\kappa(M^{2})$. 
\end{cor}

\begin{pf}
By Corollary \ref{cor-classification1}, we have $n_{0}=n(M^{1})=n(M^{2})$. Since in general, we have $n(M)=[(\kappa(M)-2)/2]$ for a surface $M$, it suffices to show that if $\kappa(M^{1})$ is even and $\kappa(M^{2})$ is odd, then the mapping class groups $\Gamma(M^{1})$ and $\Gamma(M^{2})$ are not measure equivalent.

If $\Gamma(M^{1})$ and $\Gamma(M^{2})$ were measure equivalent, then we would have two equivalence relations ${\cal R}_{1}$ on $(X_{1}, \mu_{1})$ and ${\cal R}_{2}$ on $(X_{2}, \mu_{2})$ of type ${\rm II}_{1}$ satisfying the following conditions: the relation ${\cal R}_{j}$ is generated by an essentially free, measure-preserving Borel action of $\Gamma(M^{j};m)$ for $j=1, 2$, where $m\geq 3$ is an integer. Moreover, ${\cal R}_{1}$ and ${\cal R}_{2}$ are weakly isomorphic. Let $f\colon Y_{1}\rightarrow Y_{2}$, $Y_{i}\subseteq X_{i}$ be a partial Borel isomorphism inducing the weak isomorphism between ${\cal R}_{1}$ and ${\cal R}_{2}$.

It follows from Remark \ref{rem-dagger-wrt-restriction} and Theorem \ref{thm-kappa-invariant-odd} that we have a family  
\[\{ Y_{2}, \{ {\cal S}_{1}^{1}, {\cal S}_{1}^{2}\}, \ldots, \{ {\cal S}_{n_{0}}^{1}, {\cal S}_{n_{0}}^{2}\} \}\]
of the Borel set $Y_{2}$ and subrelations of $({\cal R}_{2})_{Y_{2}}$ satisfying $(\dagger)_{n_{0}}$ and the following condition: for any $i\in \{ 1, \ldots, n_{0}\}$ and $s\in \{ 1, 2\}$, there exists a family
\[\{ Y_{2}, \{ {\cal T}_{1}^{1}, {\cal T}_{1}^{2}\}, \ldots, \{ {\cal T}_{n_{0}}^{1}, {\cal T}_{n_{0}}^{2}\} \}\]
satisfying $(\dagger)_{n_{0}}$ with respect to the relation ${\cal S}_{i}^{s}$.

On the other hand, since the family 
\[\{ Y_{1}, \{ f^{-1}({\cal S}_{1}^{1}), f^{-1}({\cal S}_{1}^{2})\}, \ldots, \{ f^{-1}({\cal S}_{n_{0}}^{1}), f^{-1}({\cal S}_{n_{0}}^{2})\} \}\]
satisfies $(\dagger)_{n_{0}}$, it follows from Theorem \ref{thm-even-n-1} that there exists $s\in \{ 1, 2\}$ satisfying the following condition: if we have a family of the Borel subset $f^{-1}(Y_{2})$ and subrelations of $({\cal R}_{1})_{f^{-1}(Y_{2})}$ satisfying $(\dagger)_{n}$ with respect to the relation $f^{-1}({\cal S}_{1}^{s})$, then $n\leq n_{0}-1$. It is a contradiction since we have the family
\[\{ Y_{1}, \{ f^{-1}({\cal T}_{1}^{1}), f^{-1}({\cal T}_{1}^{2})\}, \ldots, \{ f^{-1}({\cal T}_{n_{0}}^{1}), f^{-1}({\cal T}_{n_{0}}^{2})\} \}\]
satisfying $(\dagger)_{n_{0}}$ with respect to the relation $f^{-1}({\cal S}_{1}^{s})$ by Remark \ref{rem-dagger-wrt-restriction}.
\end{pf}


\section[Application II]{Application II (The case where complexity is odd)}\label{sec-app2}

In this section, we apply Theorem \ref{thm-odd-n-1} and give some classification result of the mapping class groups of surfaces with odd complexity. As in Section \ref{sec-app1}, we find a concrete family of a Borel subset and subrelations and apply Theorem \ref{thm-odd-n-1} to classify the mapping class groups.

For a compact orientable surface $M$ of type $(g, p)$, we denote 
\begin{equation*}
g_{0}=g_{0}(M)= \begin{cases}
           2 & {\rm if} \ g\leq 2, \\
           g & {\rm if} \ g>2.
          \end{cases}
\end{equation*}

\subsection{Classification result for surfaces with odd complexity}

\begin{thm}\label{thm-kappa-odd-class}
Let $M$ be a surface such that $\kappa(M)\geq 0$ is odd and suppose that $M$ is not of type $(g, 0)$ with $g\geq 3$. Let us denote by $n_{0}=n(M)$ and $g_{0}=g_{0}(M)$ and assume $n_{0}\geq 2$. 

Then we can find $\sigma \in S(M)$ with $n(\sigma)=n_{0}$ and $\beta_{i}^{1}, \beta_{i}^{2}\in V(C(Q_{i}))$ with $i(\beta_{i}^{1}, \beta_{i}^{2})\neq 0$ for $i\in \{ 1, \ldots, n_{0}-g_{0}+2\}$ satisfying the following conditions, where $Q_{1},\ldots, Q_{n_{0}-g_{0}+2}$ are some $n_{0}-g_{0}+2$ components of $M_{\sigma}$ not of type $(0, 3)$: 
\begin{enumerate}
\item[(i)] For any distinct $i, j\in \{ 1,\ldots, n_{0}-g_{0}+2\}$ and any $s, t\in \{ 1, 2\}$, there exists $\tau \in S(M)$ such that $\beta_{i}^{s}, \beta_{j}^{t}\in \tau$ and $n(\tau)=n_{0}$;
\item[(ii)] If we have $i\in \{ 1, \ldots, n_{0}-g_{0}+2\}$ such that $Q_{i}$ is either of type $(0, 5)$ or $(1, 2)$, then there exists $\gamma \in V(C(Q_{i}))$ with $i(\gamma, \beta_{i}^{1})=i(\gamma, \beta_{i}^{2})=0$. 
\end{enumerate} 
\end{thm}  

The condition (ii) implies that for any $i\in \{ 1, \ldots, n_{0}-g_{0}+2\}$, we have a subsurface $Q_{i}'$ of $Q_{i}$ either of type $(0, 4)$ or $(1, 1)$ with $\beta_{i}^{1}, \beta_{i}^{2}\in V(C(Q_{i}'))$. 

We prove this theorem in the following subsections. First, we give the next corollary of the theorem:

\begin{cor}\label{cor-kappa-odd-class}
\begin{enumerate}
\item[(i)] Let $M$ be a surface satisfying the hypothesis in Theorem \ref{thm-kappa-odd-class}. With the assumption $(\diamond)'$, suppose that $F=X$, $\Gamma =\Gamma(M;m)$ and the action of $\Gamma$ on $X$ is measure-preserving. Then there exists a family
\[\{ X, \{ {\cal S}_{1}^{1}, {\cal S}_{1}^{2}\}, \ldots, \{ {\cal S}_{n_{0}}^{1}, {\cal S}_{n_{0}}^{2}\} \}\]
of the Borel set $X$ and subrelations of ${\cal R}$ satisfying $(\dagger)_{n_{0}}$ and the following condition: for any distinct $i, j\in \{ 1,\ldots, n_{0}-g_{0}+2\}$ and any $s, t\in \{ 1, 2\}$, we can find a family
\[\{ X, \{ {\cal T}_{1}^{1}, {\cal T}_{1}^{2}\}, \ldots, \{ {\cal T}_{n_{0}}^{1}, {\cal T}_{n_{0}}^{2}\} \}\]
satisfying $(\dagger)_{n_{0}}$ with respect to the relations ${\cal S}_{i}^{s}$ and ${\cal S}_{j}^{t}$. 

\item[(ii)] Let $M^{1}$ and $M^{2}$ be surfaces with $\kappa(M^{1}), \kappa(M^{2})\geq 0$. Suppose that the mapping class groups $\Gamma(M^{1})$ and $\Gamma(M^{2})$ are measure equivalent and that $\kappa(M^{1})=\kappa(M^{2})$ is odd. Then $g_{0}(M^{1})=g_{0}(M^{2})$. 
\end{enumerate}
\end{cor}

\begin{pf}
First, we show the assertion (i). Let $N\in {\Bbb N}$ be the natural number in Lemma \ref{lem-free-generate-curves}. Let $\sigma \in S(M)$, $Q_{1}, \ldots, Q_{n_{0}-g_{0}+2}$ and $\beta_{i}^{1}, \beta_{i}^{2}\in V(C(Q_{i}))$ for $i\in \{ 1, \ldots, n_{0}-g_{0}+2\}$ be as in Theorem \ref{thm-kappa-odd-class}. Let $Q_{n_{0}-g_{0}+3}, \ldots, Q_{n_{0}}$ be the remaining components of $M_{\sigma}$ not of type $(0, 3)$. For $i\in \{ n_{0}-g_{0}+3, \ldots, n_{0}\}$, choose two intersecting curves $\beta_{i}^{1}, \beta_{i}^{2}\in V(C(Q_{i}))$ such that we have a subsurface $Q_{i}'$ of $Q_{i}$ either of type $(0, 4)$ or $(1, 1)$ with $\beta_{i}^{1}, \beta_{i}^{2}\in V(C(Q_{i}'))$.

We denote by ${\cal S}_{i}^{s}$ the subrelation generated by the action of $t(\beta_{i}^{s})^{N}$ for $i\in \{ 1, \ldots, n_{0}\}$ and $s\in \{ 1, 2\}$. Then the family
\[\{ X, \{ {\cal S}_{1}^{1}, {\cal S}_{1}^{2}\}, \ldots, \{ {\cal S}_{n_{0}}^{1}, {\cal S}_{n_{0}}^{2}\} \}\]
satisfies $(\dagger)_{n_{0}}$ by Lemma \ref{lem-family-good} (i). 

For any distinct $i, j\in \{ 1, \ldots, n_{0}-g_{0}+2\}$ and any $s, t\in \{ 1, 2\}$, we can find a family 
\[\{ X, \{ {\cal T}_{1}^{1}, {\cal T}_{1}^{2}\}, \ldots, \{ {\cal T}_{n_{0}}^{1}, {\cal T}_{n_{0}}^{2}\} \}\]
satisfying $(\dagger)_{n_{0}}$ with respect to the relations ${\cal S}_{i}^{s}$ and ${\cal S}_{j}^{t}$, using the existence of $\tau \in S(M)$ in Theorem \ref{thm-kappa-odd-class} and applying the same argument as in the proof of Theorem \ref{thm-kappa-invariant-odd}.

Next, we show the assertion (ii). If $\Gamma(M^{1})$ and $\Gamma(M^{2})$ are measure equivalent, then we have two equivalence relations ${\cal R}_{1}$ on $(X_{1}, \mu_{1})$ and ${\cal R}_{2}$ on $(X_{2}, \mu_{2})$ of type ${\rm II}_{1}$ satisfying the following conditions: the relation ${\cal R}_{j}$ is generated by an essentially free, measure-preserving Borel action of $\Gamma(M^{j};m)$ for $j=1, 2$, where $m\geq 3$ is an integer. Moreover, ${\cal R}_{1}$ and ${\cal R}_{2}$ are weakly isomorphic. Let $f\colon Y_{1}\rightarrow Y_{2}$, $Y_{i}\subseteq X_{i}$ be a partial Borel isomorphism inducing the weak isomorphism between ${\cal R}_{1}$ and ${\cal R}_{2}$.

Assume that $g_{0}(M^{1})<g_{0}(M^{2})$. We deduce a contradiction. Since $\Gamma(M^{1})$ and $\Gamma(M^{2})$ are measure equivalent, we see that $n_{0}=n(M^{1})=n(M^{2})$ by Corollary \ref{cor-classification1}. Remark that $g_{0}(M^{2})$ is more than $2$ and equal to the genus of $M^{2}$. It follows from Theorem \ref{thm-odd-n-1} that for any family
\[{\cal F}=\{ Y, \{ {\cal S}_{1}^{1}, {\cal S}_{1}^{2}\},\ldots, \{ {\cal S}_{n_{0}}^{1}, {\cal S}_{n_{0}}^{2}\} \}\]
of a Borel subset $Y$ of $X_{2}$ and subrelations of $({\cal R}_{2})_{Y}$ satisfying $(\dagger)_{n_{0}}$, after permuting the indices $1, \ldots, n_{0}$, we see that for any distinct $i, j\in \{ 1, \ldots, g_{0}(M^{2})-1\}$, there exist $s, t\in \{ 1, 2\}$ satisfying the following condition: if a family 
\[\{ Y, \{ {\cal T}_{1}^{1}, {\cal T}_{1}^{2}\}, \ldots, \{ {\cal T}_{n}^{1}, {\cal T}_{n}^{2}\} \}\] 
of the Borel set $Y$ and subrelations of $({\cal R}_{2})_{Y}$ satisfies $(\dagger)_{n}$ with respect to the relations ${\cal S}_{i}^{s}$ and ${\cal S}_{j}^{t}$, then $n\leq n_{0}-1$.

Let $i\in \{ 1,\ldots, n_{0}\}$ and  
\[{\cal F}=\{ Y, \{ {\cal S}_{1}^{1}, {\cal S}_{1}^{2}\},\ldots, \{ {\cal S}_{n_{0}}^{1}, {\cal S}_{n_{0}}^{2}\} \}\]
be a family of a Borel subset $Y$ of $X_{2}$ and subrelations of $({\cal R}_{2})_{Y}$ satisfying $(\dagger)_{n_{0}}$. We denote by $F(i, {\cal F})$ the set of all elements $j\in \{ 1, \ldots, n_{0}\}\setminus \{ i\}$ such that for any $s, t\in \{ 1, 2\}$, there exists a family of the Borel set $Y$ and subrelations of $({\cal R}_{2})_{Y}$ satisfying $(\dagger)_{n_{0}}$ with respect to ${\cal S}_{i}^{s}$ and ${\cal S}_{j}^{t}$. Let us denote by
\[E({\cal F}, p)=\{ i\in \{ 1,\ldots, n_{0}\} : |F(i, {\cal F})|\geq n_{0}-p+1\} \]
for $p\in {\Bbb N}$. Remark that for any natural numbers $p\leq q$, we have $E({\cal F}, p)\subseteq E({\cal F}, q)$.

Then it follows from the above argument that given a family ${\cal F}$, after permuting the indices $1, \ldots, n_{0}$, we see that for any $i\in \{ 1, \ldots, g_{0}(M^{2})-1\}$, the set $F(i, {\cal F})$ contains no numbers in $\{ 1, \ldots, g_{0}(M^{2})-1\}$. Therefore, we see that $|F(i, {\cal F})|\leq n_{0}-g_{0}(M^{2})+1$. It follows that $i\notin E({\cal F}, g_{0}(M^{2})-1)$ for any $i\in \{ 1, \ldots, g_{0}(M^{2})-1\}$ and that $|E({\cal F}, g_{0}(M^{2})-1)|\leq n_{0}-g_{0}(M^{2})+1$ for any family ${\cal F}$.

On the other hand, note that the number of boundary components of $M^{1}$ is non-zero by the inequality $g_{0}(M^{1})<g_{0}(M^{2})$ and the equality $\kappa(M^{1})=\kappa(M^{2})$. Moreover, $n_{0}\geq 3$ since the genus of $M^{2}$ is more than $2$. It follows from the assertion (i) that there exists a family 
\[{\cal F}_{1}=\{ Y_{1}, \{ {\cal U}_{1}^{1}, {\cal U}_{1}^{2}\}, \ldots, \{ {\cal U}_{n_{0}}^{1}, {\cal U}_{n_{0}}^{2}\} \}\]
of the Borel subset $Y_{1}$ and subrelations of $({\cal R}_{1})_{Y_{1}}$ satisfying $(\dagger)_{n_{0}}$ and the following condition: for any $i\in \{ 1,\ldots, n_{0}-g_{0}(M^{1})+2\}$, we have
\[(\{ 1, \ldots, n_{0}-g_{0}(M^{1})+2\}\setminus \{ i\})\subseteq F(i, f({\cal F}_{1}))\]
and thus, $|F(i, f({\cal F}_{1}))|\geq n_{0}-g_{0}(M^{1})+1$, where $f({\cal F}_{1})$ denotes the family of the images by $f$ of the Borel subset and subrelations in ${\cal F}_{1}$. It follows that 
\begin{align*}
|E(f({\cal F}_{1}), g_{0}(M^{2})-1)|&\geq |E(f({\cal F}_{1}), g_{0}(M^{1}))| \\
                                 &\geq n_{0}-g_{0}(M^{1})+2 \\
                                 &> n_{0}-g_{0}(M^{2})+1,
\end{align*}
which is a contradiction.  
\end{pf}


\subsection{Geometric lemmas for the proof of Theorem \ref{thm-kappa-odd-class}}

In what follows, we prove Theorem \ref{thm-kappa-odd-class} by constructing a concrete family of curves on the surface $M$ and using Lemma \ref{lem-family-good}. For this, we need several geometric lemmas.

\begin{lem}\label{lem-kappa-sum}
Let $M$ be a surface with $\kappa(M)\geq 0$ and $\alpha \in V(C(M))$ be a separating curve. If $Q_{1}$ and $Q_{2}$ are two components of $M_{\alpha}$, then $\kappa(M)=\kappa(Q_{1})+\kappa(Q_{2})+2$.
\end{lem}  

This lemma comes from the fact that in general, $\kappa(M)+1$ is equal to the maximal number of distinct non-trivial isotopy classes of non-peripheral simple closed curves on $M$ which can be realized disjointly on $M$.

\begin{lem}\label{lem-sep-even-odd}
Let $M$ be a surface with $\kappa(M)\geq 0$ odd. 
\begin{enumerate}
\item[(i)] Let $\alpha, \beta \in V(C(M))$ be distinct separating curves with $i(\alpha, \beta)=0$. We denote by $Q_{1}$, $Q_{2}$, $Q_{3}$ the components of $M_{\{ \alpha, \beta \}}$. Then the following three conditions are equivalent:
\begin{enumerate}
\item[(a)] There exists $\tau \in S(M)$ such that $\alpha, \beta \in \tau$ and $n(\tau)=n(M)$;
\item[(b)] There exists distinct $i, j\in \{ 1, 2, 3\}$ such that both $\kappa(Q_{i})$ and $\kappa(Q_{j})$ are even;
\item[(c)] There exists $i\in \{ 1, 2, 3\}$ such that $\kappa(Q_{i})$ is even.
\end{enumerate}
\item[(ii)] Let $\alpha, \beta \in V(C(M))$ be curves with $i(\alpha, \beta)=0$ such that $\alpha$ is separating and $\beta$ is non-separating. We denote by $Q_{1}$ and $Q_{2}$ the components of $M_{\alpha}$ and suppose $\beta \in V(C(Q_{1}))$. Then the following two conditions are equivalent:
\begin{enumerate}
\item[(a)] There exists $\tau \in S(M)$ such that $\alpha, \beta \in \tau$ and $n(\tau)=n(M)$;
\item[(b)] The number $\kappa(Q_{1})$ is odd.
\end{enumerate}  
\end{enumerate}
\end{lem}

\begin{pf}
First, we show the assertion (i). By Lemma \ref{lem-kappa-sum}, we have $\kappa(M)=\kappa(Q_{1})+\kappa(Q_{2})+\kappa(Q_{3})+4$. Thus, the equivalence of the conditions (b) and (c) is clear. If we denote by $g_{k}$ and $p_{k}$ the genus and number of boundary components of $Q_{k}$ for $k=1, 2, 3$, respectively, then we have
\[g_{1}+g_{2}+g_{3}=g,\ \ p_{1}+p_{2}+p_{3}=p+4.\]

If the condition (b) is not true, then all $\kappa(Q_{k})$ $(k=1, 2, 3)$ are odd. Then
\begin{align*}
n(Q_{1})+n(Q_{2})+n(Q_{3})&= \sum_{k=1}^{3}(g_{k}+[(g_{k}+p_{k}-2)/2])\\
                          &= g+(g+p-5)/2=n(M)-1.
\end{align*}
Hence, there exist no $\tau \in S(M)$ such that $\alpha, \beta \in \tau$ and $n(\tau)=n(M)$, that is, the condition (a) is not true.

Conversely, if the condition (b) is true, then we may assume that $\kappa(Q_{1})$ and $\kappa(Q_{2})$ are even. Then 
\begin{align*}
n(Q_{1})+n(Q_{2})+n(Q_{3})&= \sum_{k=1}^{3}(g_{k}+[(g_{k}+p_{k}-2)/2])\\
                          &=g+\sum_{k=1}^{2}(g_{k}+p_{k}-2)/2 +(g_{3}+p_{3}-3)/2\\ 
                          &=g+(g+p-3)/2=n(M).
\end{align*}
Thus, the condition (a) is true.

Next, we show the assertion (ii). If we denote by $g_{k}$ and $p_{k}$ the genus and number of boundary components of $Q_{k}$ for $k=1, 2$, respectively, then we have
\[g_{1}+g_{2}=g,\ \ p_{1}+p_{2}=p+2. \]
When we cut $Q_{1}$ along the curve $\beta$, the resulting surface $(Q_{1})_{\beta}$ is of type $(g_{1}-1, p_{1}+2)$. For any $\tau \in S(M)$ with $\alpha, \beta \in \tau$, we have 
\begin{align*}
n(\tau)&\leq n((Q_{1})_{\beta})+n(Q_{2}) \\
       &=g-1+[(g_{1}+p_{1}-1)/2]+[(g_{2}+p_{2}-2)/2].
\end{align*}
If $\kappa(Q_{1})$ is even, then $\kappa(Q_{2})$ is odd and the right hand side of the above inequality is equal to $g+(g+p-5)/2=n(M)-1$. Thus, the condition (a) cannot happen. Conversely, if $\kappa(Q_{1})$ is odd, then $\kappa(Q_{2})$ is even and 
\[n((Q_{1})_{\beta})+n(Q_{2})=g+(g+p-3)/2=n(M). \] 
Therefore, the condition (a) is satisfied.
\end{pf}

\begin{lem}\label{lem-curves-m12}
Let $Q$ be a surface of type $(1, 2)$. Then we can find a non-separating curve $\alpha \in V(C(Q))$ and two separating curves $\beta_{1}, \beta_{2}\in V(C(Q))$ with $i(\alpha, \beta_{1})=i(\alpha, \beta_{2})=0$ and $i(\beta_{1}, \beta_{2})\neq 0$.
\end{lem}

Remark that in this case, the cut surface $Q_{\alpha}$ is of type $(0, 4)$ and $\beta_{1}, \beta_{2}\in V(C(Q_{\alpha}))$. Before the proof, we give the following remark:

\begin{rem}\label{rem-curves-m04}
For each $\delta \in V(C(M_{0, 4}))$, if we cut $M_{0, 4}$ along $\delta$, then we obtain two pairs of pants and the boundary components of each of the pairs of pants consists of two of the boundary components of $M_{0, 4}$ and one corresponding to $\delta$. Conversely, for any two boundary components $P_{1}$, $P_{2}$ of $M_{0, 4}$, there exists $\delta \in V(C(M_{0, 4}))$ such that one of the components of $(M_{0, 4})_{\delta}$ has $P_{1}$ and $P_{2}$ as boundary components. In this case, we say that $\delta \in V(C(M_{0, 4}))$ {\it bounds}\index{bound} the boundary components $P_{1}$ and $P_{2}$. If $P_{3}$ and $P_{4}$ are the other boundary components of $M_{0, 4}$, then it goes without saying that $\delta$ bounds $P_{3}$ and $P_{4}$.
\end{rem}

\begin{pf*}{{\sc Proof of Lemma \ref{lem-curves-m12}}}
Choose a non-separating curve $\alpha \in V(C(Q))$. Let $p\colon Q_{\alpha}\rightarrow Q$ be the canonical map. Then $p^{-1}(\alpha)$ has two components $P_{1}$, $P_{2}$ both of which are boundary components of $Q_{\alpha}$. Let $\beta_{1}\in V(C(Q_{\alpha}))$ be a curve which bounds $P_{1}$ and $P_{2}$. Let $\gamma \in V(C(Q_{\alpha}))$ be a curve which bounds $P_{1}$ and another boundary component of $Q_{\alpha}$ other than $P_{2}$. Since $\beta_{1}\neq \gamma$ and $Q_{\alpha}$ is of type $(0, 4)$, we have $i(\beta_{1}, \gamma)\neq 0$.

In what follows, we regard $\beta_{1}$ and $\gamma$ as elements in $V(C(Q))$. It follows from Lemma \ref{lem-dehn-criterion} that $i(t_{\gamma}^{n}\beta_{1}, \beta_{1})\neq 0$ for any $n\geq 4$, where $t_{\gamma}\in \Gamma(Q)$ denotes the Dehn twist about $\gamma$. On the other hand, since $i(\alpha, \gamma)=0$, we have
\[i(t_{\gamma}^{4}\beta_{1}, \alpha)=i(\beta_{1}, t_{\gamma}^{-4}\alpha)=i(\beta_{1}, \alpha)=0.\]
Note that $t_{\gamma}^{4}\beta_{1}\in V(C(Q))$ is separating since $\beta_{1}$ is separating. Thus, the curves $\alpha$, $\beta_{1}$, $\beta_{2}=t_{\gamma}^{4}\beta_{1}\in V(C(Q))$ are desired ones.     
\end{pf*}


\subsection{The proof of Theorem \ref{thm-kappa-odd-class}}

Before the proof, we give a few remarks about notation.

\begin{rem}\label{rem-notation-class}
\begin{enumerate}
\item[(i)] For a surface $Q$ with $\kappa(Q)\geq 0$, let $\alpha \in V(C(Q))$ be a separating curve and let $Q_{1}$ and $Q_{2}$ be the components of the cut surface $Q_{\alpha}$. Let $p\colon Q_{1}\sqcup Q_{2}\rightarrow Q$ be the natural map. We denote also by $\alpha$ the boundary component of $Q_{1}$ (resp. $Q_{2}$) which is a component of $p^{-1}(\alpha)$ as far as there is no possibility of confusion. 
\item[(ii)] For a surface $Q$ with $\kappa(Q)\geq 0$, we denote by $P_{1}, \ldots, P_{p}$ the boundary components of $Q$. Let $Q_{1}$ be a subsurface of $Q$ and $q\colon Q_{1}\rightarrow Q$ be the natural embedding (see Definition \ref{defn-geometric-things} (i)). Let us denote by $\alpha_{1}, \ldots, \alpha_{j}\in V(C(Q))$ all the isotopy classes of components in $\partial_{Q}Q_{1}$. Suppose that $P_{i_{1}}, \ldots, P_{i_{k}}$ are all the boundary components which are also contained in the image by $q$ of the boundary of $Q_{1}$. Then we denote 
\[\partial Q_{1}=\{ \alpha_{1}, \ldots, \alpha_{j}, P_{i_{1}}, \ldots, P_{i_{k}}\}\]  for simplicity.  
\item[(iii)] For a surface $Q$ of type $(0, p)$ with $p\geq 4$, let $P_{1},\ldots, P_{p}$ be the boundary components of $Q$. Note that any $\alpha \in V(C(Q))$ is separating. For $\alpha \in V(C(Q))$, we say that $\alpha$ {\it bounds} $P_{i_{1}},\ldots, P_{i_{k}}$ if one of the components $Q_{1}$ of $Q_{\alpha}$ satisfies $\partial Q_{1}=\{ \alpha, P_{i_{1}}, \ldots, P_{i_{k}}\}$. It goes without saying that $\alpha$ bounds also all the boundary components of $Q$ except for $P_{i_{1}}, \ldots , P_{i_{k}}$. This definition is a generalization of that in Remark \ref{rem-curves-m04}.
\end{enumerate}  
\end{rem}

We suppose that $M$ is a compact orientable surface of type $(g, p)$ and that $\kappa(M)\geq 0$ is odd. Let us denote by $n_{0}=n(M)$. Moreover, we assume that $n_{0}\geq 2$ and $M$ is not of type $(g, 0)$ with $g\geq 3$. We denote by $P_{1},\ldots, P_{p}$ all the boundary components of $M$. 

In what follows, we give a family 
\[\sigma =\{ \alpha_{1},\ldots, \alpha_{n_{0}-1}\}\in S(M),\ \ \{ \beta_{i}^{1}, \beta_{i}^{2}, Q_{i}\}_{i=1,\ldots, n_{0}-g_{0}+2}\]
satisfying the conditions in Theorem \ref{thm-kappa-odd-class} for the surface $M$ in the cases of $g=0$, $g=1$, $g\geq 2$ with $g$ even and $g\geq 3$ with $g$ odd, respectively. 
\vspace{1em}


\subsubsection{The case of $g=0$ (Figure \ref{picture-kappa-odd-g-0})}\label{sub-odd-g0}

In this case, $p\in \{ 7, 9, 11,\ldots \}$ and $n_{0}=(p-3)/2$. First, we give $\sigma =\{ \alpha_{1}, \ldots, \alpha_{n_{0}-1}\} \in S(M)$ inductively as follows: choose $\alpha_{1}\in V(C(M))$ as a curve such that if we denote $M_{\alpha_{1}}=Q_{1}\sqcup Q_{1}'$, then $Q_{1}$ is of type $(0, 4)$ and satisfies $\partial Q_{1}=\{ \alpha_{1}, P_{1}, P_{2}, P_{3}\}$. 

For $i=2, 3,\ldots, n_{0}-1$, we choose $\alpha_{i}\in V(C(Q_{i-1}'))$ as a curve such that if we denote $(Q_{i-1}')_{\alpha_{i}}=Q_{i}\sqcup Q_{i}'$, then $Q_{i}$ is of type $(0, 4)$ and satisfies $\partial Q_{i}=\{ \alpha_{i-1}, \alpha_{i}, P_{2i}, P_{2i+1}\}$. 

Then $Q_{n_{0}-1}'$ is of type $(0, 5)$ and $\partial Q_{n_{0}-1}'=\{ \alpha_{n_{0}-1}, P_{p-3}, P_{p-2}, P_{p-1}, P_{p}\}$. We write $Q_{n_{0}}=Q_{n_{0}-1}'$. Define $\sigma =\{ \alpha_{1}, \ldots, \alpha_{n_{0}-1}\}\in S(M)$.

\begin{figure}[hp]
\centering
\vspace*{2\baselineskip}
\epsfbox{pkog0.1}
\caption{The case of $g=0$}
\label{picture-kappa-odd-g-0}
\end{figure}

\begin{figure}[hp]
\centering
\vspace*{2\baselineskip}
\epsfbox{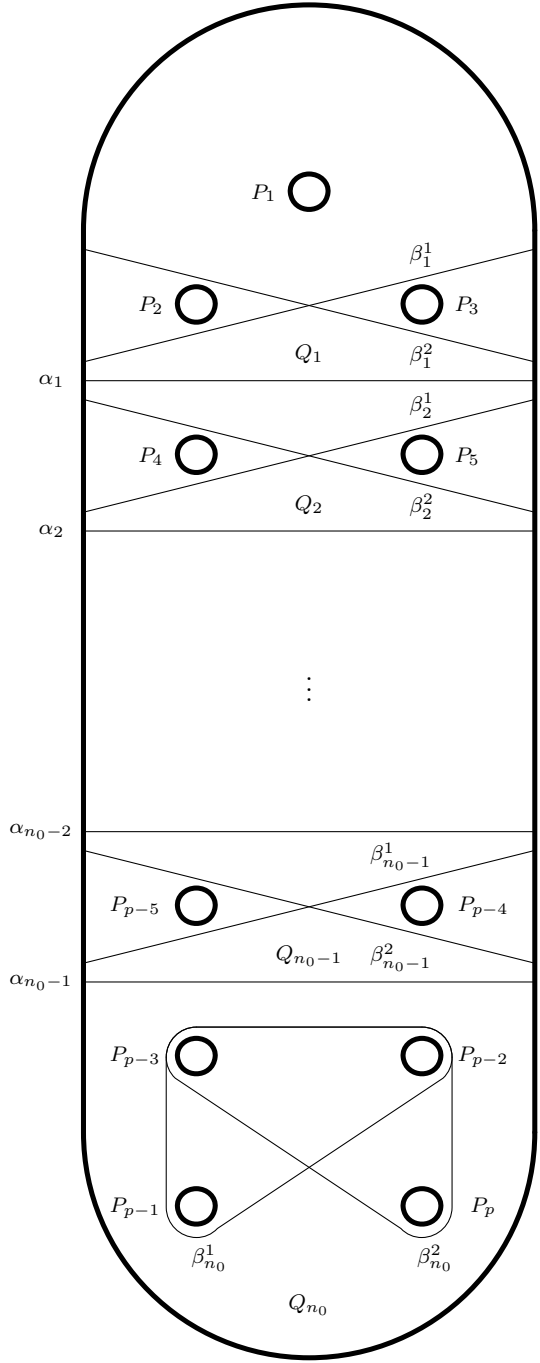}
\caption{The case of $g=1$}
\label{picture-kappa-odd-g-1}
\end{figure}

Next, we choose two intersecting curves $\beta_{i}^{1}, \beta_{i}^{2}\in V(C(Q_{i}))$ satisfying the following conditions:
\begin{itemize}
\item The curves $\beta_{1}^{1}$ and $\beta_{1}^{2}$ are arbitrary two intersecting ones;
\item For $i\in \{ 2, \ldots, n_{0}-1\}$, the curve $\beta_{i}^{1}$ (resp. $\beta_{i}^{2}$) bounds $\alpha_{i-1}$, $P_{2i}$ (resp. $\alpha_{i-1}$, $P_{2i+1}$);
\item The curve $\beta_{n_{0}}^{1}$ (resp. $\beta_{n_{0}}^{2}$) bounds $P_{p-3}$, $P_{p-2}$, $P_{p-1}$ (resp. $P_{p-3}$, $P_{p-2}$, $P_{p}$).
\end{itemize}

We show Theorem \ref{thm-kappa-odd-class} for this choice. Take any $i, j\in \{ 1,\ldots, n_{0}\}$ with $i<j$ and any $s, t\in \{ 1, 2\}$. If we cut $M$ along $\beta_{i}^{s}$, then we have two components $R_{1}$ and $R_{1}'$ of types $(0, 2i+1)$ and $(0, p-2i+1)$, respectively. Moreover, $\beta_{j}^{t}\in V(C(R_{1}'))$. If we cut $R_{1}'$ along $\beta_{j}^{t}$, then we have two components $R_{2}$ and $R_{3}$ of types $(0, 2(j-i)+2)$ and $(0, p-2j+1)$, respectively. Since $\kappa(R_{2})$ and $\kappa(R_{3})$ are even, we can apply Lemma \ref{lem-sep-even-odd} (i) and show the condition (i) in Theorem \ref{thm-kappa-odd-class}. For the condition (ii) in the theorem, let $\gamma \in V(C(Q_{n_{0}}))$ be a curve bounding $P_{p-3}$, $P_{p-2}$ with $i(\gamma, \beta_{n_{0}}^{1})=i(\gamma, \beta_{n_{0}}^{2})=0$. This completes the proof of Theorem \ref{thm-kappa-odd-class} in the case of $g=0$.     
\vspace{1em}


\subsubsection{The case of $g=1$ (Figure \ref{picture-kappa-odd-g-1})} 

In this case, $p\in \{ 4, 6, 8, \ldots \}$ and $n_{0}=p/2$. Let $\bar{M}$ be a surface of type $(0, p+3)$ with boundary components $\bar{P}_{1}, \ldots, \bar{P}_{p+3}$. Note $n(\bar{M})=n_{0}$. We apply the argument in \ref{sub-odd-g0} to $\bar{M}$ and give a family
\[\bar{\sigma}=\{ \bar{\alpha}_{1},\ldots, \bar{\alpha}_{n_{0}-1}\}\in S(\bar{M}),\ \ \{ \bar{\beta}_{i}^{1}, \bar{\beta}_{i}^{2}, \bar{Q}_{i}\}_{i=1,\ldots, n_{0}}.\] 
We construct $M$ by replacing $\bar{Q}_{n_{0}}$ with $M_{1, 2}$, that is, remove $\bar{Q}_{n_{0}}$ from $\bar{M}$ and glue one of the boundary components of $M_{1, 2}$ to $\bar{\alpha}_{n_{0}-1}$ which is a boundary component of $\bar{M}\setminus \bar{Q}_{n_{0}}$. We denote by $Q_{n_{0}}$ the new glued subsurface of $M$ of type $(1, 2)$ and by $P_{p}$ the boundary component of $Q_{n_{0}}$ other than the one corresponding to $\bar{\alpha}_{n_{0}-1}$. Let us denote by $\alpha_{n_{0}-1}\in V(C(M))$ the corresponding curve to $\bar{\alpha}_{n_{0}-1}$.  

Let us denote 
\begin{align*}
P_{i}=\bar{P}_{i} & \ \ \ {\rm for}\ i\in \{ 1, 2, \ldots, p-1\},\\
\alpha_{i}=\bar{\alpha}_{i} & \ \ \ {\rm for}\ i\in \{ 1, 2, \ldots, n_{0}-2\},\\
\beta_{i}^{s}=\bar{\beta}_{i}^{s}, \ \ Q_{i}=\bar{Q}_{i} & \ \ \ {\rm for}\ i\in \{ 1, 2,\ldots, n_{0}-1\} {\rm \ and \ } s\in \{ 1, 2\}.
\end{align*}
Define $\sigma =\{ \alpha_{1}, \ldots, \alpha_{n_{0}-1}\}\in S(M)$. We choose two intersecting curves $\beta_{n_{0}}^{1}, \beta_{n_{0}}^{2}\in V(C(Q_{n_{0}}))$ as two separating curves in Lemma \ref{lem-curves-m12}.

We show Theorem \ref{thm-kappa-odd-class} for this choice. Take any $i, j\in \{ 1,\ldots, n_{0}\}$ with $i<j$ and any $s, t\in \{ 1, 2\}$. If we cut $M$ along $\beta_{i}^{s}$, then we have two components $R_{1}$ and $R_{1}'$ of types $(0, 2i+1)$ and $(1, p-2i+1)$, respectively, with $\beta_{j}^{t}\in V(C(R_{1}'))$. If we cut $R_{1}'$ along $\beta_{j}^{t}$, then we have two components $R_{2}$ and $R_{3}$ of types $(0, 2(j-i)+2)$ and $(1, p-2j+1)$, respectively. Since $\kappa(R_{2})$ and $\kappa(R_{3})$ are even, we can apply Lemma \ref{lem-sep-even-odd} (i) and show the condition (i) in Theorem \ref{thm-kappa-odd-class}. The condition (ii) in the theorem follows from the choice of $\beta_{n_{0}}^{1}, \beta_{n_{0}}^{2}\in V(C(Q_{n_{0}}))$ in Lemma \ref{lem-curves-m12}.  
\vspace{1em}


\begin{figure}[hp]
\centering

\vspace*{2\baselineskip}
\epsfbox{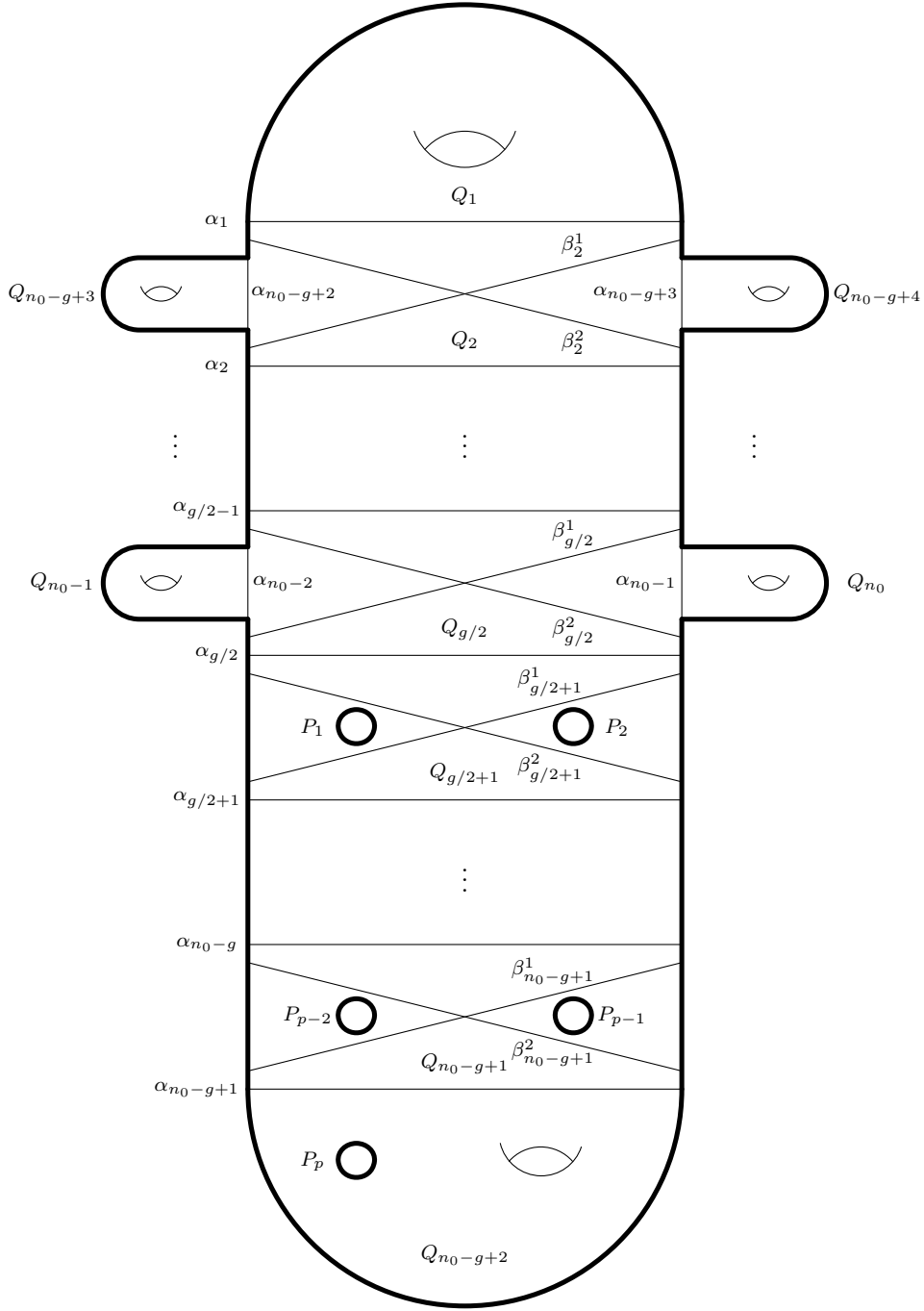}
\caption{The case where $g\geq 2$ is even}
\label{picture-kappa-odd-g-even}
\end{figure}

\subsubsection{The case where $g\geq 2$ is even (Figure \ref{picture-kappa-odd-g-even})}\label{sub-odd-g-even}

In this case, $p\in \{ 1, 3, 5,\ldots \}$. First, we give $\sigma =\{ \alpha_{1}, \ldots, \alpha_{n_{0}-1}\} \in S(M)$ as follows: we choose disjoint curves  $\alpha_{1}\in V(C(M))$ and $\alpha_{n_{0}-g+2}, \alpha_{n_{0}-g+3}, \ldots, \alpha_{n_{0}-1}\in V(C(M))$ so that they separate components $Q_{1}$ and $Q_{n_{0}-g+3}$, $Q_{n_{0}-g+4},\ldots, Q_{n_{0}}$ of type $(1, 1)$, respectively, from $M$. The remaining surface $Q_{1}'$ is of type $(1, g+p-1)$ and satisfies
\[\partial Q_{1}'=\{ \alpha_{1}, \alpha_{n_{0}-g+2}, \alpha_{n_{0}-g+3}, \ldots, \alpha_{n_{0}-1}, P_{1}, \ldots, P_{p}\}.\]

We define $\alpha_{i}\in V(C(M))$ for $i\in \{ 2, 3, \ldots, g/2\}$ inductively as follows: the curve $\alpha_{i}\in V(C(Q_{i-1}'))$ is chosen so that if we denote $(Q_{i-1}')_{\alpha_{i}}=Q_{i}\sqcup Q_{i}'$, then $Q_{i}$ is of type $(0, 4)$ and 
\[\partial Q_{i}=\{ \alpha_{i-1}, \alpha_{i}, \alpha_{n_{0}-g+2i-2}, \alpha_{n_{0}-g+2i-1}\}.\] 
Then $Q_{g/2}'$ is of type $(1, p+1)$ and 
\[\partial Q_{g/2}'=\{ \alpha_{g/2}, P_{1}, \ldots, P_{p}\}.\]

We define $\alpha_{i}\in V(C(M))$ for $i\in \{ g/2+1, g/2+2,\ldots, n_{0}-g+1\}$ inductively as follows: the curve $\alpha_{i}\in V(C(Q_{i-1}'))$ is chosen so that if we denote $(Q_{i-1}')_{\alpha_{i}}=Q_{i}\sqcup Q_{i}'$, then $Q_{i}$ is of type $(0, 4)$ and 
\[\partial Q_{i}=\{ \alpha_{i-1}, \alpha_{i}, P_{2i-g-1}, P_{2i-g}\}.\]
Then $Q_{n_{0}-g+1}'$ is of type $(1, 2)$ and 
\[\partial Q_{n_{0}-g+1}'=\{ \alpha_{n_{0}-g+1}, P_{p}\}.\]
Define $\sigma =\{ \alpha_{1}, \ldots, \alpha_{n_{0}-1}\}\in S(M)$ and denote $Q_{n_{0}-g+2}=Q_{n_{0}-g+1}'$.

Next, we choose two intersecting curves $\beta_{i}^{1}, \beta_{i}^{2}\in V(C(Q_{i}))$ satisfying the following conditions:
\begin{itemize}
\item The curves $\beta_{1}^{1}$ and $\beta_{1}^{2}$ are arbitrary two intersecting ones;
\item For $i\in \{ 2, 3,\ldots, g/2\}$, the curve $\beta_{i}^{1}$ (resp. $\beta_{i}^{2}$) bounds $\alpha_{i-1}$, $\alpha_{n_{0}-g+2i-2}$ (resp. $\alpha_{i-1}$, $\alpha_{n_{0}-g+2i-1}$);
\item For $i\in \{ g/2+1, g/2+2, \ldots, n_{0}-g+1\}$, the curve $\beta_{i}^{1}$ (resp. $\beta_{i}^{2}$) bounds $\alpha_{i-1}$, $P_{2i-g-1}$ (resp. $\alpha_{i-1}$, $P_{2i-g}$);
\item The curves $\beta_{n_{0}-g+2}^{1}$ and $\beta_{n_{0}-g+2}^{2}$ are two separating curves in Lemma \ref{lem-curves-m12};
\item For $i\in \{ n_{0}-g+3, n_{0}-g+4,\ldots, n_{0}\}$, the curves $\beta_{i}^{1}$ and $\beta_{i}^{2}$ are arbitrary two intersecting ones. (Although this choice is not needed for the proof of Theorem \ref{thm-kappa-odd-class}, it will be used in the next section.)
\end{itemize}

We show Theorem \ref{thm-kappa-odd-class} for this choice. Take any $i, j\in \{ 1,\ldots, n_{0}-g+2\}$ with $i<j$ and any $s, t\in \{ 1, 2\}$. We cut $M$ along $\beta_{j}^{t}$. 

If $j\in \{ g/2+1, g/2+2, \ldots, n_{0}-g+2\}$, then we have two components $R'$ of type $(g-1, p-2n_{0}+2g+2j-3)$ and $R$ of type $(1, 2n_{0}-2g-2j+5)$ with $\beta_{i}^{s}\in V(C(R'))$. Since $\kappa(R)$ is even, we can apply Lemma \ref{lem-sep-even-odd} (i) if $i\neq 1$ and (ii) if $i=1$.

If $j\in \{ 2, 3, \ldots, g/2\}$, then we have two components $R'$ of type $(2j-2, 1)$ and $R$ of type $(g-2j+2, p+1)$ with $\beta_{i}^{s}\in V(C(R'))$. Since $\kappa(R)$ is even, we can apply Lemma \ref{lem-sep-even-odd} (i) if $i\neq 1$ and (ii) if $i=1$ and show the condition (i) in Theorem \ref{thm-kappa-odd-class}. The condition (ii) in the theorem follows from the choice of $\beta_{n_{0}-g+2}^{1}, \beta_{n_{0}-g+2}^{2}\in V(C(Q_{n_{0}-g+2}))$ in Lemma \ref{lem-curves-m12}. 
\vspace{1em}


\subsubsection{The case where $g\geq 3$ is odd}

In this case, $p\in \{ 2, 4, 6,\ldots \}$. Let $\bar{M}$ be a surface of type $(g-1, p+1)$ with boundary components $\bar{P}_{1},\ldots, \bar{P}_{p+1}$. Note $n(\bar{M})=n_{0}-1$. We apply the argument in \ref{sub-odd-g-even} to $\bar{M}$ and give a family
\[\bar{\sigma}=\{ \bar{\alpha}_{1},\ldots, \bar{\alpha}_{n_{0}-2}\}\in S(\bar{M}),\ \ \{ \bar{\beta}_{i}^{1}, \bar{\beta}_{i}^{2}, \bar{Q}_{i}\}_{i=1,\ldots, n_{0}-1}.\]
We construct $M$ by gluing the boundary component of a new surface $Q_{n_{0}}$ of type $(1, 1)$ to $\bar{P}_{1}$. Let $\alpha_{n_{0}-1}\in V(C(M))$ be the corresponding curve to $\bar{P}_{1}$. Let us denote
\begin{align*}
P_{i-1}=\bar{P}_{i} & \ \ \ {\rm for}\ i\in \{ 2, 3, \ldots, p+1\},\\
\alpha_{i}=\bar{\alpha}_{i} & \ \ \ {\rm for}\ i\in \{ 1, 2,\ldots, n_{0}-2\},\\
\beta_{i}^{s}=\bar{\beta}_{i}^{s},\ \  Q_{i}=\bar{Q}_{i} & \ \ \ {\rm for}\ i\in \{ 1, 2,\ldots, n_{0}-1\} {\rm \ and \ } s\in \{ 1, 2\}.
\end{align*} 
Define $\sigma =\{ \alpha_{1}, \ldots, \alpha_{n_{0}-1}\}\in S(M)$. We choose arbitrary two intersecting curves $\beta_{n_{0}}^{1}, \beta_{n_{0}}^{2}\in V(C(Q_{n_{0}}))$. 

We show Theorem \ref{thm-kappa-odd-class} for this choice. Take any $i, j\in \{ 1,\ldots, n_{0}-g+2\}$ with $i<j$ and any $s, t\in \{ 1, 2\}$. As in \ref{sub-odd-g-even}, if we cut $M$ along $\beta_{j}^{t}$, then we have two components $R'$ and $R$ such that $\beta_{i}^{s}\in V(C(R'))$ and $\kappa(R)$ is even. Thus, we can apply Lemma \ref{lem-sep-even-odd} (i) if $i\neq 1$ and (ii) if $i=1$ and show the conditions in Theorem \ref{thm-kappa-odd-class} as well as in \ref{sub-odd-g-even}.

We complete the proof of Theorem \ref{thm-kappa-odd-class}.  \hfill $\square$


\section[Application III]{Application III (The case where complexity is even)}\label{sec-app3}

In this section, we apply Theorem \ref{thm-even-n-2} and give some classification result of the mapping class groups of surfaces with even complexity. We adapt the notation in Section \ref{sec-app2}.

\subsection{Classification result for surfaces with even complexity}

\begin{thm}\label{thm-kappa-even-class1}
Let $M$ be a surface of type $(g, p)$ such that $\kappa(M)\geq 0$ is even and denote $n_{0}=n(M)$. Suppose that $g\in \{ 0, 1, 2\}$ and $n_{0}\geq 3$. 

Then we can find $\sigma \in S(M)$ with $n(\sigma)=n_{0}$ and $\beta_{i}^{1}, \beta_{i}^{2}\in V(C(Q_{i}))$ with $i(\beta_{i}^{1}, \beta_{i}^{2})\neq 0$ for $i\in \{ 1, \ldots, n_{0}\}$ satisfying the following condition, where we denote by $Q_{1},\ldots, Q_{n_{0}}$ all the components of $M_{\sigma}$: for any distinct $i, j, k\in \{ 1,\ldots, n_{0}\}$ and any $s, t, u\in \{ 1, 2\}$, we have $\tau \in S(M)$ such that $\beta_{i}^{s}, \beta_{j}^{t}, \beta_{k}^{u}\in \tau$ and $n(\tau)=n_{0}-1$.
\end{thm}

\begin{cor}\label{cor-kappa-even-class1}
Let $M$ be a surface satisfying the hypothesis in Theorem \ref{thm-kappa-even-class1}. With the assumption $(\diamond)'$, suppose that $F=X$, $\Gamma =\Gamma(M;m)$ and the action of $\Gamma$ on $X$ is measure-preserving. Then there exists a family
\[\{ X, \{ {\cal S}_{1}^{1}, {\cal S}_{1}^{2}\}, \ldots, \{ {\cal S}_{n_{0}}^{1}, {\cal S}_{n_{0}}^{2}\} \}\]
of the Borel set $X$ and subrelations of ${\cal R}$ satisfying $(\dagger)_{n_{0}}$ and the following condition: for any distinct $i, j, k\in \{ 1,\ldots, n_{0}\}$ and any $s, t, u\in \{ 1, 2\}$, we can find a family
\[\{ X, \{ {\cal T}_{1}^{1}, {\cal T}_{1}^{2}\}, \ldots, \{ {\cal T}_{n_{0}-1}^{1}, {\cal T}_{n_{0}-1}^{2}\} \}\]
satisfying $(\dagger)_{n_{0}-1}$ with respect to the relations ${\cal S}_{i}^{s}$, ${\cal S}_{j}^{t}$ and ${\cal S}_{k}^{u}$. 
\end{cor}

We can prove this corollary as well as Corollary \ref{cor-kappa-odd-class} (i), using Theorem \ref{thm-kappa-even-class1}.

\begin{thm}\label{thm-kappa-even-class2}
Let $M$ be a surface of type $(g, p)$ such that $\kappa(M)\geq 0$ is even and suppose that $g\geq 2$ and $M$ is not of type $(2, 0)$. Let us denote $n_{0}=n(M)$. 

Then we can find $\sigma \in S(M)$ with $n(\sigma)=n_{0}$ and $\beta_{i}^{1}, \beta_{i}^{2}\in V(C(Q_{i}))$ with $i(\beta_{i}^{1}, \beta_{i}^{2})\neq 0$ for $i\in \{ 1, \ldots, n_{0}\}$ satisfying the following condition, where we denote by $Q_{1},\ldots, Q_{n_{0}}$ all the components of $M_{\sigma}$: for any subset $F\subseteq \{ 1,\ldots, n_{0}\}$ with $|F|\geq g+1$, there exist distinct three elements $i, j, k\in F$ satisfying that for any $s, t, u\in \{ 1, 2\}$, we have $\tau \in S(M)$ such that $\beta_{i}^{s}, \beta_{j}^{t}, \beta_{k}^{u}\in \tau$ and $n(\tau)=n_{0}-1$.
\end{thm}

Note that any surface $M$ in the assumption of the above theorem satisfies $n_{0}=n(M)\geq g+1\geq 3$ . 

\begin{cor}\label{cor-kappa-even-class2}
Let $M$ be a surface satisfying the hypothesis in Theorem \ref{thm-kappa-even-class2}. With the assumption $(\diamond)'$, suppose that $F=X$, $\Gamma =\Gamma(M;m)$ and the action of $\Gamma$ on $X$ is measure-preserving. Then there exists a family
\[\{ X, \{ {\cal S}_{1}^{1}, {\cal S}_{1}^{2}\}, \ldots, \{ {\cal S}_{n_{0}}^{1}, {\cal S}_{n_{0}}^{2}\} \}\]
of the Borel set $X$ and subrelations of ${\cal R}$ satisfying $(\dagger)_{n_{0}}$ and the following condition: for any subset $F\subseteq \{ 1,\ldots, n_{0}\}$ with $|F|\geq g+1$, there exist distinct three elements $i, j, k\in F$ satisfying that for any $s, t, u\in \{ 1, 2\}$, we can find a family
\[\{ X, \{ {\cal T}_{1}^{1}, {\cal T}_{1}^{2}\}, \ldots, \{ {\cal T}_{n_{0}-1}^{1}, {\cal T}_{n_{0}-1}^{2}\} \}\]
satisfying $(\dagger)_{n_{0}-1}$ with respect to the relations ${\cal S}_{i}^{s}$, ${\cal S}_{j}^{t}$ and ${\cal S}_{k}^{u}$. 
\end{cor}

We can prove this corollary as well as Corollary \ref{cor-kappa-odd-class} (i), using Theorem \ref{thm-kappa-even-class2}. We give the proof of Theorems \ref{thm-kappa-even-class1} and \ref{thm-kappa-even-class2} in the following subsections. 

\begin{cor}\label{cor-kappa-even-class-mcg}
Let $M^{1}$ and $M^{2}$ be surfaces with $\kappa(M^{1}), \kappa(M^{2})\geq 0$. Suppose that the mapping class groups $\Gamma(M^{1})$ and $\Gamma(M^{2})$ are measure equivalent and that $\kappa(M^{1})=\kappa(M^{2})$ is even. Then $g_{0}(M^{1})=g_{0}(M^{2})$. 
\end{cor}

\begin{pf}
If $\Gamma(M^{1})$ and $\Gamma(M^{2})$ are measure equivalent, we have two equivalence relations ${\cal R}_{1}$ on $(X_{1}, \mu_{1})$ and ${\cal R}_{2}$ on $(X_{2}, \mu_{2})$ of type ${\rm II}_{1}$ satisfying the following conditions: the relation ${\cal R}_{j}$ is generated by an essentially free, measure-preserving Borel action of $\Gamma(M^{j};m)$ for $j=1, 2$, where $m\geq 3$ is an integer. Moreover, ${\cal R}_{1}$ and ${\cal R}_{2}$ are weakly isomorphic. Let $f\colon Y_{1}\rightarrow Y_{2}$, $Y_{i}\subseteq X_{i}$ be a partial Borel isomorphism inducing the weak isomorphism between ${\cal R}_{1}$ and ${\cal R}_{2}$. Let us denote $n_{0}=n(M^{1})=n(M^{2})$ (see Corollary \ref{cor-classification1}).

First, assume that $2=g_{0}(M^{1})<g_{0}(M^{2})$. We deduce a contradiction. Remark that $n_{0}=n(M^{2})\geq 3$ since $g_{0}(M^{2})>2$. We apply Corollary \ref{cor-kappa-even-class1} to $M^{1}$. Using $f\colon Y_{1}\rightarrow Y_{2}$ and Remark \ref{rem-dagger-wrt-restriction}, we have a family 
\[{\cal F}=\{ Y_{2}, \{ {\cal S}_{1}^{1}, {\cal S}_{1}^{2}\}, \ldots, \{ {\cal S}_{n_{0}}^{1}, {\cal S}_{n_{0}}^{2}\} \}\]
of the Borel set $Y_{2}$ and subrelations of $({\cal R}_{2})_{Y_{2}}$ satisfying $(\dagger)_{n_{0}}$ and the following condition: for any distinct $i, j, k\in \{ 1, \ldots, n_{0}\}$ and any $s, t, u\in \{ 1, 2\}$, we can find a family of the Borel subset $Y_{2}$ and subrelations of $({\cal R}_{2})_{Y_{2}}$ satisfying $(\dagger)_{n_{0}-1}$ with respect to ${\cal S}_{i}^{s}$, ${\cal S}_{j}^{t}$ and ${\cal S}_{k}^{u}$.

On the other hand, we apply Theorem \ref{thm-even-n-2} to $M^{2}$ and the family ${\cal F}$. After permuting the indices $1, \ldots, n_{0}$, the family ${\cal F}$ satisfies that for any distinct $i, j, k\in \{ 1, \ldots, g_{0}(M^{2})\}$, there exist $s, t, u\in \{ 1, 2\}$ satisfying the following condition: if we have a family of the Borel subset $Y_{2}$ and subrelations of $({\cal R}_{2})_{Y_{2}}$ satisfying $(\dagger)_{n}$ with respect to ${\cal S}_{i}^{s}$, ${\cal S}_{j}^{t}$ and ${\cal S}_{k}^{u}$, then $n\leq n_{0}-2$, which is a contradiction.

Next, we assume that $3\leq g_{0}(M^{1})<g_{0}(M^{2})$ and deduce a contradiction. Note that $g_{0}(M^{1})$ (resp. $g_{0}(M^{2})$) is equal to the genus of $M^{1}$ (resp. $M^{2}$). We apply Corollary \ref{cor-kappa-even-class2} to $M^{1}$. Using $f\colon Y_{1}\rightarrow Y_{2}$ and Remark \ref{rem-dagger-wrt-restriction}, we have a family 
\[{\cal F}=\{ Y_{2}, \{ {\cal S}_{1}^{1}, {\cal S}_{1}^{2}\}, \ldots, \{ {\cal S}_{n_{0}}^{1}, {\cal S}_{n_{0}}^{2}\} \}\]
of the Borel set $Y_{2}$ and subrelations of $({\cal R}_{2})_{Y_{2}}$ satisfying $(\dagger)_{n_{0}}$ and the following condition: for any subset $F\subseteq \{ 1,\ldots, n_{0}\}$ with $|F|\geq g_{0}(M^{1})+1$, there exist distinct three elements $i, j, k\in F$ satisfying that for any $s, t, u\in \{ 1, 2\}$, we can find a family of the Borel subset $Y_{2}$ and subrelations of $({\cal R}_{2})_{Y_{2}}$ satisfying $(\dagger)_{n_{0}-1}$ with respect to the relations ${\cal S}_{i}^{s}$, ${\cal S}_{j}^{t}$ and ${\cal S}_{k}^{u}$. 

On the other hand, we apply Theorem \ref{thm-even-n-2} to $M^{2}$ and the family ${\cal F}$. After permuting the indices $1, \ldots, n_{0}$, the family ${\cal F}$ satisfies the same property as above, which is a contradiction since $g_{0}(M^{1})+1\leq g_{0}(M^{2})$.
\end{pf}

Combining Corollary \ref{cor-kappa-invariant}, Corollary \ref{cor-kappa-odd-class} (ii) and Corollary \ref{cor-kappa-even-class-mcg}, we obtain Theorem \ref{thm-main-class-final}.

\subsection{Geometric lemmas for the proof of Theorems \ref{thm-kappa-even-class1} and \ref{thm-kappa-even-class2}}

In this subsection, we denote by $M$ a surface of type $(g, p)$ with $\kappa(M)\geq 0$ even.

\begin{lem}\label{lem-sep-even1}
Let $\alpha_{1}, \alpha_{2}, \alpha_{3}\in V(C(M))$ be distinct three separating curves satisfying the following conditions: 
\begin{enumerate}
\item[(i)] If $M_{\alpha_{1}}=Q_{1}\sqcup Q_{1}'$, then $\alpha_{2}\in V(C(Q_{1}'))$;
\item[(ii)] If $(Q_{1}')_{\alpha_{2}}=Q_{2}\sqcup Q_{2}'$, then $\alpha_{3}\in V(C(Q_{2}'))$ and $Q_{2}$ contains $\alpha_{1}$ as a boundary component;
\item[(iii)] We denote $(Q_{2}')_{\alpha_{3}}=Q_{3}\sqcup Q_{4}$ and suppose that $Q_{3}$ contains $\alpha_{2}$ as a boundary component.  
\end{enumerate}
Suppose that both $\kappa(Q_{1})$ and $\kappa(Q_{4})$ are odd and that both $\kappa(Q_{2})$ and $\kappa(Q_{3})$ are even. Then there exists $\tau \in S(M)$ such that $\alpha_{1}, \alpha_{2}, \alpha_{3}\in \tau$ and $n(\tau)=n(M)-1$.
\end{lem}

\begin{figure}[h]
\centering

\epsfbox{plse1.1}
\caption{}
\label{picture-lem-sep-even1}
\end{figure}

\begin{pf}
Let us denote by $g_{i}$ and $p_{i}$ the genus and the number of boundary components of $Q_{i}$ for $i=1, 2, 3, 4$, respectively. Then we have two equalities $g_{1}+g_{2}+g_{3}+g_{4}=g$ and $p_{1}+p_{2}+p_{3}+p_{4}=p+6$. Hence, 
\begin{align*}
n(Q_{1})&+n(Q_{2})+n(Q_{3})+n(Q_{4})\\
        &=\sum_{i=1}^{4}(g_{i}+[(g_{i}+p_{i}-2)/2])\\
        &=\sum_{i=1, 4}(g_{i}+(g_{i}+p_{i}-3)/2)+\sum_{i=2, 3}(g_{i}+(g_{i}+p_{i}-2)/2)\\
        &=g+(g+p-2)/2-1=n(M)-1.
\end{align*} 
It completes the proof.
\end{pf}

\begin{lem}\label{lem-sep-even2}
Suppose $g\geq 2$. Let $\alpha_{1}, \alpha_{2}, \alpha_{3}\in V(C(M))$ be distinct three separating curves satisfying the following conditions: 
\begin{enumerate}
\item[(i)] If $M_{\alpha_{1}}=Q_{1}\sqcup Q_{1}'$, then $\alpha_{2}\in V(C(Q_{1}'))$;
\item[(ii)] If $(Q_{1}')_{\alpha_{2}}=Q_{2}\sqcup Q_{2}'$, then $\alpha_{3}\in V(C(Q_{2}'))$ and $Q_{2}$ contains $\alpha_{1}$ as a boundary component;
\item[(iii)] We denote $(Q_{2}')_{\alpha_{3}}=Q_{3}\sqcup Q_{4}$ and suppose that $Q_{3}$ contains $\alpha_{2}$ as a boundary component. 
\end{enumerate}
Suppose that both $Q_{1}$ and $Q_{4}$ are of type $(1, 1)$ and that both $\kappa(Q_{2})$ and $\kappa(Q_{3})$ are odd. Let $\beta_{1}\in V(C(Q_{1}))$ and $\beta_{4}\in V(C(Q_{4}))$. Then there exists $\tau \in S(M)$ such that $\beta_{1}, \alpha_{2}, \beta_{4}\in \tau$ and $n(\tau)=n(M)-1$.
\end{lem}

\begin{figure}[h]
\centering

\epsfbox{plse2.1}
\caption{}
\label{picture-lem-sep-even2}
\end{figure}

\begin{pf}
Let $(g_{1}, p_{1})$ (resp. $(g_{4}, p_{4})$) be the type of the surface $R_{1}$ (resp. $R_{4}$) constructed by gluing $Q_{1}$ and $Q_{2}$ on $\alpha_{1}$ (resp. $Q_{3}$ and $Q_{4}$ on $\alpha_{3}$) together. Remark that both $\kappa(R_{1})$ and $\kappa(R_{4})$ are odd by Lemma \ref{lem-kappa-sum}. Then we have two equalities $g_{1}+g_{4}=g$ and $p_{1}+p_{4}=p+2$. Hence,
\begin{align*}
n((R_{1})_{\beta_{1}})&+n((R_{4})_{\beta_{4}})\\
                      &=\sum_{i=1, 4}(g_{i}-1+[((g_{i}-1)+(p_{i}+2)-2)/2])\\
                      &=\sum_{i=1, 4}(g_{i}-1+(g_{i}+p_{i}-1)/2)\\
                      &=g+(g+p-2)/2-1=n(M)-1,
\end{align*}
which proves the lemma.
\end{pf}

\begin{lem}\label{lem-sep-even3}
Suppose $g\geq 2$. Let $\alpha_{1}, \alpha_{2}, \alpha_{3}, \alpha_{4}\in V(C(M))$ be distinct four separating curves satisfying the following conditions: 
\begin{enumerate}
\item[(i)] If $M_{\alpha_{1}}=Q_{1}\sqcup Q_{1}'$, then $\alpha_{2}\in V(C(Q_{1}'))$;
\item[(ii)] If $(Q_{1}')_{\alpha_{2}}=Q_{2}\sqcup Q_{2}'$, then $\alpha_{3}\in V(C(Q_{2}'))$ and $Q_{2}$ contains $\alpha_{1}$ as a boundary component;
\item[(iii)] If $(Q_{2}')_{\alpha_{3}}=Q_{3}\sqcup Q_{3}'$, then $\alpha_{4}\in V(C(Q_{3}'))$ and $Q_{3}$ contains $\alpha_{2}$ as a boundary component;
\item[(iv)] We denote $(Q_{3}')_{\alpha_{4}}=Q_{4}\sqcup Q_{5}$ and suppose that $Q_{4}$ contains $\alpha_{3}$ as a boundary component. 
\end{enumerate}
Suppose that both $\kappa(Q_{2})$ and $\kappa(Q_{4})$ are odd and that both $Q_{1}$ and $Q_{5}$ are of type $(1, 1)$. Let $\beta_{1}\in V(C(Q_{1}))$ and $\beta_{5}\in V(C(Q_{5}))$. Then there exists $\tau \in S(M)$ such that $\beta_{1}, \alpha_{2}, \alpha_{3}, \beta_{5}\in \tau$ and $n(\tau)=n(M)-1$.
\end{lem}

\begin{figure}[h]
\centering

\epsfbox{plse3.1}
\caption{}
\label{picture-lem-sep-even3}
\end{figure}

\begin{pf}
Note that $\kappa(Q_{3})$ is even by Lemma \ref{lem-kappa-sum}. Let $(g_{1}, p_{1})$ (resp. $(g_{5}, p_{5})$) be the type of the surface $R_{1}$ (resp. $R_{5}$) constructed by gluing $Q_{1}$ and $Q_{2}$ on $\alpha_{1}$ (resp. $Q_{4}$ and $Q_{5}$ on $\alpha_{4}$) together. Remark that both $\kappa(R_{1})$ and $\kappa(R_{5})$ are odd by Lemma \ref{lem-kappa-sum}. Let $(g_{3}, p_{3})$ be the type of $Q_{3}$. Then we have two equalities $g_{1}+g_{3}+g_{5}=g$ and $p_{1}+p_{3}+p_{5}=p+4$. Hence,
\begin{align*}
n((R_{1})_{\beta_{1}})&+n(Q_{3})+n((R_{5})_{\beta_{5}})\\
                      &=\sum_{i=1, 5}(g_{i}-1+[((g_{i}-1)+(p_{i}+2)-2)/2])+g_{3}+[(g_{3}+p_{3}-2)/2]\\
                      &=\sum_{i=1, 5}(g_{i}-1+(g_{i}+p_{i}-1)/2)+g_{3}+(g_{3}+p_{3}-2)/2\\
                      &=g+(g+p-2)/2-1=n(M)-1.
\end{align*}
It completes the proof.
\end{pf}


\subsection{The proof of Theorem \ref{thm-kappa-even-class1}}

Remark that the case of $g=2$ can be proved in Theorem \ref{thm-kappa-even-class2}. Hence, in this subsection, we show Theorem \ref{thm-kappa-even-class1} in only the cases of $g=0$ and $g=1$. 

We suppose that $M$ is a compact orientable surface of type $(g, p)$ and that $g\in \{ 0, 1\}$ and $\kappa(M)\geq 0$ is even. Let us denote by $n_{0}=n(M)$ and assume $n_{0}\geq 3$. We denote by $P_{1},\ldots, P_{p}$ all the boundary components of $M$. 

In what follows, we give a family 
\[\sigma =\{ \alpha_{1},\ldots, \alpha_{n_{0}-1}\}\in S(M),\ \ \{ \beta_{i}^{1}, \beta_{i}^{2}, Q_{i}\}_{i=1,\ldots, n_{0}}\]
satisfying the condition in Theorem \ref{thm-kappa-even-class1} for the surface $M$ in the cases of $g=0$ and $g=1$, respectively. We adapt the notation in Remark \ref{rem-notation-class}.
\vspace{1em}

\begin{figure}[hp]
\centering

\vspace*{2\baselineskip}
\epsfbox{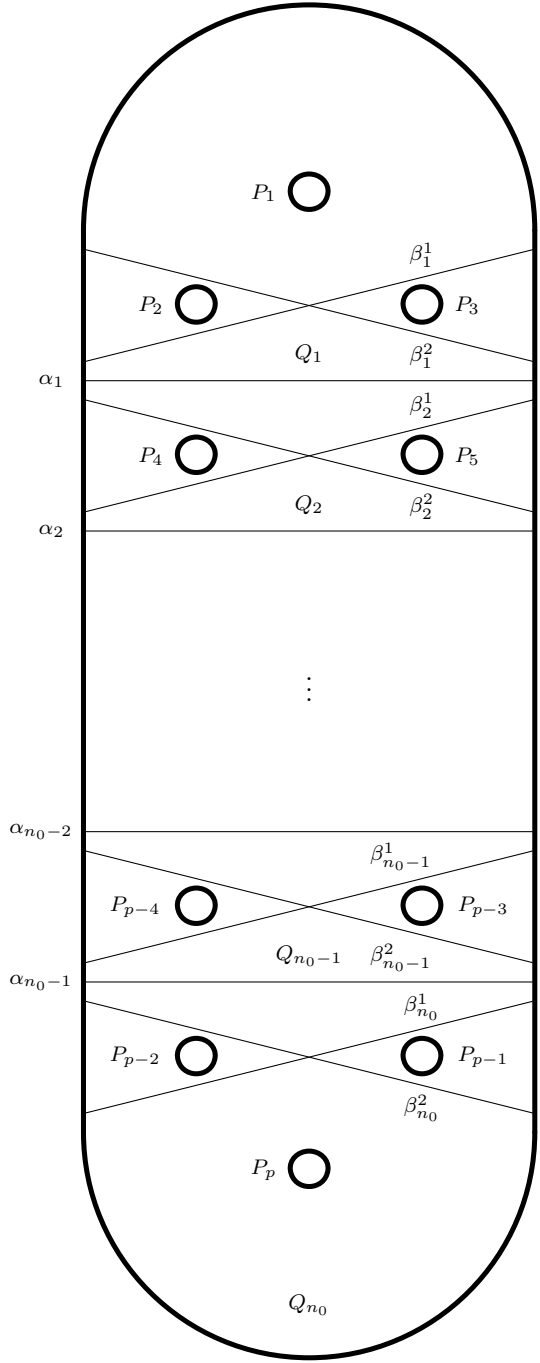}
\caption{The case of $g=0$}
\label{picture-kappa-even-g-0}
\end{figure}

\subsubsection{The case of $g=0$ (Figure \ref{picture-kappa-even-g-0})}\label{sub-even-g-0}

In this case, $p\in \{ 8, 10, 12, \ldots \}$ and $n_{0}=p/2-1$. First, we give $\sigma =\{ \alpha_{1}, \ldots, \alpha_{n_{0}-1}\} \in S(M)$ inductively as follows: choose $\alpha_{1}\in V(C(M))$ as a curve such that if we denote $M_{\alpha_{1}}=Q_{1}\sqcup Q_{1}'$, then $Q_{1}$ is of type $(0, 4)$ and satisfies $\partial Q_{1}=\{ \alpha_{1}, P_{1}, P_{2}, P_{3}\}$. 

For $i=2, 3,\ldots, n_{0}-1$, we choose $\alpha_{i}\in V(C(Q_{i-1}'))$ as a curve such that if we denote $(Q_{i-1}')_{\alpha_{i}}=Q_{i}\sqcup Q_{i}'$, then $Q_{i}$ is of type $(0, 4)$ and satisfies $\partial Q_{i}=\{ \alpha_{i-1}, \alpha_{i}, P_{2i}, P_{2i+1}\}$. 

Then $Q_{n_{0}-1}'$ is of type $(0, 4)$ and $\partial Q_{n_{0}-1}'=\{ \alpha_{n_{0}-1}, P_{p-2}, P_{p-1}, P_{p}\}$. We write $Q_{n_{0}}=Q_{n_{0}-1}'$. Define $\sigma =\{ \alpha_{1}, \ldots, \alpha_{n_{0}-1}\}\in S(M)$.   

Next, we choose two intersecting curves $\beta_{i}^{1}, \beta_{i}^{2}\in V(C(Q_{i}))$ satisfying the following conditions:
\begin{itemize}
\item The curves $\beta_{1}^{1}$ and $\beta_{1}^{2}$ are arbitrary two intersecting ones;
\item For $i\in \{ 2, \ldots, n_{0}-1\}$, the curve $\beta_{i}^{1}$ (resp. $\beta_{i}^{2}$) bounds $\alpha_{i-1}$, $P_{2i}$ (resp. $\alpha_{i-1}$, $P_{2i+1}$);
\item The curves $\beta_{n_{0}}^{1}$ and $\beta_{n_{0}}^{2}$ are arbitrary two intersecting ones.
\end{itemize}

We show Theorem \ref{thm-kappa-even-class1} for this choice. Take any $i, j, k\in \{ 1,\ldots, n_{0}\}$ with $i<j<k$ and any $s, t, u\in \{ 1, 2\}$.

\begin{figure}[hp]
\centering

\vspace*{2\baselineskip}
\epsfbox{pkeg1.1}
\caption{The case of $g=1$}
\label{picture-kappa-even-g-1}
\end{figure}

If we cut $M$ along $\beta_{i}^{s}$, then we have two components $R_{1}$ and $R_{1}'$ of types $(0, 2i+1)$ and $(0, p-2i+1)$, respectively, with $\beta_{j}^{t}\in V(C(R_{1}'))$. 

If we cut $R_{1}'$ along $\beta_{j}^{t}$, then we have two components $R_{2}$ and $R_{2}'$ of types $(0, 2(j-i)+2)$ and $(0, p-2j+1)$, respectively, with $\beta_{k}^{u}\in R_{2}'$. Moreover, $R_{2}$ has $\beta_{i}^{s}$ as a boundary component. 

If we cut $R_{2}'$ along $\beta_{k}^{u}$, then we have two components $R_{3}$ and $R_{4}$ of types $(0, 2(k-j)+2)$ and $(0, p-2k+1)$, respectively, such that $R_{3}$ has $\beta_{j}^{t}$ as a boundary component. Since both $\kappa(R_{1})$ and $\kappa(R_{4})$ are odd and both $\kappa(R_{2})$ and $\kappa(R_{3})$ are even, we can apply Lemma \ref{lem-sep-even1} and complete the proof of Theorem \ref{thm-kappa-even-class1}.
\vspace{1em}


\subsubsection{The case of $g=1$ (Figure \ref{picture-kappa-even-g-1})}

In this case, $p\in \{ 5, 7, 9, \ldots \}$ and $n_{0}=(p+1)/2$. Let $\bar{M}$ be  a surface of type $(0, p+3)$ with boundary components $\bar{P}_{1}, \ldots, \bar{P}_{p+3}$. Note $n(\bar{M})=n_{0}$. We apply the argument in \ref{sub-even-g-0} to $\bar{M}$ and give a family
\[\bar{\sigma}=\{ \bar{\alpha}_{1},\ldots, \bar{\alpha}_{n_{0}-1}\}\in S(\bar{M}),\ \ \{ \bar{\beta}_{i}^{1}, \bar{\beta}_{i}^{2}, \bar{Q}_{i}\}_{i=1,\ldots, n_{0}}.\] 
We construct $M$ by replacing $\bar{Q}_{n_{0}}$ with $M_{1, 1}$, that is, remove $\bar{Q}_{n_{0}}$ from $\bar{M}$ and glue the boundary component of $M_{1, 1}$ to $\bar{\alpha}_{n_{0}-1}$ which is a boundary component of $\bar{M}\setminus \bar{Q}_{n_{0}}$. We denote by $Q_{n_{0}}$ the new glued subsurface of $M$ of type $(1, 1)$ and by $\alpha_{n_{0}-1}\in V(C(M))$ the corresponding curve to $\bar{\alpha}_{n_{0}-1}$. 

Let us denote 
\begin{align*}
P_{i}=\bar{P}_{i} & \ \ \ {\rm for}\ i\in \{ 1, 2, \ldots, p\},\\
\alpha_{i}=\bar{\alpha}_{i} & \ \ \ {\rm for}\ i\in \{ 1, 2, \ldots, n_{0}-2\},\\
\beta_{i}^{s}=\bar{\beta}_{i}^{s},\ \ Q_{i}=\bar{Q}_{i} & \ \ \ {\rm for}\ i\in \{ 1, 2,\ldots, n_{0}-1\} {\rm \ and \ } s\in \{ 1, 2\}.
\end{align*}
Define $\sigma =\{ \alpha_{1}, \ldots, \alpha_{n_{0}-1}\}\in S(M)$. We choose arbitrary two intersecting curves $\beta_{n_{0}}^{1}, \beta_{n_{0}}^{2}\in V(C(Q_{n_{0}}))$.

We show Theorem \ref{thm-kappa-even-class1} for this choice. Take any $i, j, k\in \{ 1,\ldots, n_{0}\}$ with $i<j<k$ and any $s, t, u\in \{ 1, 2\}$. 

If we cut $M$ along $\beta_{i}^{s}$, then we have two components $R_{1}$ and $R_{1}'$ of types $(0, 2i+1)$ and $(1, p-2i+1)$, respectively, with $\beta_{j}^{t}\in V(C(R_{1}'))$. 

If we cut $R_{1}'$ along $\beta_{j}^{t}$, then we have two components $R_{2}$ and $R_{2}'$ of types $(0, 2(j-i)+2)$ and $(1, p-2j+1)$, respectively, with $\beta_{k}^{u}\in R_{2}'$. Moreover, $R_{2}$ has $\beta_{i}^{s}$ as a boundary component. 

Suppose $k=n_{0}$. Note that $\beta_{n_{0}}^{u}$ is a non-separating curve on $R_{2}'$. If we cut $R_{2}'$ along $\beta_{n_{0}}^{u}$, then the resulting surface $R_{3}$ is of type $(0, p-2j+3)$. Hence, we have
\begin{align*}
n(R_{1})&+n(R_{2})+n(R_{3})\\
        &=[(2i+1-2)/2]+[(2(j-i)+2-2)/2]+[(p-2j+3-2)/2]\\
        &=(p-1)/2=n_{0}-1.
\end{align*}
This shows the existence of $\tau \in S(M)$ in Theorem \ref{thm-kappa-even-class1}.

Next, suppose $k<n_{0}$. If we cut $R_{2}'$ along $\beta_{k}^{u}$, then we have two components $R_{3}$ and $R_{4}$ of types $(0, 2(k-j)+2)$ and $(1, p-2k+1)$, respectively, such that $R_{3}$ has $\beta_{j}^{t}$ as a boundary component. Since $\kappa(R_{1})$ and $\kappa(R_{4})$ are odd and both $\kappa(R_{2})$ and $\kappa(R_{3})$ are even, we can apply Lemma \ref{lem-sep-even1} and complete the proof of Theorem \ref{thm-kappa-even-class1}. 

Therefore, we complete the proof of Theorem \ref{thm-kappa-even-class1} in the cases of $g=0$ and $g=1$.


\begin{figure}[hp]
\centering

\vspace*{2\baselineskip}
\epsfbox{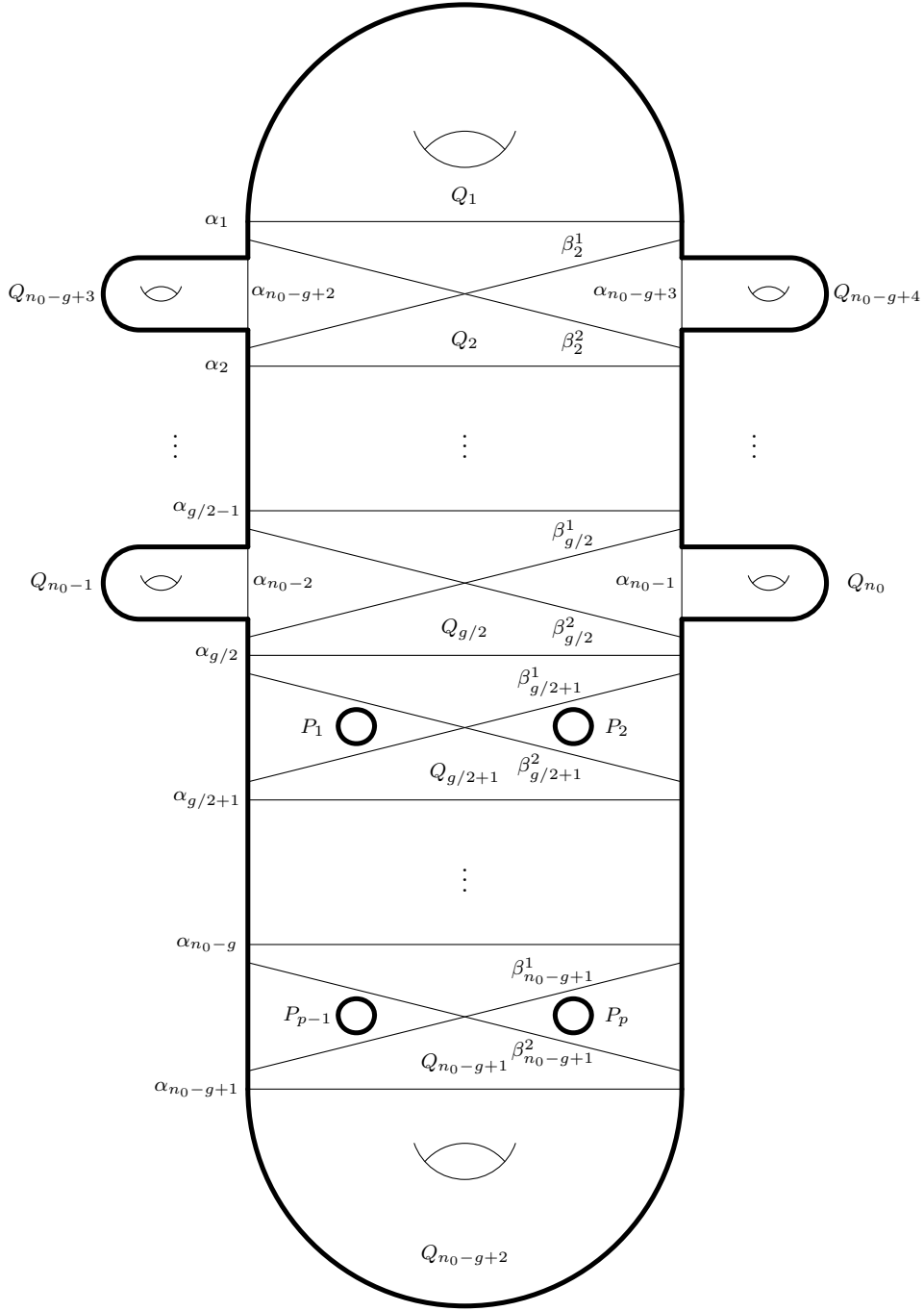}
\caption{The case where $g\geq 2$ is even}
\label{picture-kappa-even-g-even}
\end{figure}

\subsection{The proof of Theorem \ref{thm-kappa-even-class2}}

We suppose that $M$ is a compact orientable surface of type $(g, p)$ and that $g\geq 2$ and $\kappa(M)\geq 0$ is even. We assume that $M$ is not of type $(2, 0)$. Let us denote $n_{0}=n(M)$ and denote by $P_{1},\ldots, P_{p}$ all the boundary components of $M$. 

In what follows, we give a family 
\[\sigma =\{ \alpha_{1},\ldots, \alpha_{n_{0}-1}\}\in S(M),\ \ \{ \beta_{i}^{1}, \beta_{i}^{2}, Q_{i}\}_{i=1,\ldots, n_{0}}\]
satisfying the condition in Theorem \ref{thm-kappa-even-class2} for the surface $M$ in the cases where $g$ is even and odd, respectively. We adapt the notation in Remark \ref{rem-notation-class}.
\vspace{1em}

\subsubsection{The case where $g$ is even (Figure \ref{picture-kappa-even-g-even})}\label{sub-even-g-even}

Let $\bar{M}$ be a surface of type $(g, p+1)$ with boundary components $\bar{P}_{1}, \ldots, \bar{P}_{p+1}$. Note $n(\bar{M})=n_{0}$. We apply the argument in \ref{sub-odd-g-even} to $\bar{M}$ and give a family 
\[\bar{\sigma}=\{ \bar{\alpha}_{1},\ldots, \bar{\alpha}_{n_{0}-1}\}\in S(\bar{M}),\ \ \{ \bar{\beta}_{i}^{1}, \bar{\beta}_{i}^{2}, \bar{Q}_{i}\}_{i=1,\ldots, n_{0}}.\]
We construct $M$ by replacing the component $\bar{Q}_{n_{0}-g+2}$ of type $(1, 2)$ with $M_{1, 1}$, that is, remove $\bar{Q}_{n_{0}-g+2}$ from $\bar{M}$ and glue the boundary component of $M_{1, 1}$ to $\bar{\alpha}_{n_{0}-g+1}$ which is a boundary component of $\bar{M}\setminus \bar{Q}_{n_{0}-g+2}$. We denote by $Q_{n_{0}-g+2}$ the new glued subsurface of $M$ of type $(1, 1)$ and by $\alpha_{n_{0}-g+1}\in V(C(M))$ the corresponding curve to $\bar{\alpha}_{n_{0}-g+1}$.

Let us denote 
\begin{align*}
P_{i}=\bar{P}_{i} & \ \ \ {\rm for}\ i\in \{ 1, 2, \ldots, p\},\\
\alpha_{i}=\bar{\alpha}_{i} & \ \ \ {\rm for}\ i\in \{ 1, 2, \ldots, n_{0}-1\} \setminus \{ n_{0}-g+1\},\\
\beta_{i}^{s}=\bar{\beta}_{i}^{s}, \ \ Q_{i}=\bar{Q}_{i} & \ \ \ {\rm for}\ i\in \{ 1, 2,\ldots, n_{0}\} \setminus \{ n_{0}-g+2\} {\rm \ and \ } s\in \{ 1, 2\}.
\end{align*}
Define $\sigma =\{ \alpha_{1}, \ldots, \alpha_{n_{0}-1}\}\in S(M)$. We choose arbitrary two intersecting curves $\beta_{n_{0}-g+2}^{1}, \beta_{n_{0}-g+2}^{2}\in V(C(Q_{n_{0}-g+2}))$.
 
We show Theorem \ref{thm-kappa-even-class2} for this choice. Take any subset $F$ of $\{ 1, \ldots, n_{0}\}$ with $|F|\geq g+1$. It follows from $|F|\geq g+1$ that $|F\cap \{ 1, \ldots, n_{0}-g+2\}|\geq 3$. Let $i, j, k\in F\cap \{ 1,\ldots, n_{0}-g+2\}$ with $i<j<k$. Take any $s, t, u\in \{ 1, 2\}$. Remark that for each $l\in \{ 2, \ldots, n_{0}-g+1\}$ and $v\in \{ 1, 2\}$, if we cut $M$ along $\beta_{l}^{v}$, then we have two components of $M$ both of whose complexities are odd by Remark \ref{rem-lem-even-nn}.

Suppose $i\neq 1$ and $k\neq n_{0}-g+2$. When we cut $M$ along $\beta_{j}^{t}$, let $R_{1}'$ and $R_{2}'$ be the components of the cut surface with $\beta_{i}^{s}\in V(C(R_{1}'))$, $\beta_{k}^{u}\in V(C(R_{2}'))$. It follows from the above remark that both $\kappa(R_{1}')$ and $\kappa(R_{2}')$ are odd. The component $R_{1}$ of $(R_{1}')_{\beta_{i}^{s}}$ which does not have $\beta_{j}^{t}$ as a boundary component is also a component of $M_{\beta_{i}^{s}}$. Thus, $R_{1}$ has odd complexity and another component of $(R_{1}')_{\beta_{i}^{s}}$ has even complexity by Lemma \ref{lem-kappa-sum}. Similarly, the component of $(R_{2}')_{\beta_{k}^{u}}$ which does not have $\beta_{j}^{t}$ as a boundary component has odd complexity and another component of $(R_{2}')_{\beta_{k}^{u}}$ has even complexity. Thus, we can apply Lemma \ref{lem-sep-even1} and complete the proof of Theorem \ref{thm-kappa-even-class2}.

Suppose $i=1$ and $k=n_{0}-g+2$. Let $R_{1}'$ and $R_{2}'$ be the components of $M_{\beta_{j}^{t}}$ with $\beta_{i}^{s}\in V(C(R_{1}'))$, $\beta_{k}^{u}\in V(C(R_{2}'))$. Since $j\in \{ 2, \ldots, n_{0}-g+1\}$, both $\kappa(R_{1}')$ and $\kappa(R_{2}')$ are odd. If we cut $R_{1}'$ (resp. $R_{2}'$) along $\alpha_{1}$ (resp. $\alpha_{n_{0}-g+1}$), then we have two components $R_{1}$ and $R_{2}$ (resp. $R_{3}$ and $R_{4}$) with $R_{1}=Q_{1}$ (resp. $R_{4}=Q_{n_{0}-g+2}$). It follows from $\kappa(R_{1})=\kappa(R_{4})=0$ that both $\kappa(R_{2})$ and $\kappa(R_{3})$ are odd. Thus, we can apply Lemma \ref{lem-sep-even2} and complete the proof of Theorem \ref{thm-kappa-even-class2}.

Suppose $i=1$ and $k\neq n_{0}-g+2$. Let $R_{1}'$ and $R_{2}'$ be the components of $M_{\beta_{j}^{t}}$ with $\beta_{i}^{s}\in V(C(R_{1}'))$, $\beta_{k}^{u}\in V(C(R_{2}'))$. Since $j\in \{ 2, \ldots, n_{0}-g+1\}$, both $\kappa(R_{1}')$ and $\kappa(R_{2}')$ are odd. If we cut $R_{1}'$ along $\alpha_{1}$, then we have two components $R_{1}$ and $R_{2}$ with $R_{1}=Q_{1}$. It follows from $\kappa(R_{1})=0$ that $\kappa(R_{2})$ is odd. If we cut $R_{2}'$ along $\beta_{k}^{u}$, then we have two components $R_{3}$ and $R_{4}'$ such that $R_{3}$ has $\beta_{j}^{t}$ as a boundary component and $\alpha_{n_{0}-g+1}\in V(C(R_{4}'))$. Since $R_{4}'$ is also a component of $M_{\beta_{k}^{u}}$, we see that $\kappa(R_{4}')$ is odd. Let $R_{4}$ and $R_{5}$ be the components of $(R_{4}')_{\alpha_{n_{0}-g+1}}$ with $R_{5}=Q_{n_{0}-g+2}$. Since $\kappa(R_{4})$ is odd, we can apply Lemma \ref{lem-sep-even3} and complete the proof of Theorem \ref{thm-kappa-even-class2}.

When $i\neq 1$ and $k=n_{0}-g+2$, we can show Theorem \ref{thm-kappa-even-class2} similarly, using Lemma \ref{lem-sep-even3}.
\vspace{1em}


\subsubsection{The case where $g$ is odd}

Let $\bar{M}$ be a surface of type $(g-1, p+1)$ with boundary components $\bar{P}_{1}, \ldots, \bar{P}_{p+1}$. Note $n(\bar{M})=n_{0}-1$. We apply the argument in \ref{sub-even-g-even} to $\bar{M}$ and give a family 
\[\bar{\sigma}=\{ \bar{\alpha}_{1},\ldots, \bar{\alpha}_{n_{0}-2}\}\in S(\bar{M}),\ \ \{ \bar{\beta}_{i}^{1}, \bar{\beta}_{i}^{2}, \bar{Q}_{i}\}_{i=1,\ldots, n_{0}-1}.\]
We construct $M$ by gluing the boundary component of a new surface $Q_{n_{0}}$ of type $(1, 1)$ to $\bar{P}_{1}$. Let $\alpha_{n_{0}-1}\in V(C(M))$ be the corresponding curve to $\bar{P}_{1}$. Let us denote
\begin{align*}
P_{i-1}=\bar{P}_{i} & \ \ \ {\rm for}\ i\in \{ 2, 3, \ldots, p+1\},\\
\alpha_{i}=\bar{\alpha}_{i} & \ \ \ {\rm for}\ i\in \{ 1, 2,\ldots, n_{0}-2\},\\
\beta_{i}^{s}=\bar{\beta}_{i}^{s},\ \ Q_{i}=\bar{Q}_{i} & \ \ \ {\rm for}\ i\in \{ 1, 2,\ldots, n_{0}-1\} {\rm \ and \ } s\in \{ 1, 2\}.
\end{align*} 
Define $\sigma =\{ \alpha_{1}, \ldots, \alpha_{n_{0}-1}\}\in S(M)$. We choose arbitrary two intersecting curves $\beta_{n_{0}}^{1}, \beta_{n_{0}}^{2}\in V(C(Q_{n_{0}}))$. We can show Theorem \ref{thm-kappa-even-class2} for this choice as well as in \ref{sub-even-g-even}. 

This completes the proof of Theorem \ref{thm-kappa-even-class2}. \hfill $\square$


\appendix
\chapter{Amenability of a group action}\label{general-amenable} 

In this appendix, we summarize basic facts about amenability of an action of a discrete group on a standard Borel space with a quasi-invariant $\sigma$-finite measure introduced by Zimmer \cite{zim1}. A generalization of the notion to measured groupoids appears in \cite{ana}.


\section{Notation} 

Let ${\cal G}$ be a groupoid on the unit space ${\cal G}^{(0)}$. We denote the range and source maps by $r, s\colon {\cal G}\rightarrow {\cal G}^{(0)}$, respectively and define
\[{\cal G}^{(2)}=\{ (\gamma_{1}, \gamma_{2})\in {\cal G}\times {\cal G}: s(\gamma_{1})=r(\gamma_{2})\}.\]
We denote the composition of $\gamma_{1}$ and $\gamma_{2}$ with $(\gamma_{1}, \gamma_{2})\in {\cal G}^{(2)}$ by $\gamma_{1}\gamma_{2}$ and the inverse map by
\[I\colon {\cal G}\ni \gamma \mapsto \gamma^{-1}\in {\cal G}.\]
For $x\in {\cal G}^{(0)}$, we set ${\cal G}^{x}=r^{-1}(x)$ and ${\cal G}_{x}=s^{-1}(x)$.

We say that a set $Y$ is {\it fibered}\index{fibered} over a set $X$ if a surjective map $p\colon Y\rightarrow X$, called the {\it projection map}\index{projection map}, has been specified. When two sets $Y$ and $Z$ are fibered over $X$ via the maps $p\colon Y\rightarrow X$ and $q\colon Z\rightarrow X$, we define the {\it fibered product}\index{fibered!product} by
\[Y*Z=\{ (y, z)\in Y\times Z: p(y)=q(z)\}.\]
It is also a set fibered over $X$. A {\it left} ${\cal G}$-{\it space}\index{left ${\cal G}$-space} consists of a set $X$ fibered over ${\cal G}^{(0)}$ by a map $r=r_{X}\colon X\rightarrow {\cal G}^{(0)}$, called the {\it projection map}\index{projection map}, and a map $(\gamma, x)\in {\cal G}*X\mapsto \gamma x\in X$, called the {\it action map}\index{action map}, where ${\cal G}$ is fibered over ${\cal G}^{(0)}$ by the source map, such that the following equalities hold whenever they make sense: $\gamma(\gamma 'x)=(\gamma \gamma ')x$ and $ux=x$ for a unit $u\in {\cal G}^{(0)}$. 

We say that a groupoid is {\it Borel}\index{Borel groupoid} if all the associated spaces and maps are Borel spaces and Borel maps, respectively.  Similarly, we can define a Borel ${\cal G}$-space when ${\cal G}$ is a Borel groupoid. We say that ${\cal G}$ is {\it discrete}\index{discrete!groupoid} when each ${\cal G}^{x}$ is countable. In this case, we say that a Borel ${\cal G}$-space $X$ is {\it discrete}\index{discrete!Borel ${\cal G}$-space} if $r^{-1}_{X}(x)$ is countable for any $x\in {\cal G}^{(0)}$.

Let ${\cal G}$ be a discrete Borel groupoid. We say that a $\sigma$-finite measure $\mu$ on ${\cal G}^{(0)}$ is {\it quasi-invariant}\index{quasi-invariant measure} if the two measures $\tilde{\mu}$ and $\tilde{\mu}^{-1}=I_{*}\tilde{\mu}$ on ${\cal G}$ are equivalent, where $\tilde{\mu}$ is defined by
\[\tilde{\mu}(C)=\int_{{\cal G}^{(0)}}\sum_{\gamma :r(\gamma)=x}\chi_{C}(\gamma)d\mu (x),\]
for a Borel subset $C\subseteq {\cal G}$ and $\chi_{C}$ is the characteristic function on $C$. We say that $\mu$ is {\it invariant}\index{invariant!measure} if $\tilde{\mu}=\tilde{\mu}^{-1}$. A pair $({\cal G}, \mu)$ of a discrete Borel groupoid ${\cal G}$ and its quasi-invariant measure $\mu$ is called a {\it discrete measured groupoid}\index{discrete!measured groupoid}.

Let $X$ be a Borel ${\cal G}$-space. A $\sigma$-finite measure $\nu$ on $X$ is said to be {\it quasi-invariant} for the ${\cal G}$-action if $\nu$ is quasi-invariant for the Borel groupoid $X*{\cal G}$, where we regard  
\[X*{\cal G}=\{(x, \gamma)\in X\times {\cal G}: r(x)=r(\gamma )\},\]
as a Borel groupoid over $X$, with the range and source map $r(x, \gamma)=x$, $s(x, \gamma)=\gamma^{-1}x$, the composition $(x, \gamma_{1})\cdot (\gamma^{-1}_{1}x, \gamma_{2})=(x, \gamma_{1}\gamma_{2})$ and the inverse $I\colon (x, \gamma)\mapsto (\gamma^{-1}x, \gamma^{-1})$.

\begin{ex}
Let $G$ be a discrete group and $(X, \mu)$ be a Borel space with a $\sigma$-finite measure. Suppose that we have a non-singular Borel left $G$-action on $(X, \mu)$, where the word ``{\it non-singular}''\index{non-singular!action} means that the action preserves the class of the measure $\mu$. In this case, we say also that the measure $\mu$ is {\it quasi-invariant}\index{quasi-invariant measure} for the $G$-action. Such a pair $(X, G)$ is called a {\it transformation group}\index{transformation group}. Then the set $X\times G$ has a natural groupoid structure with the unit space $X$, the range and source maps $r(x, g)=x$, $s(x, g)=g^{-1}x$, the composition $(x, g_{1})\cdot (g_{1}^{-1}x, g_{2})=(x, g_{1}g_{2})$ and the inverse $I(x, g)=(g^{-1}x, g^{-1})$. This groupoid is called the {\it semi-direct product}\index{semi-direct product} of the transformation group $(X, G)$ and denoted by $X\rtimes G$.

Moreover, two measures on $X\rtimes G$ constructed from $\mu$,
\[\tilde{\mu}(C)=\int_{X}\sum_{g\in G}\chi_{C}(x, g)d\mu (x)\]
and $\tilde{\mu}^{-1}$ are equivalent by the assumption that the $G$-action is quasi-invariant. Thus, $(X\rtimes G, \mu)$ is a discrete measured groupoid.

\end{ex}

\begin{lem}\label{quasi-invariant}
Let $G$ be a discrete group and $(X, \mu)$ be a Borel space with a $\sigma$-finite measure. Suppose that we have a non-singular Borel left $G$-action on $(X, \mu)$. Let $Z$ be a countable left $G$-space. Denote $X\rtimes G$ by ${\cal G}$. Then $X\times Z$ is naturally a Borel ${\cal G}$-space consisting of 
\[r\colon X\times Z\rightarrow X, \ \ (x, z)\mapsto x,\]
\[{\cal G}*(X\times Z)\rightarrow X\times Z, \ (x, g)\cdot (g^{-1}x, z)=(x, gz).\]
Moreover, the measure $\nu$ on $X\times Z$ given by
\[\nu(C)=\int_{X}\sum_{z\in Z}\chi_{C}(x, z)d\mu (x)\]
is quasi-invariant for the ${\cal G}$-action. 
\end{lem}

\begin{pf}
The first assertion can be shown easily. It follows from 
\[(X\times Z)*{\cal G} =\{ ((x, z), (x, g))\in (X\times Z)\times {\cal G}\}\]
that $(X\times Z)*{\cal G}$ can be identified with $X\times Z\times G$. Then we have
\[r(x, z, g)=(x, z), \ \ s(x, z, g)=(g^{-1}x, g^{-1}z),\]
\[(x, z, g)^{-1}=(g^{-1}x, g^{-1}z, g^{-1}),\]
and
\begin{align*}
\tilde{\nu}(f) &= \int_{X\times Z}\sum_{\ta{\gamma \in (X\times Z)*{\cal G}}{r(\gamma)=(x, z)}}f(\gamma)d\nu(x, z)\\
                    &= \int_{X}\sum_{z\in Z}\sum_{g\in G}f(x, z, g)d\mu(x)                   
\end{align*}
for any non-negative Borel function $f$ on $(X\times Z)*{\cal G}$. The second assertion follows from the next equality:
\begin{align*}
\tilde{\nu}^{-1}(f) &= \int_{X}\sum_{z\in Z}\sum_{g\in G}f(g^{-1}x, g^{-1}z, g^{-1})d\mu(x)\\
                    &= \sum_{z\in Z}\int_{X}\sum_{g\in G}f(g^{-1}x, z, g^{-1})d\mu(x)\\
                    &= \sum_{z\in Z}\int_{X}\sum_{g\in G}f(x, z, g)\frac{d(\tilde{\mu}^{-1})}{d\tilde{\mu}}(x, g)d\mu(x)\\
                    &= \tilde{\nu}\left( f\frac{d(\tilde{\mu}^{-1})}{d\tilde{\mu}} \right)
\end{align*}   
for any non-negative Borel function $f$ on $(X\times Z)*{\cal G}$. 
\end{pf}


\section{Existence of invariant means}\label{section-amenability}

In this section, we collect some fundamental results about existence of invariant means for the ${\cal G}$-space $Y$ satisfying the following assumption: let $G$ be a discrete group and $(X, \mu)$ be a Borel space with a $\sigma$-finite measure. Assume that $G$ acts on $X$ and $\mu$ is quasi-invariant for the $G$-action. Let ${\cal G}$ denote the discrete measured groupoid $X\rtimes G$. Suppose that we have a countable $G$-space $Z$. Denote $X\times Z$ by $Y$, which is naturally a discrete Borel ${\cal G}$-space. We have a ${\cal G}$-quasi-invariant measure $\nu$ on $Y$ given in Lemma \ref{quasi-invariant}.

We introduce invariant means for the ${\cal G}$-space $Y$ of three types and state the equivalence of their existence. Note that $L^{\infty}(Y)$ has a $L^{\infty}(X)$-module structure by the projection map $r_{Y}\colon Y\rightarrow X$.

\begin{defn}[\ci{Definition 3.1.4}{ana}]
We say that a positive unital $L^{\infty}(X)$-linear map
\[m\colon L^{\infty}(Y)\rightarrow L^{\infty}(X)\]
(which is called a {\it mean}\index{mean}) is {\it invariant}\index{invariant!mean} if 
\[g\cdot m(f)=m(g\cdot f)\]
for any $g\in G$ and $f\in L^{\infty}(Y)$, where the $G$-actions on $L^{\infty}(Y)$ and $L^{\infty}(X)$ are given by
\begin{align*}
(g\cdot \varphi)(x, z)&=\varphi(g^{-1}x, g^{-1}z)\\
(g\cdot f)(x)&=f(g^{-1}x)
\end{align*}
for $x\in X$, $z\in Z$, $g\in G$, $\varphi \in L^{\infty}(Y)$ and $f\in L^{\infty}(X)$.
\end{defn}

Define $L^{\infty}(X, \ell^{1}(Z))$ as the Banach space consisting of measurable functions $f$ on $Y$ such that the essential supremum of $\sum_{z\in Z}|f(x, z)|$ for $x\in X$ is finite and whose norm is given by this quantity.

\begin{defn}[\ci{Definition 3.1.6}{ana}]\label{app-inv-mean}
A net $\{ f_{i}\}$ consisting of non-negative functions in $L^{\infty}(X, \ell^{1}(Z))$ with norm less than or equal to $1$ is called an {\it approximate weakly invariant mean}\index{approximate weakly invariant mean} if it satisfies the following two conditions:
\begin{enumerate}
\renewcommand{\labelenumi}{\rm(\roman{enumi})}
\item for any $\varphi \in L^{\infty}(Y)$, $h\in L^{1}(X)$ and $g\in G$, we have
\[\lim_{i}\int_{Y}h(x)\varphi(x, z)(f_{i}(g^{-1}x, g^{-1}z)-f_{i}(x, z))d\nu(x, z)=0;\]
\item for any $h\in L^{1}(X)$, we have
\[\lim_{i}\int_{Y}h(x)f_{i}(x, z)d\nu(x, z)=\int_{X}h(x)d\mu(x). \]
\end{enumerate}
\end{defn}

\begin{defn}[\ci{Definition 3.1.26}{ana}]
A family $m=\{ m_{x}: x\in X\}$ of states (or means) $m_{x}$ on $\ell^{\infty}(Z)$ satisfying the following conditions is called an {\it invariant measurable system of means}\index{invariant!measurable system of means}: for every $\varphi \in L^{\infty}(Y)$, 
\begin{enumerate}
\renewcommand{\labelenumi}{\rm(\roman{enumi})}
\item the function $x\mapsto m_{x}(\varphi(x, \cdot))$ is $\mu$-measurable;
\item $m_{g^{-1}x}(\varphi(g^{-1}x, \cdot))=m_{x}((g\cdot \varphi)(x, \cdot))$ for a.e. $(x, g)\in X\times G$.
\end{enumerate}
\end{defn}

\begin{thm}[\ci{Remark 3.2.6}{ana}]\label{ame-eq-cond}
Let $G$ be a discrete group and $(X, \mu)$ be a Borel space with a $\sigma$-finite measure. Assume that $G$ acts on $X$ and $\mu$ is quasi-invariant for the $G$-action. Let ${\cal G}$ denote the discrete measured groupoid $X\rtimes G$. Suppose that we have a countable $G$-space $Z$. Denote $X\times Z$ by $Y$. Let $\nu$ be the ${\cal G}$-quasi-invariant measure on $Y$ given in Lemma \ref{quasi-invariant}. 

Then the following conditions for the ${\cal G}$-space $Y$ are all equivalent:
\begin{enumerate}
\renewcommand{\labelenumi}{\rm(\roman{enumi})}
\item there exists an invariant mean;
\item there exists an approximate weakly invariant mean;
\item there exists an invariant measurable system of means.
\end{enumerate}
Moreover, in this case, we can choose an approximate weakly invariant mean which is normalized and indexed by integers.
\end{thm} 

\begin{defn}\label{defn-amenable}
We say that the ${\cal G}$-space $Y$ is {\it amenable}\index{amenable@amenable!g-space@${\cal G}$-space} when it satisfies the equivalent conditions in the above theorem. 

We say that the groupoid ${\cal G}=X\rtimes G$ is {\it amenable}\index{amenable!groupoid} if the ${\cal G}$-space $X\times G$ is amenable. In this case, we say that the action of $G$ on $(X, \mu)$ is {\it amenable (in a measurable sense)}\index{amenable!action in a measurable sense}. 
\end{defn}

\begin{ex}
\begin{enumerate}
\renewcommand{\labelenumi}{\rm(\roman{enumi})}
\item If $X$ is a single point, then amenability of $G$ as a group and amenability of ${\cal G}$ as a groupoid are equivalent.
\item If $G$ is amenable, then its non-singular Borel action on $(X, \mu)$ is always amenable.

\item Let $H$ be a subgroup of $G$. Then the set $G/H$ of left cosets is naturally a countable $G$-space. Therefore, $Y=X\times (G/H)$ can be seen as a ${\cal G}$-space with the natural quasi-invariant measure $\nu$ as in Lemma \ref{quasi-invariant}. 

If $H=G$, then $Y$ is always amenable. If $H$ is trivial, then $Y$ is amenable if and only if ${\cal G}$ is amenable. 
\end{enumerate}
\end{ex}

\begin{ex}\label{ex-amenable-action}
\begin{enumerate}
\renewcommand{\labelenumi}{\rm(\roman{enumi})}
\item The fundamental group of a closed Riemannian manifold of negative sectional curvature acts amenably on the boundary at infinity of the universal cover for some quasi-invariant measure on it \cite{spa}. 
\item The result of the example (i) is generalized to the action of a discrete subgroup of the group of isometries on a complete simply connected Riemannian manifold with negative sectional curvature bounded below and away from $0$ for any quasi-invariant measure on the boundary \cite{spa-zim}. 
\item Any hyperbolic group acts amenably on its boundary at infinity for any quasi-invariant measure on it \cite{adams1}. 
\end{enumerate}
\end{ex}

\begin{rem}\label{topologically-amenable-action}
For a discrete group $G$, set
\[{\rm Prob}(G)=\{ \mu \in \ell^{1}(G): \mu \geq 0, \sum_{g\in G}\mu(g)=1\}.\]
We equip ${\rm Prob}(G)$ with the norm topology. Note that the norm topology and the pointwise convergence topology coincide on ${\rm Prob}(G)$. There exists a natural action of $G$ on ${\rm Prob}(G)$.  

Suppose that we have a continuous action of $G$ on a compact Hausdorff space $X$. The action is said to be {\it (topologically) amenable}\index{amenable!action in a topological sense} if for every finite subset $F\subseteq G$ and $\varepsilon >0$, there exists a continuous map
\[\mu \colon X\ni x\rightarrow \mu_{x}\in {\rm Prob}(G)\]
such that 
\[\max_{s\in F}\sup_{x\in X}\Vert s\cdot \mu_{x}-\mu_{sx}\Vert_{1}< \varepsilon.\]
In this case, note that the $G$-stabilizer of every point on $X$ must be amenable by definition. We have the following relation between the topological and measurable versions of amenability: 
  
\begin{thm}[\ci{Theorem 3.3.7}{ana}]\label{thm-ana-amenable-action}
Given a discrete group $G$ and its continuous action on a compact Hausdorff space $X$, the action is topologically amenable if and only if it is amenable for any quasi-invariant measure on $X$. 
\end{thm}

Thus, all groups in Example \ref{ex-amenable-action} have property A by Theorem \ref{exact-property-A}. 
\end{rem}




\section{The fixed point property}

There exist two ways to define amenability of a discrete group $G$, the existence of an invariant mean on $\ell^{\infty}(G)$ and the fixed point property. In Section \ref{section-amenability}, we collected some results about the existence of various invariant means for a Borel action of a discrete group. In this section, we study the fixed point property introduced by Zimmer \cite{zim1}.

Let $G$ be a discrete group and $(X, \mu)$ be a standard Borel space with a $\sigma$-finite measure. Assume that $G$ acts on $X$ and $\mu$ is quasi-invariant for the $G$-action. 

Suppose that $E$ is a separable Banach space. We denote the group of isometric automorphisms of $E$ by ${\rm Iso}(E)$\index{$Iso E$@${\rm Iso}(E)$} given the strong operator topology and denote the closed unit ball of the dual $E^{*}$ of $E$ by $E_{1}^{*}$\index{$E 1 *$@$E_{1}^{*}$} given the weak* topology. Remark that ${\rm Iso}(E)$ is a separable metrizable group and the induced Borel structure is standard \cite[Lemma 1.1]{zim1}. Assume that we have a Borel cocycle $\alpha \colon X\times G\rightarrow {\rm Iso}(E)$. It means that $\alpha$ is a Borel map and satisfies that for any $g_{1}, g_{2}\in G$, we have
\[\alpha(x, g_{1})\alpha(g_{1}^{-1}x, g_{2})=\alpha(x, g_{1}g_{2})\]
for a.e.\ $x\in X$. Remark that we can show that there exists a $G$-invariant conull Borel subset $Y$ of $X$ such that the above cocycle identity is satisfied for any $g_{1}$, $g_{2}\in G$ and any $x\in Y$ since $G$ is countable. In what follows, although we often need to take a $G$-invariant Borel subset of $X$ in this way, we do not mention it. 

Suppose that for $x\in X$, we have a compact convex subset $A_{x}\subseteq E^{*}_{1}$ such that $\{ (x, a)\in X\times E^{*}_{1}: x\in X, \ a\in A_{x}\}$ is a Borel subset of $X\times E^{*}_{1}$ and for any $g\in G$, we have 
\[\alpha^{*}(x, g)A_{g^{-1}x}=A_{x}\]
for a.e.\ $x\in X$, where $\alpha^{*}(x, g)$ means the adjoint cocycle $\alpha^{*}(x, g)=(\alpha(x, g)^{-1})^{*}$. We denote by $F(X, \{ A_{x}\})$ the set of all measurable maps $\varphi \colon X\rightarrow E_{1}^{*}$ such that $\varphi(x)\in A_{x}$ for a.e.\ $x\in X$. We call $F(X, \{ A_{x}\})$ an {\it affine $G$-space}\index{affine!${\cal G}$-space} over $X$.

\begin{defn}[\ci{Definition 4.3.1}{zim3}]
We say that the $G$-action on $(X, \mu)$ has the {\it fixed point property}\index{fixed!point property} if every affine $G$-space $F(X, \{ A_{x}\})$ over $X$ has a fixed point, that is, there exists $\varphi \in F(X, \{ A_{x}\})$ such that for any $g\in G$, we have $\alpha^{*}(x, g)\varphi(g^{-1}x)=\varphi(x)$ for a.e.\ $x\in X$.
\end{defn}

This property is invariant under an isomorphism between the actions: suppose that we have two Borel $G$-actions on standard Borel spaces $(X, \mu)$ and $(Y, \nu)$ with quasi-invariant $\sigma$-finite measures $\mu$ and $\nu$ and there exist $G$-invariant conull Borel subsets $X'\subseteq X$, $Y'\subseteq Y$ and a Borel isomorphism $\varphi \colon X'\rightarrow Y'$ such that $\nu$ and $\varphi_{*}\mu$ are equivalent and $\varphi(gx)=g\varphi(x)$ for any $x\in X'$ and $g\in G$. If the $G$-action on $(X, \mu)$ has the fixed point property, then so does the $G$-action on $(Y, \nu)$.

\begin{thm}[\ci{Theorem 4.2.7}{ana}]\label{amenable-and-fixed} 
Let $G$ be a discrete group and $(X, \mu)$ be a standard Borel space with a $\sigma$-finite measure. Assume that $G$ acts on $X$ and $\mu$ is quasi-invariant for the $G$-action. Then the following two conditions are equivalent:
\begin{enumerate}
\renewcommand{\labelenumi}{\rm(\roman{enumi})}
\item The $G$-action has the fixed point property;
\item The $G$-action is amenable in the sense of Definition \ref{defn-amenable}.
\end{enumerate}
\end{thm}

We collect a few criteria for amenability of a group action.

\begin{prop}\label{lem-zim-eq-action}
Let $G$ be a discrete group and $H$ be its subgroup. Suppose that we have a non-singular Borel $H$-action on a standard Borel space $(X, \mu)$ with a $\sigma$-finite measure. Then the $H$-action on $X$ is amenable if and only if the diagonal action of $G$ on $X\times (G/H)$ is amenable.
\end{prop}

This proposition is proved in \cite[Lemma 3.3]{zim4} for the case where the ergodicity of the $H$-action is assumed. The following proof is essentially due to Zimmer \cite[Theorem 3.3]{zim1}. 

We recall the induced $G$-action from the $H$-action on $X$ \cite[Definition 4.2.21]{zim2}. We have the $H$-action on the product space $X\times G$ defined by the formula
\[h\cdot (x, g)=(hx, hg)\]
for $x\in X$, $h\in H$ and $g\in G$, which has a fundamental domain 
\[\{ (x, g_{n})\in X\times G: x\in X, \ n\in {\Bbb N}\},\]
where $\{ g_{n}\}_{n\in {\Bbb N}}$ is a set of all representatives of the set $H\backslash G$ of right cosets. Let $S$ be the quotient space of $X\times G$ by the $H$-action. On the other hand, we have the $G$-action on $X\times G$ defined by the formula
\[g\cdot (x, g')=(x, g'g^{-1})\]
for $x\in X$ and $g, g'\in G$. Since it commutes with the $H$-action, we have the induced $G$-action on $S$. The $G$-action on $S$ is called the {\it induced $G$-action}\index{induced $G$-action} from the $H$-action on $X$ (also called the {\it Mackey range}\index{Mackey range} of the cocycle from $X\times H$ into $G$ defined by the formula $(x, h)\mapsto h$ \cite[Definition 4.2.23]{zim2}).

\begin{lem}[\ci{Proposition 4.2.22}{zim2}]\label{lem-product-range}
With the above notation, the $G$-action on $S$ and the diagonal $G$-action on $X\times (G/H)$ are isomorphic.
\end{lem}

\begin{pf*}{\sc Proof of Proposition \ref{lem-zim-eq-action}.}
We may consider the $G$-action on $S$ instead of the $G$-action on $X\times (G/H)$ by Lemma \ref{lem-product-range}. It follows from \cite[Lemma 4.3.6]{zim2} that if the $G$-action on $S$ is amenable, then so is the $H$-action on $X$. We show the converse. 

Let $F(S, \{ A_{s}\}_{s\in S})$ be an affine $G$-space over $S$ with a cocycle $\alpha \colon S\times G\rightarrow {\rm Iso}(E)$. Define a Borel map $\beta \colon X\times H\rightarrow {\rm Iso}(E)$ by the formula
\[\beta(x, h)=\alpha(p(x, e), h)\]
for $x\in X$ and $h\in H$, where $p\colon X\times G\rightarrow S$ denotes the canonical projection. The map $\beta$ satisfies the cocycle identity because
\begin{align*}
\beta(x, h_{1})\beta(h_{1}^{-1}x, h_{2})&=\alpha(p(x, e), h_{1})\alpha(p(h_{1}^{-1}x, e), h_{2})\\
      &=\alpha(p(x, e), h_{1})\alpha(h^{-1}_{1}\cdot p(x, e), h_{2})\\
      &=\alpha(p(x, e), h_{1}h_{2})\\
      &=\beta(x, h_{1}h_{2}).
\end{align*} 
Put $B_{x}=A_{p(x, e)}$ for $x\in X$. Then 
\begin{multline*}
\{ (x, e, b)\in X\times \{ e\}\times E_{1}^{*}: x\in X,\ b\in B_{x}\}\\
             =(p|_{X\times \{ e\}}\times i)^{-1}(\{ (s, a)\in S\times E_{1}^{*}: s\in p(X\times \{ e\}),\  a\in A_{s}\}),
\end{multline*}
where $i$ denotes the identity map on $E_{1}^{*}$. This set is measurable because $p(X\times \{ e\})$ is a Borel subset of $S$ by Theorem \ref{thm-standard-borel-space} (iv). Thus, the set 
\[\{ (x, b)\in X\times E_{1}^{*}: x\in X,\ b\in B_{x}\}\]
is also measurable. Since $h\cdot p(x, e)=p(hx, e)$ for any $h\in H$ and $x\in X$, we see that $F(X, \{ B_{x}\})$ defines an affine $H$-space with the cocycle $\beta$ 

It follows from the amenability of the $H$-action on $X$ that there exists an element $\varphi \in F(X, \{ B_{x}\})$ such that 
\[\beta^{*}(x, h)\varphi(h^{-1}x)=\varphi(x)\]
for any $h\in H$ and $x\in X$. Define a Borel map $\psi \colon X\times G\rightarrow E_{1}^{*}$ by the formula
\[\psi(x, g)=(\alpha^{*}(p(x, e), g))^{-1}\varphi(x)\]
for $x\in X$ and $g\in G$. Since $\varphi(x)\in B_{x}=A_{p(x, e)}$, we have $\psi(x, g)\in A_{p(x, g)}$ for any $g\in G$ and $x\in X$. Moreover, for any $h\in H$, $g\in G$ and $x\in X$, we have
\begin{align*}
\psi(hx, hg)&=(\alpha^{*}(p(hx, e), hg))^{-1}\varphi(hx)=(\alpha^{*}(h\cdot p(x, e), hg))^{-1}\varphi(hx)\\   
  &=(\alpha^{*}(p(x, e), g))^{-1}(\alpha^{*}(h\cdot p(x, e), h))^{-1}\varphi(hx)\\
  &=(\alpha^{*}(p(x, e), g))^{-1}(\alpha^{*}(p(hx, e), h))^{-1}\varphi(hx)\\
  &=(\alpha^{*}(p(x, e), g))^{-1}\varphi(x)=\psi(x, g).
\end{align*}
Thus, we have the induced Borel map $\psi'\colon S\rightarrow E_{1}^{*}$ such that
\[\psi'(p(x, g))=\psi(x, g)\in A_{p(x, g)}\]
for any $g\in G$ and $x\in X$. Finally, for $x\in X$ and $g, g_{1}\in G$, we have 
\begin{align*}
\alpha^{*}(p(x, g_{1}), g)&\psi'(g^{-1}\cdot p(x, g_{1}))=\alpha^{*}(p(x, g_{1}), g)\psi'(p(x, g_{1}g))\\
  &=\alpha^{*}(g_{1}^{-1}\cdot p(x, e), g)\psi(x, g_{1}g)\\
  &=(\alpha^{*}(p(x, e), g_{1}))^{-1}\alpha^{*}(p(x, e), g_{1}g)(\alpha^{*}(p(x, e), g_{1}g))^{-1}\varphi(x)\\
  &=(\alpha^{*}(p(x, e), g_{1}))^{-1}\varphi(x)=\psi(x, g_{1})=\psi'(p(x, g_{1})).
\end{align*}
This means that the $G$-action on $S$ is amenable.
\end{pf*}

Let $G$ be a discrete group and $(X, \mu)$ be a standard Borel space with a $\sigma$-finite measure. Assume that $G$ acts on $X$ and $\mu$ is quasi-invariant for the $G$-action. Let $H$ be a subgroup of $G$. We denote the semi-direct product of the transformation group $(X, G)$ by ${\cal G}$.

\begin{prop}\label{rel-ame}
Under the above assumption, suppose that the ${\cal G}$-space $X\times (G/H)$ is amenable and that the $H$-action on $X$ is amenable. Then the $G$-action on $X$ is amenable. 
\end{prop}

\begin{pf}
By assumption, we have an invariant mean 
\[m\colon L^{\infty}(X\times (G/H))\rightarrow L^{\infty}(X)\]
for the ${\cal G}$-space $X\times (G/H)$. It follows from Proposition \ref{lem-zim-eq-action} that we have an invariant mean 
\[M\colon L^{\infty}(X\times (G/H)\times G)\rightarrow L^{\infty}(X\times (G/H))\]
for the $G$-action on $X\times (G/H)$. We denote by $K$ the $G$-invariant subspace of $L^{\infty}(X\times (G/H)\times G)$ consisting of $\varphi \in L^{\infty}(X\times (G/H)\times G)$ such that $\varphi(x, z_{1}, g)=\varphi(x, z_{2}, g)$ for a.e.\ $x\in X$ and any $g\in G$, $z_{1}, z_{2}\in G/H$. Then $K$ can be identified with $L^{\infty}(X\times G)$ as $G$-spaces. 

Thus, the restriction to $K$ of the composition 
\[m\circ M\colon L^{\infty}(X\times (G/H)\times G)\rightarrow L^{\infty}(X)\]
is an invariant mean
\[L^{\infty}(X\times G)\rightarrow L^{\infty}(X)\]
for the $G$-action on $X$. 
\end{pf}

\begin{prop}\label{prop-ame-cri-other}
Let $G$, $H$ be two discrete groups and suppose that we have a homomorphism from $G$ to $H$ with kernel amenable and a non-singular amenable action of $H$ on a standard Borel space $(X, \mu)$ with a $\sigma$-finite measure. Then the induced action of $G$ on $(X, \mu)$ is also amenable.
\end{prop}

\begin{pf}
Let $p\colon G\rightarrow H$ be the homomorphism. We may assume that $p$ is surjective. Let us denote the kernel of $p$ by $N$. It follows from Lemma \ref{lem-zim-eq-action} that the diagonal action of $G$ on $X\times (G/N)$ is amenable since $N$ is amenable. Thus, we have an invariant mean
\[L^{\infty}(X\times (G/N)\times G)\rightarrow L^{\infty}(X\times (G/N))\]
for the $G$-action. Since the groups $G/N$ and $H$ are isomorphic and the action of $H$ on $X$ is amenable, we have an invariant mean
\[L^{\infty}(X\times (G/N))\rightarrow L^{\infty}(X)\]
for the $H$-action. As in the proof of Proposition \ref{rel-ame}, we can get an invariant mean 
\[L^{\infty}(X\times G)\rightarrow L^{\infty}(X)\]
for the $G$-action on $X$ if we compose the above two invariant means and restrict it to an appropriate subspace.
\end{pf}

The following theorem is another criterion different from Propositions \ref{rel-ame} and \ref{prop-ame-cri-other} for amenability of a group action:

\begin{thm}[\ci{Corollary C}{aeg}]\label{aeg-thm}
Let $G$ be a discrete group acting on two standard Borel spaces $(X, \mu)$ and $(Y, \nu)$ with quasi-invariant measures $\mu$ and $\nu$. Let $f\colon X\rightarrow Y$ be a $G$-equivariant Borel map and assume that $f_{*}\mu =\nu$. If the $G$-action on $Y$ is amenable, then the $G$-action on $X$ is also amenable.
\end{thm}

\chapter{Measurability of the map associating image measures}\label{app-measurability}

In this appendix, we present some results concerning measurability of maps between spaces of probability measures. Since the author could not find any appropriate reference for them, proofs are presented although he believes that they are folklore among specialists. For a locally compact Polish space $\Omega$, we denote by $M(\Omega)$\index{$M \ X $@$M(\Omega)$} the space of regular probability Borel measures on $\Omega$ equipped with the weak* topology. Remark that any finite positive Borel measure on a Polish space is regular \cite[Theorems 17.10, 17.11]{kechris}.

\begin{prop}\label{prop-first-measurable}
Let $\Omega_{1}$, $\Omega_{2}$ be two locally compact Polish spaces. If $f\colon \Omega_{1}\rightarrow \Omega_{2}$ is a Borel map, then the induced map $f_{*}\colon M(\Omega_{1})\rightarrow M(\Omega_{2})$ is measurable.
\end{prop}

\begin{prop}\label{prop-second-measurable}
Let $\Omega$ be a locally compact Polish space. Then for any non-negative valued Borel function $\varphi$ on $\Omega$, the function $F_{\varphi}$ on $M(\Omega)$ defined by the formula $F_{\varphi}(\mu)=\mu(\varphi)$ for $\mu \in M(\Omega)$ is measurable.
\end{prop}

It follows from Proposition \ref{prop-second-measurable} that the $\sigma$-algebra on $M(\Omega)$ is equal to the minimal one for which the function $F_{\varphi}$ is measurable for any non-negative Borel function $\varphi$ on $\Omega$. Thus, the collection of Borel subsets in $\Omega$ of the form $\{ \omega \in \Omega : F_{\varphi}(\omega)> a\}$ for $\varphi$ and $a\in {\Bbb R}$ generates the Borel structure on $M(\Omega)$. 
 
Using this observation, we can deduce Proposition \ref{prop-first-measurable} from Proposition \ref{prop-second-measurable} easily. Therefore, we need to show only Proposition \ref{prop-second-measurable}.

For the proof, we use the following theorem called Dynkin's $\pi$-$\lambda$ theorem. We prepare some definitions necessary for the statement. 

\begin{defn}
Let $\Omega$ be a set and ${\cal P}$ be a collection of subsets in $\Omega$.  
\begin{enumerate}
\item[(i)] The collection ${\cal P}$ is said to be a {\it $\pi$-system}\index{pi-system@$\pi$-system} if it is closed under taking an intersection, that is, if $A$, $B\in {\cal P}$, then $A\cap B\in {\cal P}$.
\item[(ii)] The collection ${\cal P}$ is said to be a {\it $\lambda$-system}\index{lambda-system@$\lambda$-system} if it satisfies the following conditions:
\begin{enumerate}
\item[(a)] The set $\Omega$ is in ${\cal P}$;
\item[(b)] If $A$, $B\in {\cal P}$ and $A\subseteq B$, then we have $B\setminus A\in {\cal P}$; 
\item[(c)] If $\{ A_{n}\}_{n\in {\Bbb N}}$ is an increasing sequence in ${\cal P}$ and we denote $A=\bigcup A_{n}$, then we have $A\in {\cal P}$.
\end{enumerate}
\end{enumerate}
\end{defn}

\begin{thm}[\ci{Theorem 10.1 (iii)}{kechris}]
Let $\Omega$ be a set and ${\cal P}$, ${\cal L}$ be collections of subsets in $\Omega$. If ${\cal P}$ is a $\pi$-system and ${\cal L}$ is a $\lambda$-system which contains ${\cal P}$, then the $\sigma$-algebra generated by ${\cal P}$ is contained in the collection ${\cal L}$.  
\end{thm}

\begin{pf*}{{\sc Proof of Proposition \ref{prop-second-measurable}.}}
It suffices to show that the function $F_{A}$ on $M(\Omega)$ defined by the formula $F_{A}(\mu)=\mu (A)$ for a Borel subset $A$ of $\Omega$ is measurable. 

Let ${\cal L}$ denote the collection of all Borel subsets $A$ such that the function $F_{A}$ is measurable. First, we show that ${\cal L}$ is a $\lambda$-system.

It is clear that $\Omega \in {\cal L}$. If $A$, $B\in {\cal L}$ and $A\subseteq
B$, then we have $F_{B\setminus A}=F_{B}-F_{A}$, which is a measurable function. If $\{ A_{n}\}_{n\in {\Bbb N}}$ is an increasing sequence in ${\cal P}$ and we denote $A=\bigcup A_{n}$, then 
\begin{align*}
\{ \mu \in M(\Omega): F_{A}(\mu)>a\} &=\{ \mu \in M(\Omega): \mu(A)>a\} \\
             &=\bigcup_{n\in {\Bbb N}}\{ \mu \in M(\Omega): \mu(A_{n})>a\}  
\end{align*}
for any $a\in {\Bbb R}$, which is a Borel subset of $M(\Omega)$.

Finally, we show that the function $F_{U}$ is measurable for any open subset $U$ of $\Omega$. It completes the proof since the collection of all open subsets in $\Omega$ is a $\pi$-system and generates the Borel structure of $\Omega$.

We can take an increasing sequence $\{ f_{n}\}_{n\in {\Bbb N}}$ of non-negative valued continuous functions on $\Omega$ vanishing at infinity which converges to the characteristic function $\chi_{U}$ on $U$ pointwise, using the separability of $\Omega$ and Urysohn's lemma. Then 
\begin{align*}
\{ \mu \in M(\Omega): F_{U}>a\} &=\{ \mu \in M(\Omega): \mu(\chi_{U})>a\} \\
              &= \bigcup_{n\in {\Bbb N}}\{ \mu \in M(\Omega): \mu(f_{n})>a\}
\end{align*}
for any $a\in {\Bbb R}$, which is open in $M(\Omega)$.
\end{pf*}

\begin{cor}\label{cor-third-measurable}
Let $\Omega$ be a locally compact Polish space and $A$ be a Borel subset of $\Omega$. Then the subset $M(A)$ of $M(\Omega)$ consisting of all probability Borel measures $\mu$ with $\mu(A)=1$ is measurable.
\end{cor}


\chapter{Exactness of the mapping class group}\label{app-exact}

In this appendix, we show that the mapping class group is exact, using the functions constructed in Chapters \ref{chapter-property-A} and \ref{chapter-indec}. This claim answers a question in \cite[Section 6, Question 5]{bell-fuji} (see also Remark \ref{rem-mcg-exact}). The following proposition was communicated to us by Ozawa \cite{ozawa5}. The author is grateful to him for allowing us to present it here. The reader should be referred to Chapter \ref{chapter-property-A}, Section \ref{sec:pro} and Remark \ref{topologically-amenable-action} for exactness of a discrete group and topological amenability of a group action.

\begin{prop}\label{prop-ozawa}
Let $G$ be a discrete group acting on $X$ and $K$, where $X$ is a compact Hausdorff space and $K$ is a countable set. Assume that for any finite subset $F\subseteq G$ and $\varepsilon >0$, there exists a Borel map 
\[\mu \colon X\ni x\rightarrow \mu_{x}\in {\rm Prob}(K)\] 
such that
\[\max_{s\in F}\sup_{x\in X}\Vert s\cdot \mu_{x}-\mu_{sx}\Vert_{1}< \varepsilon.\]
Moreover, assume that for any $k\in K$, the stabilizer
\[\{ g\in G: gk=k\} \]
is exact. Then $G$ is exact. 
\end{prop}

We use the following two lemmas for the proof:

\begin{lem}[\ci{Proposition 11}{ozawa4}]
Let $G$ be a discrete group acting on $X$, $Y$ and $K$, where $X$ and $Y$ are compact Hausdorff spaces and $K$ is a countable set. Assume that for any finite subset $F\subseteq G$ and $\varepsilon >0$, there exists a Borel map 
\[\mu \colon X\ni x\rightarrow \mu_{x}\in {\rm Prob}(K)\] 
such that
\[\max_{s\in F}\sup_{x\in X}\Vert s\cdot \mu_{x}-\mu_{sx}\Vert_{1}< \varepsilon.\]
Moreover, assume that for any $k\in K$, the action of the stabilizer
\[\{ g\in G: gk=k\} \]
on $Y$ is topologically amenable. Then the diagonal action of $G$ on $X\times Y$ is topologically amenable.
\end{lem}

Given a set $Z$, we denote by $\beta Z$\index{$\ b Z $@$\beta Z$} its {\it Stone-\v Cech compactification}\index{Stone-\v Cech compactification}. Recall that $\beta Z$ is a compact Hausdorff space equipped with an inclusion of the discrete space $Z$ as an open dense subset. It has the universal property that any map from $Z$ into a compact Hausdorff space $X$ extends to a continuous map from $\beta Z$ into $X$.

\begin{lem}
Let $G$ be a discrete group and $H$ be a subgroup of $G$. If $H$ is exact, then the action of $H$ on $\beta G$ is topologically amenable.
\end{lem}

\begin{pf}
Since $H$ is exact, the action of $H$ on $\beta H$ is topologically amenable by Theorem \ref{exact-property-A}. Construct an $H$-equivariant map $G\rightarrow H$ and extend it to the $H$-equivariant continuous map $\beta G\rightarrow \beta H$. It follows from amenability of the action of $H$ on $\beta H$ that we can show amenability of the action of $H$ on $\beta G$, using the above map.   
\end{pf}

Proposition \ref{prop-ozawa} follows from these two lemmas by taking $Y$ to be $\beta G$. Using this proposition and the following lemma, we show exactness of the mapping class group.

\begin{lem}\label{lem-exact-reducible}
Let $M$ be a compact orientable surface with $\kappa(M)\geq 0$. If the mapping class group $\Gamma(Q)$ is exact for any surface $Q$ with $\kappa(Q)<\kappa (M)$, then every reducible subgroup $\Gamma$ in $\Gamma(M)$ is also exact.
\end{lem}

\begin{pf}
Note that if $H$ is a subgroup of a discrete group $G$ of finite index, then $H$ is exact if and only if $G$ is exact \cite[Theorems 4.1, 4.5]{kirch-wass}. Hence, we may assume $\Gamma$ to be a subgroup of $\Gamma(M; m)$ for some integer $m\geq 3$. Suppose that $\sigma \in S(M)$ satisfies $g\sigma =\sigma$ for any $g\in \Gamma$. It follows from Remark \ref{rem-comp-leave} that we have a natural homomorphism 
\[p_{\sigma}\colon \Gamma \rightarrow \prod_{i}\Gamma(Q_{i}),\]
where $\{ Q_{i}\}$ is the set of all components of the surface obtained by cutting $M$ along a realization of $\sigma$. It follows from $\kappa(Q_{i})<\kappa(M)$ (see Lemma \ref{main-geometric-lem} (i), (ii)) that all the $\Gamma(Q_{i})$ are exact. Note that the kernel of $p_{\sigma}$ is contained in the free abelian group generated by the Dehn twists of curves in $\sigma$ by Proposition \ref{lem-blm} and in particular, is amenable. Remark that if a discrete group $G$ has a normal subgroup $N$ and both $N$ and $G/N$ are exact, then $G$ is exact (see \cite[Theorem 5.1]{kirch-wass}). Thus, $\Gamma$ is exact.
\end{pf}

\begin{thm}
If $M$ is a compact orientable surface with $\kappa(M)\geq 0$, then the mapping class group $\Gamma(M)$ is exact.
\end{thm}

\begin{pf}
We show this by induction on the number $\kappa(M)$. If $\kappa(M)=0$, then $\Gamma(M)$ is almost isomorphic to $SL(2, {\Bbb Z})$ and thus, is exact.

Suppose that $\kappa(M)>0$ and $\Gamma(Q)$ is exact for any surface $Q$ with $0\leq \kappa(Q)<\kappa(M)$. Note that if $Q$ is a surface with $\kappa(Q)<0$, then it is clear that $\Gamma(Q)$ is exact. Using the Borel functions $f_{n}$ on $\partial C\times V(C)$ constructed in the proof of Theorem \ref{amenable-action-cc-non-excep} and the $\Gamma(M)$-equivariant continuous map $\pi \colon {\cal MIN}\rightarrow \partial C$ in Chapter \ref{chapter:amenable-action}, Section \ref{boundary-of-curve-complex}, we can construct Borel maps
\[\nu^{n}\colon {\cal MIN}\ni a\mapsto \nu_{a}^{n}\in {\rm Prob}(V(C))\]
for $n\in {\Bbb N}$ such that
\[\lim_{n\rightarrow \infty}\sup_{a\in {\cal MIN}}\sum_{x\in V(C)}\vert \nu^{n}_{a}(x)-\nu^{n}_{g^{-1}a}(g^{-1}x)\vert =0\]   
for any $g\in \Gamma(M)$.

On the other hand, we have a $\Gamma(M)$-equivariant Borel map
\[H\colon {\cal PMF}\setminus {\cal MIN}\rightarrow S(M)\]
(see Definition \ref{construction-of-H}). Combining these two maps, we can construct Borel maps
\[\mu^{n}\colon {\cal PMF}\ni a\mapsto \mu_{a}^{n}\in {\rm Prob}(S(M))\]
for $n\in {\Bbb N}$ such that 
\[\lim_{n\rightarrow \infty}\sup_{a\in {\cal PMF}}\sum_{x\in S(M)}\vert \mu^{n}_{a}(x)-\mu^{n}_{g^{-1}a}(g^{-1}x)\vert =0\]   
for any $g\in \Gamma(M)$. Moreover, for any $\sigma \in S(M)$, the stabilizer
\[\{ g\in \Gamma(M): g\sigma =\sigma \}\]
is exact by Lemma \ref{lem-exact-reducible} and our assumption. Hence, it follows from Proposition \ref{prop-ozawa} that $\Gamma(M)$ is exact.
\end{pf}

\chapter{The cost and $\ell^{2}$-Betti numbers of the mapping class group}\label{app-cost}

Levitt \cite{lev} introduced and Gaboriau \cite{gab-cost} developed generalities of {\it cost}\index{cost} of a discrete measure-preserving equivalence relation and a discrete group. Cost of an equivalence relation is a numerical invariant and is an analogous notion to the minimal number of generators of a discrete group. Using this invariant, Gaboriau showed that free groups of different ranks can not have isomorphic essentially free, measure-preserving ergodic actions \cite[Corollaire 1]{gab-cost} and that the free group of finite rank and the free group of infinite rank are not measure equivalent \cite[Propositions VI. 6, 9]{gab-cost}. We recommend the reader to consult \cite{gab-cost} and \cite{kec} for more details about cost.

Gaboriau verified that the $\ell^{2}$-{\it Betti numbers}\index{l2-Betti number@$\ell^{2}$-!Betti number} of a discrete group are invariant under measure equivalence in the following sense: if two discrete groups $\Gamma_{1}$ and $\Gamma_{2}$ are measure equivalent, then there exists a positive constant $c$ such that
\[\beta_{n}(\Gamma_{1})=c\beta_{n}(\Gamma_{2})\]
for all $n\in {\Bbb N}$, where $\beta_{n}(\Lambda)$\index{$\ b n \ C$@$\beta_{n}(\Gamma)$} denotes the $n$-th $\ell^{2}$-Betti number of a discrete group $\Lambda$ \cite[Th\'eor\`eme 6.3]{gab-l2}. Thanks to this fact, we can obtain various results about measure equivalence \cite[Corollaires 0.3, 0.4, 0.5 and 0.6]{gab-l2}. The $\ell^{2}$-{\it Euler characteristic}\index{l2-Euler characteristic@$\ell^{2}$-!Euler characteristic} of a group is defined by the alternating sum of its $\ell^{2}$-Betti numbers as in the case of ordinary Euler characteristics. It follows from the above Gaboriau's theorem that the sign of the $\ell^{2}$-Euler characteristic is invariant under measure equivalence.

In this appendix, we compute the cost and comment on the $\ell^{2}$-Betti numbers and the $\ell^{2}$-Euler characteristic of the mapping class group. In what follows, we treat only a discrete measure-preserving equivalence relation on a standard Borel space with a non-atomic finite positive measure.

\section{The cost of the mapping class group}

\begin{defn}
Let $(X, \mu)$ be a standard Borel space with a non-atomic finite positive measure. 
\begin{enumerate}
\item[(i)] A {\it graphing}\index{graphing} $\Phi$ on $(X, \mu)$ is a countable family $\Phi =\{ \varphi_{i}\colon A_{i}\rightarrow B_{i}\}_{i\in I}$ of partial Borel isomorphisms $\varphi_{i}$ on $X$ preserving $\mu$ with domain ${\rm dom}(\varphi_{i})=A_{i}$ and range ${\rm ran}(\varphi_{i})=B_{i}$.
\item[(ii)] Given a graphing $\Phi$ on $(X, \mu)$ and a relation ${\cal R}$ on $(X, \mu)$, we say that ${\cal R}$ is {\it generated}\index{generated} by $\Phi$ if ${\cal R}$ is the minimal relation satisfying $(x, \varphi(x))\in {\cal R}$ for any $\varphi \in \Phi$ and a.e.\ $x\in {\rm dom}(\varphi)$. In this case, we write ${\cal R}_{\Phi}={\cal R}$ and call $\Phi$ a {\it graphing}\index{graphing!for a relation} for ${\cal R}$.
\item[(iii)] The {\it cost}\index{cost!of a graphing} ${\cal C}_{\mu}(\Phi)$ of a graphing $\Phi$ on $(X, \mu)$ is defined by the formula
\[{\cal C}_{\mu}(\Phi)=\sum_{\varphi \in \Phi}\mu({\rm dom}(\varphi)).\]\index{$C-l-U$@${\cal C}_{\mu}(\Phi)$}
\item[(iv)] The {\it cost}\index{cost!of a discrete measure-preserving equivalence relation} ${\cal C}_{\mu}({\cal R})$ of a relation ${\cal R}$ on $(X, \mu)$ is defined by the formula
\[{\cal C}_{\mu}({\cal R})=\inf \{ {\cal C}_{\mu}(\Phi): \Phi \ {\rm is \ a \ graphing \ for \ }{\cal R}\}.\]\index{$C-lR$@${\cal C}_{\mu}({\cal R})$}
\item[(v)] The {\it cost}\index{cost!of a discrete group} ${\cal C}(\Gamma)$ of a discrete group $\Gamma$ is defined by the formula 
\[{\cal C}(\Gamma)=\inf_{\cal R}{\cal C}_{\mu}({\cal R}),\]\index{$C-C$@${\cal C}(\Gamma)$}
where ${\cal R}$ runs through all relations generated by essentially free, measure-preserving Borel actions of $\Gamma$ on a standard probability space $(X, \mu)$.
\item[(vi)] We say that a discrete group $\Gamma$ has {\it fixed price}\index{fixed!price} if ${\cal C}_{\mu}({\cal R})={\cal C}(\Gamma)$ for any relation ${\cal R}$ generated by an essentially free, measure-preserving Borel action of $\Gamma$ on a standard probability space $(X, \mu)$.
\end{enumerate}   
\end{defn}

\begin{rem}
In the definition (iv), note that there always exists a graphing for ${\cal R}$ because we have a discrete group $G$ consisting of Borel automorphisms on $X$ preserving $\mu$ such that the relation generated by $G$ is equal to ${\cal R}$ (see \cite[Theorem 1]{fm} or Example \ref{ex-thm-fm}). 

In the definition (v), remark that there always exists such an action for any $\Gamma$. When $\Gamma$ is infinite, one of the examples is the shift action on $(F^{\Gamma}, \nu)$, where $F$ is a finite set with cardinality more than one and $\nu$ is the product measure of some probability measure on $F$ which assigns non-zero measure to each element in $F$. Moreover, the action is mixing and thus, ergodic \cite[Remarque I.7]{gab-cost}, \cite[Lemma 1.17]{sauer}.  
\end{rem}

\begin{prop}[\ci{Invariance II.2}{gab-cost}]
If ${\cal R}_{1}$ and ${\cal R}_{2}$ are two isomorphic ergodic relations on standard probability spaces $(X_{1}, \mu_{1})$ and $(X_{2}, \mu_{2})$, respectively, then ${\cal C}_{\mu_{1}}({\cal R}_{1})={\cal C}_{\mu_{2}}({\cal R}_{2})$.
\end{prop}

\begin{thm}[\cite{lev}, \ci{Proposition III.3}{gab-cost}]\label{thm-lev-amenable}
Let us denote by ${\cal R}$ a relation on a standard Borel space $(X, \mu)$ with a finite positive measure. Let us denote
\[d({\cal R})=\mu(X)-\sum_{n\in ({\Bbb N}\setminus \{ 0\})\cup \{ \infty \}}(1/n)\mu(X_{n}),\]
where $X_{n}$ denotes the Borel subset $\{ x\in X: |{\cal R}x|=n\}$ of $X$. Then
\begin{enumerate}
\item[(i)] we have ${\cal C}_{\mu}({\cal R})\geq d({\cal R})$.
\item[(ii)] the relation ${\cal R}$ is amenable if and only if there exists a graphing $\Phi$ such that ${\cal C}_{\mu}(\Phi)=d({\cal R})$. In this case, ${\cal C}_{\mu}({\cal R})=d({\cal R})$. In particular, if ${\cal R}$ is amenable and recurrent, then ${\cal C}_{\mu}({\cal R})=\mu(X)$.  
\end{enumerate} 
\end{thm}

\begin{cor}[\cite{lev}, \ci{Propri\'et\'es VI.3}{gab-cost}]
All amenable groups have fixed price. If a discrete group $\Gamma$ is finite, then ${\cal C}(\Gamma)=1-(1/|\Gamma|)$ and if $\Gamma$ is infinite and amenable, then ${\cal C}(\Gamma)=1$.
\end{cor}

\begin{defn}[\ci{D\'efinition IV.5, 9}{gab-cost}]\label{defn-free-product-relation}
Let ${\cal R}$ be a relation on a standard Borel space $(X, \mu)$ with a finite positive measure. Let ${\cal R}_{1}$ and ${\cal R}_{2}$ be subrelations of ${\cal R}$ such that ${\cal R}_{1}$ and ${\cal R}_{2}$ generate ${\cal R}$.
\begin{enumerate}
\item[(i)] A finite sequence $(x_{1}, x_{2}, \ldots, x_{n})$ of points in $X$ is said to be {\it reduced}\index{reduced} if the following conditions are satisfied:
\begin{enumerate}
\item[(a)] each $(x_{i}, x_{i+1})$ is either in ${\cal R}_{1}$ or in ${\cal R}_{2}$;
\item[(b)] two successive pairs $(x_{i}, x_{i+1})$, $(x_{i+1}, x_{i+2})$ belong to distinct relations ${\cal R}_{1}$, ${\cal R}_{2}$;
\item[(c)] $x_{i}\neq x_{i+1}$ for each $i$.
\end{enumerate}
\item[(ii)] We say that ${\cal R}$ is the {\it free product}\index{free product of relations} of ${\cal R}_{1}$ and ${\cal R}_{2}$ if for a.e.\ reduced sequence $(x_{1}, x_{2}, \ldots, x_{n})$ with $n\geq 2$, we have $x_{1}\neq x_{n}$. In this case, we write ${\cal R}={\cal R}_{1}*{\cal R}_{2}$\index{$R 1 b R 2$@${\cal R}_{1}*{\cal R}_{2}$}.
\end{enumerate}
\end{defn}

Similarly, we can define the amalgamated product, HNN extension of subrelations (see \cite[D\'efinition IV.6, 20]{gab-cost}). 
  
\begin{ex}[\ci{Exemple IV.8}{gab-cost}]   
Let $\Gamma =\Gamma_{1}*\Gamma_{2}$ be the free product of two discrete groups $\Gamma_{1}$ and $\Gamma_{2}$. If ${\cal R}$ is a relation generated by an essentially free, measure-preserving action of $\Gamma$ and ${\cal R}_{i}$ denotes the subrelation of ${\cal R}$ generated by the action of $\Gamma_{i}$ for $i=1, 2$, then ${\cal R}={\cal R}_{1}*{\cal R}_{2}$. We recommend the reader to see \cite{gab-survey} for other examples of free products of equivalence relations. 
\end{ex}

\begin{thm}[\ci{Th\'eor\`eme IV.15}{gab-cost}]\label{thm-cost-free-product}
If ${\cal R}={\cal R}_{1}*{\cal R}_{2}$ on $(X, \mu)$ and both the costs of ${\cal R}_{1}$ and ${\cal R}_{2}$ are finite, then ${\cal C}_{\mu}({\cal R})={\cal C}_{\mu}({\cal R}_{1})+{\cal C}_{\mu}({\cal R}_{2})$.
\end{thm}

\begin{cor}[\ci{Proposition VI.9}{gab-cost}]
The free group of finite rank $n$ has fixed price and its cost is equal to $n$.
\end{cor}

The following proposition will be needed to compute the cost of the mapping class group:

\begin{prop}[\ci{Crit\`eres VI.24 (1)}{gab-cost}]\label{generator-cost}
If a discrete group $\Gamma$ is generated by a finite family of elements $g_{1}, g_{2}, \ldots, g_{n}$ of infinite order such that $g_{j}$ commutes with $g_{j+1}$ for each $j$, then $\Gamma$ has fixed price and ${\cal C}(\Gamma)=1$.
\end{prop}

We can find generators of the mapping class group $\Gamma(M)$ of a non-exceptional surface $M$ satisfying the hypothesis in Proposition \ref{generator-cost} as follows: let ${\rm P}\Gamma(M)$\index{$P Gamma M$@${\rm P}\Gamma(M)$} be the subgroup of $\Gamma(M)$ consisting of all elements fixing each boundary component of $M$, which is called the {\it pure mapping class group}\index{pure!mapping class group} of $M$. Note that ${\rm P}\Gamma(M)$ is naturally isomorphic to $\Gamma(M)$ when $M$ has at most one boundary component. 

It is well-known that ${\rm P}\Gamma(M)$ is generated by a finite family $\{ g_{1}, \ldots, g_{n}\}$ of the Dehn twists around some simple closed curves on $M$ (e.g., see \cite[Section 3]{kor2}). Since the curve complex $C(M)$ of $M$ is connected by Theorem \ref{thm:mm1}, we can find a path $l\colon \{ 0, 1, \ldots, m\}\rightarrow V(C(M))$ on $C(M)$ through all the curves corresponding to $g_{i}$, where the word ``path'' means that $d(l(j), l(j+1))=1$ for each $j$. Needless to say, the set of all the Dehn twists of vertices through which $l$ passes generates ${\rm P}\Gamma(M)$, and the Dehn twists about curves $l(j)$ and $l(j+1)$ commute for each $j$ since the two curves are disjoint. Moreover, all Dehn twists have infinite order. 

Assume that $M$ has more than one boundary components. For each two distinct boundary components $P_{1}$, $P_{2}$ of $M$, let $\alpha(P_{1}, P_{2})\in V(C)$ be a separating curve such that if $M$ is cut along $\alpha(P_{1}, P_{2})$, then one component $Q$ of the resulting surface $Q\sqcup R$ is a pair of pants whose boundary components are $P_{1}$, $P_{2}$ and the corresponding one to $\alpha(P_{1}, P_{2})$. Let $h(P_{1}, P_{2})\in \Gamma(M)$ be an element satisfying the following conditions:
\begin{enumerate}
\item[(i)] $h(P_{1}, P_{2})$ preserves $\alpha(P_{1}, P_{2})$ and exchanges $P_{1}$ and $P_{2}$;
\item[(ii)] $h(P_{1}, P_{2})^{2}$ is the Dehn twist about $\alpha(P_{1}, P_{2})$;
\item[(iii)] $h(P_{1}, P_{2})$ is the identity on $R$.
\end{enumerate} 
It is clear that the union $\{ g_{1}, \ldots, g_{n}\} \cup \{ h(P_{1}, P_{2})\}_{P_{1}, P_{2}}$ generates $\Gamma(M)$. Note that $V(C(R))$ is non-empty since $\kappa(M)>0$ and the Dehn twist about any element in $V(C(R))$ commutes with $h(P_{1}, P_{2})$. Therefore, using the connectedness of $C(M)$ in the same way as above, we can find a sequence $\{ h_{1}, \ldots, h_{m}\}$ of elements in $\Gamma(M)$ of infinite order such that they generate $\Gamma(M)$ and $h_{i}h_{i+1}=h_{i+1}h_{i}$ for any $i$. Hence, we have shown the following assertion:

\begin{prop}
The mapping class group of a non-exceptional surface has fixed price and its cost is equal to $1$.
\end{prop} 

The costs of the mapping class groups of exceptional surfaces have already been calculated in \cite{gab-cost}. The mapping class groups of $M_{0, 4}$, $M_{1, 0}$ and $M_{1, 1}$ are (almost) isomorphic to $SL(2, {\Bbb Z})$, which has fixed price and cost $1+1/12$ \cite[Proposition VI.9]{gab-cost} (see also \cite[Th\'eor\`eme VI.19]{gab-cost}). The mapping class groups of other surfaces are finite.


\section{The $\ell^{2}$-Betti numbers of the mapping class group}

For the definition of $\ell^{2}$-Betti numbers of discrete groups (resp. equivalence relations), see \cite{cheeger-gromov} and \cite{luck} (resp. \cite{gab-l2}). Let us denote by $\beta_{n}(\Gamma)$, $\beta_{n}({\cal R})$\index{$\ b n \ C$@$\beta_{n}(\Gamma)$}\index{$\ b n R$@$\beta_{n}({\cal R})$} the $n$-th $\ell^{2}$-Betti numbers for a discrete group $\Gamma$ and a discrete measure-preserving equivalence relation ${\cal R}$ on a standard probability space $(X, \mu)$. If the following sums converge, then we define the {\it $\ell^{2}$-Euler characteristics}\index{l2-Euler characteristic@$\ell^{2}$-!Euler characteristic} for $\Gamma$ and ${\cal R}$ by the formulas
\[\chi^{(2)}(\Gamma)=\sum_{n=0}^{\infty}(-1)^{n}\beta_{n}(\Gamma),\ \ \chi^{(2)}({\cal R})=\sum_{n=0}^{\infty}(-1)^{n}\beta_{n}({\cal R}),\]\index{$\ v2-C$@$\chi^{(2)}(\Gamma)$}\index{$\ v2R$@$\chi^{(2)}({\cal R})$}
respectively.

\begin{thm}[\ci{Corollaire 3.16, Th\'eor\`eme 6.3}{gab-l2}]\label{thm-gab-l2}
\begin{enumerate}
\item[(i)] If a relation ${\cal R}$ on a standard probability space is generated by an essentially free, measure-preserving action of a discrete group $\Gamma$, then 
\[\beta_{n}({\cal R})=\beta_{n}(\Gamma)\]
for each $n\in {\Bbb N}$.
\item[(ii)] If two discrete groups $\Gamma_{1}$ and $\Gamma_{2}$ are measure equivalent, then there exists a positive constant $c$ such that
\[\beta_{n}(\Gamma_{1})=c\beta_{n}(\Gamma_{2}),\ \ \chi^{(2)}(\Gamma_{1})=c\chi^{(2)}(\Gamma_{2})\]
for each $n\in {\Bbb N}$.
\end{enumerate}
\end{thm}

For the first $\ell^{2}$-Betti number and cost, the following relation holds:

\begin{thm}[\ci{Corollaire 3.23}{gab-l2}] 
If $\Gamma$ is a discrete group, then 
\[{\cal C}(\Gamma)-1\geq \beta_{1}(\Gamma)-\beta_{0}(\Gamma).\]
\end{thm}

Note that $\beta_{0}(\Gamma)=1/|\Gamma |$ if $|\Gamma |<\infty$ and $\beta_{0}(\Gamma)=0$ if $|\Gamma|=\infty$. 

Given a discrete group $\Gamma$ which has a torsion-free subgroup of finite index, we define its {\it virtual Euler characteristic}\index{virtual!Euler characteristic} $\chi(\Gamma)$\index{$\ v-C$@$\chi(\Gamma)$} by the formula
\[\chi(\Gamma)=\frac{1}{[\Gamma :\Gamma']}{\rm ch}(\Gamma'),\]
where ${\rm ch}(\Gamma')$ denotes the usual Euler characteristic of $\Gamma'$ and the right hand side is independent of the choice of $\Gamma'$ (see \cite[Section 6.7]{ivanov2} and the references therein). We have the following identification between the $\ell^{2}$-Euler characteristic and the virtual Euler characteristic of $\Gamma$:

\begin{thm}[\ci{Proposition 0.4 and (0.25)}{cheeger-gromov}]\label{virtual-euler}
For a discrete group $\Gamma$ for which the $\ell^{2}$-Euler characteristic and the virtual Euler characteristic can be defined, we have the equality
\[\chi^{(2)}(\Gamma)=\chi(\Gamma).\]
In particular, the sign of the virtual Euler characteristic is invariant under measure equivalence.
\end{thm}

The virtual Euler characteristic for the mapping class group has been calculated as follows:

\begin{thm}[\cite{harer-zagier}, \ci{Theorem 3.9, Corollary 3.10}{penner}]\label{thm-chi-mcg}
Let $M$ be a compact orientable surface of type $(g, p)$ and $\Gamma(M)$ be the mapping class group of $M$. Let $B_{n}$ denote the $n$-th Bernoulli number.
\begin{enumerate}
\item[(i)] If $g>1$ and $p=0$, then 
\[\chi(\Gamma(M))=\frac{B_{2g}}{4g(g-1)}.\]
\item[(ii)] If $g\geq 0$, $p\geq 1$ and $2g-2+p>0$, then 
\[\chi(\Gamma(M))=(-1)^{p}\frac{(p+2g-3)!(2g-1)}{p!(2g)!}B_{2g}.\]
\end{enumerate}  
\end{thm}
Remark that $(-1)^{n-1}B_{2n}>0$ if $n\geq 1$ and $B_{0}=1$. 

Next, we comment on the $\ell^{2}$-Betti numbers of the mapping class group. First, we recall some general theorems about $\ell^{2}$-cohomology (see \cite{cheeger-gromov} and \cite{gromov-kahler} for the terminology mentioned below). 

Let $X=(X, g)$ be a complete simply connected manifold with K\"ahler metric $g$ whose K\"ahler form is $d$(bounded) and $m=2n$ be the real dimension of $X$. Let $\Gamma$ be a discrete group of isometries of $X$ satisfying the following conditions:
\begin{enumerate}
\item[(i)] $X/\Gamma$ has finite volume;
\item[(ii)] the sectional curvature of $X$ is bounded above and below;
\item[(iii)] the injectivity radius of $X$ is bounded away from zero.
\end{enumerate} 
In this case, we say that $(X, \Gamma)$ satisfies $(\star)$. 

Let us denote by ${\cal H}^{p}$ the space of harmonic $L^{2}$-forms on $X$ of degree $p$ and denote by $\bar{H}^{p}_{(2)}(X: \Gamma)$ the $L^{2}$-cohomology space for $X$ and $\Gamma$ of degree $p$. They are $\Gamma$-modules in the sense of \cite{cheeger-gromov}. We denote by $\dim_{\Gamma}A$ the $\Gamma$-dimension for a $\Gamma$-module $A$.

\begin{thm}[\ci{Theorem 5.1}{cheeger-gromov}]
If a pair $(X, \Gamma)$ of a complete simply connected manifold $X$ with K\"ahler metric and a discrete group $\Gamma$ of isometries of $X$ satisfies $(\star)$, then the two $\Gamma$-modules ${\cal H}^{p}$ and $\bar{H}^{p}_{(2)}(X:\Gamma)$ are canonically isomorphic and $\dim_{\Gamma}{\cal H}^{p}$ is finite for any $p\in {\Bbb N}$. In particular, if $X$ is contractible, then 
\[\beta_{p}(\Gamma)=\dim_{\Gamma}{\cal H}^{p}\] 
for any $p\in {\Bbb N}$.
\end{thm}

\begin{thm}[\ci{Theorem 2.5, Remark 2.5.A (a)}{gromov-kahler}]
Let us suppose that a pair $(X, \Gamma)$ of a complete simply connected manifold $X$ with K\"ahler metric and a discrete group $\Gamma$ of isometries of $X$ satisfies $(\star)$. Let $2n$ be the real dimension of $X$. If $p\neq n$, then ${\cal H}^{p}=0$ and ${\cal H}^{n}\neq 0$. Thus, we have $\dim_{\Gamma}{\cal H}^{n}\neq 0$.
\end{thm}

On the other hand, McMullen proved the following theorem:

\begin{thm}[\ci{Theorem 1.1}{mcmullen}]
Let $M$ be an orientable surface of genus $g$ and with $p$ punctures. Suppose $2g-2+p>0$. Then there exists a K\"ahler metric $h$ on the Teichm\"uller space ${\cal T}(M)$ of $M$ comparable to the Teichm\"uller metric such that $(({\cal T}(M), h), \Gamma(M))$ satisfies $(\star)$.
\end{thm}

Note that if $M$ is an orientable surface of genus $g$ and with $p$ punctures, then the Teichm\"uller space ${\cal T}(M)$ of $M$ is homeomorphic to ${\Bbb R}^{6g-6+2p}$ and $\Gamma(M)$ acts on ${\cal T}(M)$ properly discontinuously (see \cite[Chapter 8, Theorem 6]{gardiner}). Remark that if $M'$ is a compact orientable surface of genus $g$ and with $p$ boundary components, then the two mapping class groups $\Gamma(M)$ and $\Gamma(M')$ are isomorphic (see Remark \ref{rem-mcg-exact}). Hence, we can calculate the $\ell^{2}$-Betti numbers of the mapping class group as follows, using Theorem \ref{thm-chi-mcg}:

\begin{cor}
Let $M$ be a compact orientable surface of type $(g, p)$ and $\Gamma(M)$ be the mapping class group of $M$. Suppose $2g-2+p>0$. Let $B_{n}$ denote the $n$-th Bernoulli number.
\begin{enumerate}
\item[(i)] If $p=0$, then
\[\beta_{3g-3}(\Gamma(M))=\frac{|B_{2g}|}{4g(g-1)}\]
and $\beta_{n}(\Gamma(M))=0$ for all $n\in {\Bbb N}\setminus \{ 3g-3\}$.
\item[(ii)] If $p\geq 1$, then 
\[\beta_{3g-3+p}(\Gamma(M))=\frac{(p+2g-3)!(2g-1)}{p!(2g)!}|B_{2g}|\]
and $\beta_{n}(\Gamma(M))=0$ for all $n\in {\Bbb N}\setminus \{ 3g-3+p\}$.
\end{enumerate}
\end{cor}


\chapter{A group-theoretic argument for Chapter 5}\label{app-group}

In this appendix, we give a direct proof of the following theorem. Although we use Lemma \ref{main-geometric-lem} in the course of the proof, it is a geometric lemma and can be read independently of the contents in Chapter \ref{chapter-best}.

\begin{thm}\label{main-group-version}
Let $M$ be a compact orientable surface of type $(g, p)$ satisfying $\kappa(M)\geq 0$. Let $\Gamma(M)$ be the mapping class group of $M$. 
\begin{enumerate}
\item[(i)] If $G_{1}\times \cdots \times G_{n}$ is a subgroup of $\Gamma(M)$ with each $G_{i}$ non-amenable, then 
\[n\leq g+\left[ \frac{g+p-2}{2}\right].\]
\item[(ii)] There exists a subgroup of the form $G_{1}\times \cdots \times G_{n}$ with each $G_{i}$ non-amenable and 
\[n=g+\left[ \frac{g+p-2}{2}\right].\]
\end{enumerate}
\end{thm}

Here, for $a\in {\Bbb R}$, we denote by $[a]$ the maximal integer less than or equal to $a$. For the proof, we prepare a few results.

\begin{thm}[\ci{Corollary 2}{mc}]\label{centralizer}
The normalizer of the cyclic group generated by a pseudo-Anosov element is virtually infinite cyclic.
\end{thm}

Let $M$ be a compact orientable surface of type $(g, p)$ with $\kappa(M)\geq 0$ and $\Gamma(M)$ be the mapping class group of $M$. Let $\Gamma =\Gamma(M; m)$ with $m\geq 3$ (see Theorem \ref{pure-important}). For a reducible subgroup $G$ of $\Gamma$, let us denote by $\sigma(G)$ the canonical reduction system for $G$ (see Definition \ref{def-crs-of-group} and Theorem \ref{thm-crs-non-empty}).

\begin{lem}\label{crs-invariant}
Let $G$ be a subgroup of $\Gamma$ and $H$ be an infinite reducible normal subgroup of $G$. Then $G$ is reducible and any element of $G$ fixes each curve in $\sigma(H)$.   
\end{lem}

This lemma follows from Remark \ref{rem-main-lem2}.

Let $G$ be an infinite reducible subgroup of $\Gamma$ and 
\[M_{\sigma(G)}=\bigsqcup_{i}Q_{i}\]
be the decomposition of $M$ by (cutting $M$ along a realization of) $\sigma(G)$. By Remark \ref{rem-comp-leave}, we have homomorphisms
\[p\colon G\rightarrow \prod_{i}\Gamma(Q_{i})\]
and 
\[p_{j}\colon G\rightarrow \Gamma(Q_{j})\]
by composing $p$ with the projection $\prod_{i}\Gamma(Q_{i})\rightarrow \Gamma(Q_{j})$ for each $j$. Note that the image $p_{j}(G)$ either is trivial or contains a pseudo-Anosov element in $\Gamma(Q_{j})$ by Theorem \ref{thm-ivanov-last-hope}. 

\begin{defn}
In the above notation, 
\begin{enumerate}
\item[(i)] we call $Q_{j}$ a {\it trivial subsurface}\index{trivial subsurface} for $G$ if $p_{j}(G)$ is trivial.
\item[(ii)] we call $Q_{j}$ an {\it irreducible and amenable subsurface}\index{irreducible and amenable!subsurface} for $G$ if $p_{j}(G)$ contains a pseudo-Anosov element and is amenable.
\item[(iii)] we call $Q_{j}$ an {\it irreducible and non-amenable subsurface}\index{irreducible and non-amenable!subsurface} for $G$ if $p_{j}(G)$ contains a pair of independent pseudo-Anosov elements.
\end{enumerate}
\end{defn}

For simplicity, we call these three subsurfaces {\it T}, {\it IA} and {\it IN} ones, respectively.\index{T!subsurface}\index{IA!subsurface}\index{IN!subsurface}

\begin{rem}
This classification corresponds to one for subgroups of $\Gamma(Q_{j})$ by McCarthy-Papadopoulos \cite{mc-pa} (see Theorem \ref{subgroup-classification}). If $Q_{j}$ is a pair of pants, then only the case (i) can occur.
\end{rem}

The following lemma is easy to prove. Note that the kernel of $p$ is contained in the free abelian subgroup generated by the Dehn twists around all the curves in $\sigma$ (see Proposition \ref{lem-blm}). 

\begin{lem}\label{lem-fund-group}
\begin{enumerate}
\item[(i)] If there exists $\alpha \in V(C(Q_{j}))$ such that $p_{j}(g)\alpha =\alpha$ for any $g\in G$, then $Q_{j}$ is trivial.
\item[(ii)] If $G$ is non-amenable, then there exists $j$ such that $Q_{j}$ is an IN subsurface for $G$. 
\end{enumerate}
\end{lem} 

\begin{pf*}{{\sc Proof of Theorem \ref{main-group-version}}.}
We may assume that $G_{1}\times \cdots \times G_{n}$ is a subgroup of $\Gamma(M; m)$ for some integer $m\geq 3$ and that $n\geq 2$ since $g+[(g+p-2)/2]\geq 1$ when $3g+p-4\geq 0$. Any element in $G_{i}$ can not be pseudo-Anosov by Theorem \ref{centralizer}. Thus, each $G_{i}$ is reducible by Theorem \ref{subgroup-classification}. Let us denote by $M_{\sigma(G_{i})}$ the resulting surface by cutting $M$ along (a realization of) $\sigma(G_{i})$. 

Remark that for $i\neq j$, any curve in $\sigma(G_{i})$ is invariant for $G_{j}$ and thus, it does not intersect any curve in $\sigma(G_{j})$ by the definition of essential reduction classes. Hence, the union $\sigma(G_{i})\cup \sigma(G_{j})$ is an element in $S(M)$. When $\sigma(G_{i})\cup \sigma(G_{j})$ is realized disjointly on $M$ and we denote its realization by the same symbol, we identify each component of $M_{\sigma(G_{i})}$ (resp. $M_{\sigma(G_{j})}$) with the corresponding open component of $M\setminus \sigma(G_{i})$ (resp. $M\setminus \sigma(G_{j})$).
  
\begin{lem}\label{app-lem-disjoint}
Let $i\neq j$ and realize $\sigma(G_{i})\cup \sigma(G_{j})\in S(M)$ disjointly on $M$. Let $Q_{i}$ and $Q_{j}$ be IN subsurfaces for $G_{i}$ and $G_{j}$, respectively. 
\begin{enumerate}
\item[(i)] There exist no curves in $\sigma(G_{i})$ (resp. $\sigma(G_{j})$) contained in $V(C(Q_{j}))$ (resp. in $V(C(Q_{i}))$). Thus, $Q_{j}$ (resp. $Q_{i}$) is contained in some component of $M\setminus \sigma(G_{i})$ (resp. $M\setminus \sigma(G_{j})$).

\item[(ii)] The two open components $Q_{i}$ and $Q_{j}$ are disjoint.
\end{enumerate}
\end{lem}

\begin{pf}
Let $Q_{i}$ be an IN subsurface for $G_{i}$. Since any curve in $\sigma(G_{j})$ is invariant for $G_{i}$, we see that there exist no curves $\alpha$ in $\sigma(G_{j})$ with $\alpha \in V(C(Q_{i}))$ by the definition of IN subsurfaces. This shows the assertion (i). It follows that $Q_{i}$ is contained in some component $Q_{j}'$ of $M\setminus \sigma(G_{j})$. If $Q_{j}'$ is an IN subsurface for $G_{j}$, then the above argument proves that $Q_{i}=Q_{j}'$. 

Thus, we have shown that if $Q_{i}$ and $Q_{j}$ are IN subsurfaces for $G_{i}$ and $G_{j}$, respectively, then either $Q_{i}=Q_{j}$ or $Q_{i}\cap Q_{j}=\emptyset$. 

If we had the possibility of $Q=Q_{i}=Q_{j}$, then there would exist the natural homomorphism 
\[q\colon G_{i}\times G_{j}\rightarrow \Gamma(Q)\]
and each of the images $q(G_{i})$ and $q(G_{j})$ would contain an independent pair of pseudo-Anosov elements and thus, contain a free subgroup of rank $2$ (see the comment right after Theorem \ref{subgroup-classification}). Since any elements in $q(G_{i})$ and $q(G_{j})$ commute each other, this contradicts Theorem \ref{centralizer}.
\end{pf}

Return to the proof of Theorem \ref{main-group-version}. We realize the union $\sigma =\bigcup_{i=1}^{n}\sigma(G_{i})\in S(M)$ on $M$ disjointly and fix them. By Lemma \ref{app-lem-disjoint} (i), there exist no curves in $\sigma$ contained in $V(C(Q_{i}))$ for some $i$ and some IN subsurface $Q_{i}$ for $G_{i}$. This means that any IN subsurface $Q_{i}$ for $G_{i}$ is a component of $M\setminus \sigma$ for each $i$. It follows from Lemma \ref{app-lem-disjoint} (ii) that for each $i\neq j$, any IN subsurfaces $Q_{i}$ and $Q_{j}$ satisfy $Q_{i}\cap Q_{j}=\emptyset$. Since each $G_{i}$ has an IN subsurface by Lemma \ref{lem-fund-group} (ii), the assertion (i) in Theorem \ref{main-group-version} follows from Lemma \ref{main-geometric-lem} (iii).

For the assertion (ii), take $\tau \in S(M)$ such that the number $n(\tau)$ of components of $M_{\tau}$ not of type $(0, 3)$ is equal to $g+[(g+p-2)/2]$ (see Lemma \ref{main-geometric-lem} (iii)). For each component $Q$ of $M_{\tau}$ not of type $(0, 3)$, there exist two intersecting simple closed curves $\alpha$, $\beta$ in $V(C(Q))$. Then sufficiently large powers of their Dehn twists are in $\Gamma =\Gamma(M; m)$ and generate a group containing a free subgroup of rank $2$ by Theorems \ref{tits-alternative-mcg} and \ref{commuting-dehn}. The subgroup generated by these free groups is isomorphic to the direct product of $n(\tau)$ free groups of rank $2$.    
\end{pf*}

\begin{rem}
In Theorem \ref{main-group-version}, we compute the maximal number of factors of the direct product of non-amenable groups contained in the mapping class group as a subgroup. It follows from the next theorem that we can calculate the same number for ${\Bbb Z}$ instead of non-amenable groups.

\begin{thm}[\ci{Lemma 8.8}{ivanov1}]
Let $G$ be an abelian subgroup of $\Gamma(M; m)$, where $M$ is a surface of type $(g, p)$ with $\kappa(M)\geq 0$ and $m\geq 3$ is an integer. Then $G$ is a free abelian group of rank less than or equal to $3g-3+p$.
\end{thm}  
\end{rem}


\backmatter


\printindex


\begin{thebibliography}{99}

\bibitem{adams1}S. Adams.
\newblock Boundary amenability for word hyperbolic groups and an application to smooth dynamics of simple groups.
\newblock  {\it Topology} {\bf 33} (1994), 765--783.

\bibitem{adams2}S. Adams.
\newblock Indecomposability of equivalence relations generated by word hyperbolic groups.
\newblock {\it Topology} {\bf 33} (1994), 785--798. 

\bibitem{aeg}S. Adams, G. A. Elliott and T. Giordano.
\newblock Amenable actions of groups.
\newblock {\it Trans. Amer. Math. Soc.} {\bf 344} (1994), 803--822.

\bibitem{ana}C. Anantharaman-Delaroche and J. Renault.
\newblock {\it Amenable groupoids.} 
\newblock  Monogr. Enseign. Math., 36. Enseignement Math., Geneva, 2000.  

\bibitem{bell-fuji}G. Bell and K. Fujiwara.
\newblock The asymptotic dimension of a curve graph is finite.
\newblock Preprint 2005. math.GT/0509216. 

\bibitem{bestvina-fujiwara}M. Bestvina and K. Fujiwara.
\newblock Bounded cohomology of subgroups of mapping class groups.
\newblock {\it Geom. Topol.} {\bf 6} (2002), 69--89. 


\bibitem{bb}S. J. Bigelow and R. D. Budney.
\newblock The mapping class group of a genus two surface is linear.
\newblock {\it Algebr. Geom. Topol.} {\bf 1} (2001), 699--708. 

\bibitem{blm}J. S. Birman, A. Lubotzky and J. McCarthy.
\newblock Abelian and solvable subgroups of the mapping class groups. 
\newblock {\it Duke Math. J.} {\bf 50} (1983), 1107--1120.


\bibitem{bowditch}B. H. Bowditch.
\newblock Tight geodesics in the curve complexes.
\newblock Preprint 2003.
\newblock Available at {\tt http://www.maths.soton.ac.uk/staff/Bowditch}

\bibitem{bridson-haefliger}M. R. Bridson and A. Haefliger.
\newblock {\it Metric spaces of non-positive curvature.} 
\newblock Grundlehren Math.\ Wiss., 319. 
\newblock Springer-Verlag, Berlin, 1999.  



\bibitem{bcm}J. F. Brock, R. D. Canary and Y. N. Minsky.
\newblock The classification of Kleinian surface groups, II: The Ending Lamination Conjecture.
\newblock Preprint 2004. math.GT/0412006.


\bibitem{cheeger-gromov}J. Cheeger and M. Gromov.
\newblock $L_{2}$-cohomology and group cohomology.
\newblock {\it Topology} {\bf 25} (1986), 189--215. 


\bibitem{cfw}A. Connes, J. Feldman and B. Weiss.
\newblock An amenable equivalence relation is generated by a single transformation. 
\newblock {\it Ergodic Theory Dynam. Systems} {\bf 1} (1981), 431--450 (1982).


\bibitem{cz}M. Cowling and R. J. Zimmer.
\newblock Actions of lattices in ${\rm Sp}(1,n)$.
\newblock {\it Ergodic Theory Dynam. Systems} {\bf 9} (1989), 221--237. 

\bibitem{dye}H. A. Dye.
\newblock On groups of measure preserving transformation. I.
\newblock {\it Amer. J. Math.} {\bf 81} (1959), 119--159.


\bibitem{dye2}H. A. Dye.
\newblock On groups of measure preserving transformations. II.
\newblock {\it Amer. J. Math.} {\bf 85} (1963), 551--576.


\bibitem{dykema}K. J. Dykema.
\newblock Exactness of reduced amalgamated free product $C^{*}$-algebras. 
\newblock {\it Forum Math.} {\bf 16} (2004), 161--180.


\bibitem{ek}C. J. Earle and I. Kra.
\newblock On isometries between Teichm\"uller spaces.
\newblock {\it Duke Math. J.} {\bf 41} (1974), 583--591.


\bibitem{farb}B. Farb.
\newblock The quasi-isometry classification of lattices in semisimple Lie groups. 
\newblock {\it Math. Res. Lett.} {\bf 4} (1997), 705--717.


\bibitem{FLM}B. Farb, A. Lubotzky and Y. Minsky.
\newblock Rank-1 phenomena for mapping class groups.
\newblock {\it Duke Math. J.} {\bf 106} (2001), 581--597.


\bibitem{FLP}A. Fathi, F. Laudenbach and V. Po\'enaru.
\newblock {\it Travaux de Thurston sur les surfaces.} S\`eminaire Orsay.
\newblock Ast\`erisque, 66--67. Soc. Math. France, Paris, 1979.

\bibitem{fm}J. Feldman and C. C. Moore.
\newblock Ergodic equivalence relations, cohomology, and von Neumann algebras. I.
\newblock {\it Trans. Amer. Math. Soc.} {\bf 234} (1977), 289--324.

\bibitem{fsz}J. Feldman, C. E. Sutherland and R. J. Zimmer.
\newblock Subrelations of ergodic equivalence relations.
\newblock {\it Ergodic Theory Dynam. Systems} {\bf 9} (1989), 239--269.

\bibitem{furman1}A. Furman.
\newblock Gromov's measure equivalence and rigidity of higher rank lattices.
\newblock {\it Ann. of Math. (2)} {\bf 150} (1999), 1059--1081.


\bibitem{furman2}A. Furman.
\newblock Orbit equivalence rigidity.
\newblock {\it Ann. of Math. (2)} {\bf 150} (1999), 1083--1108.

\bibitem{gab-cost}D. Gaboriau.
\newblock Co\^ut des relations d'\'equivalence et des groupes.
\newblock {\it Invent. Math.} {\bf 139} (2000), 41--98.

\bibitem{gab2}D. Gaboriau.
\newblock On orbit equivalence of measure preserving actions.
\newblock In {\it Rigidity in dynamics and geometry} (Cambridge, 2000), 167--186, Springer, Berlin, 2002.

\bibitem{gab-l2}D. Gaboriau.
\newblock Invariants $\ell^{2}$ de relations d'\'equivalence et de groupes.
\newblock {\it Publ. Math. Inst. Hautes \'Etudes Sci.} No. {\bf 95} (2002), 93--150. 


\bibitem{gab-survey}D. Gaboriau.
\newblock Examples of groups that are measure equivalent to the free group.
\newblock Preprint 2005. math.DS/0503181. To appear in {\it Ergodic Theory Dynam. Systems}.

\bibitem{gardiner}F. P. Gardiner.
\newblock {\it Teichm\"uller theory and quadratic differentials.}
\newblock Pure Appl. Math. (N.Y.). A Wiley-Interscience Publication. John Wiley \& Sons, Inc., New York, 1987. 

\bibitem{ghys-harpe}\'E. Ghys and P. de la Harpe.
\newblock {\it Sur les groupes hyperboliques d'apr\`es Mikhael Gromov.}
\newblock Progr. in Math., 83. Birkh\"auser Boston, Inc., Boston, MA, 1990. 

\bibitem{gromov1}M. Gromov.
\newblock Hyperbolic groups. 
\newblock {\it Essays in group theory}, 75--263, Math. Sci. Res. Inst. Publ., 8, Springer, New York, 1987. 

\bibitem{gromov-kahler}M. Gromov.
\newblock K\"ahler hyperbolicity and $L\sb 2$-Hodge theory. 
\newblock {\it J. Differential Geom.} {\bf 33} (1991), 263--292.


\bibitem{gromov2}M. Gromov.
\newblock Asymptotic invariants of infinite groups.
\newblock In {\it Geometric group theory, Vol.\ 2} (Sussex, 1991), 1--295, London Math.\ Soc.\ Lecture Note Ser., 182, Cambridge Univ. Press, Cambridge, 1993.


\bibitem{ghw}E. Guentner, N. Higson and S. Weinberger.
\newblock The Novikov conjecture for linear groups.
\newblock Preprint 2003.
\newblock Available at {\tt http://www.math.hawaii.edu/\verb|~|erik}, {\tt http://www.math.psu.edu/higson} or {\tt http://www.math.uchicago.edu/\verb|~|shmuel}

\bibitem{guentner-kaminker}E. Guentner and J. Kaminker.
\newblock Geometric and analytic properties of groups.
\newblock In {\it Noncommutative geometry} (Martina Franca, 2000), 253--262, Lecture Notes in Math., 1831, Springer, Berlin, 2004. 






\bibitem{ham}U. Hamenst\"adt.
\newblock Train tracks and the Gromov boundary of the complex of curves.
\newblock Preprint 2005. math.GT/0409611.

\bibitem{ham2}U. Hamenst\"adt.
\newblock Bounded cohomology and isometry groups of hyperbolic spaces.
\newblock Preprint 2005. math.GR/0507097.

\bibitem{ham3}U. Hamenst\"adt.
\newblock Geometry of the mapping class groups I: Boundary amenability.
\newblock Preprint 2005. math.GR/0510116.


\bibitem{harer-zagier}J. Harer and D. Zagier.
\newblock The Euler characteristic of the moduli space of curves. 
\newblock {\it Invent. Math.} {\bf 85} (1986), 457--485.



\bibitem{harvey}W. J. Harvey.
\newblock Boundary structure of the modular group. 
\newblock In {\it Riemann surfaces and related topics: Proceedings of the 1978 Stony Brook Conference}, 245--251, Ann. of Math. Stud., 97, Princeton Univ.\ Press, Princeton, N.J., 1981.


\bibitem{hat-thu}A. Hatcher and W. Thurston.
\newblock A presentation for the mapping class group of a closed orientable surface.
\newblock {\it Topology} {\bf 19} (1980), 221--237. 


\bibitem{higson-roe}N. Higson and J. Roe. 
\newblock Amenable group actions and the Novikov conjectures.
\newblock {\it J. Reine Angew. Math.} {\bf 519} (2000), 143--153.

\bibitem{hjorth-kechris}G. Hjorth and A. S. Kechris.
\newblock Rigidity theorems for actions of product groups and countable Borel equivalence relations.
\newblock To appear in {\it Mem. Amer. Math. Soc.}, 2005.



\bibitem{hudson}J. F. P. Hudson.
\newblock {\it Piecewise linear topology.}
\newblock W. A. Benjamin, Inc., New York-Amsterdam 1969.


\bibitem{ivanov0}N. V. Ivanov.
\newblock Automorphisms of Teichm\"uller modular groups. 
\newblock In {\it Topology and geometry---Rohlin Seminar}, 199--270, Lecture Notes in Math., 1346, Springer, Berlin, 1988. 


\bibitem{ivanov1}N. V. Ivanov.
\newblock {\it Subgroups of Teichm\"uller modular groups.}
\newblock Transl. of Math. Monogr., 115. Amer. Math. Soc., Providence, RI, 1992.


\bibitem{ivanov3}N. V. Ivanov.
\newblock Automorphism of complexes of curves and of Teichm\"uller spaces. 
\newblock {\it Internat. Math. Res. Notices} {\bf 1997}, no. 14, 651--666. 

\bibitem{ivanov4}N. V. Ivanov.
\newblock  Isometries of Teichm\"uller spaces from the point of view of Mostow rigidity.
\newblock In {\it Topology, ergodic theory, real algebraic geometry}, 131--149, Amer. Math. Soc. Transl. Ser. 2, 202, Amer. Math. Soc., Providence, RI, 2001. 

\bibitem{ivanov2}N. V. Ivanov.
\newblock Mapping class groups.
\newblock In {\it Handbook of geometric topology}, 523--633, North-Holland, Amsterdam, 2002.


\bibitem{kai}V. A. Kaimanovich.
\newblock The Poisson formula for groups with hyperbolic properties.
\newblock {\it Ann. of Math. (2)} {\bf 152} (2000), 659--692.


\bibitem{kai-mas}V. A. Kaimanovich and H. Masur.
\newblock The Poisson boundary of the mapping class group.
\newblock {\it Invent. Math.} {\bf 125} (1996), 221--264.


\bibitem{kechris}A. S. Kechris.
\newblock {\it Classical descriptive set theory.} 
\newblock Grad. Texts in Math., 156. Springer-Verlag, New York, 1995. 



\bibitem{kec}A. S. Kechris and B. Miller. 
\newblock {\it Topics in orbit equivalence.}
\newblock Lecture Notes in Math., 1852. Springer-Verlag, Berlin, 2004. 


\bibitem{kirch}E. Kirchberg.
\newblock Exact $C^{*}$-algebras, tensor products, and the classification of purely infinite algebras.
\newblock In {\it Proceedings of the International Congress of Mathematicians, Vol. 1, 2} (Z\"urich, 1994), 943--954, Birkh\"auser, Basel, 1995.


\bibitem{kirch-wass}E. Kirchberg and S. Wassermann.
\newblock Permanence properties of $C^{*}$-exact groups.
\newblock {\it Doc. Math.} {\bf 4} (1999), 513--558



\bibitem{kla}E. Klarreich.
\newblock The boundary at infinity of the curve complex and the relative Teichm\"uller space.
\newblock Preprint 1999.
\newblock Available at {\tt http://nasw.org/users/klarreich}

\bibitem{kor}M. Korkmaz.
\newblock On the linearity of certain mapping class groups. 
\newblock {\it Turkish J. Math.} {\bf 24} (2000), 367--371. 

\bibitem{kor2}M. Korkmaz.
\newblock Low-dimensional homology groups of mapping class groups: a survey. 
\newblock {\it Turkish J. Math.} {\bf 26} (2002), 101--114.

\bibitem{kos}A. A. Kosinski.
\newblock {\it Differential manifolds.} 
\newblock Pure Appl. Math., 138. Academic Press, Inc., Boston, MA, 1993.


\bibitem{lang}S. Lang.
\newblock {\it Introduction to Diophantine approximations.} 
\newblock Second edition. Springer-Verlag, New York, 1995.


\bibitem{lev}G. Levitt.
\newblock On the cost of generating an equivalence relation.
\newblock {\it Ergodic Theory Dynam. Systems} {\bf 15} (1995), 1173--1181.



\bibitem{luck}W. L\"uck.
\newblock {\it $L\sp 2$-invariants: theory and applications to geometry and $K$-theory.}
\newblock Ergeb. Math. Grenzgeb. (3), 44. Springer-Verlag, Berlin, 2002. 

\bibitem{masur2}H. Masur.
\newblock On a class of geodesics in Teichm\"uller space.
\newblock {\it Ann. of Math. (2)} {\bf 102} (1975), 205--221.

\bibitem{masur}H. Masur.
\newblock Interval exchange transformations and measured foliations. 
\newblock {\it Ann. of Math. (2)} {\bf 115} (1982), 169--200.

\bibitem{masur-minsky1}H. A. Masur and Y. N. Minsky.
\newblock Geometry of the complex of curves. I. Hyperbolicity. 
\newblock {\it Invent. Math.} {\bf 138} (1999), 103--149.

\bibitem{masur-minsky2}H. A. Masur and Y. N. Minsky.
\newblock Geometry of the complex of curves. II. Hierarchical structure.
\newblock {\it Geom. Funct. Anal.} {\bf 10} (2000), 902--974.

\bibitem{masur-wolf}H. A. Masur and M. Wolf.
\newblock Teichm\"uller space is not Gromov hyperbolic.
\newblock {\it Ann. Acad. Sci. Fenn. Ser. A I Math.} {\bf 20} (1995), 259--267.

\bibitem{mc}J. McCarthy.
\newblock Normalizers and centralizers of pseudo-Anosov mapping classes.
\newblock Preprint 1982. To appear in {\it Algebr. Geom. Topol.}
\newblock Available at {\tt http://www.math.msu.edu/\verb|~|mccarthy}


\bibitem{mc2}J. McCarthy.
\newblock A ``Tits-alternative" for subgroups of surface mapping class groups.
\newblock {\it Trans. Amer. Math. Soc.} {\bf 291} (1985), 583--612.

\bibitem{mc-pa}J. McCarthy and A. Papadopoulos.
\newblock Dynamics on Thurston's sphere of projective measured foliations.
\newblock {\it Comment. Math. Helv.} {\bf 64} (1989), 133--166.

\bibitem{mc-pa2}J. McCarthy and A. Papadopoulos.
\newblock The visual sphere of Teichm\"uller space and a theorem of Masur-Wolf.
\newblock {\it Ann. Acad. Sci. Fenn. Math.} {\bf 24} (1999), 147--154.

\bibitem{mc-pa3}J. McCarthy and A. Papadopoulos.
\newblock The mapping class group and a theorem of Masur-Wolf.
\newblock {\it Topology Appl.} {\bf 96} (1999), 75--84.


\bibitem{mcmullen}C. T. McMullen.
\newblock The moduli space of Riemann surfaces is K\"ahler hyperbolic.
\newblock {\it Ann. of Math. (2)} {\bf 151} (2000), 327--357.



\bibitem{mms}I. Mineyev, N. Monod and Y. Shalom.
\newblock Ideal bicombings for hyperbolic groups and applications.
\newblock {\it Topology} {\bf 43} (2004), 1319--1344.




\bibitem{minsky}Y. N. Minsky.
\newblock A geometric approach to the complex of curves on a surface. 
\newblock In {\it Topology and Teichm\"uller spaces} (Katinkulta, 1995), 149--158, World Sci. Publishing, River Edge, NJ, 1996. 



\bibitem{ms}N. Monod and Y. Shalom.
\newblock Orbit equivalence rigidity and bounded cohomology. 
\newblock Preprint 2002. To appear in {\it Ann. of Math.}
\newblock Available at {\tt http://www.math.uchicago.edu/\verb|~|monod}

\bibitem{mosher}L. Mosher.
\newblock Homology and dynamics in quasi-isometric rigidity of once-punctured mapping class groups. 
\newblock Preprint 2003. math.GT/0308065.


\bibitem{mvn}F. J. Murray and J. von Neumann.
\newblock On rings of operators.
\newblock {\it Ann. of Math. (2)} {\bf 37} (1936), 116--229.



\bibitem{oka}R. Okayasu.
\newblock Type III factors arising from Cuntz-Krieger algebras.
\newblock {\it Proc. Amer. Math.\ Soc.} {\bf 131} (2003), 2145--2153.

\bibitem{ow}D. S. Ornstein and B. Weiss.
\newblock Ergodic theory of amenable group actions. I. The Rohlin lemma.
\newblock {\it Bull. Amer. Math. Soc. (N.S.)} {\bf 2} (1980), 161--164.

\bibitem{ozawa1}N. Ozawa.
\newblock Amenable actions and exactness for discrete groups.
\newblock {\it C. R. Acad. Sci. Paris S\'er. I Math.} {\bf 330} (2000), 691--695.

\bibitem{ozawa2}N. Ozawa.
\newblock Solid von Neumann algebras. 
\newblock {\it Acta Math.} {\bf 192} (2004), 111--117.

\bibitem{ozawa3}N. Ozawa.
\newblock Weakly exact von Neumann algebras.
\newblock Preprint 2004. math.OA/0411437.


\bibitem{ozawa4}N. Ozawa.
\newblock Boundary amenability of relatively hyperbolic groups.
\newblock Preprint 2005. math.GR/0501555.


\bibitem{ozawa5}N. Ozawa.
\newblock Private communication. September 2005.


\bibitem{penner}R. C. Penner.
\newblock Perturbative series and the moduli space of Riemann surfaces.
\newblock {\it J. Differential Geom.} {\bf 27} (1988), 35--53. 


\bibitem{popa}S. Popa.
\newblock On a class of type ${\rm II}_{1}$ factors with Betti numbers invariants. 
\newblock Preprint 2002. math.OA/0209130. To appear in {\it Ann. of Math.} 

\bibitem{rees}M. Rees.
\newblock An alternative approach to the ergodic theory of measured foliations on surfaces. 
\newblock {\it Ergodic Theory Dynam. Systems} {\bf 1} (1981), 461--488 (1982).


\bibitem{roe}J. Roe.
\newblock {\it Lectures on coarse geometry.}
\newblock Univ. Lecture Ser., 31. Amer. Math. Soc., Providence, RI, 2003. 


\bibitem{rordam}M. R\o rdam.
\newblock Classification of nuclear, simple $C^{*}$-algebras. 
\newblock In {\it Classification of nuclear $C^{*}$-algebras. Entropy in operator algebras}, 1--145, Encyclopaedia Math.\ Sci., 126, Springer, Berlin, 2002.


\bibitem{royden}H. L. Royden.
\newblock Automorphisms and isometries of Teichm\"uller space. 
\newblock In {\it Advances in the Theory of Riemann Surfaces} (Proc.\ Conf., Stony Brook, N.Y., 1969), 369--383, Ann. of Math. Stud., No. 66. Princeton Univ.\ Press, Princeton, N.J., 1971.


\bibitem{sauer}R. Sauer.
\newblock $L^{2}$-invariants of groups and discrete measured groupoids.
\newblock Dissertation, Universit\"{a}t M\"{u}nster, 2003.
\newblock Available at {\tt http://www.math.uni-muenster.de/u/sauerr}







\bibitem{spa}R. J. Spatzier.
\newblock An example of an amenable action from geometry.
\newblock {\it Ergodic Theory Dynam. Systems} {\bf 7} (1987), 289--293.


\bibitem{spa-zim}R. J. Spatzier and R. J. Zimmer.
\newblock Fundamental groups of negatively curved manifolds and actions of semisimple groups.
\newblock {\it Topology} {\bf 30} (1991), 591--601.

\bibitem{tak1}M. Takesaki.
\newblock {\it Theory of operator algebras. I.} Reprint of the first (1979) edition.
\newblock Encyclopaedia of Math. Sci., 124. Operator Algebras and Non-commutative Geometry, 5. Springer-Verlag, Berlin, 2002. 



\bibitem{tu}J. L. Tu.
\newblock Remarks on Yu's ``property A'' for discrete metric spaces and groups.
\newblock {\it Bull. Soc. Math. France} {\bf 129} (2001), 115--139. 

\bibitem{was}S. Wassermann.
\newblock {\it Exact $C^{*}$-algebras and related topics.}
\newblock Lecture Notes Ser., 19. Seoul National University, Research Institute of Mathematics, Global Analysis Research Center, Seoul, 1994. 


\bibitem{yu}G. Yu.
\newblock The coarse Baum-Connes conjecture for spaces which admit a uniform embedding into Hilbert space. 
\newblock {\it Invent. Math.} {\bf 139} (2000), 201--240.


\bibitem{zim0}R. J. Zimmer.
\newblock Hyperfinite factors and amenable ergodic actions.
\newblock {\it Invent. Math.} {\bf 41} (1977), 23--31.


\bibitem{zim1}R. J. Zimmer.
\newblock Amenable ergodic group actions and an application to Poisson boundaries of random walks.
\newblock {\it J. Functional Analysis} {\bf 27} (1978), 350--372. 


\bibitem{zim4}R. J. Zimmer.
\newblock Induced and amenable ergodic actions of Lie groups.
\newblock {\it Ann. Sci. \'Ecole Norm. Sup. (4)} {\bf 11} (1978), 407--428.



\bibitem{zim3}R. J. Zimmer.
\newblock Ergodic actions of semisimple groups and product relations. 
\newblock {\it Ann. of Math. (2)} {\bf 118} (1983), 9--19.

\bibitem{zim2}R. J. Zimmer.
\newblock {\it Ergodic theory and semisimple groups.}
\newblock Monogr. Math., 81. Birkh\"auser Verlag, Basel, 1984. 

\end{thebibliography}
\end{document}